\documentclass{article}

\usepackage{amsmath}
\usepackage{amsfonts}
\usepackage{amssymb}
\usepackage{amsthm}
\usepackage{comment}
\usepackage{epsfig}
\usepackage{psfrag}
\usepackage{csquotes}
\usepackage{bbm}

\usepackage{bookmark}
\usepackage{xcolor}
\usepackage[many]{tcolorbox}
\usepackage{tikz}
\usepackage{lipsum}
\usepackage{esint}
\usepackage{accents}

\usepackage{xcolor}

\usepackage{titlesec}
\titleformat{\paragraph}
{\normalfont\normalsize\bfseries}{\theparagraph}{1em}{}
\titlespacing*{\paragraph}
{0pt}{3.25ex plus 1ex minus .2ex}{1.5ex plus .2ex}

\usepackage{titlesec}
\titleformat{\section}
  {\centering\large\scshape} 
  {\thesection.}{1em}{} 

\titleformat{\subsection}
  {\bfseries} 
  {\thesubsection}{1em}{} 

\usepackage{upgreek}
\usepackage{blindtext}
\usepackage{titling}

\allowdisplaybreaks

\usepackage{mathrsfs}
\usepackage{amscd}
\usepackage[all]{xy}
\usepackage{rotating}
\usepackage{lscape}
\usepackage{amsbsy}
\usepackage{verbatim}
\usepackage{moreverb}
\usepackage{fullpage}
\usepackage{amsmath,bm}

\usepackage{mathtools}

\newtheorem{theorem}{Theorem}
\newtheorem{prop}[theorem]{Proposition}
\newtheorem{lemma}[theorem]{Lemma}
\newtheorem{cor}[theorem]{Corollary}

\theoremstyle{definition}

\newtheorem{exercise*}[exercise]{$\star$ Exercise}
\newtheorem{remark}[theorem]{Remark}  

\theoremstyle{remark}

\def\Id{\mathrm{Id}}

\newcommand\inner[2]{\langle #1, #2 \rangle}
\newcommand\quadvar[1]{\langle #1 \rangle}
\newcommand\E[1]{\mathbb{E}[#1]}
\newcommand\Ex{\mathbb{E}}

\newcommand\R{\mathbb{R}} 
\newcommand\C{\mathbb{C}} 
\newcommand\N{\mathbb{N}} 
\newcommand\T{\mathbb{T}} 
\newcommand\Z{\mathbb{Z}} 
\newcommand\Dil{\textnormal{Dil}} 
\newcommand\tr{\textnormal{tr}} 
\newcommand\Tr{\textnormal{Tr}} 
\newcommand\Prob{\mathbb{P}} 
\newcommand\Cov{\textnormal{Cov}} 

\newcommand\GUE{\textnormal{GUE}}
\newcommand\CUE{\textnormal{CUE}}
\DeclarePairedDelimiter\floor{\lfloor}{\rfloor}
\DeclarePairedDelimiter{\ceil}{\lceil}{\rceil}
\newcommand\D{\textnormal{d}}

\newcommand\HS{{\textnormal{HS}}}
\DeclareMathOperator{\Ai}{Ai}
\newcommand\tGUE{{\widetilde{\textnormal{GUE}}}}
\newcommand{\RN}[1]{
  \textup{\uppercase\expandafter{\romannumeral#1}}}    

\newcommand\hlbrt{\mathcal{H}}
\newcommand\uhlbrt{\mathcal{U}}
\newcommand\vhlbrt{\mathcal{V}}
\newcommand\thlbrt{\tilde{\mathcal{U}}}
\newcommand\supp{\textnormal{supp}}
\newcommand\sgn{{\textnormal{sgn}}}
\newcommand\tepsilon{\tilde{\epsilon}} 
\newcommand\tDelta{\tilde{\Delta}} 
\newcommand\partialbar{\bar{\partial}}
\newcommand{\rhosc}{\rho_{\rm sc}}
\newcommand{\rhoarcsin}{\rho_{\rm arcsin}}

\newcommand\mathds{\mathbbm}

\numberwithin{equation}{section}
\numberwithin{theorem}{section}

\def\Xint#1{\mathchoice
   {\XXint\displaystyle\textstyle{#1}}%
   {\XXint\textstyle\scriptstyle{#1}}%
   {\XXint\scriptstyle\scriptscriptstyle{#1}}%
   {\XXint\scriptscriptstyle\scriptscriptstyle{#1}}%
   \!\int}
\def\XXint#1#2#3{{\setbox0=\hbox{$#1{#2#3}{\int}$}
     \vcenter{\hbox{$#2#3$}}\kern-.5\wd0}}

\def\dashint{\Xint-}

\DeclareMathOperator*{\argmin}{argmin}

\usepackage{authblk}

\title{Non-intersecting Brownian Motions and Gaussian \\
Multiplicative Chaos}
\author{Ahmet Keles}
\affil{\small\textit{Courant Institute, New York University\\aak10037@nyu.edu}}
\date{}

\begin{document}

\maketitle

\begin{abstract}
 We obtain Fisher-Hartwig asymptotics with root and jump type singularities in space-time under the law of the stationary Hermitian Ornstein-Uhlenbeck process, which serve as a dynamical generalization of earlier static results obtained by Riemann-Hilbert methods. This extends previous asymptotics by Krasovsky \cite{krasovsky2007correlations}, Its–Krasovsky \cite{its2007hankel}, and Charlier \cite{charlier2019asymptotics}.

As a  consequence,  fractional powers of the absolute value of the characteristic polynomial  of this process (and the exponential eigenvalues counting process) converge to a two dimensional Gaussian multiplicative chaos measure on an infinite strip in the subcritical phase. The dynamical Fisher-Hartwig asymptotics also provide the leading order of the log-characteristic polynomial, together with optimal bulk rigidity for non-intersecting Brownian motions. These results offer (i) the second connection between random matrix theory and Liouville quantum gravity measures after \cite{bourgade2022liouville}, by proving a dynamical generalization of the single-time convergence to the GMC from \cite{berestycki2018random}, (ii) a dynamical extension of the maximum of the log-characteristic polynomial \cite{lambert2019law} and the optimal rigidity \cite{claeys2021much}.
\end{abstract}
\tableofcontents

\section{Introduction}

\subsection{Background and motivation}

A log-correlated Gaussian field on a connected open set $D\subset \R^d$ is a centered Gaussian random generalized function characterized by a positive semi-definite covariance kernel
\begin{align*}
K(x,y)=\log\frac{1}{|x-y|}+g(x,y)
\end{align*}
for all $x,y\in D$, where $g$ is a continuous, bounded function on $D\times D$ (see e.g.  \cite{duplantier2017log} for  definitions). One of the earliest hints of log-correlated fields from random matrix theory dates back to 1960s, coming from the work of Dyson and Mehta \cite{dyson1963statistical}, though the explicit language of log-correlated fields was not yet developed. In 1980s, a series of studies by Spohn (e.g. \cite{spohn1986equilibrium, spohn1987interacting}) proved the Gaussian field limit for the fluctuation fields (centered empirical spectral distributions) of various particle dynamics including the eigenvalue dynamics for some random matrix models. However, an explicit recognition of log-correlated fields in random matrix theory was still missing. A more structured perspective emerged in the 1990s and early 2000s through two tightly connected directions: (i) the study of fluctuation field limits, e.g. \cite{diaconis1994eigenvalues, johansson1997random, johansson1998fluctuations, diaconis2001linear}, and (ii) the analysis of log-characteristic polynomials, e.g. \cite{keating2000random, hughes2001characteristic}. These developments laid the groundwork for rigorous studies in the later decades, such as many CLT results for the linear statistics, individual particles, and log-characteristic polynomials which are now shown to be universal in the sense that they hold for broad classes of random matrices and particle systems with similar type of log-correlated covariance structures; e.g. see \cite{anderson2006clt, rider2007noise, pastur2011eigenvalue, dumitriu2013functional, leble2018fluctuations, duits2018mesoscopic, bourgade2019gaussian, landon2020applications, landon2022almost, mody2023log, deleporte2024universality} and references therein.

A more recent development linking random matrix theory and log-correlated Gaussian fields arises through Gaussian multiplicative chaos (GMC) measures, the rigorous construction of which was provided by  Kahane \cite{kahane1985chaos} (see \cite[Section 2.1]{rhodes2014gaussian}), that can be formally written as:
\begin{align}\label{eqn:GMC_defn}
\mu_\gamma(\D z) = e^{\gamma X(z)-\frac{\gamma^2}{2}\E{X(z)^2}}\D z, \quad \gamma\in(0,\sqrt{2d})
\end{align} 
where $X$ is a log-correlated Gaussian field on a domain $D\subset \R^d$ and $\D z$ is the Lebesgue measure on $\R^d$ (for more details, see \cite{berestycki2017elementary, shamov2016gaussian, rhodes2011kpz}). The connection between random matrix theory and GMC measure dates back to 2015, when Webb proved that under suitable normalization small powers of the absolute value of the characteristic polynomials of Circular Unitary Ensemble ($\CUE$), Haar distributed random unitary matrix, converge to GMC measure in the $L^2$-phase in \cite{webb2015characteristic}. The main machinery in their proof is the asymptotic analysis of Toeplitz determinants with Fisher-Hartwig singularities obtained via Riemann-Hilbert methods.  This result was later extended in \cite{nikula2020multiplicative} to larger powers of the characteristic polynomial, demonstrating convergence in the $L^1$-phase. Further extensions to different ensembles, such as one-cut unitary invariant ensemble with a regular potential (e.g. GUE), as well as GOE, GSE, and C$\beta$E, have been explored in several works including \cite{kivimae2020gaussian, lambert2021mesoscopic, forkel2021classical, claeys2021much, lambert2023strong, lambert2024subcritical}.

Subsequent to these static results, the emergence of the GMC measure in two dimensions from the characteristic polynomial of random matrix dynamics (unitary Brownian motion) was proven in \cite{bourgade2022liouville}. Dimension two for GMC is of particular interest since, in this case, Gaussian free field (GFF) is a log-correlated Gaussian field, and the exponential of the GFF, that is also known as the  Liouville quantum gravity measure, has wide-ranging applications in mathematical physics (see \cite{duplantier2011liouville, ding2022introduction, berestycki2024gaussian} and references therein). Moreover, as the first dynamic result of this type, \cite{bourgade2022liouville} developed a probabilistic technique to obtain Fisher-Hartwig asymptotics in space-time under the law of unitary Brownian motion.\\

This article proves the second dynamical result in this direction: In Theorem \ref{thm:FH}, we obtain the Fisher-Hartwig asymptotics with merging root and jump-type singularities under the law of the stationary Hermitian Ornstein-Uhlenbeck (OU) process. This result can be situated within the broader study of the moments of characteristic polynomials in random matrix theory. These quantities have been a central focus since the early 2000s, as established in foundational works such as \cite{keating2000random, brezin2000characteristic, mehta2001moments, fyodorov2002negative, strahov2003universal}. This focus has given rise to a wide range of techniques for studying both integer and non-integer powers of the absolute values of the characteristic polynomials. The most prominent methods include integrability, supersymmetry, and the Riemann-Hilbert approach, with further developments and applications discussed in, for example, \cite{krasovsky2007correlations, osipov2010correlations, fyodorov2018characteristic, webb2019moments, assiotis2022joint, charlier2025asymptotics}. Theorem \ref{thm:FH} contributes to the understanding of the joint moments of correlated random matrices, a topic that arises naturally when studying the energy landscapes of high-dimensional random systems \cite{fyodorov2004complexity, auffinger2013random, subag2017complexity}. In this context, the result represents a step in the direction of models of interest in continuous spin-glass theory (correlated Hermitian matrices),  moving beyond the constraints imposed by the periodic symmetries of the unitary case.

A notable consequence of the multi-time Fisher-Hartwig asymptotics is Theorem \ref{thm:GMC}, the convergence of the random measures constructed by the characteristic polynomial to a GMC measure on an infinite strip. Taken together, these results further yield Theorem \ref{thm:max_log_char}, the leading order behavior of the log-characteristic polynomial and the optimal bulk rigidity for the associated non-intersecting particle dynamics. 

Our results can be viewed as (i) an extension of log-correlated field limit result from \cite{spohn1987interacting}, (ii) a dynamical generalization of the main, quantitative theorems of \cite{berestycki2018random,charlier2019asymptotics} for the quadratic external potential, (iii) analogous results in the non-periodic counterpart to the periodic setting of \cite{bourgade2022liouville} and (iv) a dynamical extension of the results on the maximum of log-characteristic polynomial and optimal rigidity in the static settings \cite{lambert2019law,claeys2021much,bourgade2025optimal}.   \\

This paper focuses primarily on non-intersecting random walks as they are ubiquitous in many topics in mathematical physics. Many of the techniques developed in this article can be generalized to different dynamics (other symmetry classes, arbitrary inverse temperatures $\beta$ and external potential), except for a few places that rely on the determinantal point process structure. Since the decoupling based on this structure (see Section \ref{sec:corr_kernel}) serves a gluing role for the rest of the proof, a replacement using a different technique seems to be fundamental for any extensive generalization.

However, only for inverse temperature $\beta=2$ and quadratic external potential do the Dyson dynamics coincide with diffusion processes conditioned not to intersect, making it a particularly central model, for which this paper obtains moment generating functions of singular statistics.

\subsection{Main results}

The Gaussian Unitary Ensemble (GUE) is defined as an ensemble of $N\times N$ Hermitian matrices where $\sqrt{N}H_{ii}$, $\sqrt{2N}\Re(H_{ij})$ and $\sqrt{2N}\Im(H_{ij})$ are independent standard Gaussian random variables, subject to the symmetry constraints on $i$ and $j$. With this convention, the limit of the normalized eigenvalue distribution is given by Wigner's semicircle law on $[-2,2]$, whose probability density function is $\rhosc(x):=\mathds{1}_{(-2,2)}\frac{1}{2\pi}\sqrt{4-x^2}$ and the eigenvalue distribution is 
\begin{align}\label{eqn:eig_val_dist}
\D\mu(\boldsymbol\lambda)=\frac{1}{Z_N}\prod_{1\leq i<j\leq N}|\lambda_i-\lambda_j|^2e^{-N\sum_{i=1}^{N}\frac{\lambda_i^2}{2}}\D \boldsymbol\lambda
\end{align}
where $\D \boldsymbol\lambda$ is the Lebesgue measure on $\{(\lambda_1,\dots,\lambda_N)\in\R^{N}|\lambda_1\leq\dots\leq \lambda_N\}$ and $Z_N$ is the normalization constant. In this article, we study the stationary Hermitian Ornstein-Uhlenbeck (OU) process, which evolves according to the matrix valued OU dynamics
\begin{equation}\label{matrix_valued_OU_dynamics}
\D H_t=\frac{1}{\sqrt{N}}\D B_t-\frac{1}{2}H_t\D t
\end{equation}
starting from its equilibrium measure (i.e. $\GUE$), where $B_t$ is $N\times N$ Hermitian matrix and for $i<j$, $\sqrt{2}\Re (B_{ij}), \ \sqrt{2}\Im (B_{ij})$ and $B_{ii}$ are independent Brownian motions. Under this dynamics, the eigenvalues $\lambda_1(t)\leq \dots \leq \lambda_N(t)$ of $H_t$ can be viewed as non-colliding OU particles and they satisfy an autonomous system of SDEs, which is called the Dyson Brownian motion \cite{dyson1962brownian}, given by
\begin{align} \label{eqn:DB_motion}
\D\lambda_i(t)=\frac{1}{\sqrt{N}}\D (B_i)_t+\Big(\frac{1}{N}\sum_{j:j\neq i}\frac{1}{\lambda_i(t)-\lambda_j(t)}-\frac{\lambda_i(t)}{2}\Big)\D t
\end{align}
where $B_i$'s are independent standard Brownian motions. 

\begin{remark}
By a time-change and rescaling of the dynamics, the matrix-valued OU process $H_t$ can be equivalently viewed as a time-rescaled Hermitian Brownian motion:  $H_t \stackrel{d}{=} e^{-t/2} B_{e^t}$ where $B_t$ is Hermitian Brownian motion with initial condition $B_0=0$. Consequently, the non-colliding OU eigenvalue dynamics are equivalent in law to non-colliding Brownian motion particles under exponential time change and exponential damping: The results stated below can be immediately rephrased in this setting.
\end{remark}

We consider two GMC measures on $\R\times(-2,2)$: $\mu_{\gamma}$ for $\gamma\in[0,2\sqrt{2})$ and $\nu_{\beta}$ for $\beta\in(-2\sqrt{2},2\sqrt{2})$, which are formally given by $e^{\gamma X(t,x)-\frac{\gamma^2}{2}\E{X(t,x)^2}}\D t\D x$ and $e^{\beta Y(t,x)-\frac{\beta^2}{2}\E{Y(t,x)^2}}\D t\D x$ where $X$ and $Y$ are centered log-correlated Gaussian fields with covariances
\begin{gather*} 
\E{X(t,x)X(s,y)}=\sum_{n=1}^{\infty}\frac{T_n(x/2)T_n(y/2)}{n}e^{-|t-s|n/2}, 
\\
 \E{Y(t,x)Y(s,y)}=\frac{\sqrt{4-x^2}}{2}\frac{\sqrt{4-y^2}}{2}\ \sum_{n=1}^{\infty}\frac{U_{n-1}(x/2)U_{n-1}(y/2)}{n}e^{-|t-s|n/2}
\end{gather*}
expressed in terms of the first and second kind of Chebyshev polynomials $T_n$ and $U_n$ respectively. We refer to \cite{kahane1985chaos,rhodes2014gaussian} for the precise meaning of the GMC measures via local averaging.

\begin{remark}\label{rmk:1/2log_corr}
The logarithmic singularity of the covariances may not be immediately apparent. However, as will be discussed in Lemma \ref{lemma:time-space_log_sing}, the singularity of the either sum, as $(t,x)\to(s,y)$, matches that of $\frac{1}{2}\log(\frac{1}{|(t,x)-(s,y)|})$, where the distance is given by $|(t,x)-(s,y)|=\sqrt{|t-s|^2+|x-y|^2}$, i.e. the Euclidean distance. This contrasts with random matrix dynamics for non-Hermitian matrices,  as it was recently proved that log-correlated fields for the {\it parabolic} distance emerge \cite{bourgade2024fluctuations}.
\end{remark}

The following theorem states that both GMC measures $\mu_{\gamma}$ and $\nu_{\beta}$ arise as scaling limits of powers of the absolute value of the characteristic polynomial and of the exponential of the trace involving a jump-type singularity, given by the function $\Xi^{x}:\R\to\R$ defined by 
\begin{align*} 
\Xi^{x}(y):=\begin{cases} \pi/2, &y<x\\ -\pi/2, &y\geq x\end{cases}.
\end{align*}

\begin{theorem}\label{thm:GMC}
Let $(H_t)_t$ be an $N\times N$ stationary Hermitian OU process. Then, for every $\gamma\in[0,2\sqrt{2})$,
\begin{align} \label{eqn:GMC_root_type}
\frac{|\det(H_t-x)|^{\gamma}}{\E{|\det(H_t-x)|^{\gamma}}}\D t\D x\xrightarrow[N\to\infty]{}\mu_{\gamma}(\D t,\D x)
\end{align}
on $\R\times(-2,2)$, where the convergence is in distribution with respect to weak topology of measures. 

Similarly, for every $\beta\in(-2\sqrt{2},2\sqrt{2})$,
\vspace{-.2cm}\begin{align} \label{eqn:GMC_jump_type}
\frac{e^{\beta\cdot \Tr\Xi^{x}(H_t)}}{\E{e^{\beta\cdot \Tr\Xi^{x}(H_t)}}}\D t\D x\xrightarrow[N\to\infty]{}\nu_{\beta}(\D t,\D x)
\end{align}
on $\R\times(-2,2)$ with the same mode of convergence.
\end{theorem}

As noted in Remark \ref{rmk:1/2log_corr}, the Gaussian fields generating the GMC measures are $1/2$-log correlated. Consequently, the ranges of the parameters $\gamma$ and $\beta$ for the $L^1$ phase are determined by $2\sqrt{d}$, rather than $\sqrt{2d}$ as in \eqref{eqn:GMC_defn}.\\

The proof of Theorem \ref{thm:GMC} relies on a dynamical generalization of the asymptotic results concerning Hankel determinants with Fisher-Hartwig singularities. In general, Fisher-Hartwig singularities arise in the asymptotic analysis of Toeplitz and Hankel determinants when their symbols exhibit root-type or jump-type singularities. These singularities lead to (i) deviations from classical asymptotic formulas, such as the strong Szeg\"o theorem, (ii) the emergence of special functions and (iii) intricate structural modifications. These deviations are captured by the Fisher-Hartwig conjectures \cite{fisher1969toeplitz}, which describe how the asymptotics of Toeplitz and Hankel determinants depend on the singularities in the symbol. The conjectures have been rigorously proven in many cases; see for example \cite{deift2011asymptotics} for most advanced results and \cite[Section 6]{deift2012toeplitz} for a comprehensive exposition, including a historical account of its development.  Closer to our setting is the $N$-dimensional Hankel determinant
\begin{align*} 
D_N(w)=\det\Big(\int_{\R} x^{i+j}w(x)\D x\Big)_{i,j=0}^{N-1}
\end{align*}
with symbol $w(x)=e^{f(x)-NV(x)}\prod_{j=1}^{J}|x-E_j|^{\gamma_j}e^{\beta_j\arg^{E_j}(x)}$ where the potential is $V(x)=\frac{x^2}{2}$, $f:\R\to\R$ is a continuous function satisfying regularity and growth conditions, 
$E_1,\dots,E_J\in(-2,2)$ and the jump-singularity $\arg^{E}$ is given by 
\begin{align*} 
\arg^{E}(x):=\frac{x-E}{2}+\Xi^{E}(x)=\frac{x-E}{2}+\frac{\pi}{2}\mathds{1}_{x<E}-\frac{\pi}{2}\mathds{1}_{x\geq E}
\end{align*}
for any $E\in\R$. By Heine's formula, this Hankel determinant admits an integral representation which leads to
\begin{align} \label{eqn:heine}
D_N(w)=Z_N\cdot \Ex_{\GUE}\Big[\prod_{j=1}^{J}|\det(H-E_j)|^{\gamma_j}e^{\beta_j\Tr(\arg^{E_j}(H))} e^{\Tr(f(H))}\Big]
\end{align}
where $Z_N$ is the normalization constant in \eqref{eqn:eig_val_dist}.

In this setting, without the function $f$ in the symbol, the asymptotic formula as $N\to\infty$ was established (i) by Garoni for root-type singularities with positive integer powers ($f=0$, $\gamma_j\in 2\N$, $\beta_j=0$) \cite{garoni2005asymptotics}; (ii) by Krasovsky for root-type singularities with real powers ($f=0$, $\Re\gamma_j>-1$, $\beta_j=0$) \cite{krasovsky2007correlations}; (iii) by Its and Krasovsky for a single jump singularity ($f=0$, $J=1$, $\gamma_1=0$, $\Im\beta_1\in(-\frac{1}{2},\frac{1}{2})$) \cite{its2007hankel}. Further generalizations include a smooth function $f$ that is real-analytic on $[-2,2]$ in the symbol, along with a general one-cut regular potential $V$. In this vein, Berestycki, Webb, and Wong proved the asymptotic formula for the root-type singularities (i.e., $\beta_j=0$) with real $\gamma_j\geq 0$ and, by incorporating the main theorem of \cite{claeys2016random}, established the convergence of small powers of the absolute value of the characteristic polynomial to the GMC measure in the $L^2$-phase, i.e. $\gamma\leq\sqrt{2}$, \cite{berestycki2018random}. Charlier further extended the asymptotic result to include both root-type and jump-type singularities with $\Re\gamma_j>-1$ and $\Im\beta_j\in(-\frac{1}{2},\frac{1}{2})$, \cite{charlier2019asymptotics}. These works provide a complete analysis when the singularities remain in a compact subset of $\{(x_1,\dots,x_j)\in(-2,2)^{J}:x_i\neq x_j \textnormal{ if }i\neq j\}$.

The analysis becomes more subtle when singularities approach each other. Claeys, Its, and Krasovsky  analyze the case with two merging Fisher-Hartwig singularities for Toeplitz determinants, identifying three regimes \cite{claeys2011emergence,claeys2015toeplitz}:  mesoscopic rate where the classical Fisher-Hartwig asymptotics remain valid; microscopic rate where the asymptotics are expressed in terms of a Painlev\'{e} transcendent;  submicroscopic rate, where the asymptotics reduce to those of the single-singularity case. An analogous result for Hankel determinants was obtained by Claeys and Fahs, who computed the asymptotics for two merging root-type singularities \cite{claeys2016random}. Moreover, Fahs derived estimates for Toeplitz determinants up to constants,  with any fixed number of singularities  \cite{fahs2021uniform}.\\

Under the restriction $V(x)=\frac{x^2}{2}$ and $\gamma_j\geq 0$, $\beta\in\R$, our main result, Theorem \ref{thm:FH} below, generalizes  results in \cite{krasovsky2007correlations,its2007hankel,charlier2019asymptotics,berestycki2018random} in three directions,  most importantly the first one which allows the computation of new random matrix statistics for a 2d ambient space:
\begin{enumerate}
\item[(i)] Dynamical Generalization: Replacing the static GUE matrix with the Hermitian OU process, introducing a time-dependent setting for the asymptotic analysis of the expressions of type \eqref{eqn:heine}.
\item[(ii)] Merging Singularities: Allowing any fixed number of singularities to merge up to any mesoscopic scale.
\item[(iii)] Relaxed Regularity Conditions: Extending the asymptotics to a larger class of functions $f$ in the symbol, such as a bump function supported in the bulk of the spectrum and a sum of $O(\log N)$ many compactly supported order $1$ functions.
\end{enumerate}

 Before stating our theorem concerning multi-time Fisher-Hartwig asymptotics, we introduce some related notation. We denote the Fourier-Chebyshev coefficients of function $f:\R\to\R$ by
\begin{align*} 
\hat{f}_n:=(1+\mathds{1}_{n\neq 0})\cdot\int_{-2}^2 f(x)T_n(x/2)\frac{1}{\pi\sqrt{4-x^2}}\D x
\end{align*}
for every $n\in\N$. We define, for $t,s\in\R$ and suitable functions $f,h:\R\to\R$, 
\begin{align*} 
\mathcal{C}(f(H_t),h(H_s)):=\frac{1}{4}\sum_{n=0}^{\infty} e^{-|t-s|n/2}n\hat{f}_n\hat{h}_n
=\int_{-2}^{2}\int_{-2}^{2} f(x) g(x,y;t-s)h(y)\D x\D y
 \end{align*}
where $g(x,y;t):=\frac{-1}{32\pi^2}\frac{1}{\sin\theta_1\sin\theta_2}\Re\Big(\frac{1}{\sin^2\frac{\theta_1-\theta_2+it}{2}}+\frac{1}{\sin^2\frac{\theta_1+\theta_2+it}{2}}\Big)$ with $x=2\cos\theta_1$, $y=2\cos\theta_2$; and we extend the definition bilinearly. We also write $\mathcal{C}(\sum_i f_i(H_{t_i}))$ as a shorthand for $\mathcal{C}(\sum_i f_i(H_{t_i}),\sum_i f_i(H_{t_i}))$.

 A more detailed discussion, along with equivalent formulations, will be given in Section \ref{subsec:multitime_loop_asymp_exploration}. The key point is that this bilinear form describes the covariance structure of the limiting Gaussian field associated with the centered empirical spectral distribution. Specifically, for sufficiently regular test functions $f$ and $h$ the centered linear statistics $\sum_{k=1}^{N}f(\lambda_k(t))-N\int f \rhosc$ and $\sum_{k=1}^{N}h(\lambda_k(s))-N\int h \rhosc $ converge to centered Gaussian random variables with covariance $\mathcal{C}(f(H_t),h(H_s))$.

The bilinear form $\mathcal{C}(\cdot,\cdot)$ is not well-defined when both arguments simultaneously contain log-type or jump-type singularities at the same space-time coordinate. This signals the necessity of Fisher-Hartwig asymptotics to properly handle such interactions. For convenience, to formalize expressions involving singular ($\log$ and $\arg$) and smooth components ($f_i$'s), we introduce the notation $\log^Ex:=\log|x-E|$ and $\mathcal{C}^{\circ}\Big(\sum_{j=1}^J(\gamma_j\log^{E_j}+\beta_j\arg^{E_j})(H_{t_j})+\sum_{i=1}^If_i(H_{s_i})\Big)$ defined as
\begin{align*} 
\mathcal{C}\big(\sum_{i=1}^{I}f_i(H_{s_i})\big)&+2\cdot\mathcal{C}\Big(\sum_{i=1}^{I}f_i(H_{s_i}),\sum_{j=1}^{J}\big(\gamma_j \log^{E_j}+\beta_j\arg^{E_j}\big)(H_{t_j})\Big)
\\
&+\sum_{1\leq j<k\leq J}2\cdot\mathcal{C}\Big(\big(\gamma_j \log^{E_j}+\beta_j\arg^{E_j}\big)(H_{t_j}),\big(\gamma_k \log^{E_k}+\beta_k\arg^{E_k}\big)(H_{t_k})\Big)
\end{align*}
which corresponds to the full expansion of $\mathcal{C}\Big(\sum_{j}\big(\gamma_j \log^{E_j}+\beta_j\arg^{E_j}\big)(H_{t_j})+\sum_{i}f_i(H_{s_i})\Big)$ with all the diagonal terms involving singularities are being removed.

We begin by stating the theorem in terms of implicit quantities to emphasize the structural form of the Fisher-Hartwig asymptotics. Explicit expressions will follow in the subsequent remark.

\begin{theorem}\label{thm:FH} Let $C>0$, $I,J\in\N$, $\Upsilon\in(0,1)$ and $\kappa\in(0,\frac{1}{10C})$. Assume that $\gamma_1,\dots,\gamma_J\in[0,C]$, $\beta_1,\dots,\beta_J\in[-C,C]$; $E_1,\dots,E_J\in[-2+\Upsilon,2-\Upsilon]$ and $t_1,\dots,t_J,s_1,\dots,s_i\in[0,C]$ with separation condition $\min_{1\leq j_1< j_2\leq J}(|(t_{j_1},E_{j_1})-(t_{j_2},E_{j_2})|)>N^{-1+C\kappa}$; $f_1,\dots,f_I\in\mathscr{S}_{C,\kappa/150}$ (see \eqref{eqn:defn_mathscr_S} for definition). Then
\begin{align*} 
\log \Ex\Big[&\prod_{j=1}^{J}\Big(|\det(H_{t_j}-E_j)|^{\gamma_j}e^{\beta_j\Tr(\arg^{E_j}(H_{t_j}))}\Big)\prod_{i=1}^{I} e^{\Tr(f_i(H_{s_i}))}\Big]=N\int \Big(\sum_{j=1}^{J}\big(\gamma_j\log^{E_j}+\beta_j\arg^{E_j}\big)+\sum_{i=1}^{I}f_i\Big)\rhosc 
\\
&+\sum_{j=1}^{J}\log \Ex\Big[|\det(H-E_j)|^{\gamma_j}e^{\beta_j\Tr(\arg_j^{E}(H))}\Big]+\frac{1}{2}\mathcal{C}^{\circ}\Big(\sum_{j}(\gamma_j\log^{E_j}+\beta_j\arg^{E_j})(H_{t_j})+\sum_{i}f_i(H_{s_i})\Big)+O(N^{-\kappa/10^5})
\end{align*}
where,
\begin{multline*} 
\log \Ex\Big[|\det(H-E)|^{\gamma}e^{\beta\Tr(\arg^{E}(H))}\Big]=\log\Big(N^{\frac{\gamma^2+\beta^2}{4}}\frac{G(1+\frac{\gamma}{2}+i\frac{\beta}{2})G(1+\frac{\gamma}{2}-i\frac{\beta}{2})}{G(1+\gamma)}\Big)
\\
+\Big(\frac{\gamma^2}{8}\log(4-E^2)+\frac{\beta\gamma}{4}(\pi-E-2\arccos\frac{E}{2})+\frac{\beta^2}{8}\big(1-2\sqrt{4-E^2}+3\log(4-E^2)\big)\Big)+O(N^{-\kappa/10^5})
\end{multline*}
with $G$ being the Barnes $G$-function. Moreover, the error term is uniform over the choice of $\boldsymbol\gamma$, $\boldsymbol\beta$, $\boldsymbol E$, $\boldsymbol t$, $\boldsymbol s$ and $\boldsymbol f$ satisfying the given conditions.
\end{theorem}
As seen from the structure of the expression above, the first term with the multiplicative factor $N$ serves merely as a normalization. Moreover, for mesoscopic functions or for interactions between singularities that are mesoscopically separated, the classical theory applies, whereas the singularities themselves give rise to Fisher-Hartwig type contributions. 

\begin{remark}\label{rmk:arg_log_C} Denoting $x=2\cos\theta$, $y=2\cos\omega$ for $\theta,\omega\in(0,\pi)$ and using the principal branch for logarithm, i.e. the branch cut $[1,\infty)$ for $\log(1-z)$, we have the following explicit expressions,
\begin{align*}
&\mathcal{C}(\log^x(H_t),\log^y(H_s))=\frac{-1}{2}\Big(\log\big|1-e^{-|t-s|/2+i(\theta-\omega)}\big|+\log\big|1-e^{-|t-s|/2+i(\theta+\omega)}\big|\Big),
\\
&\mathcal{C}(\log^x(H_t),\arg^y(H_s))= \frac{-1}{2i}\Big(\log\frac{1-e^{-|t-s|/2+i(\theta+\omega)}}{|1-e^{-|t-s|/2+i(\theta+\omega)}|}+\log\frac{1-e^{-|t-s|/2+i(\omega-\theta)}}{|1-e^{-|t-s|/2+i(\omega-\theta)}|}\Big)-e^{-|t-s|/2}\frac{x}{4},
\\
&\mathcal{C}(\arg^x(H_t),\arg^y(H_s))=\frac{-1}{2}\Big(\log\big|1-e^{-|t-s|/2+i(\theta-\omega)}\big|-\log\big|1-e^{-|t-s|/2+i(\theta+\omega)}\big|\Big)+e^{-|t-s|/2}\frac{1-\sqrt{4-x^2}-\sqrt{4-y^2}}{4}.
\end{align*}
Explicit simplifications for $\mathcal{C}(\log^x(H_t),f(H_s))$ and $\mathcal{C}(\arg^x(H_t),f(H_s))$ are not available in general, we refer the reader to Section \ref{subsec:multitime_loop_asymp_exploration} for the integral and Chebyshev series representations.
\end{remark}

In the single-time setting, the expressions in Theorem \ref{thm:FH} simplify considerably leading to an explicit asymptotic formula for Hankel determinants with Fisher-Hartwig singularities. A key novelty compared to \cite[Theorem 1.1]{charlier2019asymptotics} is that we allow singularities to merge at a mesoscopic rate (i.e., at rate $N^{-1+\kappa}$ for fixed $\kappa>0$). On the other hand, Charlier considers a general potential $V$, whereas here we focus on quadratic potential.

 We obtain the corollary below by substituting the following identities in addition to Remark \ref{rmk:arg_log_C}. For any $E\in(-2,2)$, $\int \log^E (x)\rhosc (x)\D x=\frac{E^2-2}{4}$ and $\int \arg^E (x)\rhosc (x)\D x=\frac{\pi-E}{2}+\frac{E\sqrt{4-E^2}}{4}-\arccos\frac{E}{2}$. Moreover, 
$$\mathcal{C}(f(H_t),h(H_t))=-\int\big(\dashint\frac{f(y)}{x-y}\rhoarcsin (y)\D y\big) h'(x)\rhosc (x)\D x$$
where $\dashint$ stands for the principal value integral and $\rhoarcsin (x)=\mathds{1}_{(-2,2)}\frac{1}{\pi\sqrt{4-x^2}}$. In particular, $\mathcal{C}(\log^E(H_t),\allowbreak g(H_t))=\frac{\hat{g}_0}{2}-\frac{g(E)}{2}$ and $\mathcal{C}(\arg^{E}(H_t),g(H_t))=\frac{\hat{g}_1}{4}+\frac{\sqrt{4-E^2}}{2}\cdot \dashint \frac{g(y)}{E-y}\rhosc (y)\D y$ for any suitable function $g$. We refer to the beginning of Section \ref{sec:loop_eqn} for the details on these identities.

\begin{cor} Assume that $E_1<E_2<\cdots<E_J$ lie in the bulk and are separated on a mesoscopic scale, and $f$ is a smooth function as in Theorem \ref{thm:FH}. Therefore, for some $\delta>0$ depending on the mesoscopic separation scale:
\begin{align*}
\log &\ \Ex\Big[\prod_{j=1}^{J}\Big(|\det(H-E_j)|^{\gamma_j}e^{\beta_j\Tr(\arg^{E_j}(H))}\Big) e^{\Tr(f(H))}\Big]
\\
=& N\cdot\Big(\sum_{j=1}^{J}\big(\gamma_j\frac{E_j^2-2}{4}+\beta_j\big(\frac{\pi-E_j}{2}+\frac{E_j\sqrt{4-E_j^2}}{4}-\arccos\frac{E_j}{2}\big)\big)+\int f\rhosc\Big)
+\log N\cdot \Big(\sum_{j=1}^{J}\frac{\gamma_j^2+\beta_{j}^2}{4}\Big)+\Gamma+O(N^{-\delta})
\end{align*}
where the constant term decomposes as $\Gamma=\Gamma_1+\Gamma_2+\Gamma_3+\Gamma_4+\Gamma_5+\Gamma_6$ with:
\begin{align*} 
\Gamma_1&=\sum_{j=1}^{J}\log \Big(\frac{G(1+\frac{\gamma_j}{2}+i\frac{\beta_j}{2})G(1+\frac{\gamma_j}{2}-i\frac{\beta_j}{2})}{G(1+\gamma_j)}\Big),
\\
\Gamma_2&=\sum_{j=1}^{J}\Big(\frac{\gamma_j^2}{8}\log(4-E_j^2)+\frac{\beta_j\gamma_j}{4}(\pi-E_j-2\arccos\frac{E_j}{2})+\frac{\beta_j^2}{8}\big(1-2\sqrt{4-E_j^2}+3\log(4-E_j^2)\big)\Big),
\\
\Gamma_3&=-\frac{1}{2}\int\big(\dashint\frac{f(y)}{x-y}\rhoarcsin (y)\D y\big) f'(x)\rhosc (x)\D x,
\\
\Gamma_4&=\sum_{j=1}^J \gamma_j\Big(\frac{1}{2}\int f(x) \rhoarcsin(x)\D x -\frac{f(E_j)}{2}\Big)+\sum_{j=1}^J\beta_j\Big(\frac{1}{4}\int f(x)x\rhoarcsin(x)\D x+\frac{\sqrt{4-E_j^2}}{2}\dashint \frac{f(x)}{E_j-x}\rhosc(x)\D x\Big),
\\
\Gamma_5&=-\sum_{1\leq j<k\leq J}\frac{\gamma_j\gamma_k}{2} \log|E_j-E_k|+\sum_{1\leq j<k\leq J}\gamma_j\beta_k\Big(\frac{\pi}{4}-\frac{E_j+E_k}{4}-\frac{\arccos E_k}{2}\Big)+\beta_j\gamma_k\Big(\frac{3\pi}{4}-\frac{E_j+E_k}{4}-\frac{\arccos E_j}{2}\Big),
\\
\Gamma_6&=\sum_{1\leq j<k\leq J}\frac{\beta_j\beta_k}{4}\Big(1-\sqrt{4-E_j^2}-\sqrt{4-E_k^2}+\log\Big|\frac{4-E_jE_k+\sqrt{(4-E_j^2)(4-E_k^2)}}{4-E_jE_k-\sqrt{(4-E_j^2)(4-E_k^2)}}\Big|\Big).
\end{align*}
\end{cor}

\begin{remark} Compared to \cite{charlier2019asymptotics}, there are a few notational differences:
(i) the random matrix scaling we use differs from that in \cite{charlier2019asymptotics}; (ii) our $\beta$ parameter corresponds to $i\beta$ in \cite{charlier2019asymptotics} since they use the function $e^{i\beta\arg}$ in the symbol; (iii) our $\arg$ function is not purely a step function, but rather the sum of a step function and a linear term. These differences lead to small variations in the final expressions, although the individual terms can still be matched one-to-one.
\end{remark}

A corollary of Theorem  \ref{thm:FH} is the leading order behavior for the maximum of the log-characteristic polynomial, along with the optimal bulk rigidity for the non-intersecting particle system.  To state this result,  we denote $S_N(f)(H)={\rm Tr}\,f(H)-N\int f\rhosc$ the centered linear statistics of a $N
\times N$ matrix $H$.

\begin{theorem}\label{thm:max_log_char}
For any $\Upsilon,\varepsilon>0$ and any compact set $\mathcal A\subset[-2+\Upsilon,2-\Upsilon]\times\R$ with ${\rm area}(\mathcal A)>0$,
\begin{align*} 
\lim_{N\to\infty}\Prob\Big(\max_{(E,t)\in \mathcal A}\frac{S_N(\log^E)(H_t)}{\log N}\in[\sqrt{2}-\varepsilon,\sqrt{2}+\varepsilon]\Big)=1,
\\
\lim_{N\to\infty}\Prob\Big(\max_{(E,t)\in \mathcal A}\Big|\frac{S_N(\arg^E)(H_t)}{\log N}\Big|\in[\sqrt{2}-\varepsilon,\sqrt{2}+\varepsilon]\Big)=1.
\end{align*}
As a consequence of the latter convergence, we obtain the following maximum deviation for the particles holds, where $\gamma_k$'s, the typical locations of the eigenvalues, are determined by $\int_{-2}^{\gamma_k}\rhosc=k/N$:
\begin{align*} 
\lim_{N\to\infty}\Prob\Big(\max_{k:\gamma_k\in[-2+\Upsilon,2-\Upsilon],t\in[0,1]}\frac{N\rhosc(\gamma_k)|\lambda_k(t)-\gamma_k|}{\log N}\in[\frac{\sqrt{2}}{\pi}-\varepsilon,\frac{\sqrt{2}}{\pi}+\varepsilon]\Big)=1.
\end{align*}
\end{theorem}

\begin{remark}
If $t=N^{-1+\alpha}$ with $\alpha\in[0,1]$ then the above results hold with $\sqrt{2}$ substituted with $\sqrt{1+\alpha}$. The proof is the same.
\end{remark}

\subsection{Organization and techniques}

The article is organized around the proof of Theorem \ref{thm:FH} which roughly follows the diagram below: Using the multi-time loop equation, i.e. Theorem \ref{thm:multitime_loopp}, we localize the singularities by cutting out the long range of the singularities and smooth functions. Proof of this theorem requires both the discussion in the relatively short Section \ref{sec:rig_under_biased} and mainly Lemma \ref{lemma:two_times_resolvent_mult} from Section \ref{sec:res_est}. Then we apply Theorem \ref{thm:decoupling} for decoupling the singularities far apart from each other. Using the $\GUE$-$\CUE$ comparison, i.e. Proposition \ref{prop:general_GUE_CUE_comparison}, we turn the single-time $\GUE$ problem on local functions to a single-time $\CUE$ problem. Although the transition from $\GUE$ to $\CUE$ may seem technically unmotivated at first glance,  it plays an essential role in the structure of the argument. A detailed explanation is provided in Appendix \ref{subsec:needforGUE-CUE}, since its necessity involves some technical subtleties. Finally, we recover the long-range singularities by the single-time loop equation. Note that, in the single-time setting, the $\GUE$-$\CUE$ comparison step provides the key with turning a problem concerning Hankel determinants into a problem on Toeplitz determinants. This explains clearly why the asymptotics involve the same structure and special functions.

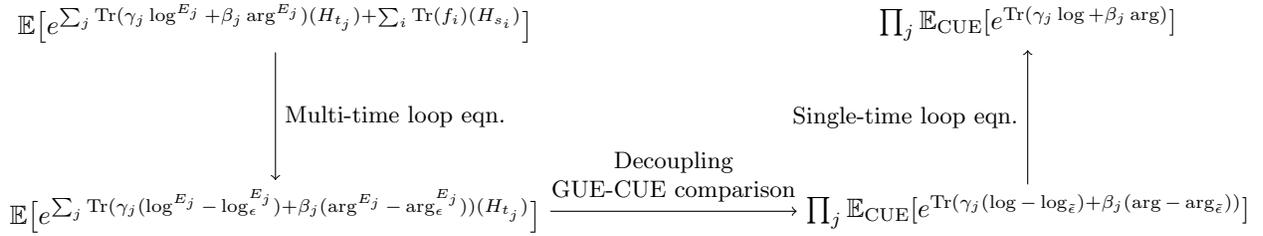
\begin{figure}[h]
    \centering
    \begin{tikzpicture}
        \node (A) at (0,0) {$\Ex\big[e^{\sum_j\Tr(\gamma_j\log^{E_j}+\beta_j\arg^{E_j})(H_{t_j})+\sum_i \Tr(f_i)(H_{s_i})}\big]$};
        \node (B) at (0,-2.5) {$\Ex\big[e^{\sum_j\Tr(\gamma_j(\log^{E_j}-\log_{\epsilon}^{E_j})+\beta_j(\arg^{E_j}-\arg_{\epsilon}^{E_j}))(H_{t_j})}\big]$};
        \node (C) at (10,-2.5) {$\prod_j\Ex_{\CUE}[e^{\Tr(\gamma_j(\log-\log_{\tilde{\epsilon}})+\beta_j(\arg-\arg_{\tepsilon}))}]$};
        \node (D) at (10,0) {$\prod_j\Ex_{\CUE}[e^{\Tr(\gamma_j\log+\beta_j\arg)}]$};
        
        \draw[->] (A) -- node[right] {\small \shortstack{Multi-time loop eqn.}} (B);
        \draw[->] (B) -- node[above] {\shortstack{\small Decoupling \\\small GUE-CUE comparison}} (C);
        \draw[->] (C) -- node[left] {\small Single-time loop eqn.} (D);
    \end{tikzpicture}
    \caption{Successive reductions leading to the proof of Theorem \ref{thm:FH}.}
\end{figure}

 With this picture in mind, we briefly explain each section in order. Section \ref{sec:corr_kernel} builds upon two correlation kernel estimations and their applications: (i) Lemma \ref{lemma:GUE_sine_CUE} provides the necessary tool to compare local $\GUE$ and $\CUE$ statistics and (ii) Proposition \ref{prop:kernel_est} provides a decorrelation result for the extended kernel when space-time points are separated by a mesoscopic scale. Section \ref{sec:rig_under_biased} discusses some sufficient conditions for the rigidity under biased measures on the law of the Hermitian OU process. It particularly addresses the biases constructed by the Laplace transform of linear statistics of smooth functions or functions with logarithmic or jump-type singularities. Throughout Section \ref{sec:res_est}, we demonstrate the stability of rigidity under Dyson Brownian motion and establish a  multi-time local law which plays a pivotal role in proving multi-time loop equation asymptotics in Section \ref{sec:loop_eqn}. If the reader wishes to skip the details regarding the SDEs on resolvent, reading the main definitions introduced in that section and Lemma \ref{lemma:two_times_resolvent_mult} would suffice for the rest of the article. Finally, we complete the proof of Theorem \ref{thm:FH} with the help of the tools developed throughout the paper and various regularization techniques. Then, continue with the convergence to the GMC measure, i.e. the proof of Theorem \ref{thm:GMC}, following the treatment of \cite{webb2015characteristic} and \cite{berestycki2018random}. As a corollary, we obtain the leading order asymptotics for the maximum of log-characteristic polynomial and maximum deviations of the particles, proof of Theorem \ref{thm:max_log_char}, following the argument in \cite{bourgade2025optimal}. In order to avoid interrupting the flow, if a result requires meticulous and lengthy calculations, we moved their proof into the Appendix.\\
 
This article obtains exponential generating functions of singular statistics in the non-periodic setting, which  introduces a number of difficulties as is customary in (S)PDEs.  Compared to the dynamics on the circle as in \cite{bourgade2022liouville},  we highlight a few representative obstacles:
 \begin{enumerate}
 \item[(i)] For the rigidity under biased measure, unlike in \cite[Section 5.2]{bourgade2022liouville}, we cannot invoke a key lemma by Johansson to a-priori bound  Laplace transforms of linear statistics. Instead, we rely on techniques based on local laws and the Helffer-Sj\"ostrand formula.  This motivates 
Section \ref{sec:rig_under_biased}.  
 \item[(ii)] As our setting lacks rotational symmetry,  an a-priori submicroscopic smoothing (Proposition \ref{prop:submic_smoothing}) requires a more technically involved approach (cf. \cite[Lemma 2.3]{bourgade2022liouville}), which relies on  a multidimensional Girsanov theorem.

 \item[(iii)] For the unitary Brownian motion, the extended kernel is expressed directly by trigonometric functions with simple frequencies \cite[Proposition 3.1]{bourgade2022liouville}, making the quantitative analysis readily tractable \cite[Lemma 3.3]{bourgade2022liouville}. In contrast, the kernel for the non-periodic case \eqref{eqn:ext_kern} is constructed by Hermite polynomials, which exhibit only asymptotically trigonometric behavior,  making it challenging to evaluate the effects of the oscillations (see Proposition \ref{prop:kernel_est} and Appendix \ref{app:pf_kernel_est}).

 \item[(iv)] In the CUE case,  the freedom to choose a constant shift parameter allows one to handle one singularity using the single-time loop equation (see Appendix of \cite{bourgade2022liouville}). This flexibility, however, is absent in the $\GUE$ setting, motivating the strategy of translating the problem from the $\GUE$ framework to the $\CUE$ framework (see Section \ref{subsec:needforGUE-CUE}). 
 
 \item[(v)] Finally, as will be clear along the proof, many technical difficulties occur due to the test functions possibly overlapping the edge (control of outliers, lack of contour integral representation of some observables, etc.), difficulties which disappear under rotational symmetry.
\end{enumerate}

\subsection{Notations}\label{subsec:not_and_defn}

When it is more convenient, we use the notations $A_N\lesssim B_N$, $A_N\asymp B_N$ and $A_N\ll B_N$ to describe the asymptotic relations $|A_N|=O(B_N)$ (i.e. there exists a constant $C>0$ such that $|A_N|\leq CB_N$ for all $N\in\N$), $A_N=\Theta(B_N)$ (i.e. there exists a $C>0$ such that $\frac{1}{C}B_N\leq A_N\leq CB_N$ for all $N\in\N$) and $A_N=o(B_N)$ (i.e. $\lim_{N\to\infty} \frac{A_N}{B_N}=0$), respectively. Most of the time, the dependence on $N$, which stands for the size of the matrices, will be omitted in notation. Also, we use $x\wedge y$ and $x\vee y$ for the $\min(x,y)$ and $\max(x,y)$. By $\mathscr{C}^{k}$ and $\mathscr{C}^{k}_{c}$ we mean $k$ times continuously differentiable functions from $\R$ to itself and compactly supported $\mathscr{C}^{k}$ functions respectively. When referring to a finite sequence of functions or parameters $a_1,\dots,a_n$ we use boldface $\boldsymbol a$ to denote the entire sequence.

 In the subsequent sections, we use $\mathbb{E}[\cdot]$ and $\Prob(\cdot)$ without subindices when the underlying measure is $\GUE$ or the Hermitian OU process defined in \eqref{matrix_valued_OU_dynamics}, if further emphasis is not necessary. Whenever we write expectation or probability for a different random matrix ensemble, such as Circular Unitary Ensemble ($\CUE$), it will be indicated by a subindex such as $\mathbb{E}_{\CUE}[\cdot]$. Also, we reserve the letter $H$ (and $H_t$) for $\GUE$ (and stationary Hermitian OU) and the letter $U$ for $\CUE$ distributed random matrix.
 
 Typical spacing between consecutive eigenvalues of $\GUE$ in the bulk are $\asymp 1/N$. Throughout the paper we use the letter $\Delta=N^{-1-\alpha}$ (with $\alpha>0$ small) for submicroscopic scaling and $\epsilon=N^{-1+\kappa}$ (with $\kappa>0$ small) for mesoscopic scaling. We explicitly distinguish  $\epsilon$ and $\varepsilon$: the former will always denote the mesoscopic scale, while the latter will be reserved for arbitrary small constants. Moreover, we use $\varphi=(\log N)^{\log\log N}$ for large, but subpolynomial scale. The symbol $\Upsilon$ will be used exclusively to denote the bulk of the spectrum. $[-2+\Upsilon,2-\Upsilon]$.

 When the underlying measure is $\GUE$ or Hermitian OU, we use the following notations for the centered linear statistics
\begin{align*} 
S_N(f)=S_N(f(H))=S_N(f)(H)=\sum_{i=1}^N f(\lambda_i)-N\int f(x)\rhosc (x)\D x
\end{align*}
and use a similar notation for the trace, $\Tr(f)=\Tr(f(H))=\Tr(f)(H)=\sum_{i=1}^N f(\lambda_i)
$. The $N$ index in the sequence of functions will typically be omitted, for example when $f=f_N:\R\to\R$ is a sequence of functions, then we write $S_N(f)$ for $S_N(f_N)$. When the measure is $\CUE$, the centered linear statistics is defined similarly: For $f:\T\to\R$, $S_N(f):=S_N(f)(U)=\sum_{k=1}^N f(\theta_k)-N\int_\T f(\theta)\frac{1}{2\pi}\D \theta$ where $\theta_k$'s are eigenangles.

Throughout the article, we study Fisher-Hartwig asymptotics for submicroscopically regularized root and jump-type singularities, addressing the relation between the formulas for the singularities and their regularizations only at the final stage. For consistency and clarity, we adhere to the following notations for bump and step functions as well as $\log$ and $\arg$ regularizations. Let $\chi:\R\to\R$ be a fixed smooth, even bump function such that $\chi(x)=1$ on $[-1,1]$, $\chi(x)=0$ on $[-2,2]^c$ and $2\geq\chi'(x)\geq 0$ for all $x\in[-2,-1]$. Define scaled bump function by $\chi_{\varepsilon}(x):=\chi(x/\varepsilon)$ for every $\varepsilon\geq 0$ (as convention, set $\chi_{0}=0$). In addition, denote the bump function centered at $E\in \R$ by $\chi_\varepsilon^E$, i.e. $\chi_\varepsilon^E(x):=\chi_\varepsilon(x-E)$. Also, we define an order one smooth right and left regularizations of step function $\frac{\pi}{2}\mathds{1}_{x<0}-\frac{\pi}{2}\mathds{1}_{x\geq0}$ by $\Xi_r,\Xi_\ell:\R\to\R$ such that
\begin{equation*}
\Xi_r(x)=\begin{cases}\pi/2,& \textnormal{if }x\in(-\infty,0]\\ -\pi/2,&\textnormal{if } x\in[1,\infty)\end{cases},\quad \Xi_\ell(x)=\begin{cases}\pi/2,& \textnormal{if }x\in(-\infty,-1]\\ -\pi/2,&\textnormal{if } x\in[0,\infty)\end{cases}
\end{equation*}
with smooth transitions in between so that $\|\Xi_r'\|_{\infty},\|\Xi_\ell'\|_{\infty}\in [0,2\pi]$. Similarly, define $\Xi_{r,\varepsilon}(x):=\Xi_{r}(x/\varepsilon)$, and $\Xi_{r,\varepsilon}^E(x):=\Xi_{r,\varepsilon}(x-E)$ for every $\varepsilon\geq 0$ and $E\in\R$ (and the same for the left regularization). 

We always use the symmetrized version of logarithm function, i.e. $\log (-x):=\log x$ for all $x\in \R^+$; and similar to the bump function, we denote the logarithm centered at $E$ by $\log^E(x):=\log(x-E)$. By the regularization of logarithm of scale $\varepsilon$ around $E$ we mean the smooth function given by
\begin{align*}
\log_\varepsilon^E(x):=\log^E(x)(1-\chi_\varepsilon^E(x))+\log(2\varepsilon)\chi_\varepsilon^E(x)=\log^E(x)-\chi_\varepsilon^E(x)\log(\frac{x-E}{2\varepsilon})
\end{align*}
for which we have $\log\leq\log_{\varepsilon_1}\leq\log_{\varepsilon_2}$ for every $0\leq \varepsilon_1\leq\varepsilon_2$. In addition, we define the right and left regularizations of $\arg$ function similarly,
\begin{align*} 
\arg_{r,\varepsilon}(x):=\frac{x}{2}+\Xi_{r,\varepsilon}(x),\quad \arg_{\ell,\varepsilon}(x):=\frac{x}{2}+\Xi_{\ell,\varepsilon}(x)
\end{align*}
with the same notation $\arg_{r,\varepsilon}^{E}(x):=\arg_{r,\varepsilon}(x-E)$ for any $\varepsilon\geq0$ and $E\in\R$ (and the same for the left regularization).

We now define a function space that captures a fundamental property shared by both regularized logarithms and step functions, enabling us to obtain the asymptotic results for these type of singularities. For every $C>0$ and $\varepsilon>0$, we define a function space $\mathscr{A}_{C,\varepsilon}$ as follows
\begin{equation*}
\mathscr{A}_{C,\varepsilon}:=\Big\{f\in\mathscr{C}^{4}\Big|\supp(f)\subseteq[A-C\varepsilon,A+C\varepsilon] \textnormal{ for some }A\in[-2+\Upsilon,2-\Upsilon], \|f^{(k)}\|_{\infty}\leq C\varepsilon^{-k},\ k=0,1,2,3,4 \Big\}.
\end{equation*}
Notice that a regularized log-singularity can be decomposed into functions in the $\mathscr{A}$ space as,
\begin{align*}
\log_{\varepsilon}=(\log_{\varepsilon}-\log_{2\varepsilon})+(\log_{2\varepsilon}-\log_{4\varepsilon})+\cdots+(\log_{2^{J-1}\varepsilon}-\log_{2^{J}\varepsilon})+\log_{2^{J}\varepsilon}
\end{align*}
and similarly for the regularized jump function,
\begin{align*} 
\Xi_{r,\varepsilon}=(\Xi_{r,\varepsilon}-\Xi^{-4\varepsilon}_{r,2\varepsilon})+(\Xi_{r,2\varepsilon}^{-4\varepsilon}-\Xi^{-16\varepsilon}_{r,4\varepsilon})+\cdots+(\Xi_{r,2^{J-1}\varepsilon}^{-4^{J-1}\varepsilon}-\Xi^{-4^{J}\varepsilon}_{r,2^{J}\varepsilon})+\Xi^{-4^{J}\varepsilon}_{r,2^{J}\varepsilon}.
\end{align*}
To capture these decompositions, we define
 \begin{align} 
\mathscr{S}_{C,\kappa,N}:=\Big\{f\in\mathscr{C}^{4}\Big| f=\sum_{j=1}^{\lfloor C\log N\rfloor}h_j \textnormal{ where } h_j\in \mathscr{A}_{C,\epsilon_{j}} \textnormal{ for some } \epsilon_{j}\in[N^{-1+\kappa},1]  \Big\} \label{eqn:defn_mathscr_S}
 \end{align}
for $C>0$ and $\kappa\in(-\infty,1]$. We say, a sequence of functions $f=f_N$ is in $\mathscr{S}_{C,\kappa}$ if for every positive integer $N$ we have, $f_N\in\mathscr{S}_{C,\kappa,N}$. For example, for any $\kappa<1$, $\epsilon=N^{-1+\kappa}$ and $E\in[-2+\Upsilon,2-\Upsilon]$, we have $\chi_{10}\log_{\epsilon}^E,\chi_{10}\arg_{r,\epsilon}^E,\chi_{10}\arg_{\ell,\epsilon}^E\in\mathscr{S}_{C,\kappa}$ for some $C>0$ by the decompositions described above.\\

\noindent\textbf{Acknowledgement.} The author is deeply grateful for Paul Bourgade's invaluable support during the preparation and revision of this work.

\section{Correlation kernel estimates and applications}\label{sec:corr_kernel}

A simple point process on a Polish space is called determinantal if there exists a correlation kernel $K$ such that for every $k\in\N$, the $k$-point correlation function $\rho_k$ is given by
\begin{align*} 
\rho_k(x_1,\dots,x_k)=\det[K(x_i,x_j)]_{i,j=1}^{k}
\end{align*}
where $\rho_k$'s are defined by the identity
\begin{align*} 
\mathbb{E}\Big[\sum_{\substack{i_1,\dots,i_k\in[\![N]\!]\\\textnormal{and distict}}} \delta_{(\lambda_{i_1},\dots,\lambda_{i_k})}(B_1\times\cdots\times B_k) \Big]=\int_{B_1\times\cdots\times B_k} \rho_k(x_1,\dots,x_k)\D x_1\dots \D x_k
\end{align*}
for all Borel sets $B_i$'s in the Polish space. Heuristically, $\rho_k$ can be thought of as
\begin{align*} 
\rho_k(x_1,\dots,x_k):=\lim_{\varepsilon\to0}\frac{\Prob(\textnormal{there is a particle in }[x_i,x_i+\varepsilon) \textnormal{ for all } i\in[\![k]\!])}{\varepsilon^k}
\end{align*}
when such a limit exists.
 
It is well known that the eigenvalues of GUE and CUE form simple determinantal point processes. The CUE correlation kernel is
\begin{equation}\label{eqn:CUE_corr_kern}
K_\CUE(x,y)=\frac{1}{2\pi}\frac{\sin (N\frac{x-y}{2})}{\sin(\frac{x-y}{2})}
\end{equation}
where $x$ and $y$ stands for the eigenangles (see \cite{dyson1970correlations}). It is evident that, in the limit $N\to\infty$, after rescaling around any point, this kernel converges to the sine kernel. On the other hand, for GUE, in the bulk and edge it converges to sine and Airy point processes respectively as $N\to\infty$ where the bulk limit behaviour of the eigenvalues can be quantified via the following asymptotics from \cite[Corollary 1]{delyon2006spectral}:
\begin{equation}\label{eqn:GUE_correlation_kernel_approximation}
K_\GUE(x,y)=\frac{1}{4\pi\sqrt{\sin\theta\sin\omega}}\Big(\frac{\sin(\frac{N}{2}(a_\theta-a_\omega))}{\sin(\frac{1}{2}(\omega-\theta))}+O(1)\Big)
\end{equation}
uniformly for $x$, $y$ in the bulk; where $x=2\cos\omega$, $y=2\cos\theta$ with $\theta,\omega\in (0,\pi)$ and $a_\omega=\sin2\omega-2\omega$, $a_\theta=\sin2\theta-2\theta$. 

From these expressions it is not hard to obtain that, when properly scaled, the correlation kernels for $\GUE$ and $\CUE$ are locally very close. More explicitly, we have the following lemma, proof of which follows directly from Taylor expansions of the necessary terms.
\begin{lemma}\label{lemma:GUE_sine_CUE}
Let $E$ be a point in the bulk, $\kappa\in[0,1/2]$ and $C>0$ be constants. Denoting $\epsilon=N^{-1+\kappa}$, the following holds uniformly in $x,y\in[E-C\epsilon,E+C\epsilon]$:
\begin{align*} 
K_\GUE&(x,y)=K_\CUE(2\pi\rhosc (E)x,2\pi\rhosc (E)y) \cdot 2\pi\rhosc (E)+O(N^\kappa).
\end{align*}
\end{lemma}
\begin{proof}\renewcommand{\qedsymbol}{}
The proof follows from \eqref{eqn:CUE_corr_kern} and Taylor expanding \eqref{eqn:GUE_correlation_kernel_approximation}. See Appendix \ref{app:pf_GUE_CUE} for details.
\end{proof}

Furthermore, the Karlin-McGregor theorem implies that both Hermitian Brownian motion starting from $\GUE$ and unitary Brownian motion starting from $\CUE$ maintain the determinantal point process structure for their eigenvalue dynamics. Explicit expressions for the corresponding correlation kernels can be found in \cite{tracy2004differential,johansson2005non,katori2007noncolliding} and \cite{bourgade2022liouville}. The correlation kernel for the eigenvalues of the stationary OU dynamics in \eqref{matrix_valued_OU_dynamics}, that is called the extended kernel, can be written as:
\begin{align}\label{eqn:ext_kern}
K(t,x;s,y)=\begin{dcases}
e^{-(N-\frac{1}{2})(s-t)/2}\sqrt{N}\sum_{k=0}^{N-1}(e^{(s-t)/2})^k\psi_k(x\sqrt{N})\psi_k(y\sqrt{N}), &\textnormal{if} \ t\leq s \\
-e^{-(N-\frac{1}{2})(s-t)/2}\sqrt{N}\sum_{k=N}^{\infty}(e^{(s-t)/2})^k\psi_k(x\sqrt{N})\psi_k(y\sqrt{N}), &\textnormal{if} \ t> s
\end{dcases}
\end{align}
where $\psi_k$'s are the Hermite functions (see appendix \ref{app:pf_kernel_est} for the definition and some significant properties of Hermite functions). We remark that thanks to the Christoffel–Darboux formula the summation significantly simplifies when $s=t$ making the kernel particularly useful for various calculations. However, to the best of our knowledge, when $s\neq t$, an analogous simplification does not exist and even the most straight-forward looking calculations demand meticulous analysis.

The diagonal terms of all three correlation kernels ($\CUE$, $\GUE$ and Hermition OU in the bulk) are $\asymp N$. The decorrelation of mesoscopically distant particles is evident in the sine kernel, i.e. the $\CUE$ case. More explicitly, the corresponding correlation kernel decays on a polynomial scale when evaluated at points separated by $N^{-1+\kappa}$. Although a similar decay is naturally expected in the multi-time setting, it is not directly visible from the expression of the extended kernel. The following proposition makes this claim precise by providing a quantitative estimate for the polynomial-scale decay.

\begin{prop}\label{prop:kernel_est} Uniformly in $\kappa\in [0,\frac{1}{12}]$, and uniformly in the pairs $(t,x),(s,y)\in\R\times[-2+\Upsilon,2-\Upsilon]$ for which $|x-y|\vee |s-t|>N^{-1+\kappa}$ holds we have 
\begin{align*}
K(t,x;s,y)=O(N^{1-\kappa/8})
\end{align*}
where $K$ is the extended kernel.
\end{prop}

\begin{proof}\renewcommand{\qedsymbol}{}
The proof (in Appendix \ref{app:pf_kernel_est}) relies on a detailed study of the asymptotically oscillatory behavior of Hermite polynomials, for which we directly quote the Riemann-Hilbert results from \cite{deift1999strong} presented in Appendix \ref{app:pf_kernel_est}.
\end{proof}

\subsection{$\GUE-\CUE$ comparison for local functions}\label{subsec:GUE_CUE_comp}
 
In the final stage of the proof of Theorem \ref{thm:FH}, we encounter a technical difficulty arising from a Hilbert transform condition on the loop equation for $\GUE$ due to submicroscopically regularized logarithms. This issue, however, can be resolved easily when the underlying measure is $\CUE$ instead. The purpose of the Proposition \ref{prop:general_GUE_CUE_comparison} below is to provide the method for translating the study of linear statistics of locally supported functions under the law of $\GUE$ to $\CUE$. The technique, which is based on the comparison between the Fredholm determinants, utilizes the kernel comparison for $\GUE$ and $\CUE$, i.e. Lemma \ref{lemma:GUE_sine_CUE}. 

Before presenting the main comparison result, we introduce a few notational conventions to simplify the expressions. When $f:\R\to\R$ is supported on a small neighborhood around a point $E\in\R$ and the underlying measure is $\CUE$, we restrict the domain of $f$ to $[E-\pi,E+\pi)$ and define $S_N(f)$ and $\Tr(f)$ via the eigenangles chosen in $[E-\pi,E+\pi)$, i.e. $\Tr(f)=\sum_{k=1}^N f(\theta_k)$ and $S_N(f):=\sum_{i=k}^N f(\theta_k)-N\int_{E-\pi}^{E+\pi} f(\theta)\frac{1}{2\pi}\D \theta$. 

Secondly, to compare the linear statistics of a local function in $\GUE$ and $\CUE$, we first need to match the local densities by rescaling one of them. For this purpose, we define dilation of a function $f:\R\to\R$ as $\Dil_E f(x):=f(\frac{x}{2\pi\rhosc (E)})$, and dilated $\GUE$ distribution as $\Dil_E\GUE\overset{\textnormal{(d)}}{=} 2\pi\rhosc (E)\cdot\GUE$ for $E\in(-2,2)$.

\begin{prop}\label{prop:general_GUE_CUE_comparison}
Let $C>0$ large and $\Upsilon>0$ small be fixed. Take $\kappa\in[0,1/30]$ and $f:\R\to\R$ be a continuous functions satisfying:
\begin{itemize}
\item[•] $\supp(f)\subseteq[E-CN^{-1+\kappa},E+CN^{-1+\kappa}]$ for some $E\in[-2+\Upsilon,2-\Upsilon]$ and,
\item[•] $0\geq f(x)\geq -(1-10\kappa)\log N$ for all $x\in\R$.
\end{itemize}
Then (the error term below is uniform over the choice of $E$ and $f$ satisfying the given properties)
\begin{align*} 
\log \mathbb{E}_{\GUE}[e^{\Tr(f)}]=\log \mathbb{E}_{\CUE}[e^{\Tr(\Dil_E f)}]+O(N^{-\kappa})
\end{align*}
 As a corollary, $\log \mathbb{E}_{\GUE}[e^{S_N(f)}]=\log \mathbb{E}_{\CUE}[e^{S_N(\Dil_E f)}]+O(N^{-\kappa})$.
\end{prop}

\begin{proof} Note that, by the definition of dilation on function and the matrix distribution,  $\mathbb{E}_{\GUE}[e^{\Tr(f)}]=\mathbb{E}_{\Dil_E\GUE}[e^{\Tr(\Dil_E f)}]$. For convenience, we will denote $\Dil_E\GUE$ by $\tGUE$ for the rest of the proof.

 First, we define an integral operator and restate the problem in terms of Fredholm determinants. Define function $k$ on $\R$ (or $\T$) by $k(x)=\sqrt{1-e^{\Dil_E f(x)}}$ and integral operator on $L^2$, $\mathcal{K}_{\tGUE}$ (similarly $\mathcal{K}_\CUE$) with convolution kernel $\mathcal{K}_{\tGUE}(x,y):=k(x)K_{\tGUE}(x,y)k(y)$ where $K_{\tGUE}$ is the correlation kernel of the determinantal process for the eigenvalues of $\tGUE$. For any $x\in\R$, using Fredholm determinant theory for integral operators (see, e.g., \cite{widom2011integral}) we have, 
\begin{equation}\label{eqn:fredholm_det}
\prod (1-x\lambda_j^{\tGUE})=\det(1-x\mathcal{K}_{\tGUE})=\mathbb{E}_{\tGUE}[\prod(1+x(e^{\Dil_E f(\lambda_i)}-1))]
\end{equation}
where $\lambda_j^{\tGUE}$'s are the eigenvalues of the integral operator $\mathcal{K}_{\tGUE}$ on $L^2$ and $\lambda_i$'s in the right-hand side are eigenvalues of the underlying matrix distribution, i.e. $\tGUE$. First, note that $K_{\tGUE}(x,y)$ is a positive semi-definite kernel because the $k$-point correlation function is always non-negative for any $k\in\N$. So, by the definition of positive definite kernel, it is easy to see that this implies positive-semi definiteness of the kernel $k(x)K_{\tGUE}(x,y)k(y)$. Thus, when evaluated as an integral operator it will be a positive-semi definite operator due to all Riemann sum approximations being non-negative. That means $\lambda_j^{\tGUE}\geq 0$ for all $j$. On the other hand, we define $\varepsilon:=\min_{x\in\R} e^{\Dil_E f(x)}\geq N^{-1+10\kappa}$, because the range of $e^{\Dil_E f}$ is in $[\varepsilon,1]$, the right hand side of the equation above is strictly positive for all $x<\frac{1}{1-\varepsilon}$. This leads $\lambda_j^{\tGUE}\in [0,1-\varepsilon]$ for all $j$. Similarly, the same holds for the integral operator $\mathcal{K}_\CUE$, i.e. the spectrum $\sigma(\mathcal{K}_\CUE)\subset [0,1-\varepsilon]$.

Plugging $x=1$ into \eqref{eqn:fredholm_det} we obtain that it suffices to prove:
\begin{equation}\label{log tGUE CUE ratio}
|\log \det(1-\mathcal{K}_{\tGUE})-\log \det(1-\mathcal{K}_\CUE)| =O(N^{-\kappa}).
\end{equation}
Using $\log \det (\Id-\mathcal{K}_\tGUE)=\sum_j \log (1-\lambda_j^\tGUE)=-\sum_{j,k\geq 1} \frac{(\lambda_j^\tGUE)^k}{k}$, we obtain that for any $m\geq 2$,
\begin{align}\label{eqn:difference_of_fredholms}
\hspace{-1cm}|\log \det&(\Id-\mathcal{K}_\tGUE) - \log \det(\Id-\mathcal{K}_\CUE)|=\Big| \sum_{k=1}^m \frac{\sum_j(\lambda_j^\tGUE)^k-\sum_j(\lambda_j^\CUE)^k}{k} + \sum_j \sum_{k=m+1}^\infty \frac{(\lambda_j^\tGUE)^k}{k}-\sum_j \sum_{k=m+1}^\infty \frac{(\lambda_j^\CUE)^k}{k}\Big|\nonumber
\\
&\leq |\tr(\mathcal{K}_\tGUE)-\tr(\mathcal{K}_\CUE)|+\sum_{k=2}^m |\frac{\tr((\mathcal{K}_\tGUE)^k)-\tr((\mathcal{K}_\CUE)^k)}{k}|+ \frac{1}{m} \Big(\sum_j \frac{|\lambda_j^\tGUE|^2}{1-|\lambda_j^\tGUE|}
+\sum_j \frac{|\lambda_j^\CUE|^2}{1-|\lambda_j^\CUE|}\Big)\nonumber
\\
&\leq |\tr(\mathcal{K}_\tGUE)-\tr(\mathcal{K}_\CUE)|+\sum_{k=2}^m |\frac{\tr((\mathcal{K}_\tGUE)^k)-\tr((\mathcal{K}_\CUE)^k)}{k}|+ \frac{1}{m\varepsilon} (\|\mathcal{K}_\tGUE\|_\HS^2+\|\mathcal{K}_\CUE\|_\HS^2)
\end{align}
where we have used that the spectra of both operators lie in $[0,1-\varepsilon]$.

In order to bound the terms of the remaining sum, we use an infinite dimensional version of Hoffman-Wielandt inequality (see \cite[Theorem 1]{bhatia1994hoffman}), and say that there exist enumerations of eigenvalues of $\mathcal{K}_\tGUE$ and $\mathcal{K}_\CUE$ such that $\sum_j |\lambda_j^\tGUE-\lambda_j^\CUE|^2\leq \|\mathcal{K}_\tGUE-\mathcal{K}_\CUE\|_\HS^2+N^{-1}$. Then, for any $k\geq2$:
\begin{align*}
\frac{|\tr((\mathcal{K}_\tGUE)^k)-\tr((\mathcal{K}_\CUE)^k)|}{k} \leq & \sum_j |\lambda_j^\tGUE-\lambda_j^\CUE|(|\lambda_j^\tGUE|^{k-1}+|\lambda_j^\CUE|^{k-1})
\\
\leq & \Big( \sum_j |\lambda_j^\tGUE-\lambda_j^\CUE|^2 \Big)^{1/2}\Big[\big( \sum_j |\lambda_j^\tGUE|^{2k-2}\big)^{1/2}+\big( \sum_j |\lambda_j^\CUE|^{2k-2}\big)^{1/2}\Big]
\\
\leq & \big(\|\mathcal{K}_\tGUE-\mathcal{K}_\CUE\|_\HS+N^{-1}\big)\Big[\big( \sum_j |\lambda_j^\tGUE|^{2}\big)^{1/2}+\big( \sum_j |\lambda_j^\CUE|^{2}\big)^{1/2}\Big]
\\
= &\big(\|\mathcal{K}_\tGUE-\mathcal{K}_\CUE\|_\HS+N^{-1}\big) (\|\mathcal{K}_\tGUE\|_\HS+\|\mathcal{K}_\CUE\|_\HS).
\end{align*}

Substituting this into \eqref{eqn:difference_of_fredholms} we obtain
\begin{equation}\label{eqn:log tGUE CUE ratio bounded}
\begin{aligned}
|\log \det(1-\mathcal{K}_\tGUE)-\log \det(1-\mathcal{K}_\CUE)| \leq& |\tr(\mathcal{K}_\tGUE)-\tr(\mathcal{K}_\CUE)|+\frac{1}{m\varepsilon}(\|\mathcal{K}_\tGUE\|_\HS^2+\|\mathcal{K}_\CUE\|_\HS^2) 
\\
&+m \big(\|\mathcal{K}_\tGUE-\mathcal{K}_\CUE\|_\HS+N^{-1}\big) (\|\mathcal{K}_\tGUE\|+\|\mathcal{K}_\CUE\|_\HS).
\end{aligned}
\end{equation}

Thus, using Lemma \ref{Lemma:GEU_CUE_Kernel_comparison} for the remaining terms in \eqref{eqn:log tGUE CUE ratio bounded}, and substituting $\varepsilon\geq N^{-1+10\kappa}$ we obtain:
\begin{align*}
|\log \det(1-\mathcal{K}_\tGUE)-\log \det(1-\mathcal{K}_\CUE)| = O\big(N^{-1+\kappa}+\frac{1}{m}N^{1-8\kappa}+ m  N^{-1+3\kappa}\big).
\end{align*}
 Taking $m=\lfloor N^{1-5\kappa}\rfloor$, we get \eqref{log tGUE CUE ratio}.
\end{proof}

\begin{lemma}\label{Lemma:GEU_CUE_Kernel_comparison}
Define $\mathcal{K}_\tGUE$ and $\mathcal{K}_\CUE$ as in the proof of Proposition \ref{prop:general_GUE_CUE_comparison}. Then,
\begin{gather*}
\|\mathcal{K}_\CUE\|_\HS =O(N^{\kappa}), \quad \|\mathcal{K}_\tGUE\|_\HS=O(N^{\kappa}),
\\
\|\mathcal{K}_\tGUE-\mathcal{K}_\CUE\|_\HS=O(N^{-1+2\kappa}),\quad
|\tr(\mathcal{K}_\tGUE)-\tr(\mathcal{K}_\CUE)|=O(N^{-1+\kappa}).
\end{gather*} 
\end{lemma}

\begin{proof} Set $\alpha=2\pi\rhosc (E)$. We have a trivial bound $|K_\CUE|\leq \frac{N}{2\pi}$, which leads to an easy estimation for $\|\mathcal{K}_\CUE\|_\HS$ due to the support size of $k$:
\begin{align*}
\|\mathcal{K}_\CUE\|_\HS^2= \int\int |k(x)K_{\CUE}(x,y)k(y)|^2\D x\D y \leq  \int\int_{[-\alpha C N^{-1+\kappa},\alpha CN^{-1+\kappa}]^2} N^2\D x\D y=O(N^{2\kappa}).
\end{align*}
We also have $K_\tGUE(x\alpha,y\alpha)-K_\CUE(x\alpha,y\alpha)=O(N^\kappa)$ by Lemma \ref{lemma:GUE_sine_CUE} whenever $x,y\in[E-CN^{-1+\kappa},E+CN^{-1+\kappa}]$. We already have $K_\CUE(x\alpha,y\alpha)=O(N)$, so is $K_\tGUE(x\alpha,y\alpha)=O(N)$ in the domain of integration. So, similarly we have $\|\mathcal{K}_\tGUE\|_\HS=O(N^{\kappa})$. Using Lemma \ref{lemma:GUE_sine_CUE}, we obtain
\begin{align*}
\|\mathcal{K}_\tGUE-\mathcal{K}_\CUE\|_\HS^2\leq &\int\int_{[\alpha E-\alpha CN^{-1+\kappa},\alpha E+\alpha C N^{-1+\kappa}]^2}|K_\tGUE(x,y)-K_\CUE(x,y)|^2\D x\D y=O(N^{-2+4\kappa})
\end{align*}
and by the asymptotics in \eqref{eqn:GUE_correlation_kernel_approximation},
\begin{align*}
|\tr(\mathcal{K}_\tGUE)-\tr(\mathcal{K}_\CUE)|=&\Big|\int_{\alpha E-\alpha CN^{-1+\kappa}}^{\alpha E+\alpha CN^{-1+\kappa}} k(x)K_\tGUE(x,x)k(x)-k(x)K_\CUE(x,x)k(x) \D x\Big|
\\
\lesssim&\int_{E-CN^{-1+\kappa}}^{E+CN^{-1+\kappa}} |\big(\frac{N}{2\pi}+O(1)\big)-\frac{N}{2\pi}| \D x=O(N^{-1+\kappa}),
\end{align*}
concluding the proof.
\end{proof}

\subsection{Space-time decoupling for local functions}

We now present a general multi-time decoupling theorem for the Laplace transform of linear statistics of local functions that are separated in space-time on a mesoscopic scale. The proof follows a similar approach to the decoupling theorem in \cite[Proposition 3.4]{bourgade2022liouville}. Although this is a fairly general decoupling statement, it represents the fundamental obstacle to any further generalization to other random matrix dynamics as it relies on the determinantal point process structure of the particle dynamics we study.

\begin{theorem}\label{thm:decoupling}
Let $\kappa\in(0,\frac{1}{1200})$, $\Upsilon>0$ and $J\in\mathbb{N}$ be fixed. Assume that $t_1,\dots,t_J\in\R$, $E_1,\dots,E_J\in[-2+\Upsilon,2-\Upsilon]$ with separation condition $\min_{i\neq j}|(t_i,E_i)-(t_j,E_j)|\geq N^{-1+100\kappa}$, and $f_1,\dots,f_J:\R\to\R$ are continuous functions satisfying the following properties
\begin{itemize}
\item[•] $\supp(f_j)\subseteq[E_j-N^{-1+\kappa},E_j+N^{-1+\kappa}]$ and,
\item[•] $0\geq f_j(x)\geq -\kappa\log N$ for all $x\in\R$.
\end{itemize}
We have the following decoupling, where the error term is uniform over the choice of $\boldsymbol t$, $\boldsymbol E$ and $\boldsymbol f$ satisfying the given properties:
\begin{align*}
\mathbb{E}\Big[\prod_{j=1}^J e^{ \Tr(f_j)(H_{t_j})}\Big]=\prod_{j=1}^J\mathbb{E}\big[e^{ \Tr(f_j)}\big]\big(1+O(N^{-\kappa})\big).
\end{align*}
\end{theorem}

\begin{proof}
We discuss the case $J=2$ for convenience, $J\geq3$ can be done by the same way. We can assume that $t_1\neq t_2$; if they are equal, following the same steps for $t_2=t_1+N^{-10}$ we obtain the decoupling and it is easy to see that the deviation of $\Tr(f_2)(H_{t_1})-\Tr(f_2)(H_{t_1+N^{-10}})$ is negligible. Define $k_i:=\sqrt{1-e^{f_{i}}}$ and denote the associated convolution kernel (the integral operator) for the extended kernel $K$ by $\mathcal{K}(i,x;j,y)=k_i(x)K(t_i,x;t_j,y)k_j(y)$ for $i,j\in\{1,2\}$. Then, as shown in \cite[Section 2]{tracy2004differential},
$
\E{e^{\Tr(f_{1})(H_{t_1})+\Tr(f_{2})(H_{t_2})}}=\det(\Id-\mathcal{K})=\det(\Id-\mathcal{K}^*)
$, 
therefore
\begin{align*}
\E{e^{\Tr(f_{1})(H_{t_1})+\Tr(f_{2})(H_{t_2})}}^2=\det((\Id-\mathcal{K})(\Id-\mathcal{K}^*)).
\end{align*}
In addition, let $\tilde{\mathcal K}(i,x;j,y)$ be the convolution kernel defined as $k_i(x)K(i,x;i,y)k_i(y)$ if $i=j$, and $0$ if $i\neq j$. We may thus write
\begin{align*}
\E{e^{\Tr(f_{1})}}^2\E{e^{\Tr(f_{2})}}^2=\det((\Id-\tilde{\mathcal{K}})(\Id-\tilde{\mathcal{K}}^*)).
\end{align*}
Note that $\tilde{\mathcal{K}}$ is already self-adjoint since $k_i(x)K(i,x;i,y)k_i(y)=k_i(y)K(i,y;i,x)k_i(x)$ for $i=1,2$. Define $\mathbf{K}=\mathcal{K}+\mathcal{K}^*-\mathcal{K}\mathcal{K}^*$ and $\tilde{\mathbf{K}}=\tilde{\mathcal{K}}+\tilde{\mathcal{K}}^*-\tilde{\mathcal{K}}\tilde{\mathcal{K}}^*$ and denote the eigenvalues of these self-adjoint operators by $\lambda^{\mathbf{K}}_j$ and $\tilde{\lambda}^{\mathbf{K}}_j$. As we showed in the proof of Proposition \ref{prop:general_GUE_CUE_comparison}, the spectrum of the integral operator $\tilde{\mathcal K}$, $\sigma(\tilde{\mathcal K})\subseteq[0,1-\varepsilon]$ where $\varepsilon:=\min_{i=1,2}\min_{x\in\R} e^{f_i(x)}\geq N^{-\kappa}$. Due to the compactness of the integral operator $\tilde{\mathcal K}$ (as it is already an Hilbert-Schmidt integral operator), we get that $\sigma(\tilde{\mathbf{K}})=\{2\lambda-\lambda^2|\lambda\in\sigma(\tilde{\mathcal K})\}\subset [0,1-\varepsilon^2]$.

Returning to the main question, the quantity we want to bound can be expressed in terms of Fredholm determinants as follows:
\begin{align*}
\log \Big(\frac{\E{e^{S_N(f_{1})(H_{t_1})+S_N(f_{2})(H_{t_2})}}^2}{\E{e^{S_N(f_{1})}}^2\E{e^{S_N(f_{2})}}^2}\Big)=\log \det(\Id-\mathbf{K}) - \log \det(\Id-\tilde{\mathbf{K}}).
\end{align*}
This difference can be bounded similar to \eqref{eqn:difference_of_fredholms}, 
\begin{align*}
|\log \det(\Id-\mathbf{K}) - \log \det(\Id-\tilde{\mathbf{K}})| \leq \sum_{k=1}^m |\frac{\tr(\mathbf{K}^k)-\tr(\tilde{\mathbf{K}}^k)}{k}|+ \frac{1}{m+1}\sum_j \frac{|\lambda^{\mathbf{K}}_j|^2}{1-|\lambda^{\mathbf{K}}_j|}
+\frac{1}{m+1}\sum_j \frac{|\lambda^{\tilde{\mathbf{K}}}_j|^2}{1-|\lambda^{\tilde{\mathbf{K}}}_j|},
\end{align*}
but here we lack an a priori range for the spectrum of $\mathbf{K}$ which is necessary to bound the sum $\sum_j \frac{|\lambda^{\mathbf{K}}_j|^2}{1-|\lambda^{\mathbf{K}}_j|}$ and obtain \eqref{eqn:log tGUE CUE ratio bounded}. Nevertheless, by comparing the eigenvalues of $\mathbf{K}$ and $\tilde{\mathbf{K}}$ via infinite dimensional Hoffman-Wielandt inequality \cite[Theorem 1]{bhatia1994hoffman} we can get sufficiently small range for the spectrum of $\mathbf{K}$ as follows: By Hoffman-Wielandt inequality there exist enumerations of eigenvalues of $\mathbf{K}$ and $\tilde{\mathbf{K}}$ such that
\begin{align*}
\sum_j |\lambda^{\mathbf{K}}_j-\lambda^{\tilde{\mathbf{K}}}_j|^2&\leq \|\mathbf{K}-\tilde{\mathbf{K}}\|_\HS^2+N^{-1}=\|(\mathcal{K}+\mathcal{K}^*-\mathcal{K}\mathcal{K}^*)-(\tilde{\mathcal{K}}+\tilde{\mathcal{K}}^*-\tilde{\mathcal{K}}\tilde{\mathcal{K}}^*) \|_\HS^2+N^{-1}
\\
&\lesssim \|\mathcal{K}-\tilde{\mathcal{K}}\|_\HS^2+(\|\mathcal{K}\|_\HS^2+\|\tilde{\mathcal{K}}\|_\HS^2)\|\mathcal{K}-\tilde{\mathcal{K}}\|_\HS^2+N^{-1}=O(N^{-21\kappa})
\end{align*}
where we have used the bound on $K$ in Proposition \ref{prop:kernel_est} with the fact that $k_i$'s have support sizes $O(N^{-1+\kappa})$ to get $\|\mathcal{K}-\tilde{\mathcal{K}}\|_\HS^2=O(N^{2\kappa-25\kappa})$ and as in Lemma \ref{Lemma:GEU_CUE_Kernel_comparison} we have $\|\mathcal{K}\|_\HS^2=O(N^{2\kappa})$, $\|\tilde{\mathcal{K}}\|_\HS^2=O(N^{2\kappa})$. On the other hand, using $\lambda^{\mathbf{K}}_j \geq \lambda^{\tilde{\mathbf{K}}}_j-|\lambda^{\mathbf{K}}_j-\lambda^{\tilde{\mathbf{K}}}_j|\geq -|\lambda^{\mathbf{K}}_j-\lambda^{\tilde{\mathbf{K}}}_j|$, we get
\begin{align*}
\sum_{j:\lambda^{\mathbf{K}}_j<0} |\lambda^{\mathbf{K}}_j|^2 <\sum _j|\lambda^{\mathbf{K}}_j-\lambda^{\tilde{\mathbf{K}}}_j|^2=O(N^{-21\kappa})
\end{align*}
and so, for all sufficiently large $N$, $-1/2<\lambda^{\mathbf{K}}_j \leq \lambda^{\tilde{\mathbf{K}}}_j+|\lambda^{\mathbf{K}}_j-\lambda^{\tilde{\mathbf{K}}}_j|\leq 1-\varepsilon^2+O(N^{-21\kappa/2})=1-\varepsilon^2+o(\varepsilon^2)$. Therefore, for all sufficiently large $N$, $\sigma(\mathbf{K})\subset(-1/2,1-\varepsilon^2/2)$.

As a result,  with exactly the same method as in the proof of Proposition \ref{prop:general_GUE_CUE_comparison},
\begin{align*}
|\log \det(\Id-\mathbf{K}) - \log \det(\Id-\tilde{\mathbf{K}})|&\lesssim m\|\mathbf{K}-\tilde{\mathbf{K}}\|_{\HS}(\|\mathbf{K}\|_{\HS}+\|\tilde{\mathbf{K}}\|_{\HS})+\frac{1}{m\varepsilon^2}(\|\mathbf{K}\|_{\HS}^2+\|\tilde{\mathbf{K}}\|_{\HS}^2)
\\
&=O(mN^{-17\kappa/2}+\frac{1}{m}N^{6\kappa}).
\end{align*}
where we have substituted $\|\mathbf{K}\|_{\HS}\leq 2\|\mathcal{K}\|_{HS}+\|\mathcal{K}\|_{HS}^2=O(N^{2\kappa})$, (similarly) $\|\tilde{\mathbf{K}}\|_{\HS}=O(N^{2\kappa})$ and $\|\mathbf{K}-\tilde{\mathbf{K}}\|_{HS}=O(N^{-21\kappa/2})$. Choosing $m=\lfloor N^{7\kappa}\rfloor$, we get $|\log \det(\Id-\mathbf{K}) - \log \det(\Id-\tilde{\mathbf{K}})|= N^{-\kappa}$ which finishes the proof. 
\end{proof}

\section{Rigidity under a biased measure}\label{sec:rig_under_biased}

The eigenvalues of $\GUE$ are typically located within the interval $(-2,2)$ and $\frac{1}{N}\sum_k\delta_{\lambda_k}$ converges to the semicircle distribution $\rhosc (x)\D x=\frac{1}{2\pi}\sqrt{(4-x^2)_+}\D x$ as $N\to\infty$. In the microscopic scale, it has been shown that each eigenvalue $\lambda_k$ is concentrated around its typical location, denoted by $\gamma_k$, which are determined by $\int_{-2}^{\gamma_k}\rhosc (x)\D x=\frac{k}{N}$. This concentration, called eigenvalue rigidity, has been shown to be universal in random matrix theory in the sense that the eigenvalues of generalized Wigner matrices are close to their typical locations with overwhelming probability (see \cite{erdHos2012rigidity}).

The proof Theorem \ref{thm:FH} begins with localizing singularities by cutting of the long ranges via loop equation techniques.  Since this requires the rigidity phenomenon to hold under certain biased measures, we introduce the relevant notation and discuss general conditions for rigidity.

For any $J\in\N$, real numbers $t_1,\dots,t_J$ and functions $f_1,\dots,f_J:\R\to\R$, we denote the biased probability measure with bias proportional to the exponential of $\sum_{j=1}^J S_N(f_j(H_{t_j}))$ by
\begin{align} \label{eqn:bias_defn}
\Ex_{\sum_{j=1}^{J}f_j(H_{t_j})}[\ \cdot\  ]:=\Ex\Big[\ \cdot \ \frac{e^{\sum_{j=1}^{J}S_N(f_j(H_{t_j}))}}{\E{e^{\sum_{j=1}^{J}S_N(f_j(H_{t_j}))}}}\Big]
\end{align}
We denote the probability under this biased measure by $\Prob_{\sum_{j=1}^{J}f_j(H_{t_j})}(\cdot)$ as well.   In the case of a single time problem, we simply write $\Ex_{f}[\cdot]$ for $\E{\cdot \frac{e^{S_N(f)}}{\E{e^{S_N(f)}}}}$. We first describe two simple conditions for the bias to preserve rigidity of eigenvalues at a time and then provide sufficient conditions that are easier to verify.

\begin{lemma}\label{lemma:rigidity_conds}
Let $J\in\N$ and $f_1,\dots,f_J:\R\to\R$ be sequences of functions satisfying the following: There exists a constant $C>0$ such that for all $N>1$,
\begin{align*} 
\mathbb{E}_{\GUE}[e^{2J S_N(f_j)}]\leq e^{(\log N)^{C}}\ \textnormal{ and }\  \mathbb{E}_{\GUE}[|S_N(f_j)|]\leq (\log N)^{C}
 \end{align*} 
for each $j=1,\dots, J$. Then, for all $D>0$ there exists an $N_0$ such that for all $N\geq N_0$ and $t,t_1,\dots,t_J\in\R$:
\begin{align*} 
\Prob_{\sum_{j}f_i(H_{t_j})}(\boldsymbol \lambda(t)\in\tilde{\mathcal{G}})\geq 1-e^{-(\log N)^{D}}
\end{align*}
where (we denote $\hat{k}:=\min(k,N+1-k)$ and $\varphi:=(\log N)^{\log\log N}$)
\begin{equation*}
\mathcal{\tilde G}:=\Big\{\Big(\lambda_1,\dots,\lambda_N)\in\R^N: |\gamma_k-\lambda_k|<\frac{\varphi}{N^{2/3}\hat{k}^{1/3}} \textnormal{ for all } k=1,\dots,N\}.
\end{equation*}
When these conditions for $f_1,\dots,f_J$ hold, we say these functions satisfy ``the rigidity conditions".
\end{lemma}
\begin{proof}
For convenience let us write $\bold{f}$ for $\sum_{j}f_i(H_{t_j})$ and fix a $D>0$. By H\"older and Jensen inequalities,
\begin{align*}
\Prob_{\bold f} (\boldsymbol \lambda(t)\in\tilde{\mathcal{G}}^c)&= \mathbb{E}\Big[\mathds{1}_{\boldsymbol \lambda(t)\in\tilde{\mathcal{G}}^c}\frac{e^{\sum_{j=1}^J S_N(f_j(H_{t_j}))}}{\E{e^{\sum_{j=1}^J S_N(f_j(H_{t_j}))}}}\Big]
\\
&\leq\Prob( \boldsymbol \lambda(t)\in\tilde{\mathcal{G}}^c)^{1/2} \prod_{j=1}^J \Big(\mathbb{E}[e^{2JS_N(f_j)}]^{1/2J}  e^{-\mathbb{E}[S_N(f_j)]}\Big)\leq \Prob_{\GUE}(\boldsymbol \lambda\in\tilde{\mathcal{G}}^c)^{1/2} e^{(\log N)^{C+1}}.
\end{align*}
On the other hand, by \cite[Lemma 3.8]{bourgade2022optimal}, for all sufficiently large $N$, $\Prob_{\GUE}(\boldsymbol \lambda\in\tilde{\mathcal{G}}^c)\leq  e^{-(\log N)^{C+D+2}}$ which completes the proof.
\end{proof}

Yet both conditions are very simple, verifying them is not most of the time. This section will provide some easy-to-verify sufficient conditions on the bias to hold the rigidity conditions. The techniques in this section are relatively standard in random matrix theory literature (e.g. \cite{erdos2010universality,erdHos2012bulk,bourgade2014universality,bourgade2016fixed}) and they rely on two main ingredients: local law and Helffer-Sj\"ostrand formula. We briefly review these below.

Denote the Cauchy-Stieltjes transform of the semicircle distribution and the normalized empirical spectral distribution by $m_{\rm sc}(z)=\int_{-2}^{2}\frac{1}{s-z}\rhosc (s)\D s=\frac{\sqrt{z^2-4}-z}{2}$ (where the branch cut $[-2,2]$ is chosen for the square root so that $\Im m_{\rm sc}(z) >0$ when $\Im z>0$) and $m(z)$ which is simply the normalized trace of the resolvent, i.e. $m(z)=\frac{1}{N}\Tr\frac{1}{H-z}=\frac{1}{N}\sum\frac{1}{\lambda_k-z}$ respectively. Although we have the explicit form for $m_{\rm sc}$, the following well-known bounds will be helpful in computations:
\begin{align}\label{eqn:basic_properties_of_m_sc}
\Im m_{\rm sc}(z)\asymp \begin{cases} \sqrt{\upkappa}+\sqrt{\eta}, &\textnormal{if }E\in[-2,2] \\
\frac{\eta}{\sqrt{\upkappa}+\sqrt{\eta}}, &\textnormal{if }E\in[-10,10]\setminus[-2,2]
\end{cases}
\end{align}
uniformly in $\eta\in(0,10]$ where $\eta=\Im z$, $E=\Re z$, $\upkappa:=|E-2|\wedge|E+2|$. The following optimal version of local law for $\GUE$ --which is a combination of Theorem 1.1, Remark 1.3, Proposition 2.5 in \cite{bourgade2022optimal}-- will be used throughout this section.
\begin{theorem}[Local law, \cite{bourgade2022optimal}] \label{thm:local_law} There exist universal constants $C,\tilde{\eta}>0$ such that if we denote the trapezoidal region $\{z=E+i\eta:\eta\in[0,\tilde\eta],-2-\eta\leq E\leq 2+\eta\}$ by $T$, then for every $N\in\N$ and $k\in \mathbb{N}$,
\begin{gather}\label{eqn:local_law}
\E{|m(z)-m_{\rm sc}(z)|^k}^{1/k} \leq  C\frac{k^{1/2}}{N\eta},\quad \textnormal{for all } z\in T,
\\
\E{|\Im(m(z)-m_{\rm sc}(z))|^k}^{1/k} \leq   C\Big(\frac{k^{1/2}}{N\eta}+\frac{k}{N\upkappa^{1/2}}\Big)\leq C\frac{k}{N\eta}, \quad \textnormal{for all } z\in(\R\times[0,\tilde\eta])\setminus T \label{eqn:local_law_outside}
\end{gather}
where $\upkappa:=|E-2|\wedge|E+2|$.
\end{theorem}
Given this theorem on Stieltjes transforms, the Helffer-Sjöstrand formula serves as a powerful tool for deriving estimates on the linear statistics for a wide class of functions.
\begin{prop} \label{prop:HS}\textnormal{(Helffer-Sj\"ostrand formula, \cite[Appendix C]{benaych2016lectures})} Given $\varepsilon>0$, $n\in\N$ and $f\in\mathscr{C}^{n}_{c}(\R)$, define $n^{th}$ order quasi-analytic (aka almost-analytic) extension of $f$ by $\tilde{f}(x)=(\sum_{k=0}^n \frac{(iy)^k}{k!} f^{(k)}(x))\chi_{\varepsilon}(y)$. Then, for every $\lambda\in\R$
\begin{align*} 
f(\lambda)=\frac{1}{\pi}\int_{\C} \frac{\partialbar\tilde{f}(z)}{\lambda-z}\D z
\end{align*}
where the anti-holomorphic derivative is given by $\partialbar=\frac{1}{2}(\partial_x+i\partial_y)$ and $\D z$ stands for the Lebesgue measure on $\C$.
\end{prop}
So, the centered linear statistics can be expressed as
\begin{align*} 
S_N(f)=\frac{N}{\pi}\int_{\C} \partialbar\tilde{f}(z) \big(m(z)-m_{\rm sc}(z)\big)\D z.
\end{align*}
Using this formulation, the propositions below provide aforementioned general conditions for the rigidity conditions. Some of them can be relaxed easily, but for our purposes this will be sufficient. 

\begin{prop}\label{rig_under_biased_prop1} There exists an $N_0$ such that for all sequence of $\mathscr{C}^2$ functions $f=(f_N)_N$ supported in $[-2,2]$ satisfying $\|f\|_{L^1}\leq (\log N)^{10}$, $\|f'\|_{L^1}\leq (\log N)^{10}$, and $\|f''\|_{L^1} \leq N^{(\log N)^{10}}$, we have
\begin{align}\label{boundfor_laplace_transform_of_lin_statistics}
\mathbb{E}[e^{ S_N(f)}] \leq e^{(\log N)^{50}},\quad \mathbb{E}[|S_N(f)|]\leq (\log N)^{25}
\end{align}
for every $N\geq N_0$.
\end{prop}

\begin{proof} Let $\varepsilon=\frac{1}{(\log N)^{10}}$, $\theta=N^{-(\log N)^{10}}$. Applying Helffer-Sj\"ostrand formula to the second order quasi-analytic extension of $f$ with bump function $\chi_{\varepsilon}$ gives
\begin{align*}
S_N(f)=\frac{N}{2\pi}\int\int_{\R^2}\Big(iyf''(x)\chi_{\varepsilon}(y)+i(f(x)+iyf'(x))\chi_{\varepsilon}'(y)\Big) \big(m(z)-m_{\rm sc}(z)\big)\D x\D y.
\end{align*}
Applying integration by parts a few times and using the fact that $f$ is a real-valued function (see \cite[(B.13)-(B.17)]{erdos2010universality}), we obtain
\begin{alignat*}{2}
\frac{1}{N}|S_N(f)| \lesssim &\ \Big|\int_{\R}\int_{y>0} f(x)(2\chi_{\varepsilon}'(y)+y\chi_{\varepsilon}''(y)) \Im\big(m(z)-m_{\rm sc}(z)\big)\D y\D x\Big| && \qquad (\RN1)
\\
&+\Big|\int_{\R}\int_{0<y\leq \theta}f''(x)y\chi_{\varepsilon}(y)\Im\big(m(z)-m_{\rm sc}(z)\big)\D y \D x\Big| && \qquad (\RN2)
\\
&+\Big|\int_{\R} f'(x)\theta\chi_{\varepsilon}(\theta)\Im\big(m(x+i\theta)-m_{\rm sc}(x+i\theta)\big)\D x\Big|&& \qquad (\RN3)
\\
&+\Big|\int_{\R}\int_{y> \theta} f'(x)(\chi_{\varepsilon}(y)+y\chi_{\varepsilon}'(y))\Im\big(m(z)-m_{\rm sc}(z)\big)\D y \D x\Big|.  &&\qquad (\RN4)
\end{alignat*}
We now bound the moments of $|S_N(f)|$ by controlling the moments of $(\RN1)-(\RN4)$ via the local law (Theorem \ref{thm:local_law}). Applying H\"older’s inequality, as in the proof of Lemma 4.3 in \cite{bourgade2022optimal}, we get:
\begin{align*}
\E{|(\RN1)|^k} &\leq \int_{(\R\times[\varepsilon,2\varepsilon])^k}\Ex\Big[\prod_{i=1}^{k}|f(x_i)(2\chi_{\varepsilon}'(y_i)+y_i\chi_{\varepsilon}''(y_i))\Im(m(z_i)-m_{\rm sc}(z_i))|\Big] \D x_1\dots\D x_k\D y_1 \dots \D y_k\nonumber
\\
&\leq \Big(\int\int_{\R\times[\varepsilon,2\varepsilon]}|f(x)||2\chi_{\varepsilon}'(y)+y\chi_{\varepsilon}''(y)|\cdot\E{|\Im(m(z)-m_{\rm sc}(z))|^k}^{1/k} \D x\D y\Big)^k\nonumber
\\
&\leq c^k \Big(\int\int_{\R\times[\varepsilon,2\varepsilon]}|f(x)||2\chi_{\varepsilon}'(y)+y\chi_{\varepsilon}''(y)|\cdot \frac{k^{1/2}}{Ny} \D x\D y\Big)^k
\leq  \frac{c^k}{N^k}\frac{k^{k/2}}{\varepsilon^k}(\log N)^{10}
\end{align*}
for some absolute constant $c$. The bounds for the moments of the other terms follows similarly:
\begin{gather*}
\E{|(\RN2)|^k} \leq   c^k \Big(\int\int_{\R\times[0,\theta]} |f''(x)|y\chi_{\varepsilon}(y)\cdot \frac{k^{1/2}}{Ny}\D x\D y\Big)^k \leq \frac{c^k}{N^k}k^{k/2}
\\
\E{|(\RN3)|^k} \leq   c^k \Big(\int_{\R} |f'(x)|\theta\cdot \frac{k^{1/2}}{N\theta} \D x\Big)^k  \leq \frac{c^k}{N^k} k^{k/2}(\log N)^{10k}
\\
\E{|(\RN4)|^k} \leq   c^k\Big(\int_\R\int_{y> \theta} |f'(x)||(\chi_{\varepsilon}(y)+y\chi_{\varepsilon}'(y))|\cdot \frac{k^{1/2}}{Ny}\D y\D x\Big)^k
\leq \frac{c^k}{N^k}k^{k/2}(\log N)^{21k}
\end{gather*}
Plugging the bounds we obtained for $(\RN1)-(\RN4)$ yields that for every $k\in\N$:
\begin{align*}
\E{|S_N(f)|^k} \leq c^k k^{k/2} (\log N)^{21k}
\end{align*}
which proves the second inequality in \eqref{boundfor_laplace_transform_of_lin_statistics}. Moreover, similar to the last step of the proof of Lemma 4.3 in \cite{bourgade2022optimal}, the first inequality in \eqref{boundfor_laplace_transform_of_lin_statistics} follows by Taylor expanding the Laplace transform, and substituting Stirling's approximation,
\begin{align}\label{taylorexpansionnew}
|\log \E{e^{ S_N(f)}}|\leq \log \E{e^{|\gamma S_N(f)|}} \leq  \log \sum_{k=0}^\infty \frac{1}{k!} \E{|S_N(f)|^k}\leq \log \Big(1+\sum_{k=1}^\infty c^k \frac{k^{k/2}}{k!} (\log N)^{21k}\Big)\leq (\log N)^{50}
\end{align}
which holds for all sufficiently large $N$.
\end{proof}

The previous proposition established the rigidity conditions for functions supported within the typical range of the spectrum. Next, we extend this result to functions with domains extended beyond the typical spectrum with the help of a classical concentration inequality for the linear statistics that can be derived by Bakry-Emery Theorem and Herbst Lemma with the help of Hoffman-Wielandt inequality (see \cite[Section 2.3.2]{anderson2010introduction} for details).

\begin{prop}\label{prop:rig_for_mathscr}
For any fixed $\Upsilon>0$, $C>0$, there exists an $N_0$ such that for all $\kappa\in[-C,1]$, $f\in\mathscr{S}_{C,\kappa}$ (recall \eqref{eqn:defn_mathscr_S}) $\gamma\in[-C,C]$, and $N\geq N_0$,
\begin{align*}
\Ex[|S_N(f)|]\leq (\log N)^{25},\quad \E{e^{\gamma S_N(f)}}\leq e^{(\log N)^{50}}.
\end{align*}
\end{prop}
\begin{proof}
By the definition of $\mathscr{S}_{C,\kappa}$, we know that $f$ can be decomposed as $f=h+g$ where $h$ is a sum of at most $\lfloor C\log N\rfloor$ many functions in $\bigcup_{\varepsilon\in[N^{-1+\kappa},\frac{\Upsilon}{2C}]}\mathscr{A}_{C,\varepsilon}$ and $g$ is a sum of at most $\lfloor C\log N\rfloor$ many order $1$ functions. By the choice of the decomposition, $h$ is supported in $[-2+\frac{\Upsilon}{2},2-\frac{\Upsilon}{2}]$, so by Proposition \ref{rig_under_biased_prop1}, $\Ex[|S_N(h)|]\leq (\log N)^{25}/4$ and $\E{e^{4\gamma S_N(h)}}\leq e^{(\log N)^{50}}$.

 On the other hand, for the function $g$, we know that there is a constant $c$ depending only on $C$ and $\Upsilon$ such that $\|g^{(k)}\|_{\infty}\leq c\log N$ for all $k=0,1,2,3,4$. Define two bump functions 
\begin{align*} 
 \chi_{in}(x)=\begin{cases}1 &,x\in[-2+2(\log N)^{-4},2-2(\log N)^{-4}] \\ 0&,x\in [-2+(\log N)^{-4}/2,2-(\log N)^{-4}]^c \end{cases},\quad \chi_{out}(x)=\begin{cases}1 &,x\in[-2-(\log N)^{-4},2+(\log N)^{-4}] \\ 0&,x\in [-2-2(\log N)^{-4},2+2(\log N)^{-4}]^c \end{cases}
 \end{align*} 
 with smooth transitions in between the given domains. By Proposition \ref{rig_under_biased_prop1} again, $\Ex[|S_N(g\cdot\chi_{in})|]\leq (\log N)^{25}/4$ and $\E{e^{4\gamma S_N(g\cdot\chi_{in})}}\leq e^{(\log N)^{50}}$. Moreover, by rigidity (e.g \cite[Lemma 3.8]{bourgade2022optimal}), it is easy to see that $ \Ex[|S_N(g\cdot(1-\chi_{out}))|]=o(1)$. Combining this with the Herbst lemma \cite[Lemma 2.3.3]{anderson2010introduction}, we also get $\E{e^{4\gamma S_N(g(1-\chi_{out}))}}\leq e^{O((\log N)^{10})}$ because $\|g(1-\chi_{out})\|_{Lip}=O((\log N)^5)$.
 
  Hence, we are only left with proving the linear statistics and Laplace estimates for the function $r(x):=g(x)(\chi_{out}(x)-\chi_{in}(x))$. Following the same steps with the proof of Proposition \ref{rig_under_biased_prop1} (with the same choice of $\varepsilon$ and $\theta$), and using \eqref{eqn:local_law_outside} instead of \eqref{eqn:local_law} for the local law estimation we get
\begin{align*} 
|S_N(r)|^k\leq c^kN^{k}((\RN1)^k+(\RN2)^k+(\RN3)^k+(\RN4)^k)
\end{align*}
where 
\begin{gather*}
\E{|(\RN1)|^k} \leq c^k \Big(\int\int_{\R\times[\varepsilon,2\varepsilon]}|r(x)||2\chi_{\varepsilon}'(y)+y\chi_{\varepsilon}''(y)|\cdot \big(\frac{k^{1/2}}{Ny}+\frac{k}{N\upkappa^{1/2}}\big) \D x\D y\Big)^k
\leq  \frac{c^k}{N^k}\frac{k^{k/2}}{\varepsilon^k(\log N)^{3k}}+\frac{c^k}{N^k}\frac{k^k}{(\log N)^{k}}
\\
\E{|(\RN2)|^k} \leq   c^k \Big(\int\int_{\R\times[0,\theta]} |r''(x)|y\chi_{\varepsilon}(y)\cdot \big(\frac{k^{1/2}}{Ny}+\frac{k}{N\upkappa^{1/2}}\big)\D x\D y\Big)^k \leq \frac{c^k}{N^k}k^{k/2}+\frac{c^k}{N^k}\frac{k^k}{(\log N)^{k}}
\\
\E{|(\RN3)|^k} \leq   c^k \Big(\int_{\R} |r'(x)|\theta\cdot \big(\frac{k^{1/2}}{N\theta}+\frac{k}{N\upkappa^{1/2}}\big) \D x\Big)^k  \leq \frac{c^k}{N^k} k^{k/2}(\log N)^{k}+\frac{c^k}{N^k}\frac{k^k}{(\log N)^{k}}
\\
\E{|(\RN4)|^k} \leq   c^k\Big(\int_\R\int_{y> \theta} |r'(x)||(\chi_{\varepsilon}(y)+y\chi_{\varepsilon}'(y))|\cdot \big(\frac{k^{1/2}}{Ny}+\frac{k}{N\upkappa^{1/2}}\big)\D y\D x\Big)^k
\leq \frac{c^k}{N^k}k^{k/2}(\log N)^{12k}+\frac{c^k}{N^k}\frac{k^k}{(\log N)^{k}}
\end{gather*}
Substituting these bounds into equation \eqref{taylorexpansionnew}, we obtain the necessary estimate.
\end{proof}

Note that $\mathscr{S}_{C,\kappa}$ contains only compactly supported functions, but it is very standard to extend the result to functions with suitable growth conditions. Below, we provide a simple argument from log-Sobolev inequality for such an extension; even though it is far from optimal, it is sufficient for our purposes.

\begin{prop}\label{rig_under_biased_prop2}
For any $C>0$, there exists a constant $c>0$ such that for all sequence of $\mathscr{C}^2$ functions $f=(f_N)_N:\R\to\R$ satisfying $\|f'\|_{\infty}+\|f''\|_{\infty}\leq C$, and any $N\in\N$,
\begin{align}\label{eqn:order1_log_sob}
\Ex[|S_N(f)|]\leq c,\quad \Ex[e^{\gamma S_N(f)}] \leq c.
\end{align}
\end{prop}
\begin{proof}
As $S_N(f)$ is invariant under constant shifts of $f$, without loss of generality we can assume that $f(0)=0$. Decompose $f$ as $f=f\cdot\chi_{3}+f\cdot(1-\chi_{3})$. By the tail bound on the extreme eigenvalue (e.g. \cite{aubrun2005sharp}) and the given growth condition on $f$, it is easy to see that $\E{|\Tr(f\cdot(1-\chi_{3}))|}=o(1)$. So, by the Herbst lemma \cite[Lemma 2.3.3]{anderson2010introduction} we get $\E{e^{\Tr(f\cdot(1-\chi_{3}))}}\leq e^{C^2/2}e^{ \E{S_N(f\cdot(1-\chi_{3}))}}= O(1)$. On the other hand, repeating the proof of the previous proposition by setting $\varepsilon$ and $\theta$ to some order $1$ constants we obtain $\mathbb{E}_{\GUE}[|S_N(f\cdot\chi_{3})|]=O(1)$. So, applying the Herbst lemma again, we get $\E{e^{S_N(f\cdot\chi_{3})}}=O(1)$. By Cauchy-Schwarz inequality, we obtain \eqref{eqn:order1_log_sob}.
\end{proof}

\begin{lemma}\label{lemma:rigidity_lemma_for_log_and_regularizations}
For any fixed $\Upsilon>0$ and $C>0$, there exists an $N_0$ such that for all $N\geq N_0$, $\iota\geq 0$, $E\in[-2+\Upsilon,2-\Upsilon]$, and $\gamma\in[0,C]$,
\begin{align*}
\Ex[|S_N(\log_\iota^E)|]\leq (\log N)^{25} \textnormal{ and } \E{e^{ S_N(\gamma\log_\iota^E)}}\leq e^{(\log N)^{50}}.
\end{align*}
Moreover, for all $d\in\{r,\ell\}$, $\beta\in[-C,C]$,
\begin{align*}
\Ex[|S_N(\arg_{d,\iota}^E)|]\leq (\log N)^{25} \textnormal{ and } \E{e^{ S_N(\beta\arg_{d,\iota}^E)}}\leq e^{(\log N)^{50}}.
\end{align*}
\end{lemma}
\begin{proof} $\iota\geq 1/N$ case is clear by Propositions \ref{prop:rig_for_mathscr} and \ref{rig_under_biased_prop2}. If $\iota\leq 1/N$, since $\log_\iota^E\leq\log^{E}_{1/N}$,
\begin{align*}
\E{e^{\gamma S_N(\log_\iota^E)}}\leq & \E{e^{\gamma S_N(\log_{1/N}^E)}} e^{C N\int (\log_{1/N}^E(s)-\log_\iota^{E}(s))\rhosc (s)\D s}\lesssim  \E{e^{\gamma S_N(\log_{1/N}^E)}}\leq e^{(\log N)^{50}}.
\end{align*}
Moreover,
\begin{align*}
\E{|\Tr(\log^{E}_{\iota})-\Tr(\log^{E}_{1/N})|} &\leq \sum_{k=1}^{N}\Ex\Big[ \mathds{1}_{|\lambda_k-E|\leq 2/N}\cdot |\log^E\lambda_k-\log^E_{1/N}\lambda_k|\Big] 
\leq \sum_{j=1}^{\infty} j\log 2 \cdot \Ex\Big[ \sum_{k=1}^{N} \mathds{1}_{\{\frac{2/N}{2^{j}}\leq |\lambda_k-E|\leq \frac{2/N}{2^{j-1}}\}}\Big]
\\
&\lesssim \sum_{j=1}^{\infty} j \cdot \int_{\frac{2/N}{2^{j}}\leq |x-E|\leq \frac{2/N}{2^{j-1}}} \rho_{1}^{\GUE}(x)\D x \lesssim \sum_{j=1}^{\infty}\frac{j}{2^{j}}=O(1)
\end{align*}
where we have used the $O(N)$ bound on the $1$-point correlation function. Together with $|\int (\log_{\iota}^E-\log_{1/N}^{E})\rhosc |=O(1/N)$ and $\E{|S_N(\log^{E}_{1/N})|}=(\log N)^{25}$, this finishes the proof. For $\arg$ function, the inequalities follow from the same steps.
\end{proof}

\begin{remark}
Studying a problem under a biased measure can be more challenging than analyzing the same problem under the standard random matrix distribution. Nevertheless,  the following basic estimate on linear statistics remain available for biased measures,  from the preservation of rigidity (e.g. \cite[Proposition 1.3]{lambert2021mesoscopic}): If we define the centered empirical measure by $\D\nu_N( x):=\sum \delta_{\lambda_i}(x)-N\rhosc (x)\D x$ and counting function $g_N(x):=\sum_{k=1}^{N} \mathds{1}_{\lambda_k\in(-\infty,x)}-N\int_{\infty}^{x}\rhosc (s)\D s$; for any $f\in \mathscr{C}^1$ (or absolutely continuous), an integration by parts yields,  on the rigidity set,
\begin{align}\label{eqn:int_on_centered_emp_meas_int_by_partss}
|S_N(f)|=\big|\int_{\R} f(x)\D\nu_N(x)\big|\leq \|g_N\|_{\infty}\int_{\R}  \big|f'(x)\big| \D x
\end{align}
When the rigidity set has an overwhelming probability, the contribution from outside of this set is typically negligible by crude estimates, e.g.  $L^\infty$.
\end{remark}

\section{Resolvent estimates}\label{sec:res_est}

In this section we will first prove resolvent estimates which will provide us with the long range stability of the eigenvalue paths (up to time $(\log N)^2$) and then the necessary tools --specifically, a multi-time local law, Proposition \ref{prop:full_rank_estimate}-- for developing multi-time loop equation asymptotics. The structure of this section closely follows \cite[Section 4]{bourgade2022liouville}, with some important modifications to the argument in the second half.  As our primary focus is the case $\beta=2$ with quadratic potential $V(x)=\frac{x^2}{2}$, we present the calculations only for that setting, but the techniques should extend to Dyson Brownian motion with general $\beta$ and potential $V$ (e.g. see \cite[Section 4]{adhikari2020dyson}).

Given the matrix-valued dynamics \eqref{matrix_valued_OU_dynamics},
for any $N\times N$ matrix $A$ we define
\begin{equation*}
m_t(z):=\frac{1}{N}\Tr(\frac{1}{H_t-z}), \ m_{t,A}(z):=\Tr(\frac{1}{H_t-z}A).
\end{equation*}

\begin{lemma}
$m_{t,A}(z)$ satisfies the following SDE,
\begin{equation}\label{SDE for m_t,A}
\D m_{t,A}(z)=\frac{-1}{\sqrt{N}}\Tr(\frac{1}{H_t-z}A\frac{1}{H_t-z}\D B_t)+\Big(\frac{1}{2}m_{t,A}(z)+\frac{1}{2}z\cdot \partial_z m_{t,A}(z)+m_t(z)\cdot \partial_zm_{t,A}(z)\Big)\D t.
\end{equation}
\end{lemma}

\begin{proof}
Let $e_{ij}$ be the $N \times N$ matrix with $(i,j)^{th}$ entry 1, and the rest zero. Define $X_{ij}^R:=e_{ij}+e_{ji}$, $X_{ij}^I:=i\cdot e_{ij}-i\cdot e_{ji}$ and $X_{ii}:=e_{ii}$ for any $i<j$. So, $B_t=\sum_{i<j}\Re(B_{ij})_tX_{ij}^R+\sum_{i<j}\Im(B_{ij})_tX_{ij}^I+\sum_i (B_{ii})_tX_{ii}$. Applying It\^o's formula we get,
\begin{align*}
\D \frac{1}{H_t-z}=&-\sum_{i<j} \frac{1}{H_t-z}X_{ij}^R\frac{1}{H_t-z}(\frac{1}{\sqrt N}\D \Re(B_{ij})_t-\frac{1}{2}\Re(H_{ij})_t\D t)-\sum_{i<j} \frac{1}{H_t-z}X_{ij}^I\frac{1}{H_t-z}(\frac{1}{\sqrt N}\D \Im(B_{ij})_t-\frac{1}{2}\Im(H_{ij})_t \D t)\\
&-\sum_{i} \frac{1}{H_t-z}X_{ii}\frac{1}{H_t-z}(\frac{1}{\sqrt N}\D (B_{ii})_t-\frac{1}{2}(H_{ii})_t\D t)\\
&+\sum_{i<j}\frac{1}{H_t-z}X_{ij}^R\frac{1}{H_t-z}X_{ij}^R\frac{1}{H_t-z}\frac{1}{2N}\D t+\sum_{i<j}\frac{1}{H_t-z}X_{ij}^I\frac{1}{H_t-z}X_{ij}^I\frac{1}{H_t-z}\frac{1}{2N}\D t\\
&+\sum_{i}\frac{1}{H_t-z}X_{ii}\frac{1}{H_t-z}X_{ii}\frac{1}{H_t-z}\frac{1}{N}\D t\\
=&\frac{-1}{\sqrt N}\frac{1}{H_t-z} \D B_t \frac{1}{H_t-z}+\frac{1}{2}\frac{1}{H_t-z}H_t\frac{1}{H_t-z}\D t+\frac{1}{N} \frac{1}{H_t-z}\frac{1}{H_t-z}\Tr(\frac{1}{H_t-z})\D t
\end{align*}
where at the second equality we have used $\frac{1}{2}\sum_{i<j}X_{ij}^R M X_{ij}^R+\frac{1}{2}\sum_{i<j}X_{ij}^I M X_{ij}^I +\sum_{i}X_{ii} M X_{ii}=\Id \cdot \Tr(M)$ for any matrix $M$. Therefore, using the invariance of trace under cyclic permutations we get
\begin{align*}
\D \Big(\Tr(\frac{1}{H_t-z}A)\Big)=\frac{-1}{\sqrt{N}}\Tr(\frac{1}{H_t-z}A\frac{1}{H_t-z}\D B_t)+\frac{1}{2}\Tr(\frac{1}{H_t-z}A\frac{1}{H_t-z}H_t)\D t+\frac{1}{N}\Tr(\frac{1}{H_t-z})\Tr(\frac{1}{H_t-z}A\frac{1}{H_t-z})\D t.
\end{align*}
Thus, substituting the equality
\begin{align*}
\Tr(\frac{1}{H_t-z}A\frac{1}{H_t-z}H_t)=\Tr\big(\frac{1}{H_t-z}A(\Id+\frac{z}{H_t-z})\big)=m_{t,A}(z)+z\cdot \partial_z m_{t,A}(z)
\end{align*}
we obtain equation \eqref{SDE for m_t,A}.
\end{proof}

If we can show that the stochastic term's contribution is negligible with overwhelming probability and $m_t(z)$ can be approximated by  $m_{\rm sc}(z)= \frac{-z+\sqrt{z^2-4}}{2}$, then we would expect
\begin{align*}
\frac{\partial m_{t,A}(z)}{\partial t}\approx  \frac{m_{t,A}(z)}{2} + \frac{\sqrt{z^2-4}}{2} \partial_z m_{t,A}(z).
\end{align*}
The characteristics of this PDE starting from $(t,z,m_{t,A}(z))$ becomes $(t-s,z_s,e^{-s/2}m_{t,A}(z))$ where 
\begin{align}\label{eqn:characteristic_curve}
z_s:=\frac{e^{s/2}(z+\sqrt{z^2-4})+e^{-s/2}(z-\sqrt{z^2-4})}{2},
\end{align}
see e.g. \cite[Lemma 2.5]{bourgade2021extreme}.
So, we would expect $m_{0,A}(z_t)\approx e^{-t/2}m_{t,A}(z)$ up to some $t$. For convenience, let us define $\tilde{m}_{t,A}(z):=e^{-t/2}m_{t,A}(z)$ and similarly $\tilde m_t(z):=e^{-t/2}m_t(z)$. The SDE for $\tilde{m}_{t,A}(z)$ along the characteristic curve for any fixed $t$ and $z$ becomes
\begin{equation}\label{SDE along char for tilde m_t,A}
\D \tilde m_{u,A}(z_{t-u})=\frac{-e^{-u/2}}{\sqrt{N}}\Tr(\frac{1}{H_u-z_{t-u}}A\frac{1}{H_u-z_{t-u}}\D B_u)+\big(\tilde m_u(z_{t-u})-e^{-t/2}m_{\rm sc}(z)\big)\partial_z m_{u,A}(z_{t-u})\D u
\end{equation}
where we used $\frac{\partial z_{t-u}}{\partial u}=-\frac{1}{2}z_{t-u}-e^{-(t-u)/2}m_{\rm sc}(z)$.

Setting $A=\frac{1}{N}\Id$, the trace term in the SDE can be simplified using the invariance of GUE distribution under unitary conjugation
\begin{align*}
\Tr(\frac{1}{H_u-z_{t-u}}\frac{1}{N}\Id\frac{1}{H_u-z_{t-u}}\D B_u)&=\frac{1}{N}\Tr(P_u^*\frac{1}{H_u-z_{t-u}}P_u P_u^* \frac{1}{H_u-z_{t-u}} P_u P_u^*\D B_u P_u)
\\
&=\frac{1}{N}\Tr(\frac{1}{(D_u-z_{t-u})^2} P_u^*\D B_u P_u)=\frac{1}{N} \sum_k \frac{1}{(\lambda_k(u)-z_{t-u})^2}\D (\tilde B_k)_{u}
\end{align*}
where $P_uD_uP_u^*$ is a diagonalization of $H_u$ and $\tilde B_k$'s are independent real Brownian motions. Substituting this into \eqref{SDE along char for tilde m_t,A} gives the SDE for the rescaled Stieltjes transform of the normalized empirical spectral distribution along the characteristics
\begin{align}\label{sde for tilde m_u(z_(t-u))}
\D \tilde m_u(z_{t-u})=\frac{e^{-u/2}}{N^{3/2}} \sum_k \frac{1}{(\lambda_k(u)-z_{t-u})^2}\D (\tilde B_k)_{u}+\big(\tilde m_u(z_{t-u})-e^{-u/2}m_{\rm sc}(z_{t-u})\big)\partial_z m_{u}(z_{t-u})\D u
\end{align}
where we flipped the signs of the Brownian motions.

Before moving on to the technical part of the section, let us introduce a few extra notations/definitions that will be used frequently throughout this section. Denote the integral of the martingale term in the SDE \eqref{sde for tilde m_u(z_(t-u))} by 
\begin{align}\label{Definition of M_t in the Stieltjes transform SDE}
M_t:=\frac{1}{N^{3/2}}\sum_k\int_0^t \frac{e^{-u/2}}{(\lambda_k(u)-z_{t-u})^2}\D (\tilde{B}_k)_u
\end{align}
Define sup-microscopic level curve/function $\mathcal{C}_\alpha:\R\to\R$ by
\begin{align}\label{eqn:defn_of_mathcal_C_alpha_curve} 
\mathcal C_\alpha(x):=\begin{cases} \frac{\alpha}{(N\sqrt{\kappa(x)})\vee N^{2/3}}, &\textnormal{if } x\in[-2,2] \\
\frac{\alpha}{N^{2/3}}, &\textnormal{otherwise}
\end{cases}
\end{align}
for any $\alpha>0$. Lastly, for later convenience, we denote any region of the form $\{E+i\eta \ | \ E\in I, f(E)\leq \eta\leq g(E)\}$ by $I\times[f(E),g(E)]$. Given $z=E+i\eta$, the imaginary and real parts of the characteristic curve will be denoted by $\eta_t=\eta_{z_{t}}=\Im(z_{t})$, $E_t=E_{z_{t}}=\Re(z_{t})$, and $\min(|E-2|,|E+2|) $ will be denoted by $\upkappa=\upkappa(z)=\upkappa(E)$. The following simple identities for the characteristic curve, which are straightforward to verify, will be helpful throughout the discussion.

\begin{lemma}
For any $a>1$, uniformly over the choice of $0\leq s\leq t$ and uniformly for all $z$ in the upper half plane we have,
\begin{gather*} 
m_{\rm sc}(z_t)=e^{-t/2}m_{\rm sc}(z), \quad \eta_{t}\Im\sqrt{z_{t}^2-4}\geq e^t\eta\Im\sqrt{z^2-4}
\\
\eta_t\asymp \eta+(e^{t/2}-1)\Im\sqrt{z^2-4}, \quad \Im\sqrt{z_t^2-4}\asymp e^{t/2}\Im\sqrt{z^2-4}
\\
\int_{0}^{\infty}\frac{1}{\eta_t^a}\D t\asymp \frac{1}{\eta^{a-1}\Im\sqrt{z^2-4}},\quad \int_{0}^{s} \frac{e^{-u/2}}{\eta_{t-u}^a}\D u\lesssim \frac{e^{-t/2}}{\eta_{t-s}^{a-1}\Im\sqrt{z^2-4}}.
\end{gather*}
Moreover, if $|E|\leq 2$ and $z$ is above the curve $\mathcal{C}_{\alpha}$, then $\eta\sqrt{z^2-4}\gtrsim \frac{\alpha}{N}$.
\end{lemma}

This section consists of a sequence of lemmas and propositions leading to the estimates that will be helpful in Section \ref{sec:loop_eqn}. In the proofs of these results, we always assume that we are on the set of overwhelming probability on which the previous lemmas and propositions hold in addition to the rigidity at time $0$, i.e. $\boldsymbol\lambda(0)\in\tilde{\mathcal{G}}$ (see \eqref{eqn:rigidity sets}). Recall that, by \cite[Lemma 3.8]{bourgade2022optimal}, the rigidity holds with overwhelming probability at a single time.

\subsection{Long-range stability of eigenvalue rigidity}\label{subsec:long_stability}

We first prove that the eigenvalue rigidity is preserved over a long time period for the process $(H_t)$. Within this section, we will be interested in the following rigidity scalings/sets:
\begin{equation}\label{eqn:rigidity sets}
\begin{aligned} 
\mathcal{G}:=\Big\{\Big(\lambda_1,\dots,\lambda_N): |\gamma_k-\lambda_k|<\frac{\varphi^{10}}{N^{2/3}\hat{k}^{1/3}}, \textnormal{ for all } k=1,\dots,N\},
\\
\mathcal{\tilde G}:=\Big\{\Big(\lambda_1,\dots,\lambda_N): |\gamma_k-\lambda_k|<\frac{\varphi}{N^{2/3}\hat{k}^{1/3}}, \textnormal{ for all } k=1,\dots,N\} 
\end{aligned}
\end{equation}
where recall that $\hat{k}:=\min(k,N+1-k)$ and $\varphi:=(\log N)^{\log\log N}$. By overwhelming probability, we mean events hold with probability at least $1-e^{-(\log N)^D}$ for any arbitrary $D>0$ for all sufficiently large $N$.

\begin{lemma}\label{first step of the proof}
For any $D>0$ there exists an $N_0$ such that for all $N\geq N_0$, $t\in[0,(\log N)^2]$, and $z\in [-2-\frac{\varphi^{6}}{N^{2/3}},2+\frac{\varphi^{6}}{N^{2/3}}]\times [\mathcal{C}_{\varphi^6}(E),1]$ we have
\begin{equation*}
\Prob\Big\{|m_{u}(z_{t-u})-m_{\rm sc}(z_{t-u})|\leq \frac{\varphi^{3}}{N\eta_{t-u}}, \textnormal{ for all } u\in[0,t]\Big\}\geq 1-e^{-(\log N)^D}.
\end{equation*}
\end{lemma}

\begin{proof}
Let $t>1$ and $z$ in the given region; for the case $t\leq 1$, the same method works with simpler calculations. Define $h(u):=\tilde{m}_u(z_{t-u})-e^{-u/2}m_{\rm sc}(z_{t-u})$ and stopping time $\tau:=\inf\{u\in[0,t]:|h(u)|\geq e^{-u/2}\frac{\varphi^{3}}{N\eta_{t-u}}\}\wedge t$. Note that it suffices to prove $|h(\tau)|<e^{-\tau/2}\frac{\varphi^{3}}{N\eta_{t-\tau}}$ holds with overwhelming probability. Substituting the definition of $h$ and \eqref{Definition of M_t in the Stieltjes transform SDE} into \eqref{sde for tilde m_u(z_(t-u))} gives,
\begin{align}\label{h(s)=h(0)+M_s+int}
h(s)=h(0)+M_s+\int_0^s h(u)\partial_z m_{u}(z_{t-u})\D u.
\end{align}
On the rigidity set $\boldsymbol{\lambda}(0) \in \mathcal{\tilde G}$, by equation \eqref{eqn:int_on_centered_emp_meas_int_by_partss}, after an appropriate truncation of the function $\frac{1}{\cdot-z}$ outside the spectrum, it is easy to see that $|h(0)|\leq \frac{\varphi^{2}}{N\eta_t}$. The martingale term can be controlled by its quadratic variation. Denoting the stopped process by $(M_s^\tau)_s$, its quadratic variation can be bounded as
\begin{align*}
\quadvar{M^\tau}_s \leq& \frac{1}{N^3}\sum_k \int_0^{s\wedge\tau}\frac{e^{-u}}{|\lambda_k(u)-z_{t-u}|^4}\D u\leq \frac{1}{N^2}\int_0^{s\wedge\tau}\frac{e^{-u/2}}{\eta_{t-u}^3}\Im(\tilde m_u(z_{t-u}))\D u
\\
\leq& \frac{1}{N^2}\int_0^{s\wedge\tau}\frac{e^{-u/2}}{\eta_{t-u}^3}\Big(e^{-u/2}\Im(m_{\rm sc}(z_{t-u}))+e^{-u/2}\frac{\varphi^3}{N\eta_{t-u}}\Big)\D u \lesssim \frac{e^{-t/2}\Im \sqrt{z^2-4}}{N^2}\int_0^{s}\frac{e^{-u/2}}{\eta_{t-u}^3}\D u
\lesssim  \frac{e^{-t}}{N^2\eta_{t-s}^2}.
\end{align*}
where we have used the definition of stopping time in the third inequality. Hence, by Burkholder–Davis-Gundy (BDG) inequality, for any $s\in[0,t]$, $|M^\tau_s| \leq \frac{\varphi e^{-t/2}}{N\eta_{t-s}}$ with overwhelming probability. Extension to the uniformity in time can be done by classical grid argument (e.g. see \cite[Remark 2.7]{benaych2016lectures} or the proof of \cite[Proposition 4.2]{bourgade2022liouville}), so, for every $D>0$ there exists an $N_0=N_0(D)$  such that for all $N\geq N_0$:
\begin{align*}
\Prob\Big\{|M^\tau_s| \leq \varphi\frac{ e^{-t/2}}{N\eta_{t-s}}, \textnormal{ for all }s\in[0,t] \Big\}\geq 1-e^{-(\log N)^D}.
\end{align*}

Applying Gr\"onwall's lemma to equation \eqref{h(s)=h(0)+M_s+int} we obtain that on a set with overwhelming probability,
\begin{align}\label{result of gronwall}
|h(s\wedge \tau)| \leq \frac{\varphi^{2}}{N\eta_t}+\frac{\varphi e^{-t/2}}{N\eta_{t-s}}+\int_0^{s\wedge \tau} \Big(\frac{\varphi^{2}}{N\eta_t}+\frac{\varphi e^{-t/2}}{N\eta_{t-u}}\Big)|\partial_z m_{u}(z_{t-u})|e^{\int_u^{s\wedge \tau}|\partial_z m_{r}(z_{t-r})|\D r}\D u
\end{align}
for all $s\in[0,t]$. The integrand of the exponential term of the right hand side is bounded by
\begin{align}\label{partial_z m_r(z_(t-r)) precise bound}
|\partial_z m_{r}(z_{t-r})|\leq \frac{e^{r/2}}{\eta_{t-r}}\Im(\tilde m_{r}(z_{t-r})) \leq \frac{e^{r/2}}{\eta_{t-r}}\Big(e^{-r/2}\Im(m_{\rm sc}(z_{t-r}))+e^{-r/2}\frac{\varphi^3}{N\eta_{t-r}}\Big) \leq \frac{e^{-(t-r)/2}\Im\sqrt{z^2-4}}{2\eta_{t-r}}(1+\varphi^{-1})
\end{align}
for every $r\leq s\wedge \tau$. Thanks to this, we obtain the following estimates for the exponential term in Grönwall’s inequality:
\begin{align}\label{exponentail term obtained bounds eqn}
e^{\int_u^{s\wedge \tau}|\partial_z m_{r}(z_{t-r})|\D r}\lesssim \begin{cases}
1& , \textnormal{if} \ s\leq t-1\\
\frac{\Im\sqrt{z^2-4}}{\eta_{t-s}}&,  \textnormal{if} \ u\leq t-1\leq s \leq t \\
\frac{\eta_{t-u}}{\eta_{t-s}}&,  \textnormal{if} \  t-1 \leq u\leq s \leq t
\end{cases}
\end{align}
where these estimates follow from calculations similar to the one below, which correspond to the case $u\leq t-1\leq s$:
\begin{align*}
\int_u^{s\wedge \tau}|\partial_z m_{r}(z_{t-r})|\D r \leq & \int_{u}^{(t-1)}2e^{-(t-r)/2}\D r+\Big(1+\varphi^{-1}\Big) \int_{t-1}^{s}\frac{1}{\frac{2\eta}{\Im(\sqrt{z^2-4})}+(t-r)}\D r
\\
\leq &\ 4+ \big(1+\varphi^{-1}\big) \log \Big(\frac{\eta+\frac{1}{2}\Im\sqrt{z^2-4}}{\eta+\frac{t-s}{2}\Im\sqrt{z^2-4}}\Big).
\end{align*}

Substituting \eqref{partial_z m_r(z_(t-r)) precise bound} and \eqref{exponentail term obtained bounds eqn} into \eqref{result of gronwall} completes the proof via the stopping time argument. The steps for the case $s>t-1$ are outlined below, while the calculations for $s\leq t-1$ are similar and more straightforward,
\begin{align*}
\int_0^{s\wedge \tau} \Big(\frac{1}{\eta_t}+\frac{e^{-t/2}}{\eta_{t-u}}\Big)|\partial_z m_{u}(z_{t-u})|&e^{\int_u^{s\wedge \tau}|\partial_z m_{r}(z_{t-r})|\D r}\D u
\\
&\lesssim \int_{0}^{t-1}\frac{1}{\eta_t}e^{-(t-u)}\frac{\Im\sqrt{z^2-4}}{\eta_{t-s}}\D u+ \int_{t-1}^{s}\frac{e^{-t/2}}{\eta_{t-u}}\frac{\Im\sqrt{z^2-4}}{\eta_{t-u}} \frac{\eta_{t-u}}{\eta_{t-s}}\D u \lesssim \frac{e^{-t/2}}{\eta_{t-s}}\log N.
\end{align*}
Consequently, substituting into \eqref{result of gronwall} gives that on a set with overwhelming probability, $|h(s\wedge \tau)| \lesssim \frac{\varphi^{2}}{N\eta_t}+\frac{\varphi e^{-t/2}}{N\eta_{t-s}}+\frac{\varphi^2 e^{-t/2}}{N\eta_{t-s}}\log N \ll e^{-s/2}\frac{\varphi^{3}}{N\eta_{t-s}}$ for every $s\in[0,t]$. So, $|h(\tau)|<e^{-\tau/2}\frac{\varphi^3}{N\eta_{t-\tau}}$ with overwhelming probability.
\end{proof}

The next lemma follows from taking $u=t$ in the previous lemma and applying the classical grid argument to establish the uniformity in space and time (see the second step of the proof of \cite[Proposition 4.2]{bourgade2022liouville} for details). Since we chose the exponent of $\varphi$ generously in the $\frac{\varphi^3}{N\eta}$ upper bound in Lemma \ref{first step of the proof}, a slight adjustment suffices to establish this bound uniformly without difficulty. In the following sequence of lemmas and propositions, we omit such technical details for simplicity.

\begin{lemma}\label{uniformity lemma above the microscopic scale}
For any $D>0$ there exists an $N_0$ such that for all $N\geq N_0$:
\begin{align*}
\Prob\Big\{|m_t(z)-m_{\rm sc}(z)|\leq \frac{\varphi^{3}}{N\eta}, \textnormal{ for all } t\in[0,(\log N)^2] \textnormal{ and } z\in [-2-\frac{\varphi^{6}}{N^{2/3}},2+\frac{\varphi^{6}}{N^{2/3}}]\times [\mathcal{C}_{\varphi^{6}}(E),1] \Big\}\geq 1-e^{-(\log N)^D}.
\end{align*}
\end{lemma}

Next lemma shows that the previous estimates can be extended to submicroscopic scales.

\begin{lemma}\label{m_t(z) estimate on the rectangle}
For any $D>0$ there exists an $N_0$ such that for all $N\geq N_0$:
\begin{align*}
\Prob\Big\{|m_t(z)-m_{\rm sc}(z)|\leq \frac{\varphi^{9}}{N\eta}, \textnormal{ for all } t\in[0,(\log N)^2] \textnormal{ and } z\in [-2-\frac{\varphi^{6}}{N^{2/3}},2+\frac{\varphi^{6}}{N^{2/3}}]\times (0,1]\Big\}\geq 1-e^{-(\log N)^D}.
\end{align*}
\end{lemma}

\begin{proof}
 Let $z$ be an arbitrary point in $[-2-\frac{\varphi^{6}}{N^{2/3}},2+\frac{\varphi^{6}}{N^{2/3}}]\times(0,\mathcal{C}_{\varphi^6}(E)]$ and $z'$ be the point on the graph of $\mathcal{C}_{\varphi^{6}}$ for which $\Re z=\Re z'$. Denote the imaginary and real parts of $z'$ by $\eta'$ and $E'$. We bound the imaginary and real parts of $(m_t(z)-m_{\rm sc}(z))$ separately, beginning with the imaginary part. 
 
  It is easy to see that for any $z$ in that given region, $\Im m_{\rm sc}(z)\lesssim \frac{\varphi^9}{N\eta}$. So $\Im (m_t(z))$ is also
\[
\frac{1}{N\eta}\sum_k \frac{\eta^2}{(\lambda_k(t)-E)^2+\eta^2} \leq \frac{1}{N\eta}\sum_k \frac{(\eta')^2}{(\lambda_k(t)-E)^2+(\eta')^2}
=\frac{\eta'}{\eta}\Im(m_t(z'))\leq \frac{\eta'}{\eta} (\Im(m_{\rm sc}(z'))+\frac{\varphi^{3}}{N\eta'}) \lesssim \frac{\varphi^{9}}{N\eta}.
\]
Together with $\Im(m_{\rm sc}(z))\lesssim \frac{\varphi^9}{N\eta}$, this gives $|\Im (m_t(z))-\Im( m_{\rm sc}(z))|\lesssim \frac{\varphi^{9}}{N\eta}$.

Next, we bound the the difference of the real parts, starting with
\begin{align}\label{re part of the difference into 3 parts}
|\Re (m_t(z))-\Re( m_{\rm sc}(z))|\leq & |\Re (m_t(z))-\Re (m_t(z'))|+|\Re (m_t(z'))-\Re( m_{\rm sc}(z'))|+|\Re( m_{\rm sc}(z'))-\Re( m_{\rm sc}(z))|.
\end{align}
The main difficulty is controlling the first term of the right hand side and it can be bounded as follows,
\begin{align*}
 \frac{1}{N}&\sum_{k:|E-\lambda_k(t)|>10\eta'}\Big| \Re\Big( \frac{1}{\lambda_k(t)-z}-\frac{1}{\lambda_k(t)-z'} \Big) \Big|+\frac{1}{N}\sum_{k:|E-\lambda_k(t)|\leq 10\eta'}\Big| \Re\Big( \frac{1}{\lambda_k(t)-z}-\frac{1}{\lambda_k(t)-z'} \Big) \Big|
\\
&\lesssim  \frac{1}{N}\sum_{k:|E-\lambda_k(t)|>10\eta'}\Big| \Im\Big( \frac{\eta'-\eta}{((\lambda_k(t)-E)^2-\eta\eta')-i(\eta+\eta')(\lambda_k(t)-E)} \Big) \Big|+\frac{\eta'}{N}\sum_{k:|E-\lambda_k(t)|\leq 10\eta'}  \frac{1}{|\lambda_k(t)-z|\eta'}
\\
&\leq  \frac{1}{N}\sum_{k:|E-\lambda_k(t)|>10\eta'}  \frac{(\eta'-\eta)(\eta+\eta')|\lambda_k(t)-E|}{|\lambda_k(t)-z|^2|\lambda_k(t)-z'|^2} +\frac{1}{N\eta}\#\{k:|E-\lambda_k(t)|<10\eta'\}
\\
&\lesssim  \frac{1}{N}\sum_{k}  \frac{\eta'}{|\lambda_k(t)-E|^2} +\frac{1}{N\eta}\sum_k\frac{(\eta')^2}{(\eta')^2+(\lambda_k(t)-E)^2}
\leq  \Im(m_t(z'))+\frac{\eta'}{\eta} \Im(m_t(z'))\lesssim  \frac{\varphi^{9}}{N\eta}.
\end{align*}
The second term of the right hand side of \eqref{re part of the difference into 3 parts} is already bounded by $\frac{\varphi^{3}}{N\eta'}$ on the set of Lemma \ref{uniformity lemma above the microscopic scale}. Lastly, the third term in \eqref{re part of the difference into 3 parts} can be controlled easily as it is a deterministic expression:
\begin{align*}
|\Re( m_{\rm sc}(z'))-\Re( m_{\rm sc}(z))|\leq | m_{\rm sc}(z')- m_{\rm sc}(z)| \lesssim \frac{\eta'}{\sqrt{\upkappa}+\sqrt{\eta}}\lesssim \frac{\varphi^{9}}{N\eta}
\end{align*}
which completes the proof.
\end{proof}

To establish the stability of rigidity, we must first show the stability of the extreme eigenvalues, in other words, the absence of outliers among the particles up to time $(\log N)^2$ with overwhelming probability. The edge rigidity of Dyson Brownian motion has been previously studied in \cite{adhikari2020dyson} and here, we employ the same stopping time argument as used in the proof of Theorem 4.3 of that article.

\begin{lemma}\label{no outliers lemma}
For any $D>0$ there exists an $N_0$ such that for all $N\geq N_0$:
\begin{align*}
\Prob\Big\{\lambda_N(t)<2+\frac{\varphi^{6}}{N^{2/3}}, \textnormal{ for all } t\in[0,(\log N)^2]\Big\}\geq 1-e^{-(\log N)^D}.
\end{align*}
\end{lemma}

\begin{proof}
Fix $z=2+\frac{\varphi^{6}}{N^{2/3}}+i\frac{\varphi}{N^{2/3}}$ and sequence of times $t_j=\frac{j}{N^4}$ for $j\in\{1,2,\dots,\lceil N^4(\log N)^2\rceil\}$.  Define $h_j(u):=\tilde{m}_u(z_{t_j-u})-e^{-u/2}m_{\rm sc}(z_{t_j-u})$ as before and a new stopping time,
\begin{align*}
\tau=\bigwedge_{j=1}^{\lceil N^4(\log N)^2\rceil}\inf\Big\{u\in[0,t_j]: |h_j(u)|\geq e^{-u/2}\frac{1}{\varphi N\eta_{t_j-u}}\Big\}\wedge \inf\Big\{u\geq 0:\lambda_N(u)=2+\upkappa\Big\}\wedge (\log N)^2
\end{align*}
where recall that $\upkappa=\upkappa(z)=|E-2|\wedge|E+2|=\frac{\varphi^{6}}{N^{2/3}}$.

Fix a $j\in\{1,2,\dots,\lceil N^4(\log N)^2\rceil\}$. The equation \eqref{h(s)=h(0)+M_s+int} for $h_j$ will be evaluated by Gr\"onwall as before and the following simple estimations for the given particular $z$ will be used throughout the proof, \begin{align*}
\upkappa_t-\upkappa\asymp \begin{cases}
t^2+t\sqrt{\upkappa}, &t<1\\
e^{t/2}, &t\geq 1
\end{cases};
 \quad \eta_t\asymp  \begin{cases}
\eta+t\frac{\eta}{\sqrt{\upkappa}} , &t<1\\
e^{t/2}\frac{\eta}{\sqrt{\upkappa}}, &t\geq 1
\end{cases}; 
\quad \upkappa_t\gtrsim\varphi^{5}\eta_t, \quad \textnormal{for all } t\geq 0.
\end{align*}
By rigidity we have $|h_j(0)|\leq \frac{\varphi^{2}}{N(\eta_{t_j}+\upkappa_{t_j})}\leq \frac{1}{\varphi^3 N \eta_{t_j}}$. We also have the following bound uniformly over the time $s\in[0,t_j]$ for the martingale term (uniformity by BDG inequality and grid argument), for every $D>0$ there exists an $N_0$ such that for all $N\geq N_0$:
\begin{align*}
\Prob\Big\{\sup_{u\in[0,t_j]}|M^\tau_u| \leq \frac{1}{N\eta_{{t_j}}\varphi^{5/2}}\Big\}\geq 1-e^{-(\log N)^D}.
\end{align*}
See Lemma \ref{stieltjes - bounding quadratic variation at a particular point outside} in the appendix for the details of bounding the quadratic variation. Substituting these into the equation \eqref{h(s)=h(0)+M_s+int} gives that with overwhelming probability:
\begin{align*}
|h_j(s\wedge \tau)| \leq \frac{1}{\varphi^2 N \eta_{t_j}}+\int_0^{s\wedge\tau} |h_j(u)||\partial_z m_{u}(z_{t-u})|\D u
\end{align*}
for all $j\in \{1,2,\dots,\lceil N^4(\log N)^2\rceil\}$ and $0\leq s\leq t_j$. Applying Gr\"onwall's lemma and following the same steps as in the proof of Lemma \ref{first step of the proof} we reach to the conclusion that in our set with overwhelming probability, the stopping time $\tau$ does not stop due to $|h_j(u)|\geq e^{-u/2}\frac{1}{\varphi N\eta_{t_j-u}}$ for any $j$. 

Now we will prove that in this set with overwhelming probability $\lambda_N(t)$ is always smaller than $2+\upkappa$. Assume the contrary. Then the stopping time stops due to detection of an outlier, i.e. $\lambda_N(\tau)=2+\upkappa$. Let $\ell=\argmin_j\{t_j>\tau\}$, which gives $0\leq t_\ell-\tau<N^{-4}$. $m_\tau$ is $N^3$-Lipschitz when $\eta>1/N$. Then $|m_\tau(z)-m_\tau(z_{t_\ell-\tau})|<N^{-1}$. Moreover, we have shown that $| m_\tau(z_{t_\ell-\tau})-e^{-(t_\ell-\tau)/2}m_{\rm sc}(z)|<\frac{1}{\varphi N\eta_{t_\ell-\tau}}\ll \frac{1}{\varphi^{1/2}N\eta}$. Thus, the contradiction follows from 
\begin{align*}
\frac{1}{N\eta}=\frac{1}{N}\Im (\frac{1}{\lambda_1(\tau)-z})<\Im(m_\tau(z))<&|m_\tau(z)-m_\tau(z_{t_\ell-\tau})|+| m_\tau(z_{t_\ell-\tau})-e^{-(t_\ell-\tau)/2}m_{\rm sc}(z)|
\\
&+|e^{-(t_\ell-\tau)/2}m_{\rm sc}(z)-m_{\rm sc}(z)| +\Im(m_{\rm sc}(z))\ll \frac{1}{\varphi^{1/2}N\eta},
\end{align*}
concluding the proof.
\end{proof}

Now that we have all the necessary tools, we can establish the stability of rigidity. The following proposition demonstrates that under the given dynamics, on any time window of length $(\log N)^2$, the eigenvalues remain near to their typical values with overwhelming probability.

\begin{prop}\label{prop:rig_until_logN2} For any $D>0$ there exists an $N_0$ such that for all $N\geq N_0$ we have
\begin{equation*}
\Prob\Big\{\boldsymbol{\lambda}(t) \in \mathcal{G}, \textnormal{ for all } t\in [0,(\log N)^2] \Big\}\geq 1-e^{-(\log N)^D}.
\end{equation*}
\end{prop}

\begin{proof}
Fix an arbitrary  time $t\in[0,(\log N)^2]$ and index $j\in[\![N]\!]$ and define two smooth counting functions $g_j:\R\to[0,1]$ and $h_j:\R\to[0,1]$:
\begin{align*}
g_j(x)=\begin{cases}
1, &x\in[-2-\frac{\varphi^{6}}{2N^{2/3}},\gamma_j]\\
0, & x\in[-2-\frac{\varphi^{6}}{N^{2/3}},\gamma_j+\frac{\varphi^{6}}{N^{2/3}\hat{j}^{1/3}}]^c
\end{cases}
\end{align*}
and
\begin{align*}
h_j(x)=\begin{cases}
1, &x\in[-2-\frac{\varphi^{6}}{N^{2/3}},\gamma_j-\frac{\varphi^{6}}{N^{2/3}\hat{j}^{1/3}}]\\
0, & x\in[-2-\frac{\varphi^{6}}{2N^{2/3}},\gamma_j]^c
\end{cases}
\end{align*}
with smooth interpolations in between the given domains. Writing Helffer-Sj\"ostrand formula for $g_j$ with bump function $\chi:=\chi_{1/2}$ and using Lemma \ref{m_t(z) estimate on the rectangle} we obtain that with overwhelming probability:
\begin{align*}
|S_N(g_j)(H_t)|&\asymp \Big| \int\int_{\R^2} \big(iyg_j''(x)\chi(y)+ig_j(x)\chi'(y)-yg_j'(x)\chi'(y))\big)N(m_t(z)-m_{\rm sc}(z))\D x\D y \Big|
\\
& \lesssim  \Big| \int\int_{\R^2} iyg_j''(x)\chi(y)N(m_t(z)-m_{\rm sc}(z))\D x\D y \Big|+\varphi^{9}.
\end{align*}
In order to bound the integral, integration by parts can be used as follows. Let $\theta=\frac{\varphi^{6}}{N^{2/3}\hat{j}^{1/3}}$, when $|y|\leq \theta$ we use the estimation in Lemma \ref{m_t(z) estimate on the rectangle} again:
\begin{align*}
\Big| \int\int_{|y|\leq \theta} iyg_j''(x)\chi(y)N(m_t(z)-m_{\rm sc}(z))\D x\D y \Big|\leq  \Big| \int\int_{|y|\leq \theta} g_j''(x)\chi(y)\D x\D y \Big|\varphi^{9}\lesssim \|g'\|_\infty \theta \varphi^{9}\lesssim \varphi^{9}
\end{align*}
and when $|y|>\theta$ we use integration by parts:
\begin{multline*}
\Big| \int\int_{y> \theta}  iyg_j''(x)\chi(y)N(m_t(z)-m_{\rm sc}(z))\D x\D y \Big|=\Big|N\int_{y>\theta}y\chi(y)\int_\R g_j'(x)(m_t'(z)-m_{\rm sc}'(z))\D x \D y\Big|
\\
=\Big|N\int_\R g_j'(x)\Big[\theta\chi(\theta)(m_t(x+i\theta)-m_{\rm sc}(x+i\theta))-\int_{y>\theta} (\chi(y)+y\chi'(y))(m_t(z)-m_{\rm sc}(z))\D y\Big] \D x\Big|
\\
\leq \varphi^{9}\int_\R |g_j'(x)|\D x + \varphi^{9} \int_\R\int_{y> \theta} |g_j'(x)| |\frac{\chi(y)+y\chi'(y)}{y}| \D y \D x \lesssim  \varphi^{9} \log N.
\end{multline*}

Thus, we obtain that 
\begin{align*}
\#\{k:\lambda_k(t)\in[-2-\frac{\varphi^{6}}{N^{2/3}},\gamma_j]\}\leq \sum_k g_j(\lambda_k(t)) \in[j-\varphi^{19/2},j+\varphi^{19/2}]
\end{align*}
by the definition of $\gamma_j$. Similarly, repeating the same steps for the function $h_j$ defined in beginning of the proof gives the other direction of the inequality. Noting that we have no outlier due to Lemma \ref{no outliers lemma} we can conclude that $|\lambda_j(t)-\gamma_j|\leq \frac{\varphi^{10}}{N^{2/3}\hat{j}^{1/3}}$.
\end{proof}

\subsection{Projections of the resolvent along the characteristic curve}\label{subsec:proj_res}

In our dynamics, the absence of rotational symmetry, unlike in the circular ensemble, renders some contour argument used for uniformity results in \cite[Section 4.3]{bourgade2022liouville} inapplicable. Consequently, establishing estimates for the projections of the resolvent over the entire upper half-plane is not straightforward. Instead, we will derive these estimates along the trajectories of characteristic curves originating from a compact set of order 1.

For any $U\subset\mathbb{C}$ and $t\geq 0$ we define the set $R_t(U)$ by extending $U$ along the trajectories of the characteristic curves and then inflating it by a scale comparable to imaginary values:
\begin{align*}
R_t(U):= \bigcup_{w\in U, s\in[0,t]} B_{\Im w_{s}/100}(w_s)
\end{align*}
and we write $R_t^n(U)$ when this expansion is applied iteratively $n$ many times, i.e. $R_t^n(U):=R_t(R_t^{n-1}(U))$. We apply inflation to enable the use of the Cauchy integral formula in our iterative argument on the successively expanded sets, as will become clear later in this section.

For convenience, we also define a modified version of the curve $\mathcal{C}_{\alpha}$ outside the spectrum, denoted by $\mathcal{S}_{\alpha}:\R\to\R$ for $\alpha>0$, as follows:
\begin{align} \label{eqn:defn_of_s_curve}
\mathcal{S}_{\alpha}(x)=\begin{cases}
\mathcal{C}_{\alpha}(x) &,\ x\in[-2-N^{-2/3},2+N^{-2/3}] \\
\alpha\cdot\frac{(\upkappa(x))^{1/4}}{N^{1/2}}&,\ \textnormal{otherwise}
\end{cases}
\end{align}
where $\upkappa(x):=|x-2|\wedge|x+2|$ and $\mathcal{C}_{\alpha}$ as defined in \eqref{eqn:defn_of_mathcal_C_alpha_curve}. By its definition, it is easy to check that for any fixed $C>0$, $\eta\Im\sqrt{z^2-4}\gtrsim\frac{\alpha}{N}$ uniformly for all $\alpha\geq 1$ and $z$ with $|\eta|,|E|<C$ above the curve $\mathcal{S}_{\alpha}$.

Later in this subsection we will use two parameters $z$ and $w$. We settle the following notation for the real and imaginary values of the characteristics starting from these parameters: $z=E+i\eta$, $z_t=E_t+i\eta_t$ and $w=E_w+i\eta_w$, $w_t=E_{w_t}+i\eta_{w_t}$. The following lemma can be viewed as a stability of certain properties along characteristic curves under small perturbations. It enables us to transfer specific properties of $U$ to its extension $R_t^n(U)$.

\begin{lemma}\label{lemma:simple_asymp_along_char} 
For any $z$ in the upper half plane, $w\in B_{\frac{\eta}{100}}(z)$, and $t\geq 0$,
\begin{align}
\Im\sqrt{z^2-4}&\asymp  \Im\sqrt{w^2-4}, \label{asymp_z_and_w_eqn1}
\\
\eta_t&\asymp  \eta_{w_{t}}, \label{asymp_z_and_w_eqn2}
\\
|z_t-y|\asymp  |w_t-y|&, \textnormal{ for all } y\in\C \textnormal{ with }\Im y\leq 0 \label{asymp_z_and_w_eqn3}
\end{align}
where $z_t$ and $w_t$ are defined as in \eqref{eqn:characteristic_curve} and the relations are uniform in $z$, $w$, $t$ and $y$. Moreover, $\eta\Im\sqrt{z^2-4}$ increases along the characteristics, more explicitly for every $t\geq 0$
\begin{align} \label{eqn:eta_Imsqrt_inc}
\eta_{t}\Im\sqrt{z_{t}^2-4}\geq e^t\eta\Im\sqrt{z^2-4}.
\end{align}

\end{lemma}
\begin{proof}\renewcommand{\qedsymbol}{}
See Appendix \ref{app:stable_asymp_along_char}.
\end{proof}

Note that when the expansion $R_t$ is applied, the growth of the size of the set is exponential in time due to the formula of the characteristic curve. So, in order to preserve uniformity of our results in space, we restrict the expansion times to be of order $\log N$.

\begin{prop}\label{prop:finiterankresolvent}
For any $C>0$, $D>0$ and $n\in\N$ there exists an $N_0$ such that for all $N\geq N_0$, $t\in[0,(\log N)^2]$, $z \in R_{n\log N}^n \Big([-C,C]\times[\mathcal{S}_{\varphi^{20}}(E),C]\Big)$ and $q\in\{\omega\in\mathbb{C}^N:|\omega|=1\}$ we have
\begin{align}
\Prob\bigg\{ \Big|e^{-u/2}\inner{q}{\frac{1}{H_u-z_{t-u}}q}-\inner{q}{\frac{1}{H_0-z_t}q}\Big|\leq e^{-(t-u)/4}&\frac{\varphi}{\sqrt{N\eta_{t-u}\Im\sqrt{z^2-4}}}\Im\inner{q}{\frac{1}{H_0-z_t}q} ,\forall u\in[0,t]  \bigg\}\geq 1-e^{-(\log N)^D}. \label{finiterank_up_to_time_t}
\end{align}
Moreover, for any $C>0$ and $D>0$ there exists an $N_0$ such that for all $N\geq N_0$ and $q\in \{\omega\in \mathbb{C}^N:|\omega|=1\}$ we have
\begin{align}\label{inside_the_box_finite_rank}
\Prob\bigg\{ \Big|e^{-t/2}\inner{q}{\frac{1}{H_t-z}q}-\inner{q}{\frac{1}{H_0-z_t}q}\Big|\leq \frac{\varphi}{\sqrt{N\eta\Im\sqrt{z^2-4}}}\Im\inner{q}{\frac{1}{H_0-z_t}q}, \ \forall t\in[0,(\log N)^2], & \nonumber
\\
 \forall z \in R_{n\log N}^n \Big([-C,C]\times[\mathcal{S}_{\varphi^{20}}(E),C]\Big)   &\bigg\}\geq 1-e^{-(\log N)^D}.
\end{align}
\end{prop}

\begin{proof}

Fix an arbitrary $t\in [0,(\log N)^2]$ and $z \in R_{n\log N}^n \Big([-C,C]\times[\mathcal{S}_{\varphi^{20}}(E),C] \Big)$.  Setting $A=qq^*$, the SDE along the characteristics for $\tilde{m}_{t,A}$, i.e. equation \eqref{SDE along char for tilde m_t,A}, becomes
\begin{align*}
\D \tilde m_{u,A}(z_{t-u})=\frac{-e^{-u/2}}{\sqrt{N}}\sum_{k,j}\frac{q^*v_k(u) v_j(u)^*q}{(\lambda_k(u)-z_{t-u})(\lambda_j(u)-z_{t-u})}\D (\tilde{B}_{k,j})_u+\big(\tilde m_u(z_{t-u})-e^{-u/2}m_{\rm sc}(z_{t-u})\big)\partial_z m_{u,A}(z_{t-u})du
\end{align*}
where $v_k(t)$'s are normalized eigenvectors corresponding to the eigenvalues $\lambda_k(t)$'s and $\tilde{B}_{k,j}$'s are independent Brownian motions up to symmetry in $k,j$. Define a stopping time
 \begin{align*}
 \tau=\inf\Big\{u\in[0,t]:\big|m_{0,A}(z_t)-\tilde{m}_{u,A}(z_{t-u})\big|\geq e^{-(t-u)/4}\frac{\varphi}{\sqrt{N\eta_{t-u}\Im\sqrt{z^2-4}}}\Im(m_{0,A}(z_t))\Big\}\wedge t.
 \end{align*}
 
 Note that for any $z\in [-C,C]\times[\mathcal{S}_{\varphi^{20}}(E),C]$, $\eta\Im\sqrt{z^2-4}\gtrsim \frac{\varphi^{20}}{N}$ by the definition of $\mathcal{S}_{\varphi^{20}}$. Hence for all $z\in R_{n\log N}^n \Big([-C,C]\times[\mathcal{S}_{\varphi^{20}}(E),C]\Big)$, by equation \eqref{eqn:eta_Imsqrt_inc}, $e^{-(t-u)/4}\frac{\varphi}{\sqrt{N\eta_{t-u}\Im\sqrt{z^2-4}}}\lesssim \varphi^{-9}$. Also note that $m_{t,A}(z)=\sum_k\frac{|q^*v_k(t)|^2}{\lambda_k(t)-z}$, $\Im(\tilde{m}_{t,A}(z))=\eta e^{-t/2}\Big(\sum_k\frac{|q^*v_k(t)|^2}{|\lambda_k(t)-z|^2}\Big)$ and $\partial_z m_{u,A}(z_{t-u})=\sum_k\frac{|q^*v_k(u)|^2}{(\lambda_k(u)-z_{t-u})^2}$. With these identities, we carry out an argument similar to those used previously in this section.
 
 The quadratic variation of the martingale term is bounded by the following expression:
\begin{align*}
\frac{1}{N}\int_0^{s\wedge\tau} \frac{\big(\Im(\tilde{m}_{u,A}(z_{t-u}))\big)^2}{\eta_{t-u}^2}\D u\asymp \frac{\big(\Im(m_{0,A}(z_{t}))\big)^2}{N}\int_0^{s\wedge\tau} \frac{1}{\eta_{t-u}^2}\D u \lesssim e^{-(t-s)/2}\frac{\big(\Im(m_{0,A}(z_{t}))\big)^2}{N\eta_{t-s}\Im\sqrt{z^2-4}}.
\end{align*}
Thus, with overwhelming probability, the stochastic integral is bounded by $e^{-(t-s)/4}\frac{\varphi^{1/2}\Im(m_{0,A}(z_{t}))}{\sqrt{N\eta_{t-s}\Im\sqrt{z^2-4}}}$. On the other hand, Riemann integral term can be easily bounded, without Gr\"onwall's inequality, using rigidity as follows:
\begin{multline*}
\Big|\int_0^{s\wedge\tau} \big(\tilde m_u(z_{t-u})-e^{-t/2}m_{\rm sc}(z)\big)\partial_z m_{u,A}(z_{t-u})\D u\Big| \leq  \int_0^{s\wedge\tau} e^{-u/2}\frac{\varphi^{11}}{N\eta_{t-u}}\cdot \frac{e^{u/2}}{\eta_{t-u}}\Im(\tilde{m}_{u,A}(z_{t-u})) \D u
\\
\lesssim  \frac{\varphi^{11}\Im(m_{0,A}(z_{t}))}{N} \int_{0}^s \frac{1}{\eta_{t-u}^2}\D u \lesssim  \frac{\varphi^{11}\Im(m_{0,A}(z_{t}))}{N} \cdot \frac{e^{-(t-s)/2}}{\eta_{t-s}\Im\sqrt{z^2-4}}\ll e^{-(t-s)/4}\frac{\Im(m_{0,A}(z_{t}))}{\sqrt{N\eta_{t-s}\Im\sqrt{z^2-4}}}.
\end{multline*}

Therefore, stopping time does not stop due to the difference $\big|m_{0,A}(z_t)-\tilde{m}_{u,A}(z_{t-u})\big|$ exceeding the upper bound, on a set with overwhelming probability. Thus, we have proved \eqref{finiterank_up_to_time_t}. Evaluating it at $u=t$ and applying the classical grid argument equation \eqref{inside_the_box_finite_rank} follows. Notice that for the uniformity of the space, the logarithmic rate of the expansion is essential so that the space variable is on a set whose size is polynomial in $N$. 
\end{proof}

\begin{remark} If the proposition is started at time $s$, $q$ can be taken as $\mathcal{F}_s:=\sigma(H_u;u\leq s)$ measurable random variable valued in $\{w\in \mathbb{C}^N:|w|=1\}$ (cf. \cite[Proposition 4.3]{bourgade2022liouville}).
\end{remark}

Later in this section, we let $v_k(t)$ denote a normalized eigenvector of $H_t$ corresponding to $\lambda_k(t)$, and we write $v_k = v_k(0)$. Using the polarization identity
\begin{align*}
\hspace{-1cm}\inner{v}{\frac{1}{H_t-z}w}=\frac{1}{4} \sum_{k=0}^3 i^k\inner{v+i^kw}{\frac{1}{H_t-z}(v+i^kw)} 
\end{align*}
the corollary below follows from the previous proposition. 

\begin{cor}\label{cor: weaker bound on finite rank projections}
For any $C>0$, the following inequality holds with overwhelming probability uniformly in $t\in[0,(\log N)^2]$ and $z\in R_{n\log N}^n \Big([-C,C]\times[\mathcal{S}_{\varphi^{20}}(E),C]\Big)$, for every $i\neq j$
\begin{align*}
e^{-t/2}\Big|\inner{v_i}{\frac{1}{H_t-z}v_j}\Big|\leq \frac{\varphi}{\sqrt{N\eta\Im\sqrt{z^2-4}}}\Big(\frac{\eta_t}{|\lambda_i(0)-z_t|^2}+\frac{\eta_t}{|\lambda_j(0)-z_t|^2}\Big).
\end{align*}

\end{cor}

The next proposition improves on the previous corollary by a bootstrap argument.

\begin{prop} \label{cross terms bound}
For any $C,D,\varepsilon>0$, $n\in\N$ there exists an $N_0$ such that for all $N\geq N_0$ and $i\neq j$,
\begin{multline*}
\Prob\bigg\{ e^{-t/2}\cdot\Big|\inner{v_i}{\frac{1}{H_t-z}v_j}\Big|\leq \frac{\varphi}{\sqrt{N\eta\Im\sqrt{z^2-4}}}\frac{\eta_t}{|\lambda_i(0)-z_t||\lambda_j(0)-z_t|}, \ \forall t\in[0,n\log N],
\\
 \forall z\in R_{n\log N}^n \Big([-C,C]\times[\mathcal{S}_{N^{\varepsilon}}(E),C]\Big) \bigg\}\geq 1-e^{-(\log N)^D}.
\end{multline*}
\end{prop}

\begin{proof}
Define $m=\lceil 100/\varepsilon\rceil$ and fix $t\in[0,n\log N]$, $z\in R_{\log N}^{n+m} \Big([-C,C]\times[\mathcal{S}_{N^{\varepsilon}}(E),C]\Big)$ and $i\neq j$. Throughout this proof we work on the set with overwhelming probability for which equation \eqref{finiterank_up_to_time_t} holds when $q$ is $v_i=v_i(0)$ and $v_j=v_j(0)$.

 Define $A=v_jv_i^*$. Then the SDE for $\tilde{m}_{u,A}$ along the characteristic curve becomes:
\begin{align*}
\D \tilde m_{u,A}(z_{t-u})=\frac{-e^{-u/2}}{\sqrt{N}}\sum_{k,l}\frac{v_k(u)^* v_j(0)v_i(0)^*v_l(u)}{(\lambda_k(u)-z_{t-u})(\lambda_l(u)-z_{t-u})}\D (\tilde{B}_{k,j})_u+\big(\tilde m_u(z_{t-u})-e^{-u/2}m_{\rm sc}(z_{t-u})\big)\partial_z m_{u,A}(z_{t-u})du
\end{align*}
where $\tilde{B}_{k,j}$'s are independent Brownian motions up to symmetry in $k,j$. Notice that $m_{t,A}(z)=\sum_k\frac{v_i^*v_k(t)v_k(t)^*v_j}{\lambda_k(t)-z}$, $\Im(\tilde{m}_{t,v_iv_i^*}(z))=\eta e^{-t/2}\Big(\sum_k\frac{|v_i^*v_k(t)|^2}{|\lambda_k(t)-z|^2}\Big)$, $\Im(m_{0,v_iv_i^*}(z))= \frac{\eta}{|\lambda_i(0)-z|^2}$, and $|\partial_z m_{u,A}(z_{t-u})|\leq \sum_k\frac{|v_i^*v_k(u)||v_k(u)^*v_j|}{|\lambda_k(u)-z_{t-u}|^2}\leq e^{u/2}\frac{1}{\eta_{t-u}}\big(\Im(\tilde{m}_{u,v_iv_i^*}(z_{t-u}))+\Im(\tilde{m}_{u,v_jv_j^*}(z_{t-u}))\big)$. The quadratic variation of the martingale term is bounded by,
\begin{multline*}
\frac{1}{N}\int_0^{t}\frac{\Im(\tilde{m}_{u,v_iv_i^*}(z_{t-u}))\Im(\tilde{m}_{u,v_jv_j^*}(z_{t-u}))}{\eta_{t-u}^2}\D u \lesssim   \frac{1}{N}\int_0^{t}\frac{\Im(m_{0,v_iv_i^*}(z_{t}))\Im(m_{0,v_jv_j^*}(z_{t}))}{\eta_{t-u}^2}\D u
\\
= \frac{\eta_t^2}{N|\lambda_i(0)-z_t|^2|\lambda_j(0)-z_t|^2}\int_0^{t}\frac{1}{\eta_{t-u}^2}\D u \lesssim \frac{1}{N\eta\Im\sqrt{z^2-4}}\frac{\eta_t^2}{|\lambda_i(0)-z_t|^2|\lambda_j(0)-z_t|^2}
\end{multline*}
Then, the martingale term is less than $\frac{\varphi^{1/2}}{\sqrt{N\eta\Im\sqrt{z^2-4}}}\frac{\eta_t}{|\lambda_i(0)-z_t||\lambda_j(0)-z_t|}$ with overwhelming probability. On the other hand, the Riemann integral term is bounded with
\begin{multline*}
\int_0^{t} e^{-u/2}\frac{\varphi^{11}}{N\eta_{t-u}}\cdot e^{u/2}\frac{1}{\eta_{t-u}}\big(\Im(\tilde{m}_{u,v_iv_i^*}(z_{t-u}))+\Im(\tilde{m}_{u,v_jv_j^*}(z_{t-u}))\big) \D u
\\
\lesssim  \int_0^{t}\frac{\varphi^{11}}{N\eta_{t-u}^2}\big(\frac{\eta_t}{|\lambda_i(0)-z_t|^2}+\frac{\eta_t}{|\lambda_j(0)-z_t|^2}\big) \D u
\lesssim  \frac{\varphi^{11}}{N\eta\Im\sqrt{z^2-4}}\big(\frac{\eta_t}{|\lambda_i(0)-z_t|^2}+\frac{\eta_t}{|\lambda_j(0)-z_t|^2}\big)
\end{multline*}
on the set with overwhelming probability. Therefore, with overwhelming probability,
\begin{align*}
e^{-t/2}\cdot\Big|\inner{v_i}{\frac{1}{H_t-z}v_j}\Big|\leq \frac{\varphi^{1/2}}{\sqrt{N\eta\Im\sqrt{z^2-4}}}\frac{\eta_t}{|\lambda_i(0)-z_t||\lambda_j(0)-z_t|}+\frac{\varphi^{12}}{N\eta\Im\sqrt{z^2-4}}\big(\frac{\eta_t}{|\lambda_i(0)-z_t|^2}+\frac{\eta_t}{|\lambda_j(0)-z_t|^2}\big)
\end{align*}
And uniformity of this inequality in $t\in[0,n\log N]$, $z\in R_{n\log N}^{n+m} \Big([-C,C]\times[\mathcal{S}_{N^{\varepsilon}}(E),C]\Big)$ follows by the grid argument.

Construct hypothesis $(\textnormal{P}_k)$ as follows: For every $C,D>0$ there exists $N_0=N_0(k,C,D)$ such that for all $N\geq N_0$, and $i\neq j$:
\begin{align*}
\Prob\bigg\{ e^{-t/2}\cdot\Big|\inner{v_i}{\frac{1}{H_t-z}v_j}\Big|\leq \frac{\varphi^{1/2}}{\sqrt{N\eta\Im\sqrt{z^2-4}}}\frac{\eta_t}{|\lambda_i(0)-z_t||\lambda_j(0)-z_t|}+\frac{\varphi^{13(k+1)}}{(N\eta\Im\sqrt{z^2-4})^{k+1}}&\big(\frac{\eta_t}{|\lambda_i(0)-z_t|^2}+\frac{\eta_t}{|\lambda_j(0)-z_t|^2}\big),
\\
  \forall t\in[0,n\log N], \ \forall z\in R_{n\log N}^{n+m-k} \Big([-C,C]\times[\mathcal{S}_{N^{\varepsilon}}(E),C]\Big) &\bigg\}\geq 1-e^{-(\log N)^D}
\end{align*}
We proved $(\textnormal{P}_0)$ and will now show that $(\textnormal{P}_{k+1})$ can be deduced from $(\textnormal{P}_k)$ for every $k=0,1,\dots,m-1$. Martingale term is already bounded as required, so, it suffices to bound the Riemann integral term. Note that for any $z\in R_{n\log N}^{n+m-(k+1)} \Big([-C,C]\times[\mathcal{S}_{N^{\varepsilon}}(E),C]\Big)$ and $s\in[0,n\log N]$, we know that $B_{\eta_s/100}(z_s)\subset R_{n\log N}^{n+m-k} \Big([-C,C]\times[\mathcal{S}_{N^{\varepsilon}}(E),C]\Big)$ by the construction of the expansion sets. Hence, the required implication $(\textnormal{P}_k)\Rightarrow (\textnormal{P}_{k+1})$ follows from Cauchy's inequality: For any $z\in R_{n\log N}^{n+m-(k+1)} \Big([-C,C]\times[\mathcal{S}_{N^{\varepsilon}}(E),C]\Big)$,
\begin{align*}
&\Big|\int_0^{t} \big(\tilde m_u(z_{t-u})-e^{-t/2}m_{\rm sc}(z)\big)\partial_z m_{u,A}(z_{t-u})du\Big| \lesssim \int_0^{t} e^{-u/2}\frac{\varphi^{11}}{N\eta_{t-u}}\cdot \frac{\max_{|w-z_{t-u}|=\eta_{t-u}/\varphi}| m_{u,A}(w)|}{\eta_{t-u}/\varphi} \D u
\\
&\hspace{-1cm}\lesssim \frac{\varphi^{12}}{N}\int_0^t \frac{1}{\eta_{t-u}^2} \Big(\frac{\varphi^{1/2}}{\sqrt{N\eta_{t-u}\Im\sqrt{z_{t-u}^2-4}}}\frac{\eta_t}{|\lambda_i(0)-z_t||\lambda_j(0)-z_t|}+\frac{\varphi^{13(k+1)}}{(N\eta_{t-u}\Im\sqrt{z_{t-u}^2-4})^{k+1}}\big(\frac{\eta_t}{|\lambda_i(0)-z_t|^2}+\frac{\eta_t}{|\lambda_j(0)-z_t|^2}\big) \Big)\D u
\\
&\hspace{-1cm}\ll \frac{\varphi^{13}}{\sqrt{N\eta\Im\sqrt{z^2-4}}}\frac{\eta_t}{|\lambda_i(0)-z_t||\lambda_j(0)-z_t|}\frac{1}{N}\int_0^t\frac{1}{\eta_{t-u}^2}\D u + \frac{\varphi^{13(k+2)}}{(N\eta\Im\sqrt{z^2-4})^{k+1}}\big(\frac{\eta_t}{|\lambda_i(0)-z_t|^2}+\frac{\eta_t}{|\lambda_j(0)-z_t|^2}\big)\frac{1}{N}\int_0^t\frac{1}{\eta_{t-u}^2}\D u
\\
&\lesssim \frac{\varphi^{13}}{\sqrt{N\eta\Im\sqrt{z^2-4}}}\frac{\eta_t}{|\lambda_i(0)-z_t||\lambda_j(0)-z_t|}\frac{1}{N\eta\Im{\sqrt{z^2-4}}} + \frac{\varphi^{13(k+2)}}{(N\eta\Im\sqrt{z^2-4})^{(k+2)}}\big(\frac{\eta_t}{|\lambda_i(0)-z_t|^2}+\frac{\eta_t}{|\lambda_j(0)-z_t|^2}\big)
\end{align*}
where we have used Lemma \ref{lemma:simple_asymp_along_char}.  As for all $z\in R_{n\log N}^{n+m} \Big([-C,C]\times[\mathcal{S}_{N^{\varepsilon}}(E),C]\Big)$ we have $N\eta\Im\sqrt{z^2-4}\gtrsim N^{\varepsilon}$, this inequality shows $(\textnormal{P}_{k+1})$. Hence, we have proved $(\textnormal{P}_m)$ and due to the same fact, i.e. $N\eta\Im\sqrt{z^2-4}\gtrsim N^{\varepsilon}$, it is easy to deduce the corollary from $(\textnormal{P}_m)$.
\end{proof}

\begin{prop} \label{prop:full_rank_estimate}
For any $C,D,\varepsilon>0$ and $n\in\N$ there exists $N_0$ such that for every $N\geq N_0$ we have
\begin{multline*}
\Prob\Big\{\Big|e^{-t/2}\Tr\big(\frac{1}{H_t-z}\frac{1}{H_0-w}\big)-\Tr\big(\frac{1}{H_0-z_t}\frac{1}{H_0-w}\big)\Big|\leq    \mathcal{E}(z,w,t),\forall t\in[0,n\log N],
\\
   \forall z\in R_{n\log N}\big([-C,C]\times[\mathcal{S}_{N^{\varepsilon}}(E),C]\big), \forall w\in[-10,10]\times(0,1]   \Big| \boldsymbol\lambda(0)\in\tilde{\mathcal{G}}\Big\}\geq 1-e^{-(\log N)^D}
\end{multline*}
where
\begin{multline*} 
\mathcal{E}(z,w,t):=\frac{N^{\varepsilon/2}}{N^{1/2}\sqrt{\eta\Im\sqrt{z^2-4}}}\frac{1}{|z_t-\bar{w}|\eta_w}\mathds{1}_{\eta_{w}<\mathcal{C}_{N^{\epsilon/2}}(E_w)}
+\frac{\varphi^{12}}{\eta\Im\sqrt{z^2-4}}\big(\int_\R \frac{\eta_t}{|s-w'||s-z_t|^2}\rhosc (s)\D s\big)
\\
+\frac{\varphi^2}{\sqrt{\eta\Im\sqrt{z^2-4}}}\big(\int_\R\frac{\eta_t^2}{|s-w'|^2|s-z_t|^4}\rhosc (s)\D s\big)^{1/2}
\end{multline*}
with $w':=w+i\mathcal{C}_{N^{\varepsilon/2}}(E_w)$.

\end{prop}

\begin{proof}
Take $A=\frac{1}{H_0-w}$. Writing the SDE along the characteristics for $\tilde{m}_{t,A}$ gives:
\begin{multline*}
\D \tilde m_{u,A}(z_{t-u})=\frac{-e^{-u/2}}{\sqrt{N}}\sum_{i,j}\frac{1}{(\lambda_i(u)-z_{t-u})(\lambda_j(u)-z_{t-u})}v_i(u)^*\frac{1}{H_0-w}v_{j}(u)\D (\tilde{B}_{i,j})_u
\\
+\big(\tilde m_u(z_{t-u})-e^{-u/2}m_{\rm sc}(z_{t-u})\big)\partial_z m_{u,A}(z_{t-u})du
\end{multline*}
where $\tilde{B}_{k,j}$'s are independent Brownian motions up to symmetry in $k,j$. The quadratic variation of the martingale term is bounded with
\begin{align*}
\frac{1}{N}\sum_{i,j} &\int_{0}^{t} e^{-u}\frac{1}{|\lambda_i(u)-z_{t-u}|^2|\lambda_j(u)-z_{t-u}|^2} \Big|v_i(u)^*\frac{1}{H_0-w}v_j(u)\Big|^{2}\D u
\\
=&\frac{1}{N}\sum_{i,j} \int_{0}^{t} e^{-u}\frac{1}{|\lambda_i(u)-z_{t-u}|^2|\lambda_j(u)-z_{t-u}|^2} \Big|\sum_{k} \frac{v_i(u)^*v_k(0)v_k(0)^*v_j(u)}{\lambda_k(0)-w}\Big|^{2}\D u
\\
=&\frac{1}{N}\sum_{i,j} \int_{0}^{t} e^{-u}\frac{1}{|\lambda_i(u)-z_{t-u}|^2|\lambda_j(u)-z_{t-u}|^2} \Big(\sum_{k,l} \frac{v_i(u)^*v_k(0)v_k(0)^*v_j(u)\cdot v_l(0)^*v_i(u)v_j(u)^*v_l(0)}{(\lambda_k(0)-w)\overline{(\lambda_l(0)-w)}}\Big)\D u
\\
=&\frac{1}{N}\sum_{k,l} \int_{0}^{t} e^{-u}\frac{1}{(\lambda_k(0)-w)\overline{(\lambda_l(0)-w)}} \Big(\sum_{j} \frac{v_k(0)^*v_j(u)v_j(u)^*v_l(0)}{|\lambda_j(u)-z_{t-u}|^2}\Big)\Big(\sum_{i} \frac{v_i(u)^*v_k(0) v_l(0)^*v_i(u)}{|\lambda_i(u)-z_{t-u}|^2}\Big)\D u
\\
\leq&\frac{1}{N}\sum_{k} \int_{0}^{t} \frac{1}{|\lambda_k(0)-w|^2} \Big(e^{-u/2}\sum_{i} \frac{|v_i(u)^*v_k(0)|^2}{|\lambda_i(u)-z_{t-u}|^2}\Big)^2\D u
\\
&+\frac{1}{N}\sum_{k\neq l} \int_{0}^{t} \frac{1}{|\lambda_k(0)-w||\lambda_l(0)-w|} \cdot \Big|e^{-u/2}\sum_{i} \frac{v_i(u)^*v_k(0) v_l(0)^*v_i(u)}{|\lambda_i(u)-z_{t-u}|^2}\Big|^2\D u.
\end{align*}
Recall that $w'$ is given by $\Re w'=\Re w$ and $\Im w'=\Im w+\mathcal{C}_{N^{\varepsilon/2}}(E_w)$. The diagonal terms (i.e. $k=l$) can be bounded on a set with overwhelming probability as follows:
\begin{align*} 
\frac{1}{N}\sum_{k} &\int_{0}^{t} \frac{1}{|\lambda_k(0)-w|^2} \Big(e^{-u/2}\sum_{i} \frac{|v_i(u)^*v_k(0)|^2}{|\lambda_i(u)-z_{t-u}|^2}\Big)^2\D u
\\
&=\frac{1}{N}\sum_{k} \int_{0}^{t} \frac{1}{|\lambda_k(0)-w|^2} \Big(\frac{1}{\eta_{t-u}} e^{-u/2}  \Im\big(\inner{v_k(0)}{\frac{1}{H_u-z_{t-u}}v_k(0)}\big) \Big)^2\D u
\\
&\lesssim \frac{1}{N}\sum_{k} \int_{0}^{t} \frac{1}{|\lambda_k(0)-w|^2} \Big(\frac{1}{\eta_{t-u}} \Im(\frac{1}{\lambda_k(0)-z_t}) \Big)^2\D u
\\
& = \frac{1}{N}\int_{0}^{t}\frac{1}{\eta_{t-u}^2} \D u\cdot \sum_{k} \frac{\eta_t^2}{|\lambda_k(0)-w|^2|\lambda_k(0)-z_t|^4} 
\\
&\lesssim \frac{1}{N\eta\Im\sqrt{z^2-4}} \Big( \sum_{k:|\lambda_k(0)-w|\leq \mathcal{C}_{N^{\epsilon/2}}(E_w)} \frac{\eta_t^2}{|\lambda_k(0)-w|^2|\lambda_k(0)-z_t|^4}+\sum_{k} \frac{\eta_t^2}{|\lambda_k(0)-w'|^2|\lambda_k(0)-z_t|^4}  \Big)
\\
&\lesssim \frac{N^{2\varepsilon/3}}{N\eta\Im\sqrt{z^2-4}}\frac{\eta_t^2}{|z_t-\bar{w}|^4\eta_{w}^2}\mathds{1}_{\eta_{w}<\mathcal{C}_{N^{\epsilon/2}}(E_w)} +\frac{1}{\eta\Im\sqrt{z^2-4}}  \int_\R \frac{\eta_t^2}{|s-w'|^2|s-z_t|^4}\rhosc (s)\D s
\end{align*}
where we have used Proposition \ref{prop:finiterankresolvent}. Similarly, for the off-diagonal terms (i.e. $k\neq l$) we have the bound
\begin{align*}
&\frac{1}{N}\sum_{k\neq l} \int_{0}^{t} \frac{1}{|\lambda_k(0)-w||\lambda_l(0)-w|} \cdot \Big|\frac{1}{2\eta_{t-u}} \Big( e^{-u/2}\sum_{i} \frac{v_i(u)^*v_k(0) v_l(0)^*v_i(u)}{\lambda_i(u)-z_{t-u}}-e^{-u/2}\sum_{i} \frac{v_i(u)^*v_k(0) v_l(0)^*v_i(u)}{\lambda_i(u)-\overline{z_{t-u}}}\Big) \Big|^2\D u
\\
&\lesssim \frac{1}{N}\sum_{k\neq l} \int_{0}^{t} \frac{1}{|\lambda_k(0)-w||\lambda_l(0)-w|}  \frac{1}{\eta_{t-u}^2}  \Big( \big|e^{-u/2}\inner{v_l(0)}{\frac{1}{H_u-z_{t-u}}v_k(0)}\big|^2+\big|e^{-u/2}\inner{v_l(0)}{\frac{1}{H_u-\overline{z_{t-u}}}v_k(0)}\big|^2\Big) \D u
\\
&=\frac{1}{N}\sum_{k\neq l} \int_{0}^{t} \frac{1}{|\lambda_k(0)-w||\lambda_l(0)-w|}  \frac{1}{\eta_{t-u}^2}  \Big( \big|e^{-u/2}\inner{v_l(0)}{\frac{1}{H_u-z_{t-u}}v_k(0)}\big|^2+\big|e^{-u/2}\inner{v_k(0)}{\frac{1}{H_u-z_{t-u}}v_l(0)}\big|^2\Big) \D u
\\
&\lesssim \frac{1}{N}\sum_{k\neq l} \int_{0}^{t} \frac{1}{|\lambda_k(0)-w||\lambda_l(0)-w|}  \frac{1}{\eta_{t-u}^2} \cdot  \frac{\varphi^2}{N\eta_{t-u}\Im\sqrt{z_{t-u}^2-4}}\frac{\eta_t^2}{|\lambda_k(0)-z_t|^2|\lambda_l(0)-z_t|^2} \D u
\\
&\lesssim \frac{\varphi^2}{N^2\Im\sqrt{z^2-4}}\sum_{k\neq l} \frac{\eta_t^2}{|\lambda_k(0)-z_t|^2|\lambda_l(0)-z_t|^2|\lambda_k(0)-w||\lambda_l(0)-w|} \int_{0}^{t}   \frac{e^{-(t-u)/2}}{\eta_{t-u}^3}   \D u
\\
&\lesssim \frac{\varphi^2}{\Im\sqrt{z^2-4}}\Big(\frac{1}{N}\sum_{k} \frac{\eta_t}{|\lambda_k(0)-z_t|^2|\lambda_k(0)-w|}\Big)^2 \frac{1}{\eta^2\Im\sqrt{z^2-4}}
\\
& \lesssim \frac{N^{4\varepsilon/3}}{N^2\eta^2(\Im\sqrt{z^2-4})^2}\frac{\eta_t^2}{|z_t-\bar{w}|^4\eta_w^2}\mathds{1}_{\eta_{w}<\mathcal{C}_{N^{\epsilon/2}}(E_w)} + \frac{\varphi^2}{\eta^2(\Im\sqrt{z^2-4})^2}\Big(\int_{\R}\frac{\eta_t}{|s-w'||s-z_t|^2}\rhosc (s)\D s\Big)^{1/2}
\end{align*}
where we have used $\overline{\inner{v_l(0)}{\frac{1}{H_u-\overline{z_{t-u}}}v_k(0)}}=\inner{v_k(0)}{\frac{1}{H_u-z_{t-u}}v_l(0)} $, Proposition \ref{cross terms bound}, $\Im\sqrt{z_t^2-4}\asymp e^{t/2}\Im\sqrt{z^2-4}$, and Proposition \ref{prop:rig_until_logN2}. Combining these we obtain that with overwhelming probability martingale term is bounded by
\begin{multline*}
\frac{N^{\varepsilon/2}}{N^{1/2}\sqrt{\eta\Im\sqrt{z^2-4}}}\frac{1}{|z_t-\bar{w}|\eta_w}\mathds{1}_{\eta_{w}<\mathcal{C}_{N^{\epsilon/2}}(E_w)}+\frac{\varphi^{2}}{\sqrt{\eta\Im\sqrt{z^2-4}}}  \Big(\int_\R \frac{\eta_t^2}{|s-w'|^2|s-z_t|^4}\rhosc (s)\D s\Big)^{1/2} 
\\
+\frac{\varphi^{2}}{\eta\Im\sqrt{z^2-4}}\Big(\int_\R \frac{\eta_t}{|s-w'||s-z_t|^2}\rhosc (s)\D s\Big)
\end{multline*}
Finally, the absolute value of the Riemann integral term $\int_0^t \tilde m_u(z_{t-u})-e^{-u/2}m_{\rm sc}(z_{t-u})\partial_z m_{u,A}(z_{t-u})\D u$ is bounded above by,
\begin{align*}
\int_0^t &e^{-u/2}\frac{\varphi^{11}}{N\eta_{t-u}}\Big|\sum_{k,l}\frac{|\inner{v_k(0)}{v_l(u)}|^2}{(\lambda_l(u)-z_{t-u})^2(\lambda_k(0)-w)}\Big|\D u
\leq  \int_0^t e^{-u/2}\frac{\varphi^{11}}{N\eta_{t-u}^2}\Big(\sum_{k}\frac{1}{|\lambda_k(0)-w|}\Im\inner{v_k(0)}{\frac{1}{H_u-z_{t-u}}v_k(0)}\Big)\D u
\\
&\lesssim  \frac{\varphi^{11}}{N} \Big(\sum_{k}\frac{1}{|\lambda_k(0)-w|}\Im\inner{v_k(0)}{\frac{1}{H_0-z_{t}}v_k(0)}\Big) \int_0^t \frac{1}{\eta_{t-u}^2}\D u
\lesssim  \frac{\varphi^{11}}{N\eta\Im\sqrt{z^2-4}}\Big(\sum_{k}\frac{1}{|\lambda_k(0)-w|}\frac{\eta_t}{|\lambda_k(0)-z_t|^2}\Big)
\\
&\lesssim  \frac{N^{\varepsilon/2}}{\sqrt{N\eta\Im\sqrt{z^2-4}}}\frac{1}{|z_t-\bar{w}|\eta_w}\mathds{1}_{\eta_{w}<\mathcal{C}_{N^{\epsilon/2}}(E_w)}+\frac{\varphi^{11}}{\eta\Im\sqrt{z^2-4}} \int_\R \frac{\eta_t}{|s-w'||s-z_t|^2}\rhosc (s)\D s
\end{align*}
Together with the bound on the martingale term, this gives the required bound $\mathcal{E}(z,w,t)$. Finally, uniformity in space and time follows from the grid argument again.
\end{proof}

\begin{remark}\label{rmk:overwhelming_under_bias}
All the propositions in Section \ref{sec:res_est} which hold with overwhelming probability (i.e. $1-e^{-(\log N)^{D}}$) also  hold under the biased measure,  if the bias satisfies the rigidity conditions: as in Lemma \ref{lemma:rigidity_conds},  $\Prob_{\mathbf{h}}(A^c)=\E{\mathds{1}_{A^c} \frac{e^{S_N(\mathbf{h})}}{\E{e^{S_N(\mathbf{h})}}}}\leq \Prob(A^c)\frac{\E{e^{2S_N(\mathbf{h})}}^{1/2}}{e^{-|\E{S_N(\mathbf{h})}|}}\leq e^{-(\log N)^{D}}$ for all sufficiently large $N$. 
\end{remark}

Note that the integrals in the error term $\mathcal{E}(z,w,t)$ can be bounded via the following lemma proof of which is straightforward.
\begin{lemma}\label{lemma:integral_with_two_sing}
Let $C>1$ be a real number, then there exists a constant $c$ such that for all $z,w\in\C$ with $\eta_{z}:=\Im z>0$, $\eta_{w}:=\Im w>0$, real numbers $a,b\in(1,C]$ and interval $I\subset \R$ with $\textnormal{len}(I)\leq C$ we have
\begin{gather*} 
\int_{I}\frac{1}{|s-z|^a|s-w|^b}\D s\leq \frac{c}{\eta_{z}^{a-1}|z-\bar{w}|^b}+\frac{c}{|z-\bar{w}|^{a}\eta_{w}^{b-1}},
\\
\int_{I}\frac{1}{|s-z|^a|s-w|}\D s \leq \frac{c}{\eta_{z}^{a-1}|z-\bar{w}|}+\frac{c}{|z-\bar{w}|^{a}}\log\frac{|z-\bar{w}|}{\eta_{w}}.
\end{gather*}
\end{lemma}

\begin{lemma}\label{lemma:two_times_resolvent_mult}
Given any $\varepsilon \in(0,1)$, $C>0$, $n\in\N$, $J\in\N$, $\mathbf{h}(H)=\sum_{j=1}^{J}h_j(H_{t_j})$ with $t_1,\dots,t_J\geq 0$ for which the rigidity conditions hold (see Lemma \ref{lemma:rigidity_conds}), uniformly in $w\in[-10,10]\times(0,1]$, $z\in R_{n\log N}\big([-C,C]\times[\mathcal{S}_{N^{\varepsilon}}(E),C]\big)$, $0\leq t\leq n\log N$, we have
\begin{align*}
\mathbb{E}_{\mathbf{h}}&\Big[\Tr\Big(\frac{1}{H_{0}-w}\frac{1}{(H_t-z)^2}\Big) \Big]=\ e^{t/2}\frac{\sqrt{z_t^2-4}}{\sqrt{z^2-4}}\mathbb{E}_{\mathbf{h}}\Big[\Tr\Big(\frac{1}{H_{0}-w}\frac{1}{(H_0-z_t)^2}\Big) \Big]+O\Big(\frac{e^{t/2}\mathcal{E}(z,w,t)}{\eta}\Big)+O(\frac{N}{\eta_w\eta^2}e^{-(\log N)^{10}}).
\end{align*}
\end{lemma}
\begin{proof} By Remark \ref{rmk:overwhelming_under_bias}, Proposition \ref{prop:full_rank_estimate} holds under the biased measure as well. Applying Cauchy's inequality to Proposition \ref{prop:full_rank_estimate} and using Lemma \ref{lemma:simple_asymp_along_char}, we obtain
\begin{align*} 
\mathbb{E}_{\mathbf{h}}\Big[\Tr\Big(\frac{1}{H_{0}-w}\frac{1}{(H_t-z)^2}\Big) \Big]=\mathbb{E}_{\mathbf{h}}\Big[\Tr\Big(\frac{1}{H_{0}-w}\frac{1}{(H_0-z_t)^2}\Big) \Big]\cdot e^{t/2}\frac{\D z_t}{\D z}+O\Big(\frac{e^{t/2}\mathcal{E}(z,w,t)}{\eta}\Big)+O(\frac{N}{\eta_w\eta^2}e^{-(\log N)^{10}})
\end{align*}
where $O(\frac{N}{\eta_w\eta^2}e^{-(\log N)^{10}})$ term stands for the expectation on the set with overwhelmingly small probability from Proposition \ref{prop:full_rank_estimate}. 
\end{proof}

\section{Multi-time loop equation asymptotics}\label{sec:loop_eqn}

In the following paragraphs, we present and interpret the single-time loop equation from \cite{johansson1998fluctuations} and outline the anticipated form of its multi-time generalization using a formula from \cite{unterberger2018global}. But first, let us introduce some notations.

For the Hilbert transform of $f:\R\to\R$, we   follow the convention
\begin{align*}
\hlbrt f(x):=\dashint_{-\infty}^{\infty} \frac{f(y)}{x-y}\D y
\end{align*}
where $\dashint$ stands for the principal value integral. For all $f\in L^p(\R)$ with $1\leq p<\infty$, $\hlbrt f$ exists almost everywhere on $\R$ and if $1<p<\infty$ then $\hlbrt$ is a bounded operator from $L^p(\R)$ to itself (see \cite[Theorems 100 and 101]{titchmarsh1937introduction}). Remember the definition of semicircle law $\rhosc (x)=\frac{\mathds{1}_{(-2,2)}(x)}{2\pi}\sqrt{4-x^2}$ and define the law of arcsine distribution $\rhoarcsin (x)=\frac{\mathds{1}_{(-2,2)}(x)}{\pi\sqrt{4-x^2}}$. We define semicircle and arcsin weighted Hilbert transforms as follows,
\begin{align*}
\uhlbrt f:=2\hlbrt(f\rhosc ),\quad \vhlbrt f:=-\hlbrt(f\rhoarcsin ).
\end{align*}
It is easy to show that $\uhlbrt$ and $\vhlbrt$ are inverse transformations for the continuous functions on $[-2,2]$ up to centering. More explicitly, we have the following proposition.
\begin{prop}\label{prop:inverse_of_uhlbrt}
Let $f$ be a function in $\mathscr{C}^1([-2,2])$. Then $\uhlbrt(\vhlbrt f)(x)$ exists everywhere in $[-2,2]$ and 
\begin{align}\label{uhlbrt_vhlbrt_on_R}
\uhlbrt(\vhlbrt f)(x)=f(x)-\int f(t)\rhoarcsin (t)\D t.
\end{align}
\end{prop}
\begin{proof}
By Lemma A.4 in \cite{lambert2019quantitative}, equation \eqref{uhlbrt_vhlbrt_on_R} holds a.e. in $(-2,2)$. Continuity of the left hand side on $[-2,2]$ can be argued easily by Plemelj-Privalov Theorem (see \cite{blaya2015cauchy}). This implies the equality everywhere on $[-2,2]$.
\end{proof}

Moreover, $\uhlbrt$ and $\vhlbrt$ have the some key properties relating them to the Chebyshev polynomials (see \cite[Chapter 10.11]{bateman1953higher}, \cite[Theorem 9.1]{mason2002chebyshev}):
\begin{equation}\label{eqn:Chebyshev_polynomials_properties}
\uhlbrt\tilde{U}_{n-1}=2\tilde{T}_{n} \textnormal{ on }[-2,2]
\quad \textnormal{and}\quad
\vhlbrt \tilde{T}_n=\frac{1}{2}\tilde{U}_{n-1},\ \vhlbrt \tilde{T}_0=0 \textnormal{ on }(-2,2), \textnormal{ for all } n\geq 1 
\end{equation}
where $\tilde{T}_n(x):=T_n(\frac{x}{2})$ and $\tilde{U}_n(x):=U_n(\frac{x}{2})$ are the rescaled versions of the first and second kind of Chebyshev polynomials $T_n$ and $U_n$ for the domain $(-2,2)$. Note that, with this scaling, the orthogonality properties of the Chebyshev polynomials are given by
\begin{align}\label{eqn:cheb_orthogonal}
\int_{-2}^{2}\tilde{U}_n(x)\tilde{U}_m(x)\rhosc (x)\D x=\mathds{1}_{n=m}, \quad \int_{-2}^{2}\tilde{T}_n(x)\tilde{T}_m(x)\rhoarcsin (x)\D x=\mathds{1}_{n=m=0}+\frac{1}{2}\mathds{1}_{n=m\neq 0}.
\end{align}
Furthermore, using the Chebyshev expansion of logarithm, i.e.,
\begin{align}\label{eqn:log_cheb_exp}
\log^Ex=\sum_{n=1}^\infty \frac{-2\tilde{T}_n(E)}{n}\tilde{T}_n(x)
\end{align}
for all $x,E\in(-2,2)$.
$\vhlbrt$ transform of $\log^{E}$ can be calculated easily:
\begin{gather} 
\vhlbrt(\log^E)(x)=\begin{cases}\frac{1}{\sqrt{4-x^2}}(\arccos(\frac{x}{2})-\pi), & x\in(-2,E) \\  
\frac{1}{\sqrt{4-x^2}}\arccos(\frac{x}{2}), & x\in(E,2) \end{cases} \label{eqn:v_trans_of_log}
\end{gather}
for any $E\in(-2,2)$. Similarly, for $\arg^{E}(x)$ we have the Chebyshev expansion
\begin{align} \label{eqn:arg_cheb_exp}
\arg^{E}x&=\big(\frac{\pi}{2}-\frac{E}{2}-\arccos\frac{E}{2}\big)+\tilde{T}_1(x)-\sum_{n=1}^{\infty}\frac{\tilde{U}_{n-1}(E)\sqrt{4-E^2}}{n}\tilde{T}_n(x) 
\end{align}
for all $x,E\in(-2,2)$ and denoting $E=2\cos\alpha$ and $x=2\cos\beta$ with $\alpha,\beta\in(0,\pi)$,
\begin{align} \label{eqn:v_trans_of_arg}
\vhlbrt \arg^{E}(x)=\frac{1}{2}+\frac{1}{\sqrt{4-x^2}}\Big(\log\sin\frac{\alpha-\beta}{2}-\log\sin\frac{\alpha+\beta}{2}\Big).
\end{align}

\subsection{Single-time and multi-time loop equations: A brief exploration}\label{subsec:multitime_loop_asymp_exploration}

We now proceed to describe the single-time loop equation. If $h:\R\to\R$ is differentiable, by doing a change of variables in the integral for the partition function, the following equation (which is called loop equation) is obtained (see \cite[(2.18)]{johansson1998fluctuations}):
\begin{align}\label{eqn:loop_eqn_by_change_of_variables}
\mathbb{E}_h[S_N(\uhlbrt g)]=&\int g(x)h'(x)\rhosc (x)\D x+\frac{1}{N}\mathbb{E}_h[S_N(g h')]+\frac{1}{N}\mathbb{E}_h\Big[\int\int \frac{g(x)-g(y)}{x-y}\D\nu_N(x)\D\nu_N(y)\Big]
\end{align}
for any continuous and piecewise continuously differentiable function $g:\R\to\R$ with $\inf_{x\in\R} g'(x)>-\infty$ where $\D\nu_N( x):=\sum_i \delta_{\lambda_i}(x)-N\rhosc (x)\D x$ is the centered empirical measure and the biased measure is defined as in \eqref{eqn:bias_defn}. By Proposition \ref{prop:inverse_of_uhlbrt} and orthogonality relations in \eqref{eqn:cheb_orthogonal} with the fact that $\tilde{T}_n'(x)=\frac{n}{2}\tilde{U}_n(x)$, a rough interpretation gives
\begin{align} \label{eqn:loop_eqn_approx}
\mathbb{E}_h[S_N(f)]\approx \int_{-2}^{2} \vhlbrt f(x) h'(x) \rhosc (x)\D x=\frac{1}{4}\sum_{n=0}^\infty n\hat{f}_n\hat{h}_n
\end{align}
when $f$ and $h$ are ``nice" functions where $\hat{f}_n$ and $\hat{h}_n$ are Chebyshev coefficients (for the scaled first type of Chebyshev polynomials $\tilde{T}_n$'s), i.e. $f=\sum_{n\geq 0} \hat{f}_n\tilde{T}_n$ on $[-2,2]$ and similarly for $h$. Indeed, it is possible to make this statement rigorous and prove the approximation for some $N$-dependent functions $f$ and $h$ via rigidity under the biased measure. Furthermore, applying the techniques we will use to prove Theorem \ref{thm:FH}, it is possible to establish single-time Fisher-Hartwig asymptotics with merging root-type singularities which generalizes Theorem 1 in \cite{krasovsky2007correlations} and prove the convergence to GMC in the single time case, regenerating the result of \cite{berestycki2018random}. However, the same change of variables strategy cannot be applied in multi-time/dynamical setting. Our purpose in this section is to give a dynamical generalization to the estimation \eqref{eqn:loop_eqn_approx} in order to prove multi-time Fisher-Hartwig asymptotics with the help of the results obtained in preceding sections. Before delving into the calculations, let us first illustrate the anticipated result on a toy model with an optimistic approach, the calculations of which can serve as a useful guide for proving the main theorem of this section.

 The convergence of the fluctuation field corresponding to particle dynamics related to some random matrix models to Gaussian fields has been proven in \cite{spohn1987interacting} and generalized through a series of studies \cite{israelsson2001asymptotic, bender2008global, unterberger2018global} with explicit forms for the covariance structure. In our case, \cite[(3.25)]{unterberger2018global} provides the following explicit covariance kernel: Assume that $f$ and $h$ are some fixed (i.e. $N$-independent), compactly supported smooth functions. As $N$ goes to infinity $(S_N(f(H_t)),S_N(h(H_0)))$ converges to a centered Gaussian vector $(X_{f,t},X_{h,0})$ with
 \begin{align*} 
 \Cov(X_{f,t},X_{h,0})=\int_{-2}^{2}\int_{-2}^{2} f(x) g(x,y;t,0)h(y)\D x\D y
 \end{align*}
where $g(x,y;t,0):=\frac{-1}{32\pi^2}\frac{1}{\sin\theta_1\sin\theta_2}\Re\Big(\frac{1}{\sin^2\frac{\theta_1-\theta_2+it}{2}}+\frac{1}{\sin^2\frac{\theta_1+\theta_2+it}{2}}\Big)$ with $x=2\cos\theta_1$ and $y=2\cos\theta_2$. Simplifying the expression, using $\cos(\theta_1-i\frac{t}{2})=\frac{1}{2}x_t$ (remember the characteristic curve at \eqref{eqn:characteristic_curve}), we obtain that,
\begin{align*} 
g(x,y;t,0)=\frac{-1}{2\pi^2}\Re\Big(\frac{1}{\sqrt{4-x^2}\sqrt{4-y^2}}\frac{4-yx_t}{(y-x_t)^2}\Big)=\frac{1}{2\pi^2}\Re\Big(\frac{1}{\sqrt{4-x^2}}\partial_y\big(\frac{\sqrt{4-y^2}}{y-x_t}\big)\Big).
\end{align*}
Substituting this we get 
 \begin{multline*} 
 \Cov(X_{f,t},X_{h,0}) =  \Re\Big(\int_{-2}^{2}\int_{-2}^{2} f(x)\rhoarcsin (x) h'(y)\frac{\rhosc (y)}{x_t-y}\D x\D y\Big)
 \\
 =\frac{1}{2}\int_{-2}^{2}\Re(\uhlbrt h'(x_t)) f(x) \rhoarcsin (x)\D x=\frac{1}{4} \sum_{k=0}^\infty ke^{-tk/2} \hat{f}_k\hat{h}_k
 \end{multline*}
 where the last equality simply follows from Lemma \ref{lemma:Uhlbrt_of_tildeU_at_z} below and the orthogonality relations \eqref{eqn:cheb_orthogonal}. That means
\begin{align*}
\E{X_{f,t}\frac{e^{X_{h,0}}}{\E{e^{X_{h,0}}}}}=\Cov(X_{f,t},X_{h,0})=\frac{1}{4} \sum_{n=0}^\infty ne^{-tn/2} \hat{f}_n\hat{h}_n.
\end{align*}
Therefore, we expect $\mathbb{E}_{h(H_0)}[S_N(f(H_t))]$ to be approximated by $\frac{1}{4} \sum_{n=0}^\infty ne^{-tn/2} \hat{f}_n\hat{h}_n$. Note that, the exponential decay in the Chebyshev coefficients expansion along the time matches the circular version with Fourier coefficients established in \cite[Lemma 5.6]{bourgade2022liouville}, and can be interpreted as a sign of relaxation phenomena.

In light of this brief exploration, we define time-dependent versions of $\uhlbrt$ and $\vhlbrt$ transforms as follows:
\begin{align*} 
\uhlbrt_t f(x)=\Re\Big(2\int_{-2}^{2}\frac{f(y)}{x_t-y}\rhosc (y)\D y\Big)=2\sum_{n=0}^{\infty}e^{-tn/2}\hat{f}_n^{U} \tilde{T}_{n+1}(x)
\end{align*}
where $\hat{f}_{n}^{U}$ stands for the Chebyshev coefficients of $f$ for the second type of scaled Chebyshev polynomials and
\begin{align*}
\vhlbrt_{t}f(x)=\Re(-\int_{-2}^{2}\frac{f(y)}{x-y_{t}}\rhoarcsin (y)\D y)=\frac{1}{2}\sum_{n=1}^{\infty}  e^{-tn/2}\hat{f}_n\tilde{U}_{n-1}(x).
\end{align*}
The equalities above are proven in the lemma below and when $t=0$ the integrals should be read as principal value integrals, so that $\uhlbrt_0$ and $\vhlbrt_0$ coincides with the definitions of $\uhlbrt$ and $\vhlbrt$. Then, we establish the definition of $\mathcal{C}$, which describes the covariance structure as mentioned above, and introduce various related notations here:
\begin{align} \label{eqn:defn_mathcal_C}
\mathcal{C}(f(H_t),h(H_0)):=\frac{1}{4}\sum_{k=0}^{\infty} e^{-tk/2}k\hat{f}_k\hat{h}_k=\int \vhlbrt_tf(x)h'(x)\rhosc (x)\D x= \int \vhlbrt_th(x) f'(x)\rhosc (x)\D x
\end{align}
for $t\geq 0$ and suitable functions $f$ and $h$ (ensuring either finiteness of the Chebyshev sum or smoothness when considering the integral representations). We extend the definition bilinearly and for simplicity we write $\mathcal{C}(\sum_i f_i(H_{t_i}))$ as shorthand for $\mathcal{C}(\sum_i f_i(H_{t_i}),\sum_i f_i(H_{t_i}))$. Due to the reversibility of the process $H_t$, if more convenient, we also use the notation $\mathcal{C}_{|t-s|}(f,g):=\mathcal{C}(f(H_t),g(H_s))$ and when the time index is $0$, we write $\mathcal{C}(f,g)$ instead of $\mathcal{C}_0(f,g)$.  
 
\begin{lemma}\label{lemma:Uhlbrt_of_tildeU_at_z}
For every $n\in\N$, $x\in[-2,2]$ and $t\geq 0$,
\begin{align*} 
\uhlbrt_t\tilde{U}_{n-1}(x)=\Re(\uhlbrt\tilde{U}_{n-1}(x_t))=2e^{-tn/2}\tilde{T}_{n}(x)
\end{align*}
and for every $x\in(-2,2)$ and $t\geq 0$,
\begin{align*} 
\vhlbrt_{t}\tilde{T}_n(x)=\Re(\int_{-2}^{2}\frac{\tilde{T}_n(y)}{y_{t}-x}\rhoarcsin (y)\D y)=e^{-tn/2}\frac{1}{2}\tilde{U}_{n-1}(x).
\end{align*} 
\end{lemma}
\begin{proof}\renewcommand{\qedsymbol}{} This is a straightforward application of the Cauchy integral formula; see Appendix \ref{app:pf_U_t and V_t}.
\end{proof}

\subsection{An integration by parts formula for biased measure}

We start the section proving a useful integration by parts formula which will be handy in the proof of the multi-time loop equation asymptotics. Given a continuous semimartingale $(M_t)_{t\geq0}$, define $(H(M)_t)_{t\geq0}$ as the strong solution of SDE $$\D H(M)_t=\frac{1}{\sqrt{N}}\D M_t-\frac{1}{2}H(M)_t \D t$$ assuming its existence and uniqueness. Let $(f_t)_t$ be an adapted bounded continuous Hermitian matrix valued process and denote $F_t:=\int_0^t f_s \D s$.

\begin{lemma}\label{derivative of H}
 Let $B$ be a Hermitian Brownian motion as defined before (see equation \eqref{matrix_valued_OU_dynamics}). Then $D_F H_t:=\lim_{\varepsilon\to0}\frac{H(B+\varepsilon F)_t-H(B)_t}{\varepsilon}$ exists and equal to $\frac{1}{\sqrt{N}}\int_0^t e^{-(t-s)/2}f_s\D s$.
\end{lemma}

\begin{proof}
As $f$ is a bounded continuous process, the existence and uniqueness of the strong solution is clear. By subtracting the SDE's of $H(B+\varepsilon F)_t$ and $H(B)_t$,
\begin{align*}
\D(H(B+\varepsilon F)_t-H(B)_t)=\frac{\varepsilon}{\sqrt{N}}f_t\D t-\frac{1}{2}(H(B+\varepsilon F)_t-H(B)_t)\D t
\end{align*}
is obtained. This implies,
\begin{align*}
\D\big(e^{t/2}\frac{H(B+\varepsilon F)_t-H(B)_t}{\varepsilon}\big)=e^{t/2}\frac{f_t}{\sqrt{N}}\D t,\ \ {i.e.}\ \ 
\frac{H(B+\varepsilon F)_t-H(B)_t}{\varepsilon}=\frac{1}{\sqrt{N}}\int_0^t e^{-(t-s)/2}f_s\D s,
\end{align*}
which concludes the proof.
\end{proof}

For the next lemma, we switch between the representations of Hermitian continuous process $H_t$ as a Hermitian matrix and a vector of size $N^2$ given by $(H^{11}_t,\dots,H^{NN}_t,H^{12,R}_t,\dots,H^{(N-1)N,I}_t)^T$ where $H^{ij,R}$ and $H^{ij,I}$ stand for the real and imaginary parts of the $(i,j)^{th}$ entry of $H$ and endow the set $\mathscr{C}([0,T])^{N^2}$ with the Borel sigma algebra given by the sup norm.

\begin{lemma}\label{int by parts for B}
Let $B$ be a Hermitian Brownian motion, $T>0$ fixed and $\Phi:\mathscr{C}([0,T])^{N^2}\to \R$ measurable. If $D_F\Phi(B):=\lim_{\varepsilon\to0}\frac{\Phi(B+\varepsilon F)-\Phi(B)}{\varepsilon}$ exists in both a.s. and $L^2$ sense, then 
\begin{align*}
\E{D_F\Phi(B)}= \E{\Phi(B)\int_0^T \inner{f_s}{\D B_s}_{\Re}}
\end{align*}
where $\inner{A}{B}_{\Re}:=\Re(\Tr(A^*B))$ when $A$ and $B$ are viewed as matrices. (Notice that when $A$ and $B$ are Hermitian, $\Tr(A^*B)\in \R$.)
\end{lemma}

\begin{proof}
Let $Y^\varepsilon_t:=\varepsilon F_t+B_t$, viewed as vectors in $\R^{N^2}$. We will use the multidimensional version of Girsanov theorem given in \cite[Theorem 8.6.4]{oksendal2013stochastic}. We have,
\begin{align*}
\D Y^{\varepsilon}_t=\varepsilon f_t \D t +\theta\D W_t
\end{align*}
where $W_t$ is a standard Brownian motion in $\R^{N^2}$ and $\theta=\textnormal{diag}(\underbrace{1,\cdots,1}_{N\text{ times}},\underbrace{\frac{1}{\sqrt{2}},\cdots,\frac{1}{\sqrt{2}}}_{(N-1)N\text{ times}})$.

 Define $u_t:=(\varepsilon f^{11}_t,\dots,\varepsilon f^{NN}_t,\varepsilon\sqrt{2} f^{12,R}_t,\dots,\varepsilon\sqrt{2} f^{(N-1)N,I}_t)^T$, so that $\theta u_t=\varepsilon f_t$ and $u_t$ satisfies Novikov's condition due to the boundedness of $f_t$. Then by the Girsanov theorem, $\tilde{W}_t=W_t+\int_0^t u_s\D s$ is a Brownian motion with respect to measure $\mathbb{Q}$ defined by $\frac{\D \mathbb{Q}}{\D\Prob}=M$ where $M_t=\exp\big(-\int_0^t u_s\D W_s- \frac{1}{2}\int_0^t|u_s|^2\D s\big)$ and $\D Y^{\varepsilon}_t=\theta \D \tilde{W}_t$. Hence,
\begin{align*}
\E{\Phi(Y^{\varepsilon})}=&\mathbb{E}_\mathbb{Q}[\Phi(Y^{\varepsilon})\exp\Big( \int_0^T u_t \D \tilde{W}_t- \frac{1}{2} \int_0^T|u_t|^2\D t \Big)]
\\
=&\mathbb{E}_\mathbb{Q}[\Phi(\theta \tilde{W})\exp\Big(\varepsilon \int_0^T (f_t \theta^{-1}) \theta^{-1} \D (\theta \tilde{W}_t)- \frac{1}{2} \int_0^T|u_t|^2\D t \Big)]
\\
=&\mathbb{E}_\mathbb{Q}[\Phi(\theta \tilde{W})\exp\Big(\varepsilon \int_0^T (f_t \theta^{-1}) \theta^{-1} \D (\theta \tilde{W}_t)- \frac{1}{2} \varepsilon^2 \int_0^T|\theta^{-1} f_t|^2\D t \Big)]
\\
=&\mathbb{E}[\Phi(B)\exp\Big(\varepsilon \int_0^T (f_t \theta^{-1}) \theta^{-1} \D B_t- \frac{1}{2} \varepsilon^2 \int_0^T|\theta^{-1} f_t|^2\D t \Big)].
\end{align*}
Differentiating with respect to $\varepsilon$ gives
\begin{align*}
\E{D_F\Phi(B)}=\E{\lim_{\varepsilon\to0}\frac{\Phi(B+\varepsilon F)-\Phi(B)}{\varepsilon}}=\lim_{\varepsilon\to0} \E{\frac{\Phi(B+\varepsilon F)-\Phi(B)}{\varepsilon}}=\frac{\D}{\D\varepsilon}\Big|_{\varepsilon=0}\E{\Phi(Y^{\varepsilon})}
\\
=\E{\Phi(B)\int_0^T f_t \theta^{-2} \D B_t}=\E{\Phi(B)\int_0^T \inner{f_s}{\D B_s}_{\Re}}
\end{align*}
and completes the proof.
\end{proof}

\begin{lemma}\label{integrationbypartslemma}
Let $(B_t)_t$ denote a Hermitian Brownian motion, and let $(H_t)_t$ be the solution of the SDE \eqref{matrix_valued_OU_dynamics} driven by this Brownian motion. Then, for any matrix valued bounded continuous adapted process $(g_t)_t$, $J\in\N$, $t_1,\dots,t_J,T\in \R_+$ and $f_1,\dots,f_J\in \mathscr{C}^1_c(\R)$,
\begin{align*}
\mathbb{E}_{\sum_{j=1}^J h_j(H_{t_j})}\Big[ \int_0^T \Tr(g_t\D B_t) \Big]=\frac{1}{\sqrt{N}} \sum_{j=1}^J\mathbb{E}_{\sum_{j=1}^J h_j(H_{t_j})}\Big[\Tr\Big(h_j'(H_{t_j})\int_0^{T\wedge t_j}e^{-(t_j-s)/2}g_s\D s\Big)\Big].
\end{align*}
\end{lemma}

\begin{proof}
Let $p_S(\cdot)$ and $p_H(\cdot)$ be the projection maps into skew-Hermitian and Hermitian matrices respectively. Then
\[
\Re(\Tr(g_t\D B_t))=\Re(\Tr(p_S(g_t)\D B_t))+\Re(\Tr(p_H(g_t)\D B_t))
=-\Re(\Tr(p_S^*(g_t)\D B_t))+\Re(\Tr(p_H(g_t)^*\D B_t))
=\inner{p_H(g_t)}{\D B_t}_\Re
\]
and 
\[
\Im(\Tr(g_t\D B_t))=\Re(-i\Tr(g_t\D B_t))=-\Re(\Tr((ip_S(g_t))\D B_t))-\Re(\Tr((ip_H(g_t))\D B_t))
\\
=-\inner{ip_S(g_t)}{\D B_t}_\Re.
\]
For simplicity consider the case $J=1$ and $t_1\leq T$. Define $\Phi(M)$ as $e^{\Tr(h_1(H(M)_{t_1}))}$ and let $F^{p_H(g)}_t:=\int_0^t p_H(g_s)\D s$, $F^{ip_S(g)}_t:=\int_0^t p_S(g_s)\D s$. Therefore,
\begin{align*}
\mathbb{E}\Big[ \int_0^T \Tr(g_t\D B_t) e^{\Tr(h_1(H_{t_1}))} \Big]=&\mathbb{E}\Big[ \int_0^T \inner{p_H(g_t)}{\D B_t}_\Re e^{\Tr(h_1(H_{t_1}))} \Big]-i\cdot \mathbb{E}\Big[ \int_0^T \inner{ip_S(g_t)}{\D B_t}_\Re e^{\Tr(h_1(H_{t_1}))} \Big]
\\
=&\E{D_{F^{p_H(g)}}\Phi(B)}-i\cdot \E{D_{F^{ip_S(g)}}\Phi(B)}
\\
=& \E{\Phi(B)D_{F^{p_H(g)}}(\Tr(h_1(H_{t_1})))}-i\cdot \E{\Phi(B)D_{F^{ip_S(g)}}(\Tr(h_1(H_{t_1})))}
\\
=& \E{\Phi(B)\Tr\big((h_1'(H_{t_1}))D_{F^{p_H(g)}}H_{t_1}\big)}-i\cdot \E{\Phi(B)\Tr\big(h_1'(H_{t_1})D_{F^{ip_S(g)}}H_{t_1}\big)}
\\
=& \frac{1}{\sqrt{N}} \E{\Phi(B)\Tr\big((h_1'(H_{t_1}))\int_0^{t_1}e^{-(t_1-s)/2}p_H(g_s)\D s\big)}
\\
&-\frac{i}{\sqrt{N}}\cdot \E{\Phi(B)\Tr\big(h_1'(H_{t_1})\int_0^{t_1}e^{-(t_1-s)/2}ip_S(g_s)\D s\big)}
\\
=&\frac{1}{\sqrt{N}} \E{\Phi(B)\Tr\big((h_1'(H_{t_1}))\int_0^{t_1}e^{-(t_1-s)/2}(p_H(g_s)+p_S(g_s))\D s\big)}
\\
=&\frac{1}{\sqrt{N}} \E{\Tr\Big((h_1'(H_{t_1}))\int_0^{t_1}e^{-(t_1-s)/2}g_s\D s\Big) e^{\Tr(h_1(H_{t_1}))} }
\end{align*}
where we have used Lemma \ref{int by parts for B} at the second equality and Lemma \ref{derivative of H} at the fifth equality.

For the $t_1>T$ case, extend $g$ up to time $t_1$ by concatenating zero function. More explicitly, $\tilde{g}_s=g_s$ if $s\in[0,T]$ and $\tilde{g}_s=0$ if $s>T$. We lose the continuity of $\tilde{g}$ at $T$ but everything still works and we obtain 
\begin{align*}
\mathbb{E}\Big[ \int_0^T \Tr(g_t\D B_t) e^{\Tr(h_1(H_{t_1}))} \Big]=&\frac{1}{\sqrt{N}} \E{\Tr\Big(h_1'(H_{t_1})\int_0^{t_1}e^{-(t_1-s)/2}\tilde{g}_s\D s\Big) e^{\Tr(h_1(H_{t_1}))} }
\\
=&\frac{1}{\sqrt{N}} \E{\Tr\Big(h_1'(H_{t_1})\int_0^{T}e^{-(t_1-s)/2}g_s\D s\Big) e^{\Tr(h_1(H_{t_1}))} }.
\end{align*}

Generalization to $J>1$ follows from the same steps with $\Phi(M):=\prod_{j=1}^J e^{\Tr(h_j(H(M)_{t_j}))}$.
\end{proof}

\subsection{Martingale approximation for $m_{u}(z_{t-u})-m_{\rm sc}(z_{t-u})$}

In the asymptotic analysis of multi-time loop equation for unitary Brownian motion, as derived in \cite[Lemma 5.6]{bourgade2022liouville}, the authors' approach for examining $\Ex_{h(H_0)}[S_N(f(H_t))]$ involves first applying the Helffer-Sjöstrand formula, 
\begin{align*} 
\Ex_{h(H_0)}[S_N(f(H_t))]=\frac{N}{\pi} \int_{\C} \partialbar\tilde{f}(z) \Ex_{h(H_0)}[m_t(z)-m_{\rm sc}(z)]\D z
\end{align*}
followed by using the SDE for $m_{u}(z_{t-u})$ to substitute the term $m_t(z)-m_{sc(z)}$. This reduces the problem to controlling the martingale term, while the other terms become negligible \cite[(5.10)-(5.12)]{bourgade2022liouville}. However, the negligibility of (5.11) depends on approximating $\partial_zm_u(z_{t-u})$ by $m_{\rm sc}'(z_{t-u})$ which vanishes in their setting. In more general models, where $m_{\rm sc}'$ is typically of order one, this approach fails. The purpose of this section is to refine the argument by incorporating the non-negligible contribution of $\partial_zm_u(z_{t-u})$ into the martingale term through an adjusted SDE approximation. This modification extends the applicability of their technique to a more general setting.

Recall the equation \eqref{SDE along char for tilde m_t,A}. Defining 
\begin{align*}
p_u=p_u(z):=\tilde{m}_u(z_{t-u})-e^{-u/2}m_{\rm sc}(z_{t-u})
\end{align*}
we can write \eqref{SDE along char for tilde m_t,A} as,
\begin{align*} 
\D p_u=\Tr(r_u\D B_u)+p_u\partial_{z}m_{u}(z_{t-u})\D u
\end{align*}
where we refer to the equation \eqref{SDE along char for tilde m_t,A} for the definition of $r_u$ and $B_u$. To approximate $p_u$ in terms of $(r_u)$ and $(B_u)$, we first show that $\partial_{z}m_{u}(z_{t-u})$ term in the SDE can be replaced by $m_{\rm sc}'(z_{t-u})$ up to a negligible error due to rigidity. This substitution simplifies the solution of the SDE, allowing us to derive an estimate for $p_u$.

Define the rigidity set for the trajectories of the particles from time $-(\log N)^2$ to $(\log N)^2$:
\begin{align*} 
\boldsymbol{\mathcal{G}}:=\{\boldsymbol\lambda(s)\in\mathcal{G},\textnormal{ for all } s\in[-(\log N)^2,(\log N)^2]\}
\end{align*} 
where $\mathcal{G}$ is as defined in \eqref{eqn:rigidity sets}. Recall that Proposition \ref{prop:rig_until_logN2} implies that $\boldsymbol{\mathcal{G}}$ is an event with overwhelming probability. The next lemma offers a very fine estimate on $p_t$ on the set $\boldsymbol{\mathcal{G}}$ via the method described above.

\begin{lemma}\label{lemma:p_t estimation} Let $C>0$ be a constant and $t\in[0,C]$ be fixed. Define $p_u$ as above and recall the notation $\upkappa=\min(|E+2|,|E-2|)$. For all sufficiently large $N$ we have, on the set $\boldsymbol{\mathcal{G}}$,
\begin{align*} 
\big|p_t-\int_{-T}^t \Tr(e^{\int_{u}^t m_{\rm sc}'(z_{t-r})\D r} r_u\D B_u)\big|<\frac{\varphi^{24}}{N^2 \eta^2\Im\sqrt{z^2-4}|\sqrt{z^2-4}|}\leq \begin{cases}
\frac{\varphi^{25}}{N^2 \eta^2(\eta+\upkappa)} &, |E|<2\\
\frac{\varphi^{25}}{N^2 \eta^3} &, |E|>2
\end{cases}
\end{align*}
for all $T\in[\log N,(\log N)^2]$ and $z\in\C$ with $\eta\in(0,C)$, $|E|<C$.
\end{lemma}

\begin{proof}
 Fix an arbitrary $z\in\C$ with $|E|\in[0,C],\eta\in(0,C]$. Let $q_u$ be the strong solution of the following SDE:
\begin{align*} 
\D q_u=\Tr(r_u\D B_u)+q_u m_{\rm sc}'(z_{t-u})\D u
\end{align*}
with the same initial condition with $p_u$ at time $-1$ (i.e. $p_{-1}=q_{-1}$). The process $p_u-q_u$ satisfies
\begin{align*} 
\D(p_u-q_u)=p_u\cdot(\partial_{z}m_{u}(z_{t-u})-m_{\rm sc}'(z_{t-u}))\D u+(p_u-q_u)m_{\rm sc}'(z_{t-u})\D u.
\end{align*}
On the set of rigidity $\boldsymbol{\mathcal{G}}$ we have
\begin{align*} 
\big|p_u\cdot(\partial_{z}m_{u}(z_{t-u})-m_{\rm sc}'(z_{t-u}))\big|\lesssim \begin{cases} 
\frac{\varphi^{22}}{N^{2}\eta_{t-u}^{3}} &,u\in[t-1,t]
\\
\frac{\varphi^{22}}{N^{2}e^{3(t-u)/2}} &, u\in [(-\log N)^2, t-1]
\end{cases}.
\end{align*}
Hence, 
\begin{align*} 
&\int_{-1}^t \big|p_u\cdot(\partial_{z}m_{u}(z_{t-u})-m_{\rm sc}'(z_{t-u}))\big|\D u\leq \frac{\varphi^{23}}{N^2 \eta^2\Im\sqrt{z^2-4}},\\
&\int_{-T}^{-1} \big|p_u\cdot(\partial_{z}m_{u}(z_{t-u})-m_{\rm sc}'(z_{t-u}))\big|\D u\leq \frac{\varphi^{23}}{N^2 }.
\end{align*}
By Gr\"onwall's lemma, we obtain
\begin{align*} 
&|p_t-q_t|\leq  \frac{\varphi^{23}}{N^2 \eta^2\Im\sqrt{z^2-4}}\cdot \exp\Big(|m_{\rm sc}(z)| \big|\int_{-1}^{t}e^{-(t-u)/2}\frac{1}{|\sqrt{z_{t-u}^2-4}|}\D u\big| \Big)\\ 
&|p_{-T}-q_{-T}|\leq  \frac{\varphi^{23}}{N^2}\cdot \exp\Big(|m_{\rm sc}(z)| \big|\int_{-T}^{-1}e^{-(t-u)/2}\frac{1}{|\sqrt{z_{t-u}^2-4}|}\D u\big| \Big).
\end{align*}
The exponential terms can be evaluated as,
\begin{align*} 
\int_{-\infty}^{t-1} \frac{e^{-(t-u)/2}}{|\sqrt{z_{t-u}^2-4}|}\D u\leq\int_{1}^{\infty}\frac{e^{-s/2}}{\frac{e^{s/2}-e^{-s/2}}{2}(|z+\sqrt{z^2-4}|)}\D s<2
\end{align*}
and
\begin{align*} 
\int_{t-1}^{t} \frac{e^{-(t-u)/2}}{|\sqrt{z_{t-u}^2-4}|}\D u\leq \int_{0}^{1} \frac{1}{\frac{s}{2}|z|+|\sqrt{z^2-4}|}\D s =\frac{2}{|z|} \log\Big(1+\frac{|z|/2}{|\sqrt{z^2-4}|}\Big).
\end{align*} 
Combining these with $|m_{\rm sc}(z)|\leq 1$ we obtain that
\begin{align*} 
|p_t-q_t|\lesssim  \frac{\varphi^{23}}{N^2 \eta^2\Im\sqrt{z^2-4}|\sqrt{z^2-4}|}, \quad |p_{-T}-q_{-T}|\lesssim \frac{\varphi^{23}}{N^2}.
\end{align*}

Now we solve the SDE for $q$ and use this approximation between processes $p$ and $q$. If we define $Q_u:=e^{-\int_{0}^{u}m_{\rm sc}'(z_{t-r})\D r}q_u$, it satisfies the following SDE:
\begin{align*} 
\D Q_u=\Tr\Big(e^{-\int_{0}^{u} m_{\rm sc}'(z_{t-r})\D r}r_u\D B_u\Big).
\end{align*}
Hence,
\begin{align*} 
q_t=Q_te^{\int_{0}^{t} m_{\rm sc}'(z_{t-r})\D r}=e^{\int_{0}^{t} m_{\rm sc}'(z_{t-r})\D r}\Big(Q_{-T}+\int_{-T}^t \Tr\Big(e^{-\int_{0}^{u} m_{\rm sc}'(z_{t-r})\D r}r_u\D B_u\Big) \Big)
\\
=e^{\int_{-T}^{t} m_{\rm sc}'(z_{t-r})\D r}q_{-T}+\int_{-T}^t \Tr\Big(e^{\int_{u}^{t} m_{\rm sc}'(z_{t-r})\D r}r_u\D B_u\Big) \D u
\end{align*}
By Lemma \ref{lemma:exp_int_msc'}, $|e^{\int_{-T}^{t} m_{\rm sc}'(z_{t-r})\D r}|=e^{-(t+T)/2}|\frac{\sqrt{z_{t+T}^2-4}}{\sqrt{z^2-4}}|\lesssim \frac{1}{|\sqrt{z^2-4}|}$. Moreover, by rigidity, $|p_{-T}|< \frac{\varphi^{11}}{N^2}$
Therefore
\begin{align*} 
|p_t-\int_{-T}^t \Tr\Big(e^{\int_{u}^{t} m_{\rm sc}'(z_{t-r})\D r}r_u\D B_u\Big) \D u|\leq |p_t-q_t|+e^{\int_{-T}^{t} m_{\rm sc}'(z_{t-r})\D r}(|p_{-T}|+|p_{-T}-q_{-T}|)
\\
\leq  \frac{\varphi^{24}}{N^2 \eta^2\Im\sqrt{z^2-4}|\sqrt{z^2-4}|}.
\end{align*}
 
 And the second inequality in the statement of the lemma simply follows from equation \eqref{eqn:basic_properties_of_m_sc}.
\end{proof}

\begin{lemma}\label{lemma:exp_int_msc'} For all $z\in \C$, and $t,u\in\R$:
\begin{align*}
\exp\Big(\int_{-u}^{t}m_{\rm sc}'(z_{t-r})\D r\Big)=e^{-(t+u)/2}\frac{\sqrt{z_{t+u}^2-4}}{\sqrt{z^2-4}}.
\end{align*}
\end{lemma}
\begin{proof}
Note that 
\begin{align*} 
\frac{\D}{\D u} \Big(\frac{\exp\Big(\frac{u}{2}+\int_{-u}^{t}m_{\rm sc}'(z_{t-r})\D r\Big)}{\sqrt{z_{t+u}^2-4}}\Big)=\exp\Big(\frac{u}{2}+\int_{-u}^{t}m_{\rm sc}'(z_{t-r})\D r\Big)\frac{(\frac{1}{2}+m_{\rm sc}'(z_{t+u}))\sqrt{z_{t+u}^2-4}-\frac{z_{t+u}}{2}}{z_{t+u}^2-4}=0
\end{align*}
and evaluating the function at $u=-t$ completes the proof.
\end{proof}

\subsection{Multi-time loop equation asymptotics}\label{subsec:multitime_loop_asymp}

Having all the necessary tools, let us discuss the main theorem of this section. Below, we formulate the theorem in a general manner; however, for our purposes the key example to keep in mind is when the bias $h=\log^E_{\Delta}$ where $E$ is a point in the bulk and $\Delta=N^{-1-\alpha}$ for some small $\alpha>0$.

\begin{theorem}\label{thm:multitime_loopp}
Let $C>0$, $\kappa\in(0,1]$, $\Upsilon>0$ be fixed. Take $t\in[0,C]$, $B=B_N\in[-2+\Upsilon,2-\Upsilon]$. Let $f\in\mathscr{S}_{C,\kappa}$ (recall the definition \eqref{eqn:defn_mathscr_S}), $h=h_1+h_2$ such that $h_1\in\mathscr{S}_{C,\kappa/2}$ and $h_2\in\mathscr{S}_{C,(-\kappa/100)}$ with $\supp(h_2)\subseteq[B-CN^{-1+\kappa/2},B+CN^{-1+\kappa/2}]$. If $\vhlbrt_t f(B)=O(N^{-\kappa/2})$, then,
\begin{align*} 
\mathbb{E}_{h(H_0)}[S_N(f(H_t))]=\mathcal{C}(f(H_t),h(H_0))+O(N^{-\kappa/5})
\end{align*}
where $\mathcal{C}$ term is defined in \eqref{eqn:defn_mathcal_C}. Moreover, the error term is uniform over the choice of $t$, $B$, $f$ and $h$ satisfying the given conditions.
\end{theorem}

\begin{proof} By Proposition \ref{prop:rig_for_mathscr}, $h$ satisfies the rigidity conditions. Take an arbitrary sequence of functions $g=g_N$ such that $g_N\in \mathscr{A}_{C,\epsilon_N}$ for a sequence $\epsilon=\epsilon_N\in[N^{-1+\kappa},1]$ and let $\Delta=N^{-1-\kappa/100}$, $\varphi=(\log N)^{\log\log N}$. By Helffer-Sj\"ostrand formula (Proposition \ref{prop:HS}),
\begin{align} \label{eqn:multi-loop-HS}
\Ex_{h(H_0)}[S_N(g(H_t))]=\frac{N}{\pi} \int_{\C} \partialbar\tilde{g}(z) \Ex_{h(H_0)}[m_t(z)-m_{\rm sc}(z)]\D z
\end{align}
where we use the third order quasi-analytic extension this time, i.e.
\begin{align*} 
\tilde{g}(z)=\big(g(x)+iyg'(x)-\frac{y^2}{2}g''(x)-i\frac{y^3}{6}g'''(x)\big)\chi_{\xi}(y)
\end{align*}
 with $\varepsilon=\frac{\kappa}{100}$ and $\xi=\epsilon N^{-\varepsilon}$. Note that this order of quasi-analytic extension provides $\int_{\C}|\frac{\partialbar\tilde{g}(z)}{y^a}|\D z\lesssim  \frac{\epsilon}{\xi^{a}}$ for all $a<4$ (and if also $1\leq a$  the right hand side is bounded by $N^{a-1}N^{a\varepsilon -(a-1)\kappa}$).
 
  Recall the definition of $\mathcal{S}_{\alpha}(E)$ from equation \eqref{eqn:defn_of_s_curve}. Evaluated on the domain $|\eta|\leq \mathcal{S}_{N^{\varepsilon}}(E)$, the integral in \eqref{eqn:multi-loop-HS} becomes negligible by rigidity.  Because if $\epsilon<\frac{\Upsilon}{2C}$:
 \begin{align*} 
 N \int_{|\eta|\leq\mathcal{S}_{N^{\varepsilon}}(E)} \partialbar\tilde{g}(z) \Ex_{h(H_0)}[m_t(z)-m_{\rm sc}(z)]\D z\lesssim \varphi^{11}\int_{|\eta|\lesssim N^{-1+\varepsilon}} |\frac{\partialbar\tilde{g}(z)}{\eta}|\D z\lesssim \varphi^{11}\epsilon^{-3}(N^{-1+\varepsilon})^{3}=O(N^{-2\kappa})
 \end{align*}
 due to rigidity and if $\epsilon\geq\frac{\Upsilon}{2C}$, $g$ is an order $1$ function, so
 \begin{align*} 
 N \int_{|\eta|\leq\mathcal{S}_{N^{\varepsilon}}(E)} \partialbar\tilde{g}(z) \Ex_{h(H_0)}[m_t(z)-m_{\rm sc}(z)]\D z\lesssim \varphi^{11}\int_{|\eta|\lesssim N^{-1/2+\varepsilon}} |\frac{\partialbar\tilde{g}(z)}{\eta}|\D z\lesssim \varphi^{11}(N^{-1/2+\varepsilon})^3=O(N^{-1}).
 \end{align*}
Thus, we write
\begin{align*} 
\Ex_{h(H_0)}[S_N(g(H_t))]=\frac{N}{\pi} \int_{|\eta|> \mathcal{S}_{N^{\varepsilon}}(E)} \partialbar\tilde{g}(z) \Ex_{h(H_0)}[m_t(z)-m_{\rm sc}(z)]\D z+O(N^{-\kappa}).
\end{align*}

Let $T=4\log N$. Using the notations $p_t$, $r_t$ and $\boldsymbol{\mathcal{G}}$ from Lemma \ref{lemma:p_t estimation}, we have
\begin{align*} 
\int_{|\eta|> \mathcal{S}_{N^{\varepsilon}}(E)} \partialbar\tilde{g}(z) \Ex_{h(H_0)}[m_t(z)-m_{\rm sc}(z)]\D z=& \ e^{t/2}\int_{|\eta|> \mathcal{S}_{N^{\varepsilon}}(E)} \partialbar\tilde{g}(z) \Ex_{h(H_0)}[\int_{-T}^t \Tr(e^{\int_{u}^t m_{\rm sc}'(z_{t-r})\D r} r_u\D B_u)]\D z
\\
&+e^{t/2}\int_{|\eta|> \mathcal{S}_{N^{\varepsilon}}(E)} \partialbar\tilde{g}(z) \Ex_{h(H_0)}[\mathds{1}_{\boldsymbol{\mathcal{G}}^c}\cdot(m_t(z)-m_{\rm sc}(z))]\D z
\\
&-e^{t/2}\int_{|\eta|> \mathcal{S}_{N^{\varepsilon}}(E)} \partialbar\tilde{g}(z) \Ex_{h(H_0)}[\mathds{1}_{\boldsymbol{\mathcal{G}}^c}\cdot\int_{-T}^t \Tr(e^{\int_{u}^t m_{\rm sc}'(z_{t-r})\D r} r_u\D B_u)]\D z
\\
&+O\Big(\int_{|\eta|> \mathcal{S}_{N^{\varepsilon}}(E)} |\partialbar\tilde{g}(z)| \Ex_{h(H_0)}\big[\mathds{1}_{\boldsymbol{\mathcal{G}}}\cdot|p_t-\int_{-T}^t \Tr(e^{\int_{u}^t m_{\rm sc}'(z_{t-r})\D r} r_u\D B_u)|\big]\D z\Big)
\end{align*}
The second, third and fourth terms in the right hand side are negligible due to rigidity and the estimation from Lemma \ref{lemma:p_t estimation}. Indeed, the second term is clear and the third term can be bounded by Ito isometry:
\begin{align*} 
\Big|\Ex_{h(H_0)}[\mathds{1}_{\mathcal{G}^c}\cdot\int_{-T}^t \Tr(e^{\int_{u}^t m_{\rm sc}'(z_{t-r})\D r} r_u\D B_u)]\Big|\leq  \Ex_{h(H_0)}[\Big|\int_{-T}^t \Tr(e^{\int_{u}^t m_{\rm sc}'(z_{t-r})\D r} r_u\D B_u)\Big|^2]^{1/2}\Prob(\boldsymbol{\mathcal{G}}^c)^{1/2}
\lesssim \Prob(\boldsymbol{\mathcal{G}}^c)^{1/2}\frac{N}{\eta^{5/2}}
\end{align*}
where we have used Lemma \ref{lemma:exp_int_msc'}. For the fourth term, we consider the expression in two cases. If $\epsilon<\frac{\Upsilon}{2C}$, then the function $g$ is supported in $[-2+\Upsilon/2,2-\Upsilon/2]$. Hence, the fourth term's contribution is bounded by 
\begin{align*} 
N\int_{\C} |\partialbar\tilde{g}(z) \frac{\varphi^{25}}{N^2\eta^2}|\D z\lesssim \varphi^{25}\frac{1}{N} N N^{2\varepsilon-\kappa}\leq N^{-\kappa+3\varepsilon}.
\end{align*}
On the other hand, if $\epsilon\geq \frac{\Upsilon}{2C}$, $g$ is an order $1$ function, so the contribution from the fourth term can be easily bounded by 
\begin{align*} 
N\int_{\C} |\partialbar\tilde{g}(z) \frac{\varphi^{25}}{N^2\eta^3}|\D z\lesssim \varphi^{25}\frac{1}{N}N^{3\varepsilon}\leq N^{-\kappa+4\varepsilon}.
\end{align*}

Now we move on the main term. In order to evaluate it, we first use the integration by parts formula from Lemma \ref{integrationbypartslemma}:
\begin{align*}
\mathbb{E}_{h(H_0)}\Big[ \int_{-T}^t \Tr(\frac{e^{-u/2+\int_{u}^t m_{\rm sc}'(z_{t-r})\D r}}{N^{3/2}}\frac{1}{(H_u-z_{t-u})^2}\D B_u) \Big]=\frac{1}{N^2} \mathbb{E}_{h(H_0)}\Big[\Tr\Big(h'(H_{0})\int_{-T}^{0}\frac{e^{\int_{u}^t m_{\rm sc}'(z_{t-r})\D r}}{(H_u-z_{t-u})^2}\D u\Big)  \Big]
\end{align*}
which leads to
\begin{align*} 
\mathbb{E}_{h(H_0)}[S_N(g(H_t))]=\frac{-e^{t/2}}{N\pi}\int_{|\eta|> \mathcal{S}_{N^{\varepsilon}}(E)}\partialbar\tilde{g}(z)\mathbb{E}_{h(H_0)}\Big[\Tr\Big(h'(H_{0})\int_{-T}^{0}\frac{e^{\int_{u}^t m_{\rm sc}'(z_{t-r})\D r}}{(H_u-z_{t-u})^2}\D u\Big)  \Big]\D z+O(N^{-\kappa/2})
\end{align*}
Applying Helffer-Sj\"ostrand formula again, for $h'$ term this time, we obtain
\begin{align*} 
\mathbb{E}_{h(H_0)}[S_N(g(H_t))]=\frac{-e^{t/2}}{N\pi^2} \int_{|\eta|> \mathcal{S}_{N^{\varepsilon}}(E)}  \int_{\C} \partialbar\tilde{g}(z)  \partialbar\widetilde{h'}(w)\int_{-T}^{0}\mathbb{E}_{h(H_0)}\Big[\Tr\Big(\frac{1}{H_0-w}\frac{e^{\int_{u}^t m_{\rm sc}'(z_{t-r})\D r}}{(H_u-z_{t-u})^2}\Big)  \Big] \D u\D w\D z+O(N^{-\kappa/2}).
\end{align*}
where $\widetilde{h'}=(h'(x)+iyh''(x))\chi_{\Delta}(y)$ is the first order quasi-analytic extension of $h'$. 

Moreover, by the reversibility of the stationary OU process $(H_t)_{t\in\R}$, $\int_{-T}^{0}\mathbb{E}_{h(H_0)}\Big[\Tr\Big(\frac{1}{H_0-w}\frac{e^{\int_{u}^t m_{\rm sc}'(z_{t-r})\D r}}{(H_u-z_{t-u})^2}\Big)  \Big] \D u=\int_{0}^{T}\mathbb{E}_{h(H_0)}\Big[\Tr\Big(\frac{1}{H_0-w}\frac{e^{\int_{-u}^t m_{\rm sc}'(z_{t-r})\D r}}{(H_u-z_{t+u})^2}\Big)  \Big] \D u$ and by Lemma \ref{lemma:two_times_resolvent_mult}:
\begin{align*} 
\mathbb{E}_{h(H_0)}\Big[\Tr\Big(\frac{1}{H_0-w}\frac{1}{(H_u-z_{t+u})^2}\Big)  \Big] =&\ e^{u/2}\frac{\sqrt{z_{t+2u}^2-4}}{\sqrt{z_{t+u}^2-4}}\mathbb{E}_{h(H_0)}\Big[\Tr\Big(\frac{1}{H_{0}-w}\frac{1}{(H_0-z_{t+2u})^2}\Big) \Big]+O\Big(|\frac{e^{u/2}\mathcal{E}(z_{t+u},w,u)}{\eta_{t+u}}|\Big)
\\
&+O(|\frac{N}{\eta_w\eta_{t+u}^2}e^{-(\log N)^{10}}|).
\end{align*}
given that $1>|\eta|> \mathcal{S}_{N^{\varepsilon}}(E)$ and $|E|\leq C+2$ (note that outside of this domain $\tilde{g}$ is already $0$). It is easy to see that contribution from the last term is negligible. The contribution from the error term $\frac{e^{u/2}\mathcal{E}(z_{t+u},w,u)}{\eta_{t+u}}$ is also negligible but it will be dealt with at the end of the proof. 

Thus, we have reduced the problem into single-time, so, we can omit the time index in $H_0$. Now we continue analyzing the main term which is:
\begin{align*}
\frac{-e^{t/2}}{N\pi^2} \int_{|\eta|> \mathcal{S}_{N^{\varepsilon}}(E)}  \int_{\C}\int_{0}^{T} \partialbar\tilde{g}(z)  \partialbar\widetilde{h'}(w)\frac{e^{u/2+\int_{-u}^t m_{\rm sc}'(z_{t-r})\D r}}{\sqrt{z_{t+u}^2-4}} \sqrt{z_{t+2u}^2-4} \cdot \mathbb{E}_{h}\Big[\Tr\Big(\frac{1}{H-w}\frac{1}{(H-z_{t+2u})^2}\Big)\Big]  \D u\D w\D z
\\
=\frac{-e^{t/2}}{N\pi} \mathbb{E}_{h}\Big[\Tr\Big( h'(H)\int_{|\eta|> \mathcal{S}_{N^{\varepsilon}}(E)} \partialbar\tilde{g}(z)   \int_{0}^{T} \frac{e^{u/2+\int_{-u}^t m_{\rm sc}'(z_{t-r})\D r}}{\sqrt{z_{t+u}^2-4}} \frac{\sqrt{z_{t+2u}^2-4}}{(H-z_{t+2u})^2}\D u\D z \Big)\Big]
\\
=\frac{-1}{N\pi} \mathbb{E}_{h}\Big[\Tr\Big( h'(H)\int_{|\eta|> \mathcal{S}_{N^{\varepsilon}}(E)} \partialbar\tilde{g}(z)\frac{1}{\sqrt{z^2-4}}   \int_{0}^{T}  \frac{\sqrt{z_{t+2u}^2-4}}{(H-z_{t+2u})^2}\D u\D z \Big)\Big]
\\
=\frac{1}{N\pi} \mathbb{E}_{h}\Big[\Tr\Big( h'(H)\int_{|\eta|> \mathcal{S}_{N^{\varepsilon}}(E)} \partialbar\tilde{g}(z)\frac{1}{\sqrt{z^2-4}}  \Big(\frac{1}{H-z_t}-\frac{1}{H-z_{t+2T}}\Big)\D z \Big)\Big]
\end{align*}
where we have used Helffer-Sj\"ostrand formula again in the second line, Lemma \ref{lemma:exp_int_msc'} in the third line and $\frac{\D(z_{t+2u})}{\D u}=\sqrt{z_{t+2u}^2-4}$ in the last line. The contribution from the term with $\frac{1}{H-z_{t+2T}}$ is negligible as $T=4\log N$ is sufficiently large. Indeed, it follows by
\begin{align*} 
\frac{1}{N}\int_{\C} \Big|\partialbar\tilde{g}(z)\frac{1}{\sqrt{z^2-4}} \frac{1}{s-z_{t+2T}}\Big|\D z\lesssim \frac{1}{N}\int_{\C} \big|\frac{\partialbar\tilde{g}(z)}{\eta} e^{-T}\big|\D z\lesssim N^{-4}
\end{align*}
and by rigidity 
\begin{align*} 
\Ex_{h}[\Tr(|h'|)]=\Ex_{h}[S_N(|h'|)]+N\int |h'(s)|\rhosc (s)\D s\lesssim \varphi^{11} \int|h''(x)|\D x+\varphi N\lesssim \varphi^{11} \frac{1}{\Delta}.
\end{align*}

Therefore, we obtained that up to negligible error $\mathbb{E}_{h(H_0)}[S_N(g(H_t))]$ is equal to 
\begin{align*} 
\frac{1}{N\pi} \mathbb{E}_{h}\Big[\Tr\Big( h'(H)\int_{|\eta|> \mathcal{S}_{N^{\varepsilon}}(E)} \partialbar\tilde{g}(z)\frac{1}{\sqrt{z^2-4}} \frac{1}{H-z_t}\D z \Big)\Big].
\end{align*}
Now we can add the region $|\eta|< \mathcal{S}_{N^{\varepsilon}}(E)$ back. If $\epsilon<\frac{\Upsilon}{2C}$, for every $s\in\R$,
\begin{align*} 
\Big|\int_{|\eta|<\mathcal{S}_{N^{\varepsilon}}(E)} \partialbar\tilde{g}(z)\frac{1}{\sqrt{z^2-4}} \frac{1}{s-z_t}\D z\Big|\lesssim \int_{|\eta|\lesssim N^{-1+\varepsilon}} |\frac{\partialbar\tilde{g}(z)}{\eta}|\D z\lesssim \epsilon^{-3} (N^{-1+\varepsilon})^3=O(N^{-2\kappa})
\end{align*}
and if $\epsilon\geq\frac{\Upsilon}{2C}$, $g$ is an order $1$ function, so
 \begin{align*} 
\Big|\int_{|\eta|<\mathcal{S}_{N^{\varepsilon}}(E)} \partialbar\tilde{g}(z)\frac{1}{\sqrt{z^2-4}} \frac{1}{s-z_t}\D z\Big|\lesssim \int_{|\eta|\lesssim N^{-1/2+\varepsilon}} |\frac{\partialbar\tilde{g}(z)}{\eta^{2}}|\D z\lesssim (N^{-1/2+\varepsilon})^{2}=O(N^{-1+2\varepsilon}).
 \end{align*}
 Hence, we obtain
 \begin{align*} 
 \mathbb{E}_{h(H_0)}[S_N(g(H_t))]=\frac{1}{N\pi}\mathbb{E}_{h}\Big[\Tr\Big( h'(H)\int_{\C} \partialbar\tilde{g}(z)\frac{1}{\sqrt{z^2-4}} \frac{1}{H-z_t}\D z \Big)\Big]+O(N^{-\kappa/2}).
 \end{align*}

Before using Lemma \ref{lemma:GreentoHS} to simplify the integral term, it is more convenient to get rid off the edge at this stage of the calculations. Let $\breve{\chi}$ be an smooth bump function which is $1$ on $[-2+N^{-1/2},2-N^{-1/2}]$ and $0$ on $[-2+N^{-1/2}/2,2-N^{-1/2}/2]$. Note that for every $s\in\R$ we have $|\int_{\C} \partialbar\tilde{g}(z)\frac{1}{\sqrt{z^2-4}} \frac{1}{s-z_t}\D z|\leq \int_{\C}| \frac{\partialbar\tilde{g}(z)}{\eta}|\D z\leq N^{\varepsilon}$ if $\epsilon<\frac{\Upsilon}{2C}$ and $|\int_{\C} \partialbar\tilde{g}(z)\frac{1}{\sqrt{z^2-4}} \frac{1}{s-z_t}\D z|\leq \int_{\C}| \frac{\partialbar\tilde{g}(z)}{\eta^2}|\D z\leq N^{2\varepsilon}$ otherwise. Moreover, on the support of $1-\breve{\chi}$, $|h'|$ is bounded by $\varphi$. Hence,
\begin{align*} 
\frac{1}{N\pi} \mathbb{E}_{h}\Big[\Tr\Big( h'(H)(1-\breve{\chi}(H))\int_{\C} \partialbar\tilde{g}(z)\frac{1}{\sqrt{z^2-4}} \frac{1}{H-z_t}\D z \Big)\Big]\leq \frac{\varphi}{N^{1-2\varepsilon}\pi}\mathbb{E}_{h}\Big[\Tr( 1-\breve{\chi}(H))\Big].
\end{align*}
On the set of rigidity, the trace is $O(\varphi N^{-1/2}(N\rhosc (-2+N^{-1/2})))=O(N^{1/4})$. Hence the difference is negligible and we have
\begin{align*} 
\mathbb{E}_{h(H_0)}[S_N(g(H_t))]& = \frac{1}{N\pi} \mathbb{E}_{h}\Big[\Tr\Big( h'(H)\breve{\chi}(H)\int_{\C} \partialbar\tilde{g}(z)\frac{1}{\sqrt{z^2-4}} \frac{1}{H-z_t}\D z \Big)\Big]+O(N^{-\kappa/4})
\\
&=\frac{1}{N} \mathbb{E}_{h}\big[\Tr( h'\cdot\breve{\chi}\cdot\vhlbrt_t g)\big]+O(N^{-\kappa/4})
\\
&=\frac{1}{N} \mathbb{E}_{h}\big[S_N( h'\cdot\breve{\chi}\cdot\vhlbrt_t g)\big]+\int_{-2}^{2} h'(s)\vhlbrt_tg(s) \rhosc (s)\D s+O(N^{-\kappa/4})
\end{align*}
where we have now used Lemma \ref{lemma:GreentoHS} in the second line and the following equality in the last:
\begin{align*}
\int_{-2}^{2} h'(s)\breve{\chi}(s)\vhlbrt_tg(s) \rhosc (s)\D s=\int_{-2}^{2} h'(s)\vhlbrt_tg(s) \rhosc (s)\D s+O(N^{-\kappa/4})
\end{align*}
which can be proven as follows. By \cite[Theorem 4.2]{trefethen2008gauss} we have $|\hat{g}_k|\lesssim  \min(\epsilon,\frac{1}{\epsilon^2 k^3})$ as $g\in\mathscr{A}_{C,\epsilon}$. We also have the following bound on the second type of Chebyshev polynomial $\tilde{U}_{k-1}(x)\lesssim \min(k,\frac{1}{\sqrt{|x-2|\wedge|x+2|}})$ uniformly in $x\in[-2,2]$ which is an immediate consequence of the form $\tilde{U}_{k-1}(2\cos\theta)=\frac{\sin(k\theta)}{\sin\theta}$. Hence using $\sup_{|x|>2-N^{-1/2}}|h'(x)|<\varphi$ we can conclude,
\begin{align*} 
\Big|\int h'(s)\vhlbrt_tg(s)(1-\breve{\chi}(s))\rhosc (s)\D s\Big| \lesssim& \  \int_{-2}^{-2+N^{-1}} \Big(\sum_{k=0}^{N^{1-\kappa}}k\epsilon+\sum_{N^{1-\kappa}}^{\infty}\frac{1}{\epsilon^2 k^2}\Big) \rhosc (s)\D s  
\\
 &+   \int_{-2+N^{-1}}^{-2+N^{-1/2}} \Big(\sum_{k=0}^{N^{1-\kappa}}\epsilon N^{1/2}+\sum_{k=N^{1-\kappa}}^{\infty}\frac{1}{\epsilon^2k^3} N^{1/2}\Big) \rhosc (s)\D s
 \lesssim   N^{-1/4}.
\end{align*}

 Note that we have proved,
 \begin{align*} 
 \mathbb{E}_{h(H_0)}[S_N(g(H_t))]=\frac{1}{N} \mathbb{E}_{h}\big[S_N( h'\cdot\breve{\chi}\cdot\vhlbrt_t g)\big]+\int_{-2}^{2} h'(s)\vhlbrt_tg(s) \rhosc (s)\D s+O(N^{-\kappa/4})
 \end{align*}
 which is valid uniformly over the function sequences $g=g_N\in\mathscr{A}_{C,\epsilon}$. Therefore, we can say that uniformly in $f\in\mathscr{S}_{C,\kappa}$,
\begin{align*} 
\mathbb{E}_{h(H_0)}[S_N(f(H_t))]=\frac{1}{N} \mathbb{E}_{h}\big[S_N( h'\cdot\breve{\chi}\cdot\vhlbrt_t f)\big]+\int_{-2}^{2} h'(s)\vhlbrt_tf(s) \rhosc (s)\D s+O(N^{-\kappa/4}\log N).
\end{align*}
Now we consider an $f\in\mathscr{S}_{C,\kappa}$ satisfying the Hilbert transform condition given in the statement of the theorem, i.e. $\vhlbrt_t f(B)=O(N^{-\kappa/2})$. Then, it suffices to prove that $\frac{1}{N} \mathbb{E}_{h}[S_N( h'\cdot\breve{\chi}\cdot\vhlbrt_tf)]=O(N^{-\kappa/4})$. By rigidity we have
\begin{align}\label{eqn:final_multi_loop} 
\Big|\frac{1}{N} \mathbb{E}_{h}\Big[S_N\big( h' (\vhlbrt_tf)\breve{\chi}\big)\Big]\Big|\leq\frac{\varphi^{11}}{N} \Big(\int |h''(x) \vhlbrt_tf(x)\breve{\chi}(x)|\D x+\int |h'(x) (\vhlbrt_tf)'(x)\breve{\chi}(x)|\D x+\int |h'(x) \vhlbrt_tf(x)\breve{\chi}'(x)|\D x\Big)
\end{align}
In order to discuss these integrals we first need to establish some bounds on $\vhlbrt_tf$. Take an arbitrary $g=g_N\in\mathscr{A}_{C,\epsilon}$ again for $\epsilon=\epsilon_N\in[N^{-1+\kappa},1]$. Using the notation $\upkappa(x)=|x-2|\wedge|x+2|$ and $|\tilde{U}_k'(x)|\lesssim \frac{k}{\upkappa(x)}+\frac{1}{(\upkappa(x))^{3/2}}$ on $x\in[-2,2]$, we obtain that for all $x\in[-2+N^{-1/2},2-N^{-1/2}]$,
\begin{align*} 
|(\vhlbrt_tg)'(x)|\lesssim \sum_{k=1}^{N^{1-\kappa}}\epsilon \big(\frac{k}{\upkappa(x)}+\frac{1}{(\upkappa(x))^{3/2}}\big)+\sum_{k=N^{1-\kappa}}^{\infty}\frac{1}{\epsilon^2k^3}\big(\frac{k}{\upkappa(x)}+\frac{1}{(\upkappa(x))^{3/2}}\big) \lesssim \frac{1}{\epsilon\upkappa(x)}\leq \frac{N^{1-\kappa}}{\upkappa(x)}.
\end{align*}
So, $|(\vhlbrt_tf)'(x)|\lesssim \frac{N^{1-\kappa}}{\upkappa(x)}\log N$. Recall the decomposition $h=h_1+h_2$ and the condition $|\vhlbrt_tf(B)|=O(N^{-\kappa/2})$. Due to the bound on $|(\vhlbrt_tf)'|$, we have $|\vhlbrt_tf(x)|=O(N^{-\kappa/2})$ on the support of $h_2$. Similarly, we have
\begin{align*} 
|\vhlbrt_tg(x)|\lesssim \sum_{k=1}^{N^{1-\kappa}}\epsilon \frac{1}{\sqrt{\upkappa(x)}}+\sum_{k=N^{1-\kappa}}^{\infty}\frac{1}{\epsilon^2k^3}\frac{1}{\sqrt{\upkappa(x)}} \lesssim \frac{1}{\sqrt{\upkappa(x)}}
\end{align*}
for all $x\in[-2,2]$. So, $|\vhlbrt_tf(x)|\lesssim \frac{\log N}{\sqrt{\upkappa(x)}}$ for all $x\in[-2,2]$. Hence the first integral term in \eqref{eqn:final_multi_loop} can be bounded by:
\begin{align*} 
\int|h_1''(x)|\frac{1}{\sqrt{\upkappa(x)}}\D x+\int |h_2''(x)| N^{-\kappa/2}\D x\leq \varphi(N^{1-\kappa/2}+N^{1+\kappa/100-\kappa/2})=O(N^{1-\kappa/3}).
\end{align*}
The second integral term in \eqref{eqn:final_multi_loop} is bounded by
\begin{align*} 
\int_{-2+N^{-1/2}}^{-2+N^{-\kappa/2}} |h'(x)|\frac{1}{\epsilon\upkappa(x)}\D x+\int_{2-N^{-\kappa/2}}^{2-N^{-1/2}}|h'(x)|\frac{1}{\epsilon\upkappa(x)}\D x+\int_{-2+N^{-\kappa/2}}^{2-N^{-\kappa/2}}|h'(x)|\frac{1}{\epsilon\upkappa(x)}\D x=O(N^{1-\kappa/3})
\end{align*}
where we have used $\sup\{|h'(x)|:x\in(-2,-2+N^{-\kappa/2})\cup(2-N^{-\kappa/2},2)\}<\varphi$ and $\int|h'|\leq\varphi$. Lastly, using $|\vhlbrt_tf(x)|\lesssim \frac{\log N}{\sqrt{\upkappa(x)}}=O(N^{1/4}\log N)$ for all $x$ on the support of $\breve{\chi}'$, the third integral in \eqref{eqn:final_multi_loop} is bounded by $\int \varphi N^{1/4} \log N|\breve{\chi'}|=O(N^{1/3})$. Thus, we proved 
\begin{align*} 
\Big|\frac{1}{N} \mathbb{E}_{h}\Big[S_N\big( h'\breve{\chi} \vhlbrt_tf\big)\Big]\Big|=O(N^{-\kappa/4})
\end{align*}
and this completes the proof of
\begin{align*} 
\Ex_{h(H_0)}[S_N(f(H_t))]=\int_{-2}^{2} h'(s)\vhlbrt_t f(s) \rhosc (s)\D s+O(N^{-\kappa/4})=\frac{1}{4}\sum_{k=0}^\infty ke^{-tk/2}\hat{f}_k\hat{h}_k +O(N^{-\kappa/4}).
\end{align*}

Finally, we discuss the error term with $\mathcal{E}(z_{t+u},w,u)$, i.e. we will prove that
\begin{align*} 
\frac{e^{t/2}}{N\pi^2} \int_{\C}  \int_{\C} |\partialbar\tilde{g}(z)|  |\partialbar\widetilde{h'}(w)|\int_{0}^{T} |\frac{e^{u/2}\mathcal{E}(z_{t+u},w,u)}{\eta_{t+u}}| \D u\D w\D z=O(N^{-\kappa/3})
\end{align*}
where choosing $\varepsilon=\frac{\kappa}{100}$ in the definition of $\mathcal{E}$ and substituting the bounds given in Lemma \ref{lemma:integral_with_two_sing}, denoting $w'=w+i\sgn(\eta_w)\mathcal{C}_{N^{\kappa/200}}(E_w)$,
\begin{align*} 
\mathcal{E}(z_{t+u},w,u)\lesssim |\frac{N^{-1/2+\kappa/200}}{\eta_w\sqrt{|\eta_{t+u}|}\eta_{t+2u}}|+|\frac{\varphi^3}{\eta_{t+u}\eta_{t+2u}}|+|\frac{\varphi^2}{\sqrt{|\eta_{t+u}|}\sqrt{|\eta_{w'}|}\eta_{t+2u}}|.
\end{align*}
Hence, using $\eta_{t+u}e^{u/2} \asymp\eta_{t+2u}$ uniformly in all $u\geq0$, simply dismissing $t$'s by $|\eta_{t+u}|\geq|\eta_u|$ and $|\eta_{w'}|\gtrsim N^{-1+\kappa/200}$ we obtain
\begin{align}\label{eqn:3errorterms}
|\frac{e^{u/2}\mathcal{E}(z_{t+u},w,u)}{\eta_{t+u}}|\lesssim \frac{N^{-1/2+\kappa/200}}{|\eta_{u}|^{5/2}|\eta_w|}+\frac{\varphi^3}{|\eta_{u}|^{3}}+\frac{N^{1/2}}{|\eta_{u}|^{5/2}}.
\end{align}
Now, we plug these three in one by one. For any $a>1$, $\int_0^\infty \frac{1}{|\eta_u|^a}\D u\lesssim \frac{1}{|\eta|^{a-1}\Im\sqrt{z^2-4}}$. First, consider the case $\epsilon<\frac{\Upsilon}{2C}$. The first term in \eqref{eqn:3errorterms} gives
\begin{align*} 
N^{-3/2+\kappa/200} \int_{\C}   |\partialbar\tilde{g}(z)\frac{1}{\eta^{3/2}}|\D z\cdot   \int_{\C}|\frac{\partialbar\widetilde{h'}(w)}{\eta_w}| \D w\lesssim N^{-3/2+\kappa/200}N^{1/2}N^{3\varepsilon/2-\kappa/2}\varphi N^{1+\kappa/100}=O(N^{-\kappa/3}).
\end{align*}
Similarly, the second term can be evaluated as 
\begin{align*} 
N^{-1} \int_{\C}   |\partialbar\tilde{g}(z)\frac{1}{\eta^{2}}|\D z\cdot   \int_{\C}|\partialbar\widetilde{h'}(w)| \D w\lesssim N^{-1} \varphi N N^{2\varepsilon-\kappa} \varphi  =O(N^{-\kappa/3}).
\end{align*}
 Lastly, the third term's contribution can be bounded by
\begin{align*} 
N^{-1/2} \int_{\C}   |\partialbar\tilde{g}(z)\frac{1}{\eta^{3/2}}|\D z\cdot   \int_{\C}|\partialbar\widetilde{h'}(w)| \D w\lesssim N^{-1/2} \varphi N^{1/2}N^{3\varepsilon/2-\kappa/2} \varphi =O(N^{-\kappa/3}).
\end{align*}

On the other hand, if $\epsilon\geq\frac{\Upsilon}{2C}$, then $g$ is an order $1$ function and it suffices to use $\int_0^\infty \frac{1}{|\eta_u|^a}\D u\lesssim \frac{1}{|\eta|^{a}}$ and $\int_{\C}|\frac{\partialbar\tilde{g}(z)}{\eta^{a}}|\D z\lesssim \varphi N^{a\varepsilon}$. So the error collected from \eqref{eqn:3errorterms} will be:
\begin{align*} 
N^{-3/2+\kappa/200} N^{5\varepsilon/2}\varphi N^{1+\kappa/100}+N^{-1}  N^{3\varepsilon} \varphi+ N^{-1/2} \varphi N^{5\varepsilon/2} \varphi \ll N^{-1/3},
\end{align*}
concluding the proof.
\end{proof}

\begin{remark}
In the statement of the theorem, we assumed that both $f$ and $h$ have compact supports. However, the result naturally extends to any function satisfying suitable growth conditions.  Although we do not provide the details of this extension, we will apply the theorem to the case that both $f$ and $h$ are regularized logarithms.
\end{remark}

In order to keep the statement of the above theorem simple, both the bias and the function whose linear statistics is being calculated are given at one time ($0$ and $t$). The generalization, as stated below, follows from exactly the same proof.

\begin{theorem}\label{thm:multitime_loop_gen}
Let $C>0$, $\kappa\in(0,1]$, $\Upsilon>0$, $J_1,J_2\in\N$, $I\in\Z^{+}$ be fixed. Take $s_1,\dots,s_I, t_1,\allowbreak\dots,t_{J_1+J_2}\in[0,C]$, $B_1,\dots,B_{J_1}\in[-2+\Upsilon,2-\Upsilon]$, $f_1,\dots,f_I\in\mathscr{S}_{C,\kappa}$, $h_1,\dots,h_{J_1}\in\mathscr{S}_{C,(-\kappa/100)}$ such that $\supp(h_j)\subseteq[B_j-CN^{-1+\kappa/2},B_j+CN^{-1+\kappa/2}]$ for $j=1,\dots,J_1$ and $h_{J_1+1},\dots,h_{J_1+J_2}\in\mathscr{S}_{C,\kappa/2}$. Then the following asymptotic formula is valid provided that $\sum_{i=1}^{I} \vhlbrt_{|s_i-t_j|}f_i(B_j)=O(N^{-\kappa/2})$ for every $j=1,\dots,J_1$,
\begin{align*} 
\mathbb{E}_{\sum_{j=1}^{J_1+J_2} h_j(H_{t_j})}\big[S_N(\sum_{i=1}^{I}f_{i}(H_{s_i}))\big]=\mathcal{C}\Big(\sum_{i=1}^{I}f_{i}(H_{s_i}),\sum_{j=1}^{J_1+J_2} h_j(H_{t_j})\Big)  +O(N^{-\kappa/5}).
\end{align*}
Moreover, the error term is uniform over the choice of $\boldsymbol s$, $\boldsymbol t$, $\boldsymbol B$, $\boldsymbol f$ and $\boldsymbol h$ satisfying the given conditions.
\end{theorem}

\begin{lemma}\label{lemma:GreentoHS}
Let $k\in\N$, $f$ be a compactly supported $\mathscr{C}^k$ function on $\R$, $t\geq 0$ and $\tilde{f}$ be its $k^{th}$ order quasi-analytic extension with a bump function $\chi_{c}$ for a $c>0$. Then, for every $s\in \R$:
\begin{align*} 
\frac{1}{\pi}\int_{\C}\partialbar\tilde{f}(z)\frac{1}{\sqrt{z^2-4}(s-z_t)}\D z=\Re\Big(\int_{-2}^2\frac{f(x)}{x_t-s}\rhoarcsin (x)\D x\Big)=\vhlbrt_t f(s).
\end{align*}
\end{lemma}
\begin{proof}
Noting that $z_t$ and $\sqrt{z^2-4}$ are analytic on upper half plane and lower half plane separately, the proof becomes a simple application of Green's theorem:
\begin{align*} 
\int_{\C}\partialbar\tilde{f}(z)\frac{1}{\sqrt{z^2-4}(s-z_t)}\D z=\int_{\eta>0} \partialbar\Big(\tilde{f}(z)\frac{1}{\sqrt{z^2-4}} \frac{1}{s-z_t}\Big)\D z+\int_{\eta<0} \partialbar\Big(\tilde{f}(z)\frac{1}{\sqrt{z^2-4}} \frac{1}{s-z_t}\Big)\D z
\\
\frac{-i}{2}\int_{\R}\frac{f(x)}{\sqrt{x^2-4}(s-x_t)}\D x+\frac{i}{2}\int_{\R}\frac{f(x)}{\overline{\sqrt{x^2-4}}(s-\bar{x_t})}\D x=-\frac{1}{2}\int_{\R}\frac{f(x)}{\sqrt{4-x^2}}(\frac{1}{s-x_t}+\frac{1}{s-\bar{x_t}})\D x,
\end{align*}
concluding the proof.
\end{proof}
%
%

\section{Proofs of main theorems}\label{sec:FH}

The techniques developed so far are valid for smooth functions; however, we aim to handle logarithmic and jump singularities. To address this, we first establish two distinct types of regularizations for these singularities and prove that the original expression with the singular terms is bounded between the expressions obtained from the regularized versions.

\subsection{Submicroscopic smoothing}

\begin{prop}\label{prop:submic_smoothing}
Let $C>0$, $I,J\in\N$, $\alpha,\Upsilon\in(0,1)$; define $\Delta=N^{-1-\alpha}$. Assume that $\gamma_1,\dots,\gamma_J\in[0,C]$, $\beta_1,\dots,\beta_J\in[-C,C]$; $E_1,\dots,E_J\in[-2+\Upsilon,2-\Upsilon]$ and $t_1,\dots,t_J,s_1,\dots,s_i\in[0,C]$; $f_1,\dots,f_I\in\mathscr{S}_{C,0}$. Therefore,
\begin{gather*}
\E{e^{\sum_j \Tr (\gamma_j\log^{E_j}+\beta_j\arg^{E_j})(H_{t_j})+\sum_i \Tr f_{s_i}(H_{s_i})}}\leq \E{e^{\sum_j \Tr (\gamma_j\log^{E_j}_{\Delta}+\beta_j\arg_{d(j),\Delta}^{E_j})(H_{t_j})+\sum_i \Tr f_{s_i}(H_{s_i})}}
\\
\E{e^{\sum_j \Tr (\gamma_j\log^{E_j}+\beta_j\arg^{E_j})(H_{t_j})+\sum_i \Tr f_{s_i}(H_{s_i})}}\geq \E{e^{\sum_j \Tr \big((\gamma_j\log^{E_j}+\beta_j\arg^{E_j})\ast\tilde{\chi}_\Delta\big)(H_{t_j})+\sum_i \Tr f_{s_i}(H_{s_i})}}(1+O(N^{-\alpha/3}))
\end{gather*}
where $\tilde{\chi}_{\Delta}:=\|\chi_{\Delta}\|_{L^1}^{-1}\chi_{\Delta}$ and $d(j)$ is $r$ if $\beta_j\geq 0$, $\ell$ if $\beta_j<0$. The $O(N^{-\alpha/3})$ term is uniform over the choice of $\boldsymbol \gamma$, $\boldsymbol \beta$, $\boldsymbol E$, $\boldsymbol t$, $\boldsymbol s$ and $\boldsymbol f$ satisfying the given conditions.
\end{prop}

\begin{proof}
The first inequality is obvious by $\gamma_j\log^{E_j}+\beta_j\arg^{E_j}\leq\gamma_j\log^{E_j}_{\Delta}+\beta_j\arg_{d(j),\Delta}^{E_j}$, which holds due to the choice of right and left regularization for $\arg$ function depending on the sign of $\beta_j$. For the second inequality, we begin with proving
\begin{align}\label{eqn:EF(H)=EF(H+x)}
\E{e^{\sum_j  \Tr (\gamma_j\log^{E_j}+\beta_j\arg^{E_j})(H_{t_j})+\sum_i \Tr f_{s_i}(H_{s_i})}}=\E{e^{\sum_j \Tr (\gamma_j\log^{E_j}+\beta_j\arg^{E_j})(H_{t_j}+\varepsilon)+\sum_i \Tr f_{s_i}(H_{s_i}+\varepsilon)}}(1+O(N^{\alpha/3}))
\end{align}
uniformly in $|\varepsilon|\leq 2\Delta$ via multidimensional Girsanov's Theorem \cite[Theorem 8.6.4]{oksendal2013stochastic}. For convenience, define $F(H)=e^{\sum_i S_N (f_{i})(H_{s_i})+\sum_j  S_N (\gamma_j\log^{E_j}+\beta_j\arg^{E_j})(H_{t_j})}$. By translation invariance of the OU process we can assume that $\min_i s_i,\min_j t_j\geq 1$. Let $g:\R\to[0,1]$ be an order $1$ smooth step function such that $g(x)=0$ when $x\leq 0$ and $g(x)=1$ when $x\geq 1$. Define $\tilde{H}_t:=H_t+g(t)\varepsilon$. Note that the SDE for $\tilde{H}_t$ becomes 
\begin{align}\label{tilde_H_matrix_SDE}
\D\tilde{H}_t=\frac{1}{\sqrt{N}}\D B_t+(-\frac{1}{2}\tilde{H}_t+\frac{1}{2}g(t)\varepsilon+g'(t)\varepsilon)\D t
\end{align}

 View $H_t$ as a vector of size $N^2$ given by $(H^{11}_t,\dots,H^{NN}_t,H^{12,R}_t,\dots,H^{(N-1)N,I}_t)^T$ where $H^{ij,R}$ and $H^{ij,I}$ stands for the real and imaginary parts of the $(i,j)^{th}$ entry of $H$ (similarly for $\tilde{H}_t$). Let $\theta=\textnormal{diag}(\underbrace{1,\dots,1}_{N\text{ times}},\underbrace{\frac{1}{\sqrt{2}},\dots,\frac{1}{\sqrt{2}}}_{(N-1)N\text{ times}})$ and $W_t$ be the Brownian motion in $\R^{N^2}$ so that
\vspace{-0.8cm}\begin{align*}
\D H_t=\frac{1}{\sqrt{N}}\theta \D W_t -\frac{1}{2}H_t\D t.
\end{align*} 
 Define the vector $\phi_t:=(\underbrace{1,\dots,1}_{N\text{ times}},\underbrace{0,\dots,0}_{(N-1)N\text{ times}}\hspace{-0.3cm})^T(\frac{1}{2}g(t)\varepsilon+g'(t)\varepsilon)$. The vector form of the SDE \eqref{tilde_H_matrix_SDE} becomes
 \vspace{-0.2cm}\begin{align*}
 \D \tilde{H}_t=\frac{1}{\sqrt{N}}\theta\D W_t+(-\frac{1}{2}\tilde{H}_t+\phi_t)\D t
 \end{align*}
 Then, by Girsanov's Theorem, if we denote $u_t:=(\frac{1}{\sqrt{N}}\theta)^{-1}\phi_t$, we have that $\tilde{W}_t=W_t+\int_0^t u_s\D s$ is a Brownian motion with respect to measure $\mathbb{Q}$ for which $\frac{\D\Prob}{\D \mathbb{Q}}=M$ where $M=\exp\big(\int_0^\cdot u_s\cdot \D \tilde{W}_s- \frac{1}{2}\int_0^\cdot |u_s|^2\D s\big)$; and $\D \tilde{H}_t=\frac{1}{\sqrt{N}}\theta \D \tilde{W}_t-\frac{1}{2}\tilde{H}_t\D t$. That leads
 \begin{align*}
\E{F(\tilde{H})}=\mathbb{E}_{\mathbb{Q}}[F(\tilde{H})M_T]=\E{F(H)\exp\big(\int_0^T u_s\cdot\D W_s- \frac{1}{2}\int_0^T |u_s|^2\D s\big)}
 \end{align*}
 Notice that the martingale term is very close to $1$ with high probability, because:
 \begin{align*}
 \int_0^T u_s\cdot\D W_s=\sqrt{N}\varepsilon\sum_{i=1}^N \int_0^T(\frac{g(s)}{2}+g'(s))\D W^{ii}_s \stackrel{(d)}{=} \mathcal{N}(0,N^2 \varepsilon^2 \int_0^T |\frac{g(s)}{2}+g'(s)|\D s)
 \end{align*}
 and 
 \begin{align*}
 \int_0^T |u_s|^2\D s=N^2\varepsilon^2\int_0^T |\frac{g(s)}{2}+g'(s)|\D s\lesssim N^{-2\alpha}.
 \end{align*}
So if we denote $\mathcal S:=\Big\{\exp\big(\int_0^T u_s\cdot\D W_s- \frac{1}{2}\int_0^T |u_s|^2\D s\big)\in(1-N^{-\alpha/3},1+N^{-\alpha/3})\Big\}$, by Gaussian tail bound we obtain $\Prob(\mathcal S)\geq 1-e^{-N^{\alpha}}$. Therefore,
 \begin{align*}
 \E{F(\tilde{H})}
 &=\E{F(H)(1+O(N^{-\alpha/3}))\mathds{1}_{\mathcal{S}}}+\E{F(H)\exp\big(\int_0^T u_s\cdot\D W_s- \frac{1}{2}\int_0^T |u_s|^2\D s\big)\mathds{1}_{\mathcal{S}^c}}
 \\
 &=\E{F(H)}(1+O(N^{-\alpha/3}))-\E{F(H)(1+O(N^{-\alpha/3}))\mathds{1}_{\mathcal{S}^c}}+\E{F(H)\exp\big(\int_0^T u_s\cdot\D W_s- \frac{1}{2}\int_0^T |u_s|^2\D s\big)\mathds{1}_{\mathcal{S}^c}}
 \end{align*}
Note that $e^{-(\log N)^{60}}\leq\E{F(H)^2}\leq e^{(\log N)^{60}}$ due to Proposition \ref{prop:rig_for_mathscr} and Lemma \ref{lemma:rigidity_lemma_for_log_and_regularizations} by Jensen and H\"older inequalities. Moreover, because $\big(\exp(\int_0^\cdot 3u_s\cdot\D W_s- \frac{1}{2}\int_0^\cdot 9|u_s|^2\D s)\big)$ is a martingale, we have,
\begin{align*}
\E{\exp(\int_0^T 3u_s\cdot\D W_s- \frac{1}{2}\int_0^T 3|u_s|^2\D s)}^{1/3}=\E{\exp(\int_0^T 3u_s\cdot\D W_s- \frac{1}{2}\int_0^T 9|u_s|^2\D s)}^{1/3}e^{\int_0^T |u_s|^2\D s}=e^{\int_0^T |u_s|^2\D s}\leq e^{N^{-\alpha}}.
\end{align*}
Thus, by H\"older inequality, we obtain
\begin{align*}
 \E{F(\tilde{H})}=\E{F(H)}(1+O(N^{-\alpha/3}))
 \end{align*}
 which implies the equation \eqref{eqn:EF(H)=EF(H+x)}.

Having \eqref{eqn:EF(H)=EF(H+x)}, if we denote by $X$ a random variable with probability density function $\tilde{\chi}_{\Delta}$, and independent from the process $H$, we can write:
\begin{align}\label{eqn:log_conv_reg}
\E{e^{\sum_j  \Tr (\gamma_j\log^{E_j}+\beta_j\arg^{E_j})(H_{t_j})+\sum_i \Tr f_{s_i}(H_{s_i})}}&=\E{\mathbb{E}_X [e^{\sum_j  \Tr (\gamma_j\log^{E_j}+\beta_j\arg^{E_j})(H_{t_j}+X)+\sum_i \Tr f_{s_i}(H_{s_i}+X)}]}(1+O(N^{-\alpha/3})) \nonumber
\\
&\hspace{-3cm}\geq\E{e^{\sum_j \mathbb{E}_X[\Tr (\gamma_j \log^{E_j}+\beta_j\arg^{E_j})(H_{t_j}+X)]+\sum_i \mathbb{E}_X[\Tr f_{s_i}(H_{s_i}+X)]}}(1+O(N^{-\alpha/3})) \nonumber
\\
&\hspace{-3cm}=\E{e^{\sum_j  \Tr \big((\gamma_j\log^{E_j}+\beta_j\arg^{E_j})\ast\tilde{\chi}_\Delta)(H_{t_j})+\sum_i \Tr f_{s_i}\ast\tilde{\chi}_\Delta(H_{s_i})}}(1+O(N^{-\alpha/3})).
\end{align}

 Define rigidity set $\mathcal{R}:=\bigcap\limits_{r\in\{s_i:i\}\cup\{t_j:j\}} \{|\lambda_k(r)-\gamma_k|<\frac{(\log N)^{\log\log N}}{N^{2/3}\hat{k}^{1/3}},\textnormal{ for all } k=1,\dots,N\}$ for which we have $\Prob(\mathcal{R})\geq 1-e^{-(\log N)^{100}}$ by \cite[Lemma 3.8]{bourgade2022optimal}. It is easy to see that on $\mathcal{R}$, $ |\Tr f_{s_i}\ast\tilde{\chi}_{\Delta}(H_{s_i})-\Tr f_{s_i}(H_{s_i})|\leq N^{-\alpha/2}$. Hence,
\begin{align*}
\E{&e^{\sum_j  \Tr \big((\gamma_j\log^{E_j}+\beta_j\arg^{E_j})\ast\tilde{\chi}_\Delta)(H_{t_j})+\sum_i \Tr f_{s_i}\ast\tilde{\chi}_\Delta(H_{s_i})}}\\
&\hspace{4cm}\geq \Big(\E{e^{\sum_j \Tr \big((\gamma_j\log^{E_j}+\beta_j\arg^{E_j})\ast\tilde{\chi}_\Delta)(H_{t_j})+\sum_i \Tr f_{s_i}(H_{s_i})}}-\mathcal{E}\Big)(1+O(N^{-\alpha/2}))
\end{align*}
where
\begin{align*}
\mathcal{E}=\E{e^{\sum_j \Tr \big((\gamma_j\log^{E_j}+\beta_j\arg^{E_j})\ast\tilde{\chi}_\Delta)(H_{t_j})+\sum_i \Tr f_{s_i}(H_{s_i})}\cdot & \mathds{1}_{\mathcal{R}^c}}
\leq e^{N\int ( \sum_{j}(\gamma_j\log^{E_j}+\beta_j\arg^{E_j})\ast\tilde{\chi}_\Delta)+\sum_{i}f_{s_i}) \rhosc }e^{(\log N)^{60}}\Prob(\mathcal{R}^{C})^{1/2}
\\
&\leq N^{-1}\cdot \E{e^{\sum_j \Tr \big((\gamma_j\log^{E_j}+\beta_j\arg^{E_j})\ast\tilde{\chi}_\Delta)(H_{t_j})+\sum_i \Tr f_{s_i}(H_{s_i})}}.
\end{align*}
by Proposition \ref{prop:rig_for_mathscr}, H\"older and Jensen inequalities. Together with the equation \eqref{eqn:log_conv_reg}, this completes the proof.
\end{proof}

\subsection{Fisher-Hartwig asymptotics for regularized singularities}

We now proceed to discuss the Fisher-Hartwig asymptotics for the different regularizations of singularities introduced in the previous section. 

\begin{prop}\label{prop:FH_upper_bound}
Let $C>1$, $I,J\in\N$, $\Upsilon\in(0,1)$ and $\kappa\in(0,\frac{1}{1500C})$; define $\Delta=N^{-1-\kappa/200}$. Assume that $\gamma_1,\dots,\gamma_J\in[0,C]$, $\beta_1,\dots,\beta_J\in[-C,C]$; $E_1,\dots,E_J\in[-2+\Upsilon,2-\Upsilon]$ and $t_1,\dots,t_J,s_1,\dots,s_i\in[0,C]$ with separation condition $\min_{j_1\neq j_2}(|(t_{j_1},E_{j_1})-(t_{j_2},E_{j_2})|)>N^{-1+150C\kappa}$; $f_1,\dots,f_I\in\mathscr{S}_{C,\kappa}$. Therefore,
\begin{align*} 
\log \Ex\Big[&e^{\sum_{j=1}^{J} S_N(\gamma_j\log_{\Delta}^{E_j}+\beta_j\arg_{d(j),\Delta}^{E_j})(H_{t_j})+\sum_{i=1}^{I}S_N(f_i)(H_{s_i}) }\Big]
\\
=&\sum_{j=1}^{J}N^{\frac{\gamma_j^2+\beta_j^2}{4}}\frac{G(1+\frac{\gamma_j}{2}+i\frac{\beta_j}{2})G(1+\frac{\gamma_j}{2}-i\frac{\beta_j}{2})}{G(1+\gamma_j)}+\frac{1}{2}\mathcal{C}^{\circ}\Big(\sum_{j}(\gamma_j\log^{E_j}+\beta_j\arg^{E_j})(H_{t_j})+\sum_{i}f_i(H_{s_i})\Big)
\\
&+\sum_{j}\Big(\frac{\gamma_j^2}{8}\log(4-E_j^2)+\frac{\beta_j\gamma_j}{4}(\pi-E_j-2\arccos\frac{E_j}{2})+\frac{\beta_j^2}{8}(1-2\sqrt{4-E_j^2}+3\log(4-E_j^2))\Big)+O(N^{-\kappa/2})
\end{align*}
where $d(j)$ is $r$ if $\beta_j\geq 0$ and $\ell$ if $\beta_j<0$, the $\mathcal{C}^{\circ}$ term is defined as in Theorem \ref{thm:FH}, and $G$ stands for the Barnes $G$-function. Moreover, the error term is uniform over the choice of $\boldsymbol \gamma$, $\boldsymbol \beta$, $\boldsymbol E$, $\boldsymbol t$, $\boldsymbol s$ and $\boldsymbol f$ satisfying the given conditions.
\end{prop}

\begin{proof}
Define $\epsilon:=N^{-1+\kappa}$, $\alpha:=\kappa/200$. 
The first step of the proof is to eliminate the long-range/regular components of the functions via the multi-time loop equation and obtain a Laplace transform of non-positive the local functions which are separated in space-time. To satisfy the $\vhlbrt$ transform condition in Theorem \ref{thm:multitime_loopp}, we introduce two bump functions, left and right, $q_{\ell,\theta}$ and $q_{r,\theta}$ for any $\theta>0$ and incorporate them into the analysis:
\begin{align}\label{eqn:comp_q}
q_{\ell,\theta}(x):=\begin{cases}0,&x\leq-2\theta \\ -1, &x\in[-\theta,-\theta N^{-\kappa/2}] \\ 0, &x\geq -\theta N^{-\kappa/2}/2 \end{cases}, \quad q_{r,\theta}(x):=\begin{cases}0,&x\leq\theta N^{-\kappa/2}/2 \\ -1, &x\in[\theta N^{-\kappa/2},\theta ] \\ 0, &x\geq 2\theta \end{cases}
\end{align}
with smooth transitions in between and $q_{r,\theta}^E:=q_{r,\theta}(\cdot-E)$ for any $\theta>0$ and $E\in\R$ (similarly $q_{\ell,\theta}^E$). Also, define a function $Y(\xi)$ for $\xi\in[0,1]$ by
\begin{align*} 
Y(\xi):=\log \Ex\Big[\exp\Big(\sum_{j}S_N\big(\gamma_j(\log_{\Delta}^{E_j}-\xi\log_{\epsilon}^{E_j})+\beta_j(\arg_{d(j),\Delta}^{E_j}-\xi\arg_{d(j),\epsilon}^{E_j})+\xi \alpha_{j}q_{s(j),\epsilon}^{E_j}\big)(H_{t_j})+\sum_{i}S_N((1-\xi)f_i)(H_{s_i})\Big)\Big]
\end{align*}
for some $s(j)\in\{\ell,r\}$ and order $1$ constant $\alpha_j$'s to be determined later. Then
\begin{align*} 
Y'(\xi):= -\Ex_{h_\xi}\Big[\sum_{j} S_N(\gamma_j\log_{\epsilon}^{E_j}+\beta_j\arg_{d(j),\epsilon}^{E_j}-\alpha_jq_{s(j),\epsilon}^{E_j})(H_{t_j})+\sum_{i}S_N(f_i)(H_{s_i})\Big]
\end{align*}
where $h_\xi:=\sum_{j}\big(\gamma_j(\log_{\Delta}^{E_j}-\xi\log_{\epsilon}^{E_j})+\beta_j(\arg_{d(j),\Delta}^{E_j}-\xi\arg_{d(j),\epsilon}^{E_j})+\xi \alpha_{j}q_{s(j),\epsilon}^{E_j}\big)(H_{t_j})+\sum_{i}((1-\xi)f_i)(H_{s_i})$ and we use the biased measure notation introduced in \eqref{eqn:bias_defn}. By Propositions \ref{prop:rig_for_mathscr}-\ref{rig_under_biased_prop2} and Lemma \ref{lemma:rigidity_lemma_for_log_and_regularizations}, $h_\xi$ satisfies the rigidity conditions. Moreover, $f_i,q_{s(j),\epsilon}^{E_j},\log_{\epsilon}^{E_j},\arg_{d(j),\epsilon}^{E_j}\in\mathscr{S}_{C,\kappa/2}$ (after a suitable truncation to make $\log_{\epsilon}$ and $\arg_{\epsilon}$ compactly supported); $\log_{N^{-1+\kappa/4}}^{E_j},\arg_{d(j),N^{-1+\kappa/4}}^{E_j}\in\mathscr{S}_{C,\kappa/4}$ and $(\log_{\Delta}^{E_j}-\log_{N^{-1+\kappa/4}}^{E_j}),\allowbreak(\arg_{d(j),\Delta}^{E_j}-\arg_{d(j),N^{-1+\kappa/4}}^{E_j})\in\mathscr{S}_{C,(-\kappa/200)}$ are supported around $E_j$ with size $O(N^{-1+\kappa/4})$. Therefore, to be able to apply Theorem \ref{thm:multitime_loop_gen}, we only need to check the $\vhlbrt$-transform condition, which can be written as,
\begin{align*} 
\sum_j \gamma_j \vhlbrt_{|t_j-t_n|} \log_{\epsilon}^{E_j}(E_n)+\sum_j \beta_j \vhlbrt_{|t_j-t_n|} \arg_{d(j),\epsilon}^{E_j}(E_n)+\sum_i \vhlbrt_{|s_i-t_n|}f_i(E_n)-\sum_j \alpha_j \vhlbrt_{|t_j-t_n|} q_{s(j),\epsilon}^{E_j}(E_n) =O(N^{-\kappa/4})
\end{align*}
for all $n=1,\dots,J$.

By Lemma \ref{lemma:vhlbrt_cond_check} (d), we see that $\vhlbrt_{|t_j-t_n|} q_{s(j),\epsilon}^{E_j}(E_n)$ is $\asymp\log N$ when $j=n$ and $o(N^{-\kappa})$ otherwise. Moreover, by Lemma \ref{lemma:vhlbrt_cond_check} (a)-(c)-(e)-(f) each term in the first three summations is $O(\log N)$. Our goal is to choose $s(j)$'s and non-negative $O(1)$ constants $\alpha_j$'s satisfying this system of linear equations, which can be achieved by 
\begin{align*} 
s(j)=\begin{cases} \ell, & \textnormal{if }\sum_k \gamma_k \vhlbrt_{|t_k-t_j|} \log_{\epsilon}^{E_k}(E_j)+\sum_k \beta_k \vhlbrt_{|t_k-t_j|} \arg_{d(j),\epsilon}^{E_k}(E_j)+\sum_i \vhlbrt_{|s_i-t_j|}f_i(E_j)\geq 0 \\ r, & \textnormal{otherwise}\end{cases}
\end{align*}
and
\begin{align*} 
\alpha_j=\frac{1}{\vhlbrt q_{s(j),\epsilon}^{E_j}(E_j)} \Big(\sum_k \gamma_k \vhlbrt_{|t_k-t_j|} \log_{\epsilon}^{E_k}(E_j)+\sum_k \beta_k \vhlbrt_{|t_k-t_j|} \arg_{d(j),\epsilon}^{E_k}(E_j)+\sum_i \vhlbrt_{|s_i-t_j|}f_i(E_j)\Big).
\end{align*}
Thus, we can apply Theorem \ref{thm:multitime_loop_gen}, and obtain the following, after integrating $Y'(\xi)$ from $0$ to $1$:
\begin{align*}
Y(0)=&Y(1)+\mathcal{C}\Big(\sum_{j}\big(\gamma_j(\log_{\Delta}^{E_j}-\frac{1}{2}\log_{\epsilon}^{E_j})+\beta_j(\arg_{d(j),\Delta}^{E_j}-\frac{1}{2}\arg_{d(j),\epsilon}^{E_j})+\frac{1}{2} \alpha_{j}q_{s(j),\epsilon}^{E_j}\big)(H_{t_j})+\sum_{i}\frac{1}{2}f_i(H_{s_i})
\\
&\hspace{6cm},\sum_{j} (\gamma_j\log_{\epsilon}^{E_j}+\beta_j\arg_{d(j),\epsilon}^{E_j}-\alpha_jq_{s(j),\epsilon}^{E_j})(H_{t_j})+\sum_{i}f_i(H_{s_i})\Big).
\end{align*}
By straightforward calculations, this expression can be simplified to the following, details are provided in Lemma \ref{lemma:log_arg_index_removals},
\begin{equation}\label{eqn:Y(0)=Y(1)+...}
\begin{aligned}
Y(0)=& Y(1)+\frac{1}{2}\mathcal{C}^{\circ}\Big(\sum_{j}(\gamma_j\log^{E_j}+\beta_j\arg^{E_j})(H_{t_j})+\sum_{i}f_i(H_{s_i})\Big)+\mathcal{L}_1+\mathcal{L}_2+\mathcal{L}_3+O(N^{-\kappa/2})
\\
&+\sum_{j}\Big(\frac{\gamma_j^2}{4}\log\frac{1}{2\Delta}+\frac{\gamma_j\beta_j}{4}\big(\frac{1-d(j)}{2}\pi-E_j-2\arccos\frac{E_j}{2}\big)+\frac{\beta_j^2}{8}\big(1-2\sqrt{4-E_j^2}-2\log\frac{\Delta}{4-E_j^2}\big)\Big)
\end{aligned}
\end{equation}
where abusing the notation we assigned a numerical value to $d(j)$ as $1$ if it is $r$ and $-1$ if it is $\ell$; $\mathcal{L}_1$, $\mathcal{L}_2$ and $\mathcal{L}_3$, defined in Lemma \ref{lemma:log_arg_index_removals}, are quantities related to local functions and they will be canceled out later.

Now, we begin our discussion of $Y(1)$. By our choice of $d(j)$ and $s(j)$ --the directions of the $\arg$ regularizations and the bump functions-- for each $j$, the local functions we have are non-negative, more explicitly, $0\geq \gamma_j(\log_{\Delta}^{E_j}-\log_{\epsilon}^{E_j})+\beta_j(\arg_{d(j),\Delta}^{E_j}-\arg_{d(j),\epsilon}^{E_j})+\alpha_jq_{s(j),\epsilon}^{E_j}\geq-C\frac{101\kappa}{100}\log N $. Hence, by Theorem \ref{thm:decoupling}, $Y(1)$ is equal to $\sum_{j}\log \Ex\Big[\exp\Big(S_N\big(\gamma_j(\log_{\Delta}^{E_j}-\log_{\epsilon}^{E_j})+\beta_j(\arg_{\Delta}^{E_j}-\arg_{\epsilon}^{E_j})+\alpha_jq_{s(j),\epsilon}^{E_j}\big)\Big)\Big]$ up to an error of $o(N^{-\kappa})$. This decoupling reduces the problem to a single-time setting with regularized singularity only at a single point. At first glance, one might expect that the same technique could recover the long-range parts of the singularities in the single-time expressions. However, this approach encounters a fundamental obstacle: it requires order $1$ change in $\alpha_j$ constants. So, the asymptotics of the multi-time loop equation, or even the single-time loop equation for $\GUE$ fail to recover the long range of logarithm due to the absence of a cost-free constant shift. This key property --being able to shift by a constant without incurring any cost-- is exclusive to the single-time loop equation for $\CUE$ (this has been discussed in detail in Appendix \ref{app:loop_eqn_CUE}). To address this issue, we apply Proposition \ref{prop:general_GUE_CUE_comparison} to convert the problem into a single-time $\CUE$ setting. Then, we can either proceed using the single-time loop equation (see \cite[Lemma 2.1]{lambert2021mesoscopic}) and carry out the necessary error estimations similarly to \cite{bourgade2022liouville}; or we can directly use \cite[Theorem 1.2]{bourgade2022liouville}. Restricting $\arg$ function to the domain $[-\pi,\pi]$ in the following lines (which is simply equal to $\Im\log(1-e^{ix})$), this method provides us with the following:
\begin{align*} 
Y(1)=&\sum_{j=1}^{J}\log \Ex_{\CUE}\Big[\exp\Big( S_N\big(\Dil_{E_j}\big(\gamma_j(\log_{\Delta}^{E_j}-\log_{\epsilon}^{E_j})+\beta_j(\arg_{d(j),\Delta}^{E_j}-\arg_{d(j),\epsilon}^{E_j})+\alpha_jq_{s(j),\epsilon}^{E_j}\big)\big)\Big)\Big]+O(N^{-\kappa})
\\
=&\sum_{j=1}^{J}\log \Ex_{\CUE}\Big[\exp\Big( S_N\big(\gamma_j(\log_{\tDelta_j}-\log_{\tepsilon_j})|1-e^{i\cdot}|+\beta_j(\arg_{d(j),\tDelta_j}-\arg_{d(j),\tepsilon_j})+\alpha_jq_{s(j),\tepsilon_j}\big)\Big)\Big]+O(N^{-\kappa})
\\
=&\sum_{j=1}^{J}\log\Ex_{\CUE}\Big[|\det(\Id-U)|^{\gamma_j}e^{\beta_j\Im\log\det(\Id-U)}\Big]-\mathcal{L}_1-\mathcal{L}_2-\mathcal{L}_3
\\
& -\sum_{j}\Big(\frac{\gamma_j^2}{4}\big(\log\frac{1}{2\Delta}-\log(\sqrt{4-E_j^2})\big)-\frac{\beta_j\gamma_j}{4}\big(\frac{1+d(j)}{2}\pi\big)-\frac{\beta_j^2}{4}\big(\log\Delta+\log\sqrt{4-E_J^2}\big)\Big)+O(N^{-\kappa/2})
\end{align*}
where $\tilde{\epsilon}_j:=\epsilon2\pi\rhosc (E_j)$, $\tilde{\Delta}_j:=\Delta2\pi\rhosc (E_j)$. We have used Proposition \ref{prop:general_GUE_CUE_comparison} in the first equality, $|1-e^{ix}|=x+O(x^3)$ in the second, and Lemmas \ref{lemma:CUE_loop_app}-\ref{lemma:log_arg_index_removal_circular}-\ref{lemma:L_cancellation} in the third. After substituting this into the equation \eqref{eqn:Y(0)=Y(1)+...}, it suffices to obtain an asymptotic formula for $\Ex_{\CUE}\Big[|\det(\Id-U)|^{\gamma}e^{\beta\Im\log\det(\Id-U)}\Big]$ uniformly in $\gamma\in[0,C]$, $\beta\in[-C,C]$ to complete the proof. Indeed, this was computed in \cite[(71)]{keating2000random} using Selberg's integral formula,
\begin{align*} 
\Ex_{\CUE}\Big[|\det(\Id-U)|^{\gamma}e^{\beta\Im\log\det(\Id-U)}\Big]&=\frac{G(1+N)G(1+N+\gamma)}{G(1+N+\frac{\gamma}{2}+i\frac{\beta}{2})G(1+N+\frac{\gamma}{2}-i\frac{\beta}{2})}\frac{G(1+\frac{\gamma}{2}+i\frac{\beta}{2})G(1+\frac{\gamma}{2}-i\frac{\beta}{2})}{G(1+\gamma)}
\\
&=N^{\frac{\gamma^2+\beta^2}{4}}\frac{G(1+\frac{\gamma}{2}+i\frac{\beta}{2})G(1+\frac{\gamma}{2}-i\frac{\beta}{2})}{G(1+\gamma)}(1+O(\frac{1}{N}))
\end{align*}
for any $\gamma,\beta\in\C$ with $\Re\gamma,\frac{\Re\gamma+\Im\beta}{2},\frac{\Re\gamma-\Im\beta}{2}>-1$, where $G$ is the Barnes $G$-function; and for the second line we have used the Stirling-like asymptotic expansion $\log G(1+z)=z^2(\frac{\log z}{2}-\frac{3}{4})+\frac{z}{2}\log(2\pi)-\frac{\log z}{12}+\zeta'(-1)+O(1/z)$ as $z\to\infty$. 
\end{proof}

Note that the final expression in Proposition \ref{prop:FH_upper_bound} is independent of the submicroscopic regularization index. The next proposition restates Proposition \ref{prop:FH_upper_bound} for a slightly different regularization, yielding the same asymptotics. Since the proof is identical, we omit it here for brevity.

\begin{prop}\label{prop:FH_lower_bound}
Let $C>1$, $I,J\in\N$, $\Upsilon\in(0,1)$ and $\kappa\in(0,\frac{1}{1500C})$; define $\Delta=N^{-1-\kappa/200}$. Assume that $\gamma_1,\dots,\gamma_J\in[0,C]$, $\beta_1,\dots,\beta_J\in[-C,C]$; $E_1,\dots,E_J\in[-2+\Upsilon,2-\Upsilon]$ and $t_1,\dots,t_J,s_1,\dots,s_i\in[0,C]$ with separation condition $\min_{j_1\neq j_2}(|(t_{j_1},E_{j_1})-(t_{j_2},E_{j_2})|)>N^{-1+150C\kappa}$; $f_1,\dots,f_I\in\mathscr{S}_{C,\kappa}$. Therefore,
\begin{align*} 
\log \Ex\Big[&e^{\sum_{j=1}^{J}S_N(\gamma_j\log^{E_j}\ast\tilde{\chi}_{\Delta}+\beta_j\arg^{E_j}\ast\tilde{\chi}_{\Delta})(H_{t_j})+\sum_{i=1}^{I}S_N(f_i)(H_{s_i}) }\Big] 
\\
=&\sum_{j=1}^{J}N^{\frac{\gamma_j^2+\beta_j^2}{4}}\frac{G(1+\frac{\gamma_j}{2}+i\frac{\beta_j}{2})G(1+\frac{\gamma_j}{2}-i\frac{\beta_j}{2})}{G(1+\gamma_j)}+\frac{1}{2}\mathcal{C}^{\circ}\Big(\sum_{j}(\gamma_j\log^{E_j}+\beta_j\arg^{E_j})(H_{t_j})+\sum_{i}f_i(H_{s_i})\Big)
\\
&+\sum_{j}\Big(\frac{\gamma_j^2}{8}\log(4-E_j^2)+\frac{\beta_j\gamma_j}{4}(\pi-E_j-2\arccos\frac{E_j}{2})+\frac{\beta_j^2}{8}(1-2\sqrt{4-E_j^2}+3\log(4-E_j^2))\Big)+O(N^{-\kappa/2})
\end{align*}
where $d(j)$ is $r$ if $\beta_j\geq 0$ and $\ell$ if $\beta_j<0$, $\tilde{\chi}_{\Delta}=\frac{\chi_{\Delta}}{\|\chi_{\Delta}\|_{L^1}}$, and the $\mathcal{C}^{\circ}$ term is defined as in Theorem \ref{thm:FH}. Moreover, the error term is uniform over the choice of $\boldsymbol \gamma$, $\boldsymbol \beta$, $\boldsymbol E$, $\boldsymbol t$, $\boldsymbol s$ and $\boldsymbol f$ satisfying the given conditions.
\end{prop}

\subsection{Proofs of Theorems \ref{thm:FH},  \ref{thm:GMC}, and \ref{thm:max_log_char}}

\begin{proof}[Proof of Theorem \ref{thm:FH}]
Combining Propositions \ref{prop:submic_smoothing}, \ref{prop:FH_upper_bound} and \ref{prop:FH_lower_bound} completes the proof.
\end{proof}

Now we discuss the proof of \eqref{eqn:GMC_root_type}; the proof of \eqref{eqn:GMC_jump_type} follows in a step-by-step parallel manner.

\begin{proof}[Proof of Theorem \ref{thm:GMC}]
After proving Theorem \ref{thm:FH}, i.e. having Fisher-Hartwig asymptotics with merging root-type singularities up to some mesoscopic scale, the convergence to the GMC measure is a very standard corollary in the random matrix literature (see \cite{webb2015characteristic,berestycki2018random,lambert2018subcritical,nikula2020multiplicative,
kivimae2020gaussian,forkel2021classical,claeys2021much,bourgade2022liouville}). Here, we briefly present the key steps for $\gamma\in[0,2)$ ($L^2$-regime) with few calculations. In \cite{lambert2018subcritical}, \cite{nikula2020multiplicative}, and \cite{claeys2021much}, the extension of the result to the $L^1$-regime ($\gamma\in[2,2\sqrt{2})$ in our problem) has been discussed in detail. This is achieved using barrier estimates, an argument involving the random measures under study being supported on set of thick points. These techniques allow for the neglect of non-typical events via the first-moment methods while following an approach similar to the $L^2$-regime calculations for the typical region.

Recall that we have $\log^{x}y=-2\sum_{n=1}^{\infty}\frac{\tilde{T}_n(x)}{n}\tilde{T}_n(y)$. Let 
\begin{align*} 
X_N(t,x):=\log\det(H_t-x)=-2\sum_{n=1}^{\infty}\frac{\sum_{j=1}^{N}\tilde{T}_n(\lambda_j(t))}{n}\tilde{T}_n(x).
\end{align*}
The object of interest is 
\begin{align*} 
\mu_{N,\gamma}(\D t,\D x)=\frac{e^{\gamma X_N(t,x)}}{\E{e^{\gamma X_N(t,x)}}}\D t\D x=\frac{|\det(H_t-x)|^{\gamma}}{\E{|\det(H_t-x)|^{\gamma}}}\D t\D x.
\end{align*}
Introduce an approximating measure with smooth bias as follows. For each $n\in\N$, define $\breve{T}_n$ as a compactly supported continuous function on $\R$ such that $\breve{T}_n|_{(-2-\varepsilon,2+\varepsilon)}=\tilde{T}_n|_{(-2-\varepsilon,2+\varepsilon)}$ (for some small positive number $\varepsilon$ that can depend on $M$) and let
\begin{align*} 
X_{N}^{(M)}(t,x)=(-2)\sum_{n=1}^{M}\frac{\sum_{j=1}^N \breve{T}_n(\lambda_j(t))}{n}\tilde{T}_n(x)
\end{align*}
and similarly define the random measure
\begin{align*} 
\mu_{N,\gamma}^{(M)}(\D t,\D x)=\frac{e^{\gamma X_{N}^{(M)}(t,x)}}{\E{e^{\gamma X_N(t,x)}}}\D t\D x.
\end{align*}
Moreover, let
\begin{align*} 
G_M(t,x):=\sum_{n=1}^{M}\frac{\tilde{T}_n(x)}{\sqrt{n}}(A_n)_{nt}
\end{align*}
where $(A_n)_{t\in\R}$'s are independent stationary OU processes in $\R$ satisfying $\D (A_n)_t=\D (B_n)_t-\frac{1}{2}(A_n)_t\D t$ and define random measure 
\begin{align*} 
\mu_{\gamma}^{(M)}(\D t,\D x):=e^{\gamma G_M(t,x)-\frac{\gamma^2}{2}\E{G_M(t,x)^2}}\D t \D x
\end{align*}
on $\R\times[-2,2]$.

$\mu_{N,\gamma}^{(M)}$ converges weakly in distribution to $\mu_{\gamma}^{(M)}$ as $N$ goes to infinity (due to the convergence of the fluctuation field to the Gaussian process, see Appendix \ref{app:gaussian_field_approx} for details) and $\mu_{\gamma}^{(M)}$ converges weakly in distribution to the GMC measure $\mu_{\gamma}$ that is described in Theorem \ref{thm:GMC}. In other words, for any compactly supported continuous function $\psi:\R\times(-2,2)\to\R$, when we take $N\to\infty$ then $M\to \infty$, $\mu_{N,\gamma}^{(M)}(\psi)$ converges in distribution to $\mu^{\gamma}(\psi)$ (see \cite[Lemma 4.8]{kallenberg2017random} for the convergence of random measures). Thus, it suffices to show $\big(\mu_{N,\gamma}(\psi)-\mu_{N,\gamma}^{(M)}(\psi)\big)$ converges to $0$ in distribution as $N\to\infty$ then $M\to \infty$. Indeeed, when $\gamma\in[0,2)$ it is easy to show that $\big(\mu_{N,\gamma}(\psi)-\mu_{N,\gamma}^{(M)}(\psi)\big)$ converges to $0$ in $L^2$ once we have Theorem \ref{thm:FH}, as demonstrated below.

 Take $\Upsilon>0$ such that $\psi$ is supported in $\R\times(-2+\Upsilon,2-\Upsilon)$. Let $\kappa>0$ be a small constant to be determined later. Define the set $A:=\{(t,x,s,y):|(t,x),(s,y)|>N^{-1+\kappa}\}$. By applying Theorem \ref{thm:FH} and using Lemma \ref{lemma:time-space_log_sing} in order to bound the contributions from $A^c$ we get,
\begin{align*} 
\E{\mu_{N,\gamma}(\psi)^2}=&\int_{(\R\times(-2,2))^2}\psi(t,x)\psi(s,y)\frac{\E{|\det(H_t-x)|^{\gamma}|\det(H_s-y)|^{\gamma}}}{\E{|\det(H_t-x)|^{\gamma}}\E{|\det(H_s-y)|^{\gamma}}}\D t\D x\D s \D y
\\
=&\int_{A}\psi(t,x)\psi(s,y) e^{\gamma^2\mathcal{C}(\log^{x}(H_t),\log^{y}(H_s))+O(N^{-\delta})}\D t\D x\D s \D y
\\
&+O\Big(\int_{A^c} \frac{\E{|\det(H_t-x)|^{2\gamma}}^{1/2}\E{|\det(H_s-y)|^{2\gamma}}^{1/2}}{\E{|\det(H_t-x)|^{\gamma}}\E{|\det(H_s-y)|^{\gamma}}}\D t\D x\D s \D y\Big)
\\
=&\int_{A}\psi(t,x)\psi(s,y) \exp\big({\gamma^2\sum_{n=1}^{\infty}\frac{\tilde{T}_n(x)\tilde{T}_n(y)}{n}e^{-|t-s|n/2}}\big)\D t\D x\D s \D y+O(N^{-2+2\kappa+\frac{\gamma^2}{2}}+N^{-\delta})
\end{align*}
where $\delta>0$ is a constant depending on $\Upsilon$ and $\kappa$. This means, we can choose $\kappa\in(0,2-\frac{\gamma^2}{2})$ and the error term becomes negligible. Moreover, because $\frac{\gamma^2}{2}<2$,
\begin{align*} 
\int_{A^{C}}\psi(t,x)\psi(s,y) \exp\big({\gamma^2\sum_{n=1}^{\infty}\frac{\tilde{T}_n(x)\tilde{T}_n(y)}{n}e^{-|t-s|n/2}}\big)\D t\D x\D s \D y\xrightarrow[N\to\infty]{}0
\end{align*}
which can be seen by Lemma \ref{lemma:time-space_log_sing}. Thus we have 
\begin{align*} 
\lim_{N\to\infty}\E{\mu_{N,\gamma}(\psi)^2}=\int_{(\R\times(-2,2))^2}\psi(t,x)\psi(s,y) \exp\big({\gamma^2\sum_{n=1}^{\infty}\frac{\tilde{T}_n(x)\tilde{T}_n(y)}{n}e^{-|t-s|n/2}}\big)\D t\D x\D s \D y.
\end{align*}
Similarly we obtain,
\begin{equation*}
\begin{aligned} 
\lim_{N\to\infty}\E{\mu_{N,\gamma}(\psi)\mu_{N,\gamma}^{(M)}(\psi)}&=\lim_{N\to\infty}\E{\mu_{N,\gamma}^{(M)}(\psi)^2}
\\
&=\int_{(\R\times(-2,2))^2}\psi(t,x)\psi(s,y)\exp\big(\gamma^2\sum_{n=1}^{M}\frac{\tilde{T}_n(x)\tilde{T}_n(y)}{n}e^{-|t-s|n/2}\big)\D t\D x\D s\D y.
\end{aligned}
\end{equation*}
Therefore, $\lim_{N\to\infty}\E{(\mu_{N,\gamma}(\psi)-\mu_{N,\gamma}^{(M)}(\psi))^2}$ is equal to
\begin{align*} 
\int_{(\R\times(-2,2))^2}\psi(t,x)\psi(s,y) \Big(\exp\big({\gamma^2\sum_{n=1}^{\infty}\frac{\tilde{T}_n(x)\tilde{T}_n(y)}{n}e^{-|t-s|n/2}}\big)-\exp\big(\gamma^2\sum_{n=1}^{M}\frac{\tilde{T}_n(x)\tilde{T}_n(y)}{n}e^{-|t-s|n/2}\big)\Big)\D t\D x\D s \D y
\end{align*}
By Fatou's Lemma, 
\begin{align*} 
\int_{(\R\times(-2,2))^2}\psi(t,x)\psi(s,y) &\exp\big({\gamma^2\sum_{n=1}^{\infty}\frac{\tilde{T}_n(x)\tilde{T}_n(y)}{n}e^{-|t-s|n/2}}\big)\D t\D x\D s \D y
\\
&\leq\liminf_{M}\int_{(\R\times(-2,2))^2}\psi(t,x)\psi(s,y) \exp\big({\gamma^2\sum_{n=1}^{M}\frac{\tilde{T}_n(x)\tilde{T}_n(y)}{n}e^{-|t-s|n/2}}\big)\D t\D x\D s \D y
\end{align*}
On the other hand, because $\E{\big(\mu_{N,\gamma}(\psi)-\mu_{N,\gamma}^{(M)}(\psi)\big)^2}\geq 0$ for every $M\in\N$, we also have
\begin{align*} 
\int_{(\R\times(-2,2))^2}\psi(t,x)\psi(s,y) &\exp\big({\gamma^2\sum_{n=1}^{\infty}\frac{\tilde{T}_n(x)\tilde{T}_n(y)}{n}e^{-|t-s|n/2}}\big)\D t\D x\D s \D y
\\
&\geq\limsup_{M}\int_{(\R\times(-2,2))^2}\psi(t,x)\psi(s,y) \exp\big({\gamma^2\sum_{n=1}^{M}\frac{\tilde{T}_n(x)\tilde{T}_n(y)}{n}e^{-|t-s|n/2}}\big)\D t\D x\D s \D y
\end{align*}
which completes the proof.
\end{proof}

\begin{lemma}\label{lemma:time-space_log_sing} For every $\Upsilon\in(0,1)$, there exists a constant $C>1$ such that for all $x,y\in[-2+\Upsilon,2-\Upsilon]$ and $t\in[0,1]$,
\begin{gather*} 
\frac{1}{C}|(t,x)-(0,y)|^{-1/2}\leq e^{\mathcal{C}_t(\log^x,\log^y)}\leq C |(t,x)-(0,y)|^{-1/2},
\\
\frac{1}{C}|(t,x)-(0,y)|^{-1/2}\leq e^{\mathcal{C}_t(\Xi^x,\Xi^y)}\leq C |(t,x)-(0,y)|^{-1/2}
\end{gather*}
where $|(t,x)-(0,y)|=\sqrt{|t-0|^2+|x-y|^2}$. Furthermore, $C$ can be chosen so that the upper bounds hold for all $x,y\in(-2,2)$. On the other hand, for all $t\geq 1$ and  $x,y\in(-2,2)$,
\begin{align*} 
e^{\mathcal{C}_t(\log^x,\log^y)}\leq e, \quad e^{\mathcal{C}_t(\Xi^x,\Xi^y)}\leq e.
\end{align*}
\end{lemma}

\begin{proof}
Start with $t\in[0,1]$. Let $x=2\cos \theta$ and $y=2\cos\omega$ with $\theta,\omega\in(0,\pi)$. Using the Chebyshev expansion of $\log$ we get:
\begin{multline*}
\mathcal{C}_t(\log^x,\log^y)=\sum_{n=1}^{\infty} e^{-tn/2}n \frac{\tilde{T}_n(x)}{n}\frac{\tilde{T}_n(y)}{n}=\frac{1}{2}\sum_{n=1}^{\infty}e^{-tn/2}\frac{\cos(n(\theta-\omega))+\cos(n(\theta+\omega))}{n} 
\\
= \frac{-1}{2}\Big(\log\big|1-e^{-t/2+i(\theta-\omega)}\big|+\log\big|1-e^{-t/2+i(\theta+\omega)}\big|\Big)
\end{multline*}
where we have used $\log(1-z)=-\sum_{n=1}^{\infty}\frac{z^n}{n}$. Then
\begin{align*} 
e^{\mathcal{C}_t(\log^x,\log^y)}= \frac{1}{\big|1-e^{-t/2+i(\theta-\omega)}\big|^{1/2}\big|1-e^{-t/2+i(\theta+\omega)}\big|^{1/2}} \asymp \frac{1}{t^{1/2}+|x-y|^{1/2}}
\end{align*}
where we have used that uniformly in $s\in[0,1/2]$ and $\alpha\in[-\pi.\pi]$, $|1-e^{-s+i\alpha}|\asymp s+|\alpha|$.

 On the other hand, when $t\geq 1$, $\mathcal{C}_t(\log^x,\log^y)\leq \sum_{n=1}^{\infty}\frac{e^{-tn/2}}{n}=-\log(1-e^{-t/2})\leq 1$.  For $\Xi$ we start with  the same type of expression
 \begin{multline*}
\mathcal{C}_t(\Xi^x,\Xi^y)=\sum_{n=1}^{\infty} e^{-tn/2}n \frac{\tilde{U}_{n-1}(x)\frac{\sqrt{4-x^2}}{2}}{n}\frac{\tilde{U}_{n-1}(y)\frac{\sqrt{4-y^2}}{2}}{n}=\sum_{n=1}^{\infty} e^{-tn/2}\frac{\sin n\theta \sin n\omega}{n}
\\
= \frac{-1}{2}\Big(\log\big|1-e^{-t/2+i(\theta-\omega)}\big|-\log\big|1-e^{-t/2+i(\theta+\omega)}\big|\Big),
\end{multline*}
and the rest of the proof follows similarly.
\end{proof}

Using Theorems \ref{thm:GMC} and \ref{thm:FH}, the proof of the Theorem \ref{thm:max_log_char} follows closely the approach of \cite[Section 3]{bourgade2025optimal}. Since the overall argument is the same, we only outline the necessary changes below.

\begin{proof}[Sketch of the proof of Theorem \ref{thm:max_log_char}]
We follow \cite[Section 3]{bourgade2025optimal}, adapting their arguments to our setting. The proof of convergence in probability of the maximum of the characteristic polynomial splits into two parts: (i) the upper bound, \cite[Sections 3.1-2]{bourgade2025optimal} and (ii) the lower bound, \cite[Section 3.3]{bourgade2025optimal}. 
 
 (i) For the upper bound, we choose the grid $J$ as $\mathcal A\cap (N^{-1-\varepsilon}\Z\times N^{-1-\varepsilon}\Z)$. This requires controlling the variations of the Stieltjes transform of the empirical spectral distribution under submicroscopic time changes. This is achieved by the same computations as in Section \ref{subsec:long_stability}. The only necessary modification is that, in place of $\varphi=(\log N)^{\log\log N}$ used for the related upper bounds, we take $(\log N)^D$ for a suitable constant $D$. This choice calls for some additional care but requires no fundamental changes, since rigidity at this scale is already available from \cite[Lemma 3.8]{bourgade2022optimal}. (ii) For the lower bound, no further input is needed beyond the similar modifications mentioned for the upper bound, because \cite[Theorem 3.4]{claeys2021much} holds in any dimension, and Theorem \ref{thm:GMC} has already been established.
 
 Finally, the proof for controlling the fluctuations of the particles is completely identical to the \cite[Section 3.4]{bourgade2025optimal}.
\end{proof}

\appendix
\section*{Appendices}
\addcontentsline{toc}{section}{Appendices}
\renewcommand{\thesubsection}{\Alph{subsection}}
\numberwithin{equation}{subsection}
\renewcommand{\theequation}{\thesubsection.\arabic{equation}}
\numberwithin{theorem}{subsection}

\subsection{Kernel Estimates: Proofs of Lemma \ref{lemma:GUE_sine_CUE} and Proposition~\ref{prop:kernel_est}}

\subsubsection{$\GUE-\CUE$ kernel comparison}\label{app:pf_GUE_CUE}
\begin{proof}[Proof of Lemma \ref{lemma:GUE_sine_CUE}]
Remember the definition of $\theta,\omega,a_\theta,a_\omega$ in \eqref{eqn:GUE_correlation_kernel_approximation} and denote the scaling constant $\alpha:=2\pi\rhosc (E)$. By Taylor expansion we obtain
\begin{align*}
a_\theta-a_\omega= 4(\omega-\theta)\sin^2\omega-(2\theta-2\omega)^2\frac{\sin2\omega}{2}+O((\theta-\omega)^3)
\end{align*}
and substituting this into the expression in \eqref{eqn:GUE_correlation_kernel_approximation} yields
\begin{align*}
\frac{\sin(\frac{N}{2}(a_\theta-a_\omega))}{\sin(\frac{1}{2}(\omega-\theta))}=\frac{\sin\big(2N(\omega-\theta)\sin^2\omega\big)- \cos\big(2N(\omega-\theta)\sin^2\omega\big) N(\omega-\theta)^2\sin2\omega}{\sin(\frac{1}{2}(\omega-\theta))}+O(N^{-1+\kappa}).
\end{align*} 

On the other hand, using $y-x= 2(\omega-\theta)\sin\omega-(\omega-\theta)^2\cos\omega+O((\omega-\theta)^3)$ and $\sin(ax)=a\sin(x)+O(x^3)$ (when $a$ is of order $1$ and $x$ is small) the kernel for $\CUE$ can be expressed in terms of $\omega$ and $\theta$ as follows:
\begin{multline*}
K_\CUE(x\alpha,y\alpha)=\frac{1}{2\pi}\frac{\sin(N\alpha(\omega-\theta)\sin\omega-\frac{N}{2}\alpha(\omega-\theta)^2\cos\omega)}{\sin(\alpha2\sin\omega\cdot\frac{1}{2}(\omega-\theta))}+O(1)
\\
=\frac{1}{2\pi}\frac{\sin(N\alpha(\omega-\theta)\sin\omega-\frac{N}{2}\alpha(\omega-\theta)^2\cos\omega)}{\alpha2\sin\omega\sin(\frac{1}{2}(\omega-\theta))}+O(1)
=\frac{1}{\alpha 4\pi\sqrt{\sin\omega\sin\theta}}\frac{\sin(N\alpha(\omega-\theta)\sin\omega-\frac{N}{2}\alpha(\omega-\theta)^2\cos\omega)}{\sin(\frac{1}{2}(\omega-\theta))}+O(1)
\end{multline*}
where the last equality follows from the fact that $\frac{1}{\sqrt{\sin\omega}}-\frac{1}{\sqrt{\sin\theta}}=O(\omega-\theta)$. Thus, the difference between the correlation kernels, $K_\GUE(x,y)-K_\CUE(x\alpha,y\alpha)\cdot\alpha$ becomes
\begin{align*}
\frac{C(\omega,\theta)}{\sin(\frac{1}{2}(\omega-\theta))}\bigg( \sin\Big(2N(\omega-\theta)\sin^2\omega-N(\omega-\theta)^2\sin2\omega\Big)-\sin\Big(N\alpha(\omega-\theta)\sin\omega-\frac{N}{2}\alpha(\omega-\theta)^2\cos\omega\Big)\bigg)+O(1)
\end{align*}
where $C(\omega,\theta):=\frac{1}{4\pi\sqrt{\sin\theta\sin\omega}}$ is a term of order $1$. Define $\beta=\arccos (E/2) \in (0,2\pi)$, so that $\alpha=2\sin\beta$. The difference of the arguments of the sine functions simplifies as
\begin{align*}
\Big(2N(\omega-\theta)\sin^2\omega&-N(\omega-\theta)^2\sin2\omega\Big)-\Big(N\alpha(\omega-\theta)\sin\omega-\frac{N}{2}\alpha(\omega-\theta)^2\cos\omega\Big)=
\\
&=N(\omega-\theta)(2\sin\omega-2\sin\beta)\sin\omega-\sin(2\omega)\frac{N}{2}(\omega-\theta)^2+\frac{N}{2}(\omega-\theta)^2(2\sin\beta-2\sin\omega)\cos\omega
\\
&=NE(\omega-\theta)(\omega-\beta)\sin\omega-(\omega-\theta)^2N\sin\omega\cos\omega+O((\omega-\theta)N^{-1+2\kappa}).
\end{align*}
Plugging this back into the kernel difference and using Taylor expansion of the sine function completes the derivation of the desired result.
\end{proof}

\subsubsection{Spatial-temporal decay of the extended kernel}\label{app:pf_kernel_est}

In this section, we prove Proposition~\ref{prop:kernel_est}. Before doing so, we first present several asymptotic results for Hermite functions that will be used in the proof. We use the following notation for Hermite functions:
\begin{align*}
\psi_k(x)=\textnormal{He}_k(x)\frac{e^{-x^2/4}}{(\sqrt{2\pi}k!)^{1/2}}
\end{align*}
where $\textnormal{He}_k(x)=(-1)^ke^{x^2/2}\frac{\D^k}{\D x^k}e^{-x^2/2}$ stands for the so called probabilist's Hermite polynomial. Several asymptotic results have been established for Hermite functions via various methods. Below we collect the ones used in this paper.

For convenience, define the scaled variable 
\begin{align*} 
\tilde x_t:=\frac{x\sqrt{N}}{2\sqrt{t}}
\end{align*}
for $t>0$ and function $F(x):=\big|\int_x^1 \sqrt{|1-t^2|}\D t\big|$. 

Using steepest descent analysis on the associated Riemann–Hilbert problem, it was shown in \cite{deift1999strong} (see also \cite{gustavsson2005gaussian}) that there exists an absolute constant $\delta_0>0$ such that for any $\delta<\delta_0$ and $k\in \Z^+$:
\begin{enumerate}
\item If $\tilde x_k\in [-1+\delta,1-\delta]$, (or equivalently $k\geq (\frac{x}{2(1-\delta)})^2 N$), then
\begin{align}\label{Region 1 asymptotics}
\psi_k(x\sqrt N)=\pi^{-1/2}k^{-1/4}(1-\tilde{x}_k^2)^{-1/4}\cos\big[2kF(\tilde{x}_k)-\frac{1}{2}\arcsin(\tilde{x}_k)\big]+O(k^{-5/4}),
\end{align}

\item If $\tilde x_k\in [1-\delta,1)$, (or equivalently $ (\frac{x}{2(1-\delta)})^2 N\geq k > (\frac{x}{2})^2 N$ and $x>0$), then

\begin{align}\label{Region 2 asymptotics}
\psi_k(x\sqrt N)=2^{-1/2}\pi^{-1/2}k^{-1/4}\bigg(&(\frac{1+\tilde{x}_k}{1-\tilde{x}_k})^{1/4}\Big(\cos\big[2kF(\tilde{x}_k)-\frac{\pi}{4}\big]+O(\frac{1}{kF(\tilde{x}_k)})\Big)(1+O(k^{-1}))\nonumber
\\
&-(\frac{1-\tilde{x}_k}{1+\tilde{x}_k})^{1/4}\Big(\sin\big[2kF(\tilde{x}_k)-\frac{\pi}{4}\big]+O(\frac{1}{kF(\tilde{x}_k)})\Big)(1+O(k^{-1}))\bigg),
\end{align}

\item If $\tilde x_k\in (1,1+\delta]$, (or equivalently $ (\frac{x}{2})^2 N> k \geq (\frac{x}{2(1+\delta)})^2  N$ and $x>0$), then
\begin{align}\label{Region 3 asymptotics}
\psi_k(x\sqrt N)= O\Big( k^{-1/4}e^{-2kF(\tilde{x}_k)}(\tilde{x}_k-1)^{-1/4} \big(1+\frac{1}{kF(\tilde{x}_k)}\big)\Big),
\end{align}

\item If $\tilde x_k\in [1+\delta,\infty)$, (or equivalently $ (\frac{x}{2(1+\delta)})^2  N\geq k  $ and $x>0$), then
\begin{align}\label{Region 4 asymptotics}
\psi_k(x\sqrt N)=O(k^{-1/4}e^{-2kF(\tilde{x_k})}).
\end{align}
\end{enumerate}
In these formulas, asymptotic expansions for the Airy function and its derivative from \cite[page 394]{olver1997asymptotics} were used to simplify expressions.

\begin{cor}\label{hermitefunctionpsiasymptotic}
There is a constant $C$ depending only on $\Upsilon<1$ such that for any pair of positive integers $k$ and $N$ satisfying $k\geq (1-\frac{\Upsilon}{10})N$, we have $|\psi_k(x\sqrt{N})|<C k^{-1/4}$ when $x\in[-2+\Upsilon,2-\Upsilon]$.
\end{cor}

In addition, the following rough bound from \cite[Theorem 1]{bonan1990estimates} will be used when estimating the extended kernel:
\begin{align}\label{sup norm bound of hermite functions}
\|\psi_k\|_\infty=O(k^{-1/12}).
\end{align}
Lastly, Mehler’s formula provides the following identity which is useful when analyzing multi-time correlation kernels:
\begin{align}\label{mehler'sformula}
\sum_{k=0}^\infty u^k\psi_k(x)\psi_k(y)=&\frac{1}{\sqrt{2\pi(1-u^2)}}\exp\Big(-\frac{1-u}{1+u}\frac{(x+y)^2}{8}-\frac{1+u}{1-u}\frac{(x-y)^2}{8}  \Big)
\end{align}
for any $u\in(0,1)$.

Equipped with these estimates, we are now prepared to prove Proposition \ref{prop:kernel_est}.

\begin{proof}[Proof of Proposition \ref{prop:kernel_est}]

To illustrate the structure of the argument, we will treat several cases in detail and briefly indicate how the others follow by similar reasoning.

When the difference between the time variables exceeds $N^{-1+\kappa}$, the kernel becomes significantly easier to estimate. There are essentially two cases to consider:

\begin{itemize}
\item[•] If $s-t>N^{-1+\kappa}$:

Let $a:=1-\frac{\Upsilon}{10}$ and $\tau:=s-t$. Then, for all sufficiently large $N$ and for every $x,y\in [-2+\Upsilon,2-\Upsilon]$:
\begin{align*} 
|K(t,x;s,y)|\lesssim& \sqrt{N}e^{-(N-\frac{1}{2})\tau/2}\Big(\sum_{k=0}^{\floor{aN}}e^{k\tau/2}+\sum_{k=\ceil{aN}}^{N-1}e^{k\tau/2}\frac{1}{k^{1/2}}\Big)\lesssim 1+\frac{1}{1-e^{-\tau}}=O(N^{1-\kappa})
\end{align*}
where we have used Corollary \ref{hermitefunctionpsiasymptotic} and equation \eqref{sup norm bound of hermite functions}.

\item[•] If  $t-s>N^{-1+\kappa}$:

Let $\tau:=t-s$ again. Then,
\begin{align*}
|K(t,x;s,y)|&\lesssim\sqrt{N} e^{(N-\frac{1}{2})\tau/2}\sum_{k=N}^{\infty}e^{-k\tau/2}\frac{1}{k^{1/2}}\lesssim \frac{1}{1-e^{-\tau}} =O(N^{1-\kappa})
\end{align*}
by Corollary \ref{hermitefunctionpsiasymptotic}.
\end{itemize}

The remainder of the proof is divided into four main cases:

\begin{enumerate}
\item $N^{-1+\kappa}<|x-y|<N^{-\kappa/4}$ and $0\leq s-t\leq N^{-1+\kappa/8}$,
\item $N^{-1+\kappa}<|x-y|<N^{-\kappa/4}$ and $N^{-1-\kappa/8} \leq t-s\leq N^{-1+\kappa/8}$,
\item $N^{-1+\kappa}<|x-y|<N^{-\kappa/4}$ and $0<t-s\leq N^{-1-\kappa/8}$,
\item $N^{-\kappa/4}<|x-y|$ and $|s-t|\leq N^{-1+\kappa/8}$.
\end{enumerate}

\textit{The first case.} Without loss of generality assume that $x>0$ and $|x|\geq |y|$ and let $\delta=\Upsilon/10$ be sufficiently small for \eqref{Region 1 asymptotics}-\eqref{Region 4 asymptotics} to hold.
We partition the index set $k\in\{1,2,\dots,N\}$ into seven (possibly overlapping) regions in order to apply the asymptotics \eqref{Region 1 asymptotics}-\eqref{Region 4 asymptotics}.\\

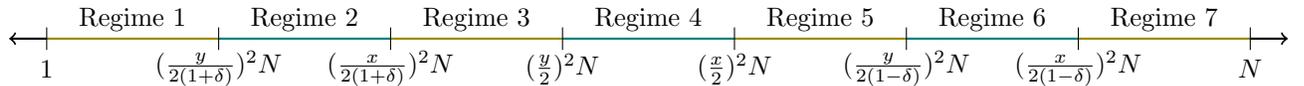
\begin{figure}[h]
    \centering
   \begin{tikzpicture}
\draw[<-,thick] (-0.5,0) -- (0,0);   
\draw[->,thick] (16,0) -- (16.5,0);  

\draw[thick, olive]        (0,0) -- (2.285714,0);        
\draw[thick, teal]       (2.285714,0) -- (4.571428,0); 
\draw[thick, olive]       (4.571428,0) -- (6.857142,0); 
\draw[thick, teal]     (6.857142,0) -- (9.142856,0); 
\draw[thick, olive]     (9.142856,0) -- (11.42857,0); 
\draw[thick, teal]      (11.42857,0) -- (13.714284,0);
\draw[thick, olive]      (13.714284,0) -- (16,0);      

\foreach \i/\dash/\regime in {
  0/$1$/Regime $1$, 
  1/$(\frac{y}{2(1+\delta)})^2N$/Regime $2$, 
  2/$(\frac{x}{2(1+\delta)})^2N$/Regime $3$, 
  3/$(\frac{y}{2})^2N$/Regime $4$, 
  4/$(\frac{x}{2})^2N$/Regime $5$, 
  5/$(\frac{y}{2(1-\delta)})^2N$/Regime $6$, 
  6/$(\frac{x}{2(1-\delta)})^2N$/Regime $7$} {
    \draw[thin] (\i*2.285714, -0.15) -- (\i*2.285714, 0.15); 
    \node at (\i*2.285714, -0.4) {\dash};                    
    \node at (\i*2.285714 + 1.142857, 0.25) {\regime};       
}

\draw[thin] (16, -0.15) -- (16, 0.15);
\node at (16, -0.4) {$N$};

\end{tikzpicture}
    \caption{Partition of the index set $\{1,2,\dots,N\}$ for the asymptotics \eqref{Region 1 asymptotics}-\eqref{Region 4 asymptotics}.}
\end{figure}

We will evaluate the sum for each regime separately and show that for each one $\sum (e^{(s-t)/2})^k\psi_k(x\sqrt{N})\psi_k(y\sqrt{N})$ is $O(N^{1/2-\kappa/8})$. For convenience, numerate and denote the intervals by $R_1,\dots,R_7$ as shown in the figure. Note that the sizes of the transition regions where $\psi_k(x\sqrt{N})$ and $\psi_k(y\sqrt{N})$ obey different asymptotics, i.e. regions $2$, $4$, and $6$, are at most of order $N^{1-\kappa/4}$ as $|x-y|<N^{-\kappa/4}$.

Starting with the regime 7 (i.e. $\tilde{x}_k,\tilde{y}_k\in(-1+\delta,1-\delta))$) substituting the asymptotic formula \eqref{Region 1 asymptotics} gives:
\begin{align}\label{region_1_summation}
&\pi\cdot \sum_{k\in R_7}(e^{(s-t)/2})^{k-(N-1)}\psi_k(x\sqrt{N})\psi_k(y\sqrt{N})=
\nonumber\\
&= O(1)+\sum_{k\in R_7}(e^{(s-t)/2})^{k-(N-1)}k^{-1/2}(1-\tilde{x}_k^2)^{-1/4}(1-\tilde{y}_k^2)^{-1/4}\cos\big[2kF(\tilde{x}_k)-\frac{1}{2}\arcsin(\tilde{x}_k)\big]\cos\big[2kF(\tilde{y}_k)-\frac{1}{2}\arcsin(\tilde{y}_k)\big]
\nonumber\\
&=O(N^{1/2-\kappa/4}\log N)+\sum_{k\in R_7}(e^{(s-t)/2})^{k-(N-1)}k^{-1/2}(1-\tilde{x}_k^2)^{-1/2}\cos\big[2kF(\tilde{x}_k)-\frac{1}{2}\arcsin(\tilde{x}_k)\big]\cos\big[2kF(\tilde{y}_k)-\frac{1}{2}\arcsin(\tilde{x}_k)\big]
\nonumber\\
&=O(N^{1/2-\kappa/4}\log N)+\frac{1}{2}\sum_{k\in R_7}(e^{(s-t)/2})^{k-(N-1)}k^{-1/2}(1-\tilde{x}_k^2)^{-1/2}\Big(\cos\big[2k(F(\tilde{x}_k)+F(\tilde{y}_k))-\arcsin(\tilde{x}_k)\big]
\nonumber\\
&\hspace{13cm}+\cos\big[2k(F(\tilde{y}_k)-F(\tilde{x}_k))\big]\Big)
\nonumber\\
&=O(N^{1/2-\kappa/4}\log N)+\frac{1}{2}\sum_{k\in R_7}(e^{(s-t)/2})^{k-(N-1)}k^{-1/2}(1-\tilde{x}_k^2)^{-1/2}\Big(\cos\big[2k(F(\tilde{x}_k)+F(\tilde{y}_k))\big](1-\tilde{x}_k^2)^{1/2}
\nonumber\\
&\hspace{6cm}+\sin\big[2k(F(\tilde{x}_k+F(\tilde{y}_k))\big]\tilde{x}_k+\cos\big[2k(F(\tilde{y}_k)-F(\tilde{x}_k))\big]\Big)
\end{align}
As every term in the sum is trivially $O(k^{-1/2})$, we can omit the sum over $R_7\cap[0,N^{1-\kappa/4}]$. From now on we will evaluate the sum over $R_7\cap[N^{1-\kappa/4},N]$ and for simplicity, we continue using the same letter $R_7$ instead of writing $R_7\cap[N^{1-\kappa/4},N]$ every time throughout the calculations of region $7$. We now evaluate the right hand side above by decomposing it into three sums.

Starting with the first sum in \eqref{region_1_summation}, i.e. $\sum_{k\in R_7}(e^{(s-t)/2})^{k-(N-1)}k^{-1/2}\cos\big[2k(F(\tilde{x}_k)+F(\tilde{y}_k))\big]$, define $b_t:=2t(F(\tilde{x}_t)+F(\tilde{y}_t))$ and $d_k=b_{k+1}-b_k$. Then
\begin{align*}
\frac{\D b_t}{\D t}=\arccos \tilde{x}_t +\arccos \tilde{y}_t=O(1), \quad \frac{\D^2b_t}{\D t^2}=\frac{1}{2t}\Big(\frac{\tilde{x}_t}{\sqrt{1-\tilde{x}_t^2}}+\frac{\tilde{y}_t}{\sqrt{1-\tilde{y}_t^2}}\Big)=O(\frac{1}{t}).
\end{align*}
Split the sum into groups of size $\asymp N^\kappa$ (so that we have at most $N^{1-\kappa}$ many groups), and denote them by $G_1,G_2,\dots$ (say that $2N^\kappa>|G_i|>N^\kappa$ for all $i$). We will evaluate the sum in a single group first. For simplicity let's say $G_i=\{m,m+1,\dots,m+A\}$, so $A\asymp N^\kappa$ and $m\geq N^{1-\kappa/4}$. Then we can conclude,
\begin{align*}
|\sum_{k=m}^{m+A}(e^{(s-t)/2})^{k-(N-1)}k^{-1/2}\cos b_k|\leq & |\sum_{k=m}^{m+A}(e^{(s-t)/2})^{-((m+A)-k)}k^{-1/2}\cos b_k|
\\
=&O(N^{-1+3\kappa})+|\sum_{k=m}^{m+A}k^{-1/2}\cos b_k|
\\
=&O(N^{-1+3\kappa})+O(m^{-3/2}N^{2\kappa})+m^{-1/2}|\sum_{k=m}^{m+A}\cos b_k|
\\
=&O(N^{-1+3\kappa})+O(m^{-3/2}N^{3\kappa})+m^{-1/2}|\sum_{k=m}^{m+A}\cos (b_m+(k-m)d_m)|
\\
=&O(N^{-1+3\kappa})+m^{-1/2}|\frac{1}{\sin(d_m/2)}\sin\Big(\frac{1}{2}(A+1)d_m\Big)\cos(b_m+\frac{Ad_m}{2})|
\\
=&O(N^{-1+3\kappa})+m^{-1/2}O(1) = O(N^{-1/2+\kappa/8})
\end{align*}
where we have used $|(e^{(s-t)/2})^{-((m+A)-k)}-1|=O(N^{-1+9\kappa/8})$ in the second line; $k^{-1/2}-m^{-1/2}=O(m^{-3/2}N^\kappa)$ in the third line; $\cos b_k=\cos (b_m+(k-m)d_m)+O(m^{-1}N^{2\kappa})$ in the fourth line; and $3\pi/2>d_m>2\arccos(1-\delta)$ in the last line. When this sum is added up for all groups $G_1,G_2,\dots$ (remember, there are $O(N^{1-\kappa})$ many groups of sizes $O(N^\kappa)$) total sum can be bounded by $O(N^{1/2-7\kappa/8})$.

The second sum in the expression \eqref{region_1_summation}, i.e. $\sum_{k\in R_1}(e^{(s-t)/2})^{k-(N-1)}k^{-1/2}\frac{\tilde{x}_k}{\sqrt{1-\tilde{x}_k^2}}\sin\big[2k(F(\tilde{x}_k)+F(\tilde{y}_k))\big]$ can be evaluated by the same way. We replace $\tilde x_k$'s by $\tilde x_m$ using $\tilde x_m-\tilde{x}_k=O(N^{-1+3\kappa/8})$ and $\frac{1}{\sqrt{1-\tilde{x}_k^2}}$ by $\frac{1}{\sqrt{1-\tilde{x}_m^2}}$ using $\frac{1}{\sqrt{1-\tilde{x}_m^2}}-\frac{1}{\sqrt{1-\tilde{x}_k^2}}=O(\tilde{x}_m-\tilde{x}_k)=O(N^{-1+3\kappa/8})$. But for the third sum in \eqref{region_1_summation}, i.e., $\sum_{k\in R_1}(e^{(s-t)/2})^{k-(N-1)}k^{-1/2}(1-\tilde{x}_k^2)^{-1/2}\cos\big[2k(F(\tilde{y}_k)-F(\tilde{x}_k))\big]$, a different approach is required, which we now develop.

Define $b_t:=2t(F(\tilde{y}_t)-F(\tilde{x}_t))>0$. Then
\begin{align*}
\frac{\D b_t}{\D t}=\arccos \tilde{y}_t -\arccos \tilde{x}_t>0, \quad \frac{\D^2b_t}{\D t^2}=\frac{1}{2t}\Big(\frac{\tilde{y}_t}{\sqrt{1-\tilde{y}_t^2}}-\frac{\tilde{x}_t}{\sqrt{1-\tilde{x}_t^2}}\Big)<0.
\end{align*}
So, $\frac{\D b_t}{\D t} \asymp (x-y)N^{1/2}t^{-1/2}$ that gives $b_{k+1}-b_k=O(k^{-1/2}N^{1/2-\kappa/4})$ as $|x-y|<N^{-\kappa/4}$. Therefore, the sum can be replaced with an integral up to negligible error:
\begin{align*}
\sum_{k\in R_7}(e^{(s-t)/2})^{k-(N-1)}k^{-1/2}(1-\tilde{x}_k^2)^{-1/2}\cos b_k=&\int_{R_7} (e^{(s-t)/2})^{r-(N-1)}r^{-1/2}(1-\tilde{x}_r^2)^{-1/2}\cos b_r \D r+O(N^{1/2-\kappa/4} \log N ).
\end{align*}

 Let $r_1$ be the minimum value in the set $R_7 \cap (2\pi\mathbb{Z}+3\pi/4)$. Define the other $r_2,r_3,\dots$ recursively by the relation $b_{r_{k+1}}-b_{r_k}=\pi$ until the $r_i$ value reaches $N$. More precisely, given $r_k<N$, if the $\hat{r}$ value satisfying $b_{\hat{r}}-b_{r_k}=\pi$ is less than $N$ define $r_{k+1}:=\hat{r}$, otherwise set $r_{k+1}:=N$. Write the obtained sequence as $r_0=\big((\frac{x}{2(1-\delta)})^2N\big) \vee N^{1-\kappa/4}<r_1<r_2<\dots<r_m<r_{m+1}=N$.

As $b_t$ is increasing with derivative of order $(x-y)N^{1/2}t^{-1/2}$, $r_{k+1}-r_k= O(\frac{1}{x-y}N^{-1/2}(r_{k+1})^{1/2})= O(N^{1-\kappa})$. Then, for each $k$, $|\int_{r_k}^{r_{k+1}} (e^{(s-t)/2})^{r-(N-1)}r^{-1/2}(1-\tilde{x}_r^2)^{-1/2}\cos b_r \D r|\leq \int_{r_k}^{r_{k+1}} r^{-1/2} \D r=O(N^{1/2-\kappa/2})$. That means the integral over $[r_k,r_{k+1}]$ is negligible for each single $k$, so a few of them can be neglected when necessary.

Define
\begin{align*}
I_k:=\int_{r_k}^{r_{k+1}} (e^{(s-t)/2})^{r-(N-1)}r^{-1/2}(1-\tilde{x}_r^2)^{-1/2}\cos b_r \D r.
\end{align*}
Observe that the integrand of $I_{2k+1}$ is always positive as $b_r\in 2\pi\mathbb{Z}+(-\pi/2,\pi/2)$ for every $r\in [r_{2k+1},r_{2k+2}]$ and similarly the integrand of $I_{2k}$ is always negative. We may therefore group the integrals into pairs and write
\begin{align*}
\int_{R_7} (e^{(s-t)/2})^{r-(N-1)}r^{-1/2}(1-\tilde{x}_r^2)^{-1/2}\cos b_r \D r=O(N^{1/2-\kappa/2})+(I_1+I_2)+(I_3+I_4)+\cdots.
\end{align*}
It can easily be verified that every integral $I_k$ in the sum can be replaced by $\hat{I}_k$ defined below, with a total error of order $N^{1/2-\kappa/4}$: 
\begin{align*}
\hat{I}_{2k+1}:=\int_{r_{2k+1}}^{r_{2k+2}} (e^{(s-t)/2})^{r-(N-1)}r_{2k+1}^{-1/2}(1-\tilde{x}_{r_{2k+1}}^2)^{-1/2}\cos b_r \D r,\\
\hat{I}_{2k+2}:=\int_{r_{2k+2}}^{r_{2k+3}} (e^{(s-t)/2})^{r-(N-1)}r_{2k+1}^{-1/2}(1-\tilde{x}_{r_{2k+1}}^2)^{-1/2}\cos b_r \D r.
\end{align*}
Now consider the sum of a pair 
\begin{align*}
\hat{I}_{2k+1}+\hat I_{2k+2}=r_{2k+1}^{-1/2}(1-\tilde{x}_{r_{2k+1}}^2)^{-1/2}\Big(\int_{r_{2k+1}}^{r_{2k+2}} (e^{(s-t)/2})^{r-(N-1)}\cos b_r \D r+\int_{r_{2k+2}}^{r_{2k+3}} (e^{(s-t)/2})^{r-(N-1)}\cos b_r \D r\Big).
\end{align*}
The integrand of the first integral is always positive, while that of the second is always negative and $b_r$ is increasing with decreasing derivative. Moreover, prefactor $(e^{(s-t)/2})^{r-(N-1)}$ with $r$. These facts together imply that each pair satisfies $\hat{I}_{2k+1}+\hat I_{2k+2}<0$. From this we can conclude that the main sum $\sum_{k\in R_7}(e^{(s-t)/2})^{k-(N-1)}k^{-1/2}(1-\tilde{x}_k^2)^{-1/2}\cos b_k$ is equal to $O(N^{1/2-\kappa/8})+(\textnormal{some strictly negative quantity})$. Applying the same argument with a parity shift, i.e. pairing $I_k$'s as $(I_2+I_3)+(I_3+I_4)+\cdots$ and defining $\hat{I}$'s accordingly and doing the same estimations tields that the same sum is equal to $O(N^{1/2-\kappa/8})+(\textnormal{some strictly positive quantity})$. Taken together, these bounds imply that the sum must be $O(N^{1/2-\kappa/8})$. Thus, we are done with the region 1.

The calculations regarding regime $5$ are similar, we will only point out a few differences here. 

Throughout the analysis for this region, we may assume without loss of generality that $y\geq 0$. If not, since Hermite polynomials with odd/even indices are odd/even, the repetition of the all arguments for the sum over odd terms and the sum over even terms separately solves the problem.

Moreover, in the asymptotic formula for regime $5$, i.e. \eqref{Region 2 asymptotics}, the error term may blow up as $\tilde{x}_k$ goes to $1$. In order to avoid this problem, we first split the interval into two phases: $k\in((\frac{x}{2})^2N,(\frac{x}{2})^2N+\sqrt{N})$ and $k\in((\frac{x}{2})^2N+\sqrt{N},(\frac{y}{2(1-\delta)})^2N)$. The first phase can be dealt with the sup norm bound \eqref{sup norm bound of hermite functions} easily. So, we consider $R_2=((\frac{x}{2})^2N+\sqrt{N},(\frac{y}{2(1-\delta)})^2N)$. Note that having $k>(\frac{x}{2})^2N+\sqrt{N}$ gives $1-\tilde{x}_k>\frac{1}{2\sqrt{N}}$.

Repeating the same steps with the calculations of regime $7$, gives three different sums again where the first two can be handled by splitting the region into groups of sizes $N^{1/3}$ and estimating the sum by turning it into a trigonometric sum of an arithmetic sequence as before. Lastly, the third sum can be approximated by an integral and the same argument works here without any change.

The calculations for other regions are much more simpler. For regime $6$, it suffices to apply the corresponding asymptotic formulas and follow analogous steps. In fact, the conclusion is even easier to obtain, since the size of $R_6$ is at most of order $N^{1-\kappa/4}$. For regions $1-4$, Hermite functions with $x$ variable exhibit rapid decay as in \eqref{Region 3 asymptotics} and \eqref{Region 4 asymptotics}. Specifically, for regimes $3$ and $4$, when the edges values of $k$'s are trimmed via \eqref{sup norm bound of hermite functions} as before, $\psi_k(x\sqrt N)$ can be bounded by $e^{-N^{1/6}}$; and for regimes $1$ and $2$, $\psi_k(x\sqrt N)=O(e^{-N^{1/3}})$. Hence, the desired estimate follows immediately in each case.\\

\textit{The second case.} The only difference in this case is that we start with an infinite sum as $t\geq s$.  As $|\tilde{x}_k|=\frac{|x|\sqrt{N}}{2\sqrt{k}}\leq 1-\delta$ for every $k\geq N$, the asymptotic \eqref{Region 1 asymptotics} is valid for every term. So, the problem can be reduced to a finite sum problem on $k\in[N,N^{1+\kappa/4}]$ because the sum over $k\geq N^{1+\kappa/4}$ is already negligible,
\begin{align*}
\Big|\sum_{k>N^{1+\kappa/4}}(e^{(t-s)/2})^{-(k-N)}&\psi_k(x\sqrt{N})\psi_k(y\sqrt{N})\Big|\leq o(1)+\sum_{k>N^{1+\kappa/4}}(e^{(t-s)/2})^{-(k-N)}k^{-1/2}
=O(N^{1/2-\kappa/8})
\end{align*}
where we have used $t-s\geq N^{-1-\kappa/8}$. On the other hand, the finite sum $\sum_{k=N}^{\lfloor N^{1+\kappa/4}\rfloor }(e^{(t-s)/2})^{-(k-N)}\psi_k(x\sqrt{N})\psi_k(y\sqrt{N})$ can be evaluated similar to the first case.\\

\textit{The third case.} Thanks to the submicroscopic bound on $t-s$, Mehler's formula (see equation \eqref{mehler'sformula}) can be used to reduce the problem into a finite sum problem as follows
\begin{multline*}
\sum_{k=N}^{\infty}(e^{-(t-s)/2})^{k-N}\psi_k(x\sqrt{N})\psi_k(y\sqrt{N})\asymp \sum_{k=0}^{\infty}e^{-k(t-s)/2}\psi_k(x\sqrt{N})\psi_k(y\sqrt{N})- \sum_{k=0}^{N}e^{-k(t-s)/2}\psi_k(x\sqrt{N})\psi_k(y\sqrt{N})
\\
=  O(N^{1/2-\kappa})- \sum_{k=0}^{N}e^{-k(t-s)/2}\psi_k(x\sqrt{N})\psi_k(y\sqrt{N})
\end{multline*}
where we have used Mehler's formula, $t-s\leq N^{-1-\kappa/8}$, and $|x-y|>N^{-1+\kappa}$. So, it suffices to bound the finite sum $\sum_{k=0}^{N}e^{-k(t-s)/2}\psi_k(x\sqrt{N})\psi_k(y\sqrt{N})$ which can be done similarly to the first case. \\

\textit{The fourth case.} In this case all calculations are similar, and mostly easier because of the large gap between $x$ and $y$. However, computations look messier due to the fact that any of $\frac{1}{2}\arcsin \tilde{y}_k$, $(1-\tilde{y}_k^2)^{-1/4}$, $(\frac{1-\tilde{y}_k}{1+\tilde{y}_k})^{1/4}$, or $(\frac{1+\tilde{y}_k}{1-\tilde{y}_k})^{1/4}$ terms cannot be replaced by their $x$ versions. The only fundamental difference is that integral approximations are not valid anymore. Nevertheless, grouping every region with appropriate group size and taking sums over the groups is sufficient for every case thanks to the large gap in between $x$ and $y$.
\end{proof}

\begin{remark}
It can also be possible to use the contour integral representation of the extended kernel (see \cite[(2.13)]{johansson2005non}) and use similar methods to \cite{johansson2005arctic} to prove the result.
\end{remark}

\subsection{Quadratic variation estimate in Lemma \ref{no outliers lemma}}\label{bounding quad var}
Here we provide the detailed calculations for the quadratic variation estimates used in Lemma \ref{no outliers lemma}, where we established the absence of outliers with overwhelming probability.

\begin{lemma}\label{stieltjes - bounding quadratic variation at a particular point outside}
Fixed $z=2+\frac{\varphi^{6}}{N^{2/3}}+i\frac{\varphi}{N^{2/3}}$, $j\in\{1,2,\dots,\lceil N^4(\log N)^2\rceil \}$, given $\boldsymbol\lambda(0)\in\tilde{\mathcal{G}}$, $t_j=\frac{j}{N^4}$ and 
\begin{align*}
\tau\leq \inf\Big\{u\in[0,t_j]: |h_j(u)|\geq e^{-u/2}\frac{1}{\varphi N\eta_{t_j-u}}\Big\}\wedge \inf\Big\{u\in[0,t_j]:\lambda_N(u)=2+\kappa\Big\}\wedge t_j
\end{align*} 
$(M_{t})_t$ is defined as in \eqref{Definition of M_t in the Stieltjes transform SDE}. Then, for all $0\leq s\leq t_j$:
\begin{align*}
\quadvar{M^{\tau}}_s \lesssim \frac{1}{\varphi^{5}N^2\eta_{t_j}^2}.
\end{align*}
\end{lemma}

\begin{proof} For convenience, we will write $t$ in place of $t_j$ throughout this proof. To streamline the argument, we divide the analysis into three separate cases:
\begin{enumerate}
\item[•] If $t-1\geq s$:

We have $\kappa_{t-u}-\kappa\asymp e^{(t-u)/2}$, and $\eta_{t-u}\asymp \frac{\eta}{\sqrt{\kappa}}e^{(t-u)/2}$ for all $u\leq s\wedge\tau$. Using these,
\begin{align*}
\quadvar{M^\tau}_s \leq& \frac{1}{N^3}\sum_k \int_0^{s\wedge\tau}\frac{e^{-u}}{|\lambda_k(u)-z_{t-u}|^4}\D u\lesssim \frac{1}{N^2} \int_0^{s\wedge\tau} \frac{e^{-u/2}}{((\kappa_{t-u}-\kappa)+\eta_{t-u})^2\eta_{t-u}}\Im(\tilde{m}_u(z_{t-u}))\D u
\\
\lesssim & \frac{1}{N^2} \int_0^{s} \frac{e^{-u/2}}{e^{3(t-u)/2}\frac{\eta}{\sqrt{\kappa}}}\Big(e^{-u/2}\Im(m_{\rm sc}(z_{t-u}))+e^{-u/2}\frac{1}{\varphi N\eta_{t-u}}\Big)\D u\lesssim \frac{e^{s-2t}}{N^2}.
\end{align*}

\item[•] If $t-\sqrt{\upkappa}\geq s\geq t-1$:

We still have $\upkappa_{t-u}-\upkappa\asymp e^{(t-u)/2}$ and $\eta_{t-u}\asymp \frac{\eta}{\sqrt{\upkappa}}e^{(t-u)/2}$ for all $u\leq t-1$. Moreover, $\upkappa_{t-u}-\upkappa\asymp (t-u)^2$ and $\eta_{t-u}\asymp (t-u)\frac{\eta}{\sqrt{\upkappa}}$ for all $s\geq  u \geq t-1$. Using these and the result we obtained in the first case, we get
\begin{align*}
\quadvar{M^\tau}_s \leq& \frac{1}{N^3}\sum_k \int_0^{(t-1)\wedge\tau}\frac{e^{-u}}{|\lambda_k(u)-z_{t-u}|^4}\D u + \frac{1}{N^3}\sum_k \int_{(t-1)\wedge\tau}^{s\wedge\tau}\frac{e^{-u}}{|\lambda_k(u)-z_{t-u}|^4}\D u
\\
\lesssim & \frac{e^{-t}}{N^2}+\frac{1}{N^2} \int_{(t-1)\wedge\tau}^{s\wedge\tau} \frac{e^{-u/2}}{((\upkappa_{t-u}-\upkappa)+\eta_{t-u})^2\eta_{t-u}}\Im(\tilde{m}_u(z_{t-u}))\D u
\\
\lesssim & \frac{e^{-t}}{N^2}+\frac{1}{N^2} \int_{(t-1)}^{s} \frac{e^{-u/2}}{(t-u)^5}\cdot e^{-t/2}\D u \lesssim \frac{e^{-t}}{N^2(t-s)^4} \leq \frac{e^{-t}}{N^2\upkappa^2} .
\end{align*}

\item[•] If $t\geq s\geq t-\sqrt{\upkappa}$:

We have $\upkappa_{t-u}-\upkappa\asymp (t-s)\sqrt{\upkappa}$ and $\eta_{t-u}\asymp\eta$ for all $s\geq u\geq t-\sqrt{\upkappa}$. Using these and the result we obtained in the second case, we get
\begin{align*}
\quadvar{M^\tau}_s \lesssim & \frac{e^{-t}}{N^2\upkappa^2}+\frac{1}{N^2} \int_{(t-\sqrt{\upkappa})\wedge\tau}^{s\wedge\tau} \frac{e^{-u/2}}{((\upkappa_{t-u}-\upkappa)+\eta_{t-u})^2\eta_{t-u}}\Im(\tilde{m}_u(z_{t-u}))\D u
\\
\lesssim & \frac{e^{-t}}{N^2\upkappa^2}+\frac{e^{-t}}{N^2} \int_{(t-\sqrt{\upkappa})}^{s} \frac{1}{((t-u)\sqrt{\upkappa}+\eta)^2\eta}\cdot \frac{1}{\varphi^2N^{1/3}}\D u
\\
\lesssim & \frac{e^{-t}}{N^2\upkappa^2}+\frac{e^{-t}}{N^2\upkappa}\Big(\frac{1}{(t-s)\sqrt{\upkappa}+\eta}-\frac{1}{\upkappa+\eta}\Big) \lesssim  \frac{e^{-t}}{N^2\upkappa\eta} .
\end{align*}
\end{enumerate}
\end{proof}

\subsection{Proof of Lemma \ref{lemma:simple_asymp_along_char}}\label{app:stable_asymp_along_char}

In the proof below, we use \eqref{eqn:basic_properties_of_m_sc} repeatedly without explicitly mentioning it each time.

\begin{proof}[Proof of Lemma \ref{lemma:simple_asymp_along_char}]
Assume wlog that $\Re(z)\geq 0$. We begin with the proof of \eqref{asymp_z_and_w_eqn1}, which we divide into two cases.

\begin{itemize}
\item If $z\notin B_{1/10}(2)$:

$B_{\frac{\eta}{100}}(z)\cap B_{\frac{1}{20}}(2)=\emptyset$ because for every $r\in B_{\frac{\eta}{100}}(z)$, $|z-2|\geq \eta>10|z-r|$ so $|r-2|\geq|z-2|-|z-r|>\frac{9|z-2|}{10}>\frac{1}{20}$. Hence, for all $r\in B_{\frac{\eta}{100}}(z)$, $|\frac{\D}{\D r}\sqrt{r^2-4}|^2=\frac{|r|^2}{|r-2||r+2|}\leq 49$ which yields $|\frac{\D}{\D r}\sqrt{r^2-4}|\leq 7$. Using this we obtain $|\sqrt{z^2-4}-\sqrt{w^2-4}|\leq 7|z-w|\leq\frac{\eta}{10}$ and as $\Im\sqrt{z^2-4}\geq\eta$, we get $\Im\sqrt{z^2-4}\asymp \Im\sqrt{w^2-4}$.

\item If $z\in B_{1/10}(2)$:

Note that if $|E|\leq2$, then $\Im m_{\rm sc}(z)\asymp \sqrt{\kappa+\eta}$, so $(\Im\sqrt{z^2-4})^2\asymp (\eta+\sqrt{\kappa+\eta})^2\asymp \eta+\kappa\asymp |z-2|$ and if $|E|\geq 2$ then $\Im m_{\rm sc}(z)\asymp \frac{\eta}{\sqrt{\kappa+\eta}}$, so $(\Im\sqrt{z^2-4})^2\asymp (\eta+\frac{\eta}{\sqrt{\kappa+\eta}})^2\asymp \frac{\eta^2}{\kappa+\eta}\asymp \frac{\eta^2}{|z-2|}$. Using these two, we can consider the four cases depending on $\Re z$ and $\Re w$ being less or greater than $2$ separately.
\begin{itemize}
\item  If $\Re z,\Re w\leq 2$: $(\Im\sqrt{z^2-4})^2\asymp|z-2|\asymp|w-2|\asymp(\Im\sqrt{w^2-4})^2$.

\item If $\Re z,\Re w\geq 2$: $(\Im\sqrt{z^2-4})^2\asymp\frac{\eta^2}{|z-2|}\asymp\frac{\eta_w^2}{|w-2|}\asymp(\Im\sqrt{w^2-4})^2$.

\item If $\Re z\geq 2\geq \Re w$ or $\Re z\leq 2\leq \Re w$: In this case $\kappa$ has to be less than $|z-w|\leq\frac{\eta}{100}$, so $\eta\asymp \eta+\kappa\asymp|z-2|$, and similarly $\eta_w\asymp|w-2|$. Therefore $(\Im\sqrt{z^2-4})^2\asymp\eta\asymp\eta_w\asymp(\Im\sqrt{w^2-4})^2$.
\end{itemize}
\end{itemize}
Hence, we have proved the equation \eqref{asymp_z_and_w_eqn1}. Using this, equation \eqref{asymp_z_and_w_eqn2} follows by:
\begin{align*}
\eta_t=\frac{e^{t/2}+e^{-t/2}}{2}\eta+\frac{e^{t/2}-e^{-t/2}}{2}\Im\sqrt{z^2-4}\asymp \frac{e^{t/2}+e^{-t/2}}{2}\eta_w+\frac{e^{t/2}-e^{-t/2}}{2}\Im\sqrt{w^2-4}\asymp \eta_{w_{t}}.
\end{align*}

The proof of inequality \eqref{eqn:eta_Imsqrt_inc} is straightforward substituting the explicit formula for $z_t$:
\begin{align*}
\eta_t\Im\sqrt{z_t^2-4}=\frac{e^{t}-e^{-t}}{4}(\eta^2+(\Im\sqrt{z^2-4})^2)+\frac{e^{t}+e^{-t}}{2}\eta\Im\sqrt{z^2-4}\geq e^t\eta\Im\sqrt{z^2-4}.
\end{align*}

Lastly, to prove the equation \eqref{asymp_z_and_w_eqn3} notice that it suffices to prove $|z-x|\asymp |w-x|$ for $x\in\R$ because $|z-y|=\inf_{x\in\R}\big(|z-x|+|x-y|\big)$ and $|w-y|=\inf_{x\in\R}\big(|w-x|+|x-y|\big)$. First, we show that $|\sqrt{z^2-4}-\sqrt{w^2-4}|\lesssim \Im\sqrt{z^2-4}$:
\begin{itemize}
\item If $z\notin B_{1/10}(2)$, then $|\sqrt{z^2-4}-\sqrt{w^2-4}|\leq \frac{\eta}{10}\leq\frac{\Im\sqrt{z^2-4}}{10}$.
\item If $z\in B_{1/10}(2)$, then
\begin{align*}
|\sqrt{z^2-4}-\sqrt{w^2-4}|=\frac{|z^2-w^2|}{|\sqrt{z^2-4}+\sqrt{w^2-4}|}\asymp \frac{|w-z||w+z|}{|\sqrt{z^2-4}|+|\sqrt{w^2-4}|}\lesssim \frac{\eta}{\sqrt{|z-2|}}\asymp \frac{\eta}{\sqrt{\eta}+\sqrt{\kappa}}
\end{align*}
as $\Re\sqrt{z^2-4},\Im\sqrt{z^2-4},\Re\sqrt{w^2-4},\Im\sqrt{w^2-4}\geq 0$. On the other hand, due to $\Im\sqrt{z^2-4}\geq \Im m_{\rm sc}(z)$ and \eqref{eqn:basic_properties_of_m_sc}, we have $\Im\sqrt{z^2-4}\gtrsim \frac{\eta}{\sqrt{\eta}+\sqrt{\kappa}}\asymp |\sqrt{z^2-4}-\sqrt{w^2-4}|$.
\end{itemize}
 Hence,
 \begin{multline*}
 |\Re(z_t-x)-\Re(w_t-x)|=|\Re(z_t-w_t)|\leq \frac{e^{t/2}+e^{-t/2}}{2}|\Re(z-w)|+\frac{e^{t/2}-e^{-t/2}}{2}|\sqrt{z^2-4}-\sqrt{w^2-4}|
 \\
 \lesssim  \frac{e^{t/2}+e^{-t/2}}{2}\eta+\frac{e^{t/2}-e^{-t/2}}{2}\Im\sqrt{z^2-4}\leq \eta_t \asymp\eta_{w_t}.
 \end{multline*}
Thus $|\Re(w_t-x)|\lesssim \eta_t +|\Re(z_t-x)|\asymp |z_t-x|$. As $|z_t-x|\gtrsim |\Re(w_t-x)|$ with $|z_t-x|\geq\eta_t\asymp\eta_{w_t}=\Im(w_t-x)$, we obtain $|z_t-x|\gtrsim|w_t-x|$. Similarly $|w_t-x|\gtrsim|z_t-x|$. Hence the equation \eqref{asymp_z_and_w_eqn3} is proven.
\end{proof}

\subsection{Proof of Lemma \ref{lemma:Uhlbrt_of_tildeU_at_z}}\label{app:pf_U_t and V_t}

\begin{proof} We already have these equalities when $t=0$, see \eqref{eqn:Chebyshev_polynomials_properties}. Assume $t>0$. Using the form of $\tilde{U}_n(2\cos\theta)=\frac{\sin((n+1)\theta)}{\sin\theta}$ and doing classical change of variables $y=2\cos\theta$ and
\begin{align}\label{eqn:uhlbrtU_N} 
\uhlbrt \tilde{U}_{n}(z)=2\int_{-2}^2 \frac{\tilde{U}_n(y)}{z-y}\rhosc (y)\D y=\frac{2}{\pi}\int_{0}^{\pi}\frac{\cos(n\theta)-\cos((n+2)\theta)}{z-2\cos\theta}\D \theta.
\end{align}
Writing cosines in exponential form and using change of variables $e^{i\theta}=w$ and $e^{-i\theta}=w$:
\begin{align*} 
\frac{2}{\pi}\int_{0}^{\pi}\frac{\cos(n\theta)}{z-2\cos\theta}\D \theta &= \int_{0}^{\pi}\frac{e^{in\theta}}{z-(e^{i\theta}+e^{-i\theta})}\D\theta+\int_{0}^{\pi}\frac{e^{-in\theta}}{z-(e^{i\theta}+e^{-i\theta})}\D\theta
\\
&=\frac{i}{\pi}\int_{\partial \mathbb{D}}\frac{w^{n}}{w^2-wz+1}\D w=\frac{2}{\sqrt{z^2-4}}\big(\frac{z-\sqrt{z^2-4}}{2}\big)^n.
\end{align*}
Substituting this into \eqref{eqn:uhlbrtU_N} we obtain $\uhlbrt \tilde{U}_{n}(z)=2\big(\frac{z-\sqrt{z^2-4}}{2}\big)^{n+1}$.
Denote $x=2\cos\alpha$. Plugging $x_t$ into $z$ finishes the proof of the first equality
\begin{multline*} 
\Re(\uhlbrt\tilde{U}_{n-1}(x_t))=\Re\Big(2e^{-tn/2}\big(\frac{x-\sqrt{x^2-4}}{2}\big)^{n}\Big)
=e^{-tn/2}\Big(\big(\frac{x-\sqrt{x^2-4}}{2}\big)^{n}+\big(\frac{x+\sqrt{x^2-4}}{2}\big)^{n}\Big)
\\
=e^{-tn/2}(e^{-i\alpha n}+e^{i\alpha n})=2e^{-tn/2}\tilde{T}_n(x)
\end{multline*}
where we have used $\tilde{T}_n(2\cos\alpha)=\cos (n\alpha)$.

The second equality in Lemma \ref{lemma:Uhlbrt_of_tildeU_at_z} can be shown similarly. Doing the same change of variables $y=2\cos\theta$ and using $y_t=e^{i\theta}e^{t/2}+e^{-i\theta}e^{-t/2}$ we write
\begin{align*} 
\Re(\int_{-2}^{2}&\frac{\tilde{T}_n(y)}{y_{t}-x}\rhoarcsin (y)\D y)=\frac{1}{2\pi}\Re\Big(\int_{0}^{\pi} \frac{e^{in\theta}+e^{-in\theta}}{e^{i\theta}e^{t/2}+e^{-i\theta}e^{-t/2}-x}\D\theta\Big)
\\
&=\frac{1}{4\pi}\int_{0}^{\pi} (e^{in\theta}+e^{-in\theta})\Big(\frac{1}{e^{i\theta}e^{t/2}+e^{-i\theta}e^{-t/2}-x}+\frac{1}{e^{-i\theta}e^{t/2}+e^{i\theta}e^{-t/2}-x}\Big)\D\theta
\\
&=\frac{-i}{4\pi}\int_{\partial\mathbb{D}} w^n\Big(\frac{1}{we^{t/2}+w^{-1}e^{-t/2}-x}+\frac{1}{w^{-1}e^{t/2}+we^{-t/2}-x}\Big)\frac{1}{w}\D w
\\
&=\frac{-i}{4\pi e^{t/2}}\int_{\partial\mathbb{D}} \frac{w^n}{(w-\frac{x+\sqrt{x^2-4}}{2e^{t/2}})(w-\frac{x-\sqrt{x^2-4}}{2e^{t/2}})}\D w+\frac{-i}{4\pi e^{-t/2}}\int_{\partial\mathbb{D}} \frac{w^n}{(w-\frac{x+\sqrt{x^2-4}}{2e^{-t/2}})(w-\frac{x-\sqrt{x^2-4}}{2e^{-t/2}})}\D w.
\end{align*}
As $e^{t/2}>1=|\frac{x+\sqrt{x^2-4}}{2}|=|\frac{x-\sqrt{x^2-4}}{2}|$, the poles of the second integral are outside $\bar{\mathbb{D}}$, so the integral is zero. Note that both poles of the first integral are in the unit disc, so we obtain
\begin{align*} 
\Re(\int_{-2}^{2}&\frac{\tilde{T}_n(y)}{y_{t}-x}\rhoarcsin (y)\D y)=\frac{1}{2\sqrt{x^2-4}}\Big((\frac{x+\sqrt{x^2-4}}{2e^{t/2}})^n-(\frac{x-\sqrt{x^2-4}}{2e^{t/2}})^n\Big)=e^{-tn/2}\frac{1}{2}\tilde{U}_{n-1}(x)
\end{align*}
where we have used $\tilde{U}_{n-1}(2\cos\alpha)=\frac{\sin(n\alpha)}{\sin\alpha}$.
\end{proof}

\subsection{Details of the proof of Proposition \ref{prop:FH_upper_bound}}\label{app:pf_of_C_approx}

\subsubsection{Single-time loop equation for $\CUE$}\label{app:loop_eqn_CUE}

Before presenting the technical lemmas that are used in the proof of Proposition \ref{prop:FH_upper_bound}, we introduce the circular analogues of the $\uhlbrt$, $\vhlbrt$ operators, $\mathcal{C}$ in the circular setting ($\CUE$) and recall some key properties. The loop equation in $\CUE$ (see Lemma 2.1 in \cite{lambert2021mesoscopic}) is given by (cf. \eqref{eqn:loop_eqn_by_change_of_variables})
 \begin{align}\label{eqn:loop_CUE}
\mathbb{E}_{\CUE,h}[S_N(\tilde\uhlbrt g)]=&\int_{\T} g(x)h'(x)\frac{1}{2\pi}\D x+\frac{1}{N}\mathbb{E}_{\CUE,h}[S_N(g h')]+\frac{1}{2N}\mathbb{E}_{\CUE,h}\Big[\int_{\T}\int_{\T} \frac{g(x)-g(y)}{\tan(\frac{x-y}{2})}\D\nu_N(x)\D\nu_N(y)\Big]
\end{align}
where $\mathbb{E}_{\CUE,h}[\cdot]=\mathbb{E}_{\CUE}[\cdot\frac{e^{S_N(h)}}{\mathbb{E}_\CUE[e^{S_N(h)}]}]$, $\D\nu_N(x):=\sum_i\delta_{\theta_i}(x)-N\frac{1}{2\pi}\D x$ and
\begin{align*} 
\thlbrt f(x):=\dashint_{\T}\frac{f(y)}{\tan(\frac{x-y}{2})}\frac{1}{2\pi}\D y
\end{align*}
(cf. the definition of $\vhlbrt$-transform). Denoting the Fourier coefficients of $f\in L^2(\T)$ by $\hat{f}_n=\int_\T f(x)e^{-inx}\frac{1}{2\pi}\D x$, it holds that $(\widehat{\thlbrt f})_n=-i\sgn(n)\hat{f}_n$ for all $n\in\Z$ where we adopt the convention $\sgn(0)=0$ (cf. equation \eqref{eqn:Chebyshev_polynomials_properties}).  Consequently, for all $f\in L^2(\T)$, we have $\thlbrt(-\thlbrt f)=f-\hat{f}_0$ a.e. Moreover, if $f\in \mathscr{C}^{1}(\T)$, the equality holds everywhere on $\T$ (cf. Proposition \ref{prop:inverse_of_uhlbrt}). Note that we use the same notation for the Fourier coefficients on $\T$ and Chebyshev coefficients on $(-2,2)$, but which one is meant will be clear  from the context.

 Additionally, we have the following Fourier expansions of $\log|1-e^{ix}|$ and $\arg$ functions, (on $\T$ we take $\arg=\Im\log(1-e^{ix})$ which is $\arg$ restricted on $[-\pi,\pi]$),
\begin{align*}
\log|1-e^{ix}|=\sum_{n\neq 0} \frac{-1}{2|n|} e^{ixn}, \quad \arg x=\sum_{n\neq 0} \frac{i}{2n} e^{ixn}
\end{align*}
from which it is easy to see $\thlbrt(\log|1-e^{i\cdot}|)=\arg$ and $\thlbrt(\arg)=-\log|1-e^{i\cdot}|$ (cf. equations \eqref{eqn:log_cheb_exp}-\eqref{eqn:v_trans_of_arg}). Lastly, the analogue of our definition $\mathcal{C}(f,h)$ at \eqref{eqn:defn_mathcal_C} is $\mathcal{\tilde{C}}(f,h):=\sum_{k=-\infty}^\infty |k| \hat{f}_k\hat{h}_{-k}=-\int_{\T} \thlbrt f(x)h'(x)\frac{1}{2\pi}\D x=-\int_{\T} \thlbrt h(x)f'(x)\frac{1}{2\pi}\D x$, which describes the covariance structure of the limit of the fluctuation field for $\CUE$. 

\subsubsection{Motivation  for the $\GUE-\CUE$ transition}\label{subsec:needforGUE-CUE}

At first glance, the transition from $\GUE$ to $\CUE$ may appear somewhat intricate and unmotivated. In the following paragraph, we provide a heuristic explanation of the need for this transition. The key point is that the loop equation for $\CUE$ allows us to handle a single singularity via a simple algebraic trick, something that does not carry over to the $\GUE$ setting. Heuristically, when $h$ has a single singularity and $g$ vanishes at that singularity, the loop equation in $\CUE$ reads,
\begin{align*} 
\mathbb{E}_{\CUE,h}[S_N(\tilde\uhlbrt g)]\approx &\int_{\T} g(x)h'(x)\frac{1}{2\pi}\D x
\end{align*}
where $\thlbrt$ satisfies $\thlbrt(-\thlbrt f)=f-\hat{f}_0$ and $\thlbrt 1=0$. This structure yields the following identity for any $c\in\R$,
\begin{align*} 
\mathbb{E}_{\CUE,h}[S_N(f)]=\mathbb{E}_{\CUE,h}[S_N(f-\hat{f}_0)]=\mathbb{E}_{\CUE,h}[S_N(\thlbrt (-\thlbrt f))]=\mathbb{E}_{\CUE,h}[S_N(\thlbrt (-\thlbrt f+c))].
\end{align*}
By choosing $c$ so that $-\thlbrt f+c$ vanishes at the singularity of $h$, we can apply the loop equation directly. 

 In contrast, for the $\GUE$, again assuming $g$ vanishes at the singularity of $h$, the analogous identity becomes (see Section \ref{subsec:multitime_loop_asymp_exploration}),
\begin{align*} 
\mathbb{E}_h[S_N(\uhlbrt g)]=&\int g(x)h'(x)\rhosc (x)\D x
\end{align*}
where the transforms $\uhlbrt$ and $\vhlbrt$ satisfy $\uhlbrt\vhlbrt f=f-\hat{f}_0$ and $\uhlbrt 1=x$. However, for $\GUE$, the corresponding identity is,
\begin{align*} 
\mathbb{E}_h[S_N(f)]=\mathbb{E}_h[S_N(f-\hat{f}_0)]=\mathbb{E}_h[S_N(\uhlbrt\vhlbrt f)]=\mathbb{E}_h[S_N(\uhlbrt(\vhlbrt f+c))]-c\cdot \mathbb{E}_h[S_N(x)]
\end{align*}
where the choice of $c$ to make $\vhlbrt f+c$ vanish at the singularity of $h$ introduces an extra term involving $\mathbb{E}_h[S_N(x)]$, which cannot be eliminated in a straightforward way. This obstruction motivates the shift from $\GUE$ to $\CUE$ in our analysis.

\subsubsection{Technical lemmas}\label{subsubsec:tech_lemmas}

\begin{lemma}\label{lemma:vhlbrt_cond_check} Let $C>0$ be large, $\Upsilon,\kappa>0$ be small fixed constants, and $\epsilon=N^{-1+\kappa}$. The following holds:
\begin{itemize}
\item[\textnormal{(a)}] $|\vhlbrt_{0}\log_{\epsilon}^{E}(E)| \lesssim 1$ uniformly in $E\in[-2+\Upsilon,2-\Upsilon]$.
\item[\textnormal{(b)}] $|\vhlbrt_{t}\log^{E}(x)|\lesssim 1$ and $|(\vhlbrt_t\log^E)'(x)|\lesssim 1+\frac{1}{t+|x-E|}$ uniformly in $t\geq 0$ and $E,x\in[-2+\Upsilon,2-\Upsilon]$ with $(t,x-E)\neq(0,0)$.
\item[\textnormal{(c)}] $|\vhlbrt_{t}\log_{\epsilon}^{E}(x)|\lesssim 1+\frac{\epsilon}{t+|E-x|}$ and $|(\vhlbrt_{t}\log_{\epsilon}^{E})'(y)|\lesssim 1+\frac{1}{t+|x-E|}+\frac{\epsilon}{(t+|x-E|)^2}$ uniformly in $t\geq0$ and $E,x\in[-2+\Upsilon,2-\Upsilon]$.
\item[\textnormal{(d)}] Let $q$ be defined as in \eqref{eqn:comp_q}. For $i\in\{r,\ell\}$, uniformly in the choice of $E,x\in[-2+\Upsilon,2-\Upsilon]$ and $t\geq0$ with $t+|E-x|\gtrsim N^{-1+2\kappa}$, we have $\vhlbrt_{t}q_{i,\epsilon}^{E}(x)=O(N^{-\kappa})$, $(\vhlbrt_{t}q_{i,\epsilon}^{E})'(x)=O( N^{1-3\kappa})$ and $|\vhlbrt_0q_{i,\epsilon}^{E}(E)|\asymp \log N$. 
\item[\textnormal{(e)}] If $f\in\mathscr{S}_{C,\kappa}$, then $|\vhlbrt_tf(x)|\lesssim\log N$ uniformly in $t\geq0$ and $x\in[-2+\Upsilon,2-\Upsilon]$.
\item[\textnormal{(f)}] Uniformly in the choice of $E,x\in[-2+\Upsilon,2-\Upsilon]$ and $t\geq0$, $\vhlbrt_{t}(\arg_{d,\epsilon}^{E})(x)=O(\log N)$ for $d\in\{r,\ell\}$, i.e. for any direction of regularization. If we also have the separation condition $t+|E-x|\gtrsim N^{-1+2\kappa}$, then $(\vhlbrt_{t}(\arg_{d,\epsilon}^{E}))'(x)=O( N^{1-2\kappa})$. 
\end{itemize}
\end{lemma}

\begin{proof} For convenience we will write $\rho$ for $\rhoarcsin $.

\noindent\textit{Proof of \textnormal{(a)}:} Using $\vhlbrt (1)=0$ on $[-2,2]$ we can write,
\begin{align*} 
|\vhlbrt \log_{\epsilon}^{E}(E)|=&\Big|\int\frac{(1-\chi_{\epsilon}^{E}(y))(\log^{E}y-\log(2\epsilon))}{E-y}\rho(y)\D y\Big|
\\
\leq &\Big|\int_{B_{\Upsilon/2}(E)\setminus B_{\epsilon}(E)} \frac{1-\chi^E_{\epsilon}(y)}{E-y}\log^Ey\rho(y)\D y\Big|+\Big|\log(2\epsilon)\cdot \dashint \frac{\chi^E_{\epsilon}(y)}{E-y}\rho(y)\D y\Big|+O(1)
\\
=&\Big|\int_{\epsilon}^{\Upsilon/2} (1-\chi_{\epsilon}(s)) \log s\frac{\rho(E-s)-\rho(E+s)}{s}\D s\Big|+|\log(2\epsilon)|\int_0^{2\epsilon} \big|\frac{\rho(E-s)-\rho(E+s)}{s}\big|\D s +O(1)=O(1).
\end{align*}

\noindent\textit{Proof of \textnormal{(b)}:}
For $E$ and $x$ in the bulk, denoting $E=2\cos\alpha$, $x=2\cos\beta$ with $\alpha,\beta\in(0,\pi)$, we have
\begin{align*} 
\vhlbrt_t\log^E(x)=-\sum_{n=1}^\infty\frac{e^{-tn/2}}{n}\frac{\sin (n\beta)\cos(n \alpha)}{\sin\beta}
\end{align*}
so, denoting $\theta=\alpha+\beta$ and $\omega=\alpha-\beta$, we get
\begin{align*} 
\frac{\D}{\D t} (\vhlbrt_t\log^E(x))=\frac{1}{4\sin\beta}\sum_{n=1}^\infty e^{-tn/2}(\sin(n\theta)+\sin(n\omega)).
\end{align*}
Note that uniformly in $s\in[0,1)$ and $\gamma\in[0,\pi]$, $|\sum_{n=1}^\infty s^n\sin(n\gamma)|=|\Im\sum_{n=1}^\infty s^ne^{in\gamma}|=\frac{|s\sin\gamma|}{|1-se^{i\gamma}|^2}\asymp \frac{s\cdot (\gamma\wedge(\pi-\gamma))}{(1-s+\gamma)^2}$. Thus, we get
\begin{align*} 
|\frac{\D}{\D t} (\vhlbrt_t\log^E(x))|\lesssim  \frac{e^{-t}\theta'}{(1-e^{-t}+\theta')^2}+ \frac{e^{-t}\omega'}{(1-e^{-t}+\omega')^2}
\end{align*}
where $\theta':=\theta\wedge|\pi-\theta|\wedge(2\pi-\theta)$ and $\omega':=(\omega+\pi)\wedge|\omega|\wedge(\pi-\omega)$. Then, for any $t\geq 0$,
\begin{align*} 
|\vhlbrt_t\log^E(x)|\leq |\vhlbrt_0\log^E(x)|+\int_{0}^{\infty}|\frac{\D}{\D s} (\vhlbrt_s\log^E(x))|\D s\lesssim 1+\int_{0}^{1} \frac{\theta}{(1-s+\theta)^2}+\frac{\omega}{(1-s+\omega)^2}\D s \lesssim 1
\end{align*}
where we have used equation \eqref{eqn:v_trans_of_log}. Furthermore, using this, we can bound the derivative of $\vhlbrt_t\log^E$ as well,
\begin{align*} 
(\vhlbrt_t\log^E)'(x)&=\sum_{n=1}^{\infty}e^{-tn/2}\frac{\tilde{T}_{n}(E)}{n}\frac{2n\tilde{T}_n(x)-x\tilde{U}_{n-1}(x)}{4-x^2}
=\frac{2}{4-x^2}\sum_{n=1}^{\infty}e^{-tn/2}\tilde{T}_{n}(E)\tilde{T}_{n}(x)-\frac{x}{4-x^2}\sum_{n=1}^{\infty}e^{-tn/2}\frac{\tilde{T}_{n}(E)}{n}\tilde{U}_{n-1}(x)
\\
&=\frac{1}{4-x^2}\sum_{n=1}^{\infty}e^{-tn/2}(\cos(n\theta)+\cos(n\omega))+\frac{x}{4-x^2}\vhlbrt_t\log^E(x)
\\
&\lesssim 1+\Big(\frac{1}{1-e^{-t/2+i\theta}}+\frac{1}{1-e^{-t/2-i\theta}}+\frac{1}{1-e^{-t/2+i\omega}}+\frac{1}{1-e^{-t/2-i\omega}}\Big)\lesssim 1+\frac{1}{t+|\omega|}\asymp 1+\frac{1}{t+|x-E|}.
\end{align*}

\noindent\textit{Proof of \textnormal{(c)}:}
Substituting the definition of $\log_{\epsilon}^E$ we get,
\begin{align*} 
\vhlbrt_t\log_{\epsilon}^{E}(x)=-\Re\int\frac{\log_{\epsilon}^{E}}{x-y_t}\rho(y)\D y=-\Re\int\frac{\log^{E}}{x-y_t}\rho(y)\D y+\Re\int\frac{\chi_{\epsilon}^{E}(y)(\log^{E}y-\log(2\epsilon))}{x-y_t}\rho(y)\D y.
\end{align*}
Then,
\begin{align*} 
|\vhlbrt_{t}\log_{\epsilon}^{E}(x)|&\leq |\vhlbrt_{t}\log^{E_j}(x)|+\Big|\int\frac{\chi_{\epsilon}^{E}(y)(\log^{E}y-\log(2\epsilon))}{x-y_t}\rho(y)\D y\Big|
\\
&\lesssim 1+\frac{1}{t+|E-x|}\int \chi_{\epsilon}^{E}(y)|\log^{E}y-\log(2\epsilon)|\rho(y)\D y\lesssim 1+\frac{\epsilon}{t+|E-x|}.
\end{align*}
Similarly, for the derivative, we obtain
\begin{align*} 
|(\vhlbrt_{t}\log_{\epsilon}^{E})'(x)|&\lesssim |(\vhlbrt_{t}\log^{E})'(x)|+\frac{1}{(t+|E-x|)^2}\int \chi_{\epsilon}^{E}(y)|\log^{E}y-\log(2\epsilon)|\rho(y)\D y
\\
&\lesssim 1+\frac{1}{t+|x-E|}+\frac{\epsilon}{(t+|x-E|)^2}.
\end{align*}

\noindent\textit{Proof of \textnormal{(d)}:} $\vhlbrt_0q_{i,\epsilon}^{E}(E)\asymp \log N$ is clear by the definition of function $q$. Moreover $|\vhlbrt_{t}q_{i,\epsilon}^{E}(x)|\leq \big|\int \frac{q_{i,\epsilon}^{E}(y)}{x-y_t}\rho(y)\D y\big|\lesssim N^{1-2\kappa}\int |q_{i,\epsilon}^{E_j}(y)|\rho(y)\D y= O(N^{-\kappa})$ and similarly $|(\vhlbrt_{t}q_{i,\epsilon}^{E})'(x)|\lesssim (N^{1-2\kappa})^2\int |q_{i,\epsilon}^{E_j}(y)|\rho(y)\D y=O(N^{1-3\kappa})$.\\

\noindent\textit{Proof of \textnormal{(e)}:} For any $f\in\mathscr{A}_{C,\varepsilon}$, $\sup_{t/\geq0,x\in[-2+\Upsilon,2-\Upsilon]}|\vhlbrt_t f(x)|\lesssim \sum_{n=1}^{1/\varepsilon}\varepsilon+\sum_{n=1/\varepsilon}^{\infty}\frac{1}{n^2\varepsilon}\lesssim 1$ by \cite[Theorem 4.2]{trefethen2008gauss}; and any $f\in\mathscr{S}_{C,\kappa}$ is sum of $O(\log N)$ many such functions.\\

\noindent\textit{Proof of \textnormal{(f)}:} Take any $d\in\{r,\ell\}$. The bound on $\vhlbrt_t\arg^{E}_{d,\epsilon}$ follows easily by the bounds on Chebyshev coefficients, $|\vhlbrt\arg_{d,\epsilon}^{E}(E)|\lesssim \sum_{n=1}^{1/\epsilon}\frac{1}{n}+\sum_{n=1/\epsilon}^{\infty}\frac{1}{n^2\epsilon}\lesssim \log N$. For the derivative when the separation condition is given, using $\vhlbrt_t \tilde{T}_1(x)=\tilde{U}_0(x)$ and $\vhlbrt_t(1)(x)=0$ we can write
\begin{align*} 
(\vhlbrt_t\arg^{E})'(x)=\frac{e^{-t/2}}{2}+2\cdot\Re\int_{-2}^{E}\frac{\pi/2}{(x-y_t)^2}\rho(y)\D y=\frac{e^{-t/2}}{2}-2\cdot\Re\int_{E}^{2}\frac{\pi/2}{(x-y_t)}\rho(y)\D y
\end{align*}
It is easy to see that if $x\geq E$, $\int_{-2}^{E}\frac{1}{(x-y_t)^2}\rho(y)\D y\lesssim N^{1-2\kappa}$ and if $x\leq E$, $\int_{E}^{2}\frac{1}{(x-y_t)^2}\rho(y)\D y\lesssim N^{1-2\kappa}$ given $t+|E-x|\gtrsim N^{-1+2\kappa}$. To obtain the same bounds on the regularization of $\arg$, it suffices to check $\vhlbrt$ transform of the difference, and since the difference is local, we have the trivial bound $(\vhlbrt_t(\arg^{E}-\arg^E_{d,\epsilon}))'(x)\lesssim \frac{\epsilon}{(N^{-1+2\kappa})^2}$ which completes the proof. 
\end{proof}

\begin{lemma}\label{lemma:log_arg_index_removals} Everything is given as in Proposition \ref{prop:FH_upper_bound}, then we have 
\begin{align*}
\mathcal{C}\Big(&\sum_{j}\big(\gamma_j(\log_{\Delta}^{E_j}-\frac{1}{2}\log_{\epsilon}^{E_j})+\beta_j(\arg_{d(j),\Delta}^{E_j}-\frac{1}{2}\arg_{d(j),\epsilon}^{E_j})+\frac{1}{2} \alpha_{j}q_{s(j),\epsilon}^{E_j}\big)(H_{t_j})+\sum_{i}\frac{1}{2}f_i(H_{s_i})
\\
&\hspace{6cm},\sum_{j} (\gamma_j\log_{\epsilon}^{E_j}+\beta_j\arg_{d(j),\epsilon}^{E_j}-\alpha_jq_{s(j),\epsilon}^{E_j})(H_{t_j})+\sum_{i}f_i(H_{s_i})\Big)
\\
=&\frac{1}{2}\mathcal{C}^{\circ}\Big(\sum_{j}(\gamma_j\log^{E_j}+\beta_j\arg^{E_j})(H_{t_j})+\sum_{i}f_i(H_{s_i})\Big)+\mathcal{L}_1+\mathcal{L}_2+\mathcal{L}_3+O(N^{-\kappa/2})
\\
&+\sum_{j}\Big(\frac{\gamma_j^2}{4}\log\frac{1}{2\Delta}+\frac{\gamma_j\beta_j}{4}\big(\frac{1-d(j)}{2}\pi-E_{j}-2\arccos\frac{E_{j}}{2}\big)+\frac{\beta_j^2}{8}\big(1-2\sqrt{4-E_j^2}-2\log\frac{\Delta}{4-E_j^2}\big)\Big)
\end{align*}
where $\mathcal{L}_1$, $\mathcal{L}_2$ and $\mathcal{L}_3$ are expressions in terms of local functions,
\begin{align*} 
\mathcal{L}_1&=\frac{-1}{2}\sum_{j}\mathcal{C}\big(\gamma_j(\log^{E_j}-\log_{\epsilon}^{E_j})+\beta_j(\arg^{E_j}-\arg_{d(j),\epsilon}^{E_j})+\alpha_jq_{s(j),\epsilon}^{E_j},\gamma_j(\log_{\Delta}^{E_j}-\log_{\epsilon}^{E_j})+\beta_j(\arg_{d(j),\Delta}^{E_j}-\arg_{d(j),\epsilon}^{E_j})+\alpha_jq_{s(j),\epsilon}^{E_j}\big)
\\
\mathcal{L}_2&=\frac{-1}{4}\sum_{j}\sqrt{4-E_j^2}\beta_j\cdot \vhlbrt\Big(\gamma_j(\log_{\Delta}^{E_j}-\log_{\epsilon}^{E_j})+\beta_j(\arg_{d(j),\Delta}^{E_j}-\arg_{d(j),\epsilon}^{E_j})+\alpha_jq_{s(j),\epsilon}^{E_j}\Big)(E_j)
\\
\mathcal{L}_3&=\frac{-1}{4}\sum_{j}\sqrt{4-E_j^2}\beta_j\cdot\vhlbrt\big(\gamma_j(-\log^{E_j}+\log_{\epsilon}^{E_j})+\beta_j(-\arg^{E_j}+\arg_{d(j),\epsilon}^{E_j})-\alpha_jq_{s(j),\epsilon}^{E_j}\big)(E_j+\Delta)
\end{align*}
\end{lemma}

\begin{proof} We have,
\begin{align*}
\mathcal{C}\Big(\sum_{j}&\big(\gamma_j(\log_{\Delta}^{E_j}-\frac{1}{2}\log_{\epsilon}^{E_j})+\beta_j(\arg_{d(j),\Delta}^{E_j}-\frac{1}{2}\arg_{d(j),\epsilon}^{E_j})+\frac{1}{2} \alpha_{j}q_{s(j),\epsilon}^{E_j}\big)(H_{t_j})+\sum_{i}\frac{1}{2}f_i(H_{s_i})
\\
&\hspace{6cm},\sum_{j} (\gamma_j\log_{\epsilon}^{E_j}+\beta_j\arg_{d(j),\epsilon}^{E_j}-\alpha_jq_{s(j),\epsilon}^{E_j})(H_{t_j})+\sum_{i}f_i(H_{s_i})\Big)
\\
=&\frac{1}{2}\sum_{i,m}\mathcal{C}_{|s_i-s_m|}(f_i,f_m)+\sum_{i\in I,j\in J} \mathcal{C}_{|s_i-t_j|}(f_i,\gamma_j\log_{\Delta}^{E_j}+\beta_j\arg_{d(j),\Delta}^{E_j}) 
\\
&+\frac{1}{2}\sum_{j}\mathcal{C}\big(\gamma_j\log_{\epsilon}^{E_j}+\beta_j\arg_{d(j),\epsilon}^{E_j}-\alpha_jq_{s(j),\epsilon}^{E_j},\gamma_j(\log_{\Delta}^{E_j}-\log_{\epsilon}^{E_j})+\beta_j(\arg_{d(j),\Delta}^{E_j}-\arg_{d(j),\epsilon}^{E_j})+\alpha_jq_{s(j),\epsilon}^{E_j}\big) 
\\
&+\frac{1}{2}\sum_{j}\mathcal{C}\big(\gamma_j\log_{\epsilon}^{E_j}+\beta_j\arg_{d(j),\epsilon}^{E_j}-\alpha_jq_{s(j),\epsilon}^{E_j},\gamma_j\log_{\Delta}^{E_j}+\beta_j\arg_{d(j),\Delta}^{E_j}\big) 
\\
&+\frac{1}{2}\sum_{j\neq k }\mathcal{C}_{|t_j-t_k|}\big(\gamma_j\log_{\epsilon}^{E_j}+\beta_j\arg_{d(j),\epsilon}^{E_j}-\alpha_jq_{s(j),\epsilon}^{E_j},\gamma_k(\log_{\Delta}^{E_k}-\log_{\epsilon}^{E_k})+\beta_k(\arg_{d(k),\Delta}^{E_k}-\arg_{d(k),\epsilon}^{E_k})+\alpha_kq_{s(k),\epsilon}^{E_k}\big) 
\\
&+\frac{1}{2}\sum_{j\neq k }\mathcal{C}_{|t_j-t_k|}\big(\gamma_j\log_{\epsilon}^{E_j}+\beta_j\arg_{d(j),\epsilon}^{E_j}-\alpha_jq_{s(j),\epsilon}^{E_j},\gamma_k\log_{\Delta}^{E_k}+\beta_k\arg_{d(k),\Delta}^{E_k}\big)
\\
=&(\RN1)+(\RN2)+(\RN3)+(\RN4)+(\RN5)+(\RN6)
\end{align*}
where $(\RN1)-(\RN6)$ are the summations. Next, we discuss every summation one-by-one, but before proceeding with that, we will list some helpful identities that will be used frequently throughout the proof.

We have the following bounds on the Chebyshev coefficients by \cite[Theorem 4.2]{trefethen2008gauss} for every natural number $n$: $|\widehat{(\log_{\epsilon}^{E_j})}_n|\lesssim \min(1,\frac{\log N}{n},\frac{1}{\epsilon n^2})$, $|\widehat{(\arg_{d(j),\epsilon}^{E_j})}_n|\lesssim \min(1,\frac{1}{n},\frac{1}{\epsilon n^2})$, $|\widehat{(q_{\alpha(j),\epsilon})}_n|\lesssim \min(N^{-1+\kappa},\frac{1}{n},\frac{N^{1-\kappa/2}}{n^2})$, $|\widehat{(\log_{\Delta}^{E_j}-\log_{\epsilon}^{E_j})}_n|\allowbreak\lesssim \min(\epsilon,\frac{\log N}{n},\frac{1}{\Delta n^2})$; by equation \eqref{eqn:log_cheb_exp}, $\widehat{(\log^{E_j})}_n\lesssim \frac{1}{n}$, so $|\widehat{(\log^{E_j}-\log_{\Delta}^{E_j})}_n|\lesssim \min(\Delta,\frac{\log N}{n})$; similarly by equation \eqref{eqn:arg_cheb_exp}, $\widehat{(\arg^{E_j})}_n\lesssim \frac{1}{n}$, so $|\widehat{(\arg^{E_j}-\arg_{d(j),\Delta}^{E_j})}_n|\lesssim \min(\Delta,\frac{1}{n})$. We assume that $f_i\in \mathscr{A}_{C,\epsilon_i}$ for some $\epsilon_i\in[N^{-1+\kappa},1]$, which gives $|\widehat{(f_i)}_n|\lesssim \min(\epsilon_i,\frac{1}{n^2\epsilon_i})$. Recall that $f_i$ is a sum of $O(\log N)$ many such functions, so the error term we obtain will be multiplied only by $\log N$ at the end. Moreover, by straightforward calculations, we have
\begin{gather} \label{eqn:C_log_and_arg}
\mathcal{C}(\log^E,g)=\frac{\hat{g}_0}{2}-\frac{g(E)}{2}, \quad \mathcal{C}_t(\arg^{E},g)=e^{-t/2}\frac{\hat{g}_1}{4}-\frac{\sqrt{4-E^2}}{2}\vhlbrt_t g(E)
\end{gather}
for all $E\in(-2,2)$ and sufficiently regular functions $g$. Using these we now simplify the sums $(\RN2)-(\RN6)$.\\

\noindent\textit{The second summation:} Using the Chebyshev expansion definition of $\mathcal{C}$ in \eqref{eqn:defn_mathcal_C},
\begin{align*} 
\mathcal{C}_{|s_i-t_j|}(f_i,\log^{E_j}-\log_{\Delta}^{E_j})\lesssim \sum_{n=1}^{1/\epsilon_i} n \epsilon_i\Delta +\sum_{n=1/\epsilon_i}^{1/\Delta} n \frac{1}{n^2\epsilon_i}\Delta +\sum_{n=1/\Delta}^{\infty} n \frac{1}{n^2\epsilon_i}\frac{\log N}{n}\ll N^{-\kappa},
\end{align*}
and similarly $\mathcal{C}_{|s_i-t_j|}(f_i,\arg^{E_j}-\arg_{d(j),\Delta}^{E_j})\ll N^{-\kappa}$. Hence, we can omit the $\Delta$ indices in $(\RN2)$ with an $O(N^{-\kappa})$ error.\\

\noindent\textit{The third summation:} By the definition of $\mathcal{L}_1$,
\begin{align*} 
(\RN3)=\mathcal{L}_1+\frac{1}{2}\sum_{j}\mathcal{C}\big(\gamma_j\log^{E_j}+\beta_j\arg^{E_j},\gamma_j(\log_{\Delta}^{E_j}-\log_{\epsilon}^{E_j})+\beta_j(\arg_{d(j),\Delta}^{E_j}-\arg_{d(j),\epsilon}^{E_j})+\alpha_jq_{s(j),\epsilon}^{E_j}\big)
\end{align*}
Substituting \eqref{eqn:C_log_and_arg} we obtain
\begin{align*} 
(\RN3)=\mathcal{L}_1+\mathcal{L}_2-\sum_{j}\frac{\gamma_j^2}{4}\log\frac{\Delta}{\epsilon}+O(N^{-1+\kappa}).
\end{align*}

\noindent\textit{The fourth summation:} First, we show that the $\Delta$ regularizations can be omitted again,
\begin{align*} 
\mathcal{C}(\gamma_j\log_{\epsilon}^{E_j}+\beta_j\arg_{d(j),\epsilon}^{E_j}-\alpha_jq_{s(j),\epsilon}^{E_j},\log^{E_j}-\log_{\Delta}^{E_j})&\lesssim \sum_{n=1}^{1/\epsilon} n \frac{\log N}{n} \Delta+\sum_{n=1/\epsilon}^{1/\Delta} n \frac{N^{1-\kappa/2}}{ n^2} \Delta+\sum_{n=1/\Delta}^{\infty}n\frac{N^{1-\kappa/2}}{n^2}\frac{\log N}{n}\ll N^{-\kappa/2}
\end{align*}
and similarly for $\arg^{E_j}-\arg_{d(j),\Delta}^{E_j}$ in the second entry. After we have $\log^{E_j}$ and $\arg^{E_j}$ in the second entry, applying equation \eqref{eqn:C_log_and_arg} gives,
\begin{align*} 
(\RN4)=&\frac{-1}{4}\sum_{j}\sqrt{4-E_j^2}\beta_j\cdot\vhlbrt(\gamma_j\log_{\epsilon}^{E_j}+\beta_j\arg_{d(j),\epsilon}^{E_j}-\alpha_jq_{s(j),\epsilon}^{E_j})(E_j)
\\
&-\sum_{j}\frac{\gamma_j^2}{4}\log(2\epsilon)+\sum_{j}\frac{\beta_j\gamma_j}{4}\big(\frac{1-d(j)}{2}\pi-E_{j}-\arccos\frac{E_{j}}{2}\big)+\sum_{j}\frac{\beta_j^2}{8}\big(1-\sqrt{4-E^2}\big)+O(N^{-\kappa/2})
\end{align*}
where, abusing the notation, we have assigned numerical values to $d(j)$ so that $d(j)=1$ when it is $r$, and $d(j)=-1$ when it is $\ell$. Next, we want to localize the $\vhlbrt$ term. Notice that in the bulk $[-2+\Upsilon,2-\Upsilon]$ we have $(\vhlbrt g)'\lesssim \sum_{n=1}^{\infty}n|\hat{g}_n|$ as $\tilde{U}_{n}'\lesssim n$. This easily leads to $\big(\vhlbrt(\gamma_j\log_{\epsilon}^{E_j}+\beta_j\arg_{d(j),\epsilon}^{E_j}-\alpha_jq_{s(j),\epsilon}^{E_j})\big)'\lesssim N^{1-\kappa/2}$ in the bulk. Hence,
\begin{align*} 
\frac{-1}{4}\sum_{j}\sqrt{4-E_j^2}\beta_j\vhlbrt&(\gamma_j\log_{\epsilon}^{E_j}+\beta_j\arg_{d(j),\epsilon}^{E_j}-\alpha_jq_{s(j),\epsilon}^{E_j})(E_j)
\\
&=\mathcal{L}_3-\frac{1}{4}\sum_{j}\sqrt{4-E_j^2}\beta_j\vhlbrt(\gamma_j\log^{E_j}+\beta_j\arg^{E_j})(E_j+\Delta)+O(N^{-\kappa/2})
\end{align*}
Substituting equations \eqref{eqn:v_trans_of_log} and \eqref{eqn:v_trans_of_arg} and denoting $E_j=2\cos\alpha$, $E_j+\Delta=2\cos\beta$ we obtain
\begin{align*} 
\vhlbrt(\gamma_j\log^{E_j}+\beta_j\arg^{E_j})(E_j+\Delta)&=\gamma_j\frac{\arccos(\frac{E_j+\Delta}{2})}{\sqrt{4-(E_j+\Delta)^2}}+\beta_j\Big(\frac{1}{2}+\frac{1}{\sqrt{4-(E_j+\Delta)^2}}\big(\log\sin\frac{\alpha-\beta}{2}-\log\sin\frac{\alpha+\beta}{2}\big)\Big)
\\
&=\gamma_j\frac{\arccos(\frac{E_j}{2})}{\sqrt{4-E_j^2}}+\beta_j\Big(\frac{1}{2}+\frac{\log\frac{\Delta}{4-E_j^2}}{\sqrt{4-E_j^2}}\Big)+O(\Delta).
\end{align*}

\noindent\textit{The fifth summation:} Applying integration by parts, $(\RN5)$ is equal to,
\begin{align*} 
- \int \Big(\vhlbrt_{|t_j-t_k|}\big(\gamma_j\log_{\epsilon}^{E_j}+\beta_j\arg_{d(j),\epsilon}^{E_j}-\alpha_jq_{s(j),\epsilon}^{E_j}\big)\cdot \rhosc \Big)'(x)\big(\gamma_k(\log_{\Delta}^{E_k}-\log_{\epsilon}^{E_k})+\beta_k(\arg_{d(k),\Delta}^{E_k}-\arg_{d(k),\epsilon}^{E_k})+\alpha_kq_{s(k),\epsilon}^{E_k}\big)(x)\D x
\\
\lesssim N^{1-2\kappa} \int \big|\gamma_k(\log_{\Delta}^{E_k}-\log_{\epsilon}^{E_k})(x)+\beta_k(\arg_{d(k),\Delta}^{E_k}-\arg_{d(k),\epsilon}^{E_k})(x)+\alpha_kq_{s(k),\epsilon}^{E_k}(x)\big|\D x\lesssim\epsilon N^{1-2\kappa}
\end{align*}
where we have used parts (c), (d) and (f) of Lemma \ref{lemma:vhlbrt_cond_check}. So, it is negligible. \\

\noindent\textit{The sixth summation:} First, note that the contribution from $q$ is negligible. The proof is similar to the what we've done for the fifth summation, i.e. integration by parts,
\begin{align*} 
\int\vhlbrt_{|t_j-t_k|}\big(\gamma_k\log_{\Delta}^{E_k}+\beta_k\arg_{d(k),\Delta}^{E_k}\big)(x)(q_{s(j),\epsilon}^{E_j})'(x)\rhosc (x)\D x \lesssim \epsilon N^{1-2\kappa}.
\end{align*}
We claim that all the regularization indices can be removed up to negligible error. We first remove the $\epsilon$ index of $\arg$, using integration by parts
\begin{align*} 
\mathcal{C}_{|t_j-t_k|}&\big(\arg_{d(j),\Delta^2}^{E_j}-\arg_{d(j),\epsilon}^{E_j},\gamma_k\log_{\Delta}^{E_k}+\beta_k\arg_{d(k),\Delta}^{E_k}\big)
\\
&=-\int \Big(\vhlbrt_{|t_j-t_k|}\big(\gamma_k\log_{\Delta}^{E_k}+\beta_k\arg_{d(k),\Delta}^{E_k}\big)\cdot \rhosc \Big)'(x)\big(\arg_{d(j),\Delta^2}^{E_j}-\arg_{d(j),\epsilon}^{E_j}\big)(x)\D x\lesssim \epsilon N^{1-2\kappa}
\end{align*}
and the bounds on Chebyshev coefficients, 
\begin{align*} 
\mathcal{C}_{|t_j-t_k|}\big(\arg^{E_j}-\arg_{d(j),\Delta^2}^{E_j},\gamma_k\log_{\Delta}^{E_k}+\beta_k\arg_{d(k),\Delta}^{E_k}\big)
\lesssim \sum_{n=1}^{1/\Delta}n \Delta^2\frac{\log N}{n}+\sum_{n=1/\Delta}^{1/\Delta^2}n \Delta^2\frac{1}{\Delta n^2}+\sum_{n=1/\Delta^2}^{\infty}n \frac{1}{n}\frac{1}{\Delta n^2}\ll N^{-1}.
\end{align*}

Secondly, we remove the $\Delta$ index of $\arg$ by the bounds on Chebyshev coefficients,
\begin{align*} 
\mathcal{C}_{|t_j-t_k|}\big(\log_{\epsilon}^{E_j},\arg^{E_k}-\arg_{d(k),\Delta}^{E_k}\big)\lesssim \sum_{n=1}^{1/\epsilon}n\frac{\log N}{n}\Delta+\sum_{n=1/\epsilon}^{1/\Delta}n\frac{1}{\epsilon n^2}\Delta+\sum_{1/\Delta}^{\infty}n\frac{1}{\epsilon n^2}\frac{1}{n}\ll N^{-\kappa}
\end{align*}
and by equation \eqref{eqn:C_log_and_arg},
\begin{align*} 
\mathcal{C}_{|t_j-t_k|}\big(\arg^{E_j},\arg^{E_k}-\arg_{d(k),\Delta}^{E_k}\big)=-\frac{\sqrt{4-E_j^2}}{2}\vhlbrt_{|t_j-t_k|}(\arg^{E_k}-\arg_{d(k),\Delta}^{E_k})(E_j)+O(\Delta)\lesssim \Delta N^{1-2\kappa}.
\end{align*}
where the second inequality is due to the separation condition between $(t_{j},E_{j})$ and $(t_k,E_k)$. 

Thirdly, the removal of the $\epsilon$ regularization of $\log$ follows similarly:
\begin{align*} 
\mathcal{C}_{|t_j-t_k|}\big(\log_{\Delta^2}^{E_j}-\log_{\epsilon}^{E_j},\log_{\Delta}^{E_k}\big)
=-\int \big((\vhlbrt_{|t_j-t_k|}\log_{\Delta}^{E_k})\cdot \rhosc \big)'(x)\big(\log_{\Delta^2}^{E_j}-\log_{\epsilon}^{E_j}\big)(x)\D x\lesssim \epsilon N^{1-2\kappa},
\end{align*}
in addition to $\mathcal{C}_{|t_j-t_k|}\big(\log^{E_j}-\log_{\Delta^2}^{E_j},\log_{\Delta}^{E_k}\big)\lesssim N^{-\kappa}$ by Chebyshev coefficients calculations as before, and
\begin{align*} 
\mathcal{C}_{|t_j-t_k|}\big(\log^{E_j}-\log_{\epsilon}^{E_j},\arg^{E_k}\big)=-\frac{-\sqrt{4-E_k^2}}{2}\vhlbrt_{|t_j-t_k|}(\log^{E_j}-\log_{\epsilon}^{E_j})(E_k)+ O(\epsilon)\lesssim \epsilon N^{1-2\kappa}
\end{align*}
again due to the separation condition between $(t_{j},E_{j})$ and $(t_k,E_k)$.

Lastly, we discuss the removal of the $\Delta$ regularization index of $\log$, i.e. we are left to prove that $\mathcal{C}_{|t_j-t_k|}\big(\gamma_j\log^{E_j}+\beta_j\arg^{E_j},\log^{E_k}-\log_{\Delta}^{E_k}\big)$ is negligible. The term with $\arg$ can be bounded by the same way, due to the separation condition. For the other term, we use the symmetry of $\log$ function as follows:
\begin{align*} 
\mathcal{C}_{|t_j-t_k|}(\log^{E_j},\log^{E_k}-\log_{\Delta}^{E_k})&=\dashint\vhlbrt_{|t_j-t_k|}\log^{E_j}(x)(\log^{E_k}-\log_{\Delta}^{E_k})'(x)\rhosc (x)\D x 
\\
\lesssim & \int_{0}^{2\Delta} \Big|\frac{(\vhlbrt_{|t_j-t_k|}\log^{E_j}\cdot\rhosc )(E_k+x)-(\vhlbrt_{|t_j-t_k|}\log^{E_j}\cdot\rhosc )(E_k-x)}{x}\Big|\D x\lesssim \Delta N^{1-2\kappa}.
\end{align*}

Combining all these simplifications of $(\RN2)-(\RN6)$ concludes the proof of the lemma.
\end{proof}

\begin{lemma}\label{lemma:CUE_loop_app}
Everything is given as in the proof of Proposition \ref{prop:FH_upper_bound} and Appendix \ref{app:loop_eqn_CUE}, then for every $j=1,\dots,J$,
\begin{align*} 
&\log \Ex_{\CUE}\Big[\exp\Big( S_N\big(\gamma_j(\log_{\tDelta_j}-\log_{\tepsilon_j})|1-e^{i\cdot}|+\beta_j(\arg_{d(j),\tDelta_j}-\arg_{d(j),\tepsilon_j})+\alpha_jq_{s(j),\tepsilon_j}\big)\Big)\Big]
\\
&=\log\Ex_{\CUE}\Big[|\det(\Id-U)|^{\gamma_j}e^{\beta_j\Im\log\det(\Id-U)}\Big]+O(N^{-\kappa})
\\
&-\tilde{\mathcal{C}}\Big(\gamma_j(\log_{\tDelta_j}-\frac{1}{2}\log_{\tepsilon_j})|1-e^{i\cdot}|+\beta_j(\arg_{d(j),\tDelta_j}-\frac{1}{2}\arg_{d(j),\tepsilon_j})+\frac{1}{2}\alpha_jq_{s(j),\tepsilon_j}   ,   \gamma_j\log_{\tepsilon_j}|1-e^{i\cdot}|+\beta_j\arg_{d(j),\tepsilon_j}-\alpha_jq_{s(j),\tepsilon_j}\Big).
\end{align*}

\end{lemma}
\begin{proof}
This is a straightforward consequence of the single-time loop equation \eqref{eqn:loop_CUE}. The details have been worked out in the appendix of \cite{bourgade2022liouville}, albeit in a slightly minor variations in constants. The same proof applies here.
\end{proof}

\begin{lemma}\label{lemma:log_arg_index_removal_circular} Everything is given as in the proof of Proposition \ref{prop:FH_upper_bound} and Appendix \ref{app:loop_eqn_CUE}, then we have
\begin{align*} 
\sum_{j}\tilde{\mathcal{C}}\big(&\gamma_j(\log_{\tDelta_j}-\frac{1}{2}\log_{\tepsilon_j})|1-e^{i\cdot}|+\beta_j(\arg_{d(j),\tDelta_j}-\frac{1}{2}\arg_{d(j),\tepsilon_j})+\frac{1}{2}\alpha_jq_{s(j),\tepsilon_j},\gamma_j\log_{\tepsilon_j}|1-e^{i\cdot}|+\beta_j\arg_{d(j),\tepsilon_j}-\alpha_jq_{s(j),\tepsilon_j}\big)
\\
=&\tilde{\mathcal{L}}_1+\tilde{\mathcal{L}}_2+\tilde{\mathcal{L}}_3+O(N^{-\kappa/2}) +\sum_{j}\Big(\frac{\gamma_j^2}{4}\big(\log\frac{1}{2\Delta}-\log(\sqrt{4-E_j^2})\big)-\frac{\beta_j\gamma_j}{4}\big(\frac{1+d(j)}{2}\pi\big)-\frac{\beta_j^2}{4}\big(\log\Delta+\log\sqrt{4-E_J^2}\big)\Big)
\end{align*}
where 
\begin{align*} 
\tilde{\mathcal{L}}_1&=\frac{-1}{2}\sum_{j}\tilde{\mathcal{C}}\big(\gamma_j(\log-\log_{\tepsilon_j})|1-e^{i\cdot}|+\beta_j(\arg-\arg_{d(j),\tepsilon_j})+\alpha_jq_{s(j),\tepsilon_j}
\\
&\hspace{5cm},\gamma_j(\log_{\tDelta_j}-\log_{\tepsilon_j})|1-e^{i\cdot}|+\beta_j(\arg_{d(j),\tDelta_j}-\arg_{d(j),\tepsilon_j})+\alpha_jq_{s(j),\tepsilon_j}\big)
\\
\tilde{\mathcal{L}}_2&=\frac{1}{4}\sum_{j}\beta_j \cdot \thlbrt\big(\gamma_j(\log_{\tDelta_j}-\log_{\tepsilon_j})|1-e^{i\cdot}|+\beta_j(\arg_{d(j),\tDelta_j}-\arg_{d(j),\tepsilon_j})+\alpha_jq_{s(j),\tepsilon_j}\big)(0)
\\
\tilde{\mathcal{L}}_3&=\frac{1}{4}\sum_{j}\beta_j\cdot \thlbrt\big(\gamma_j(-\log+\log_{\tepsilon_j})|1-e^{i\cdot}|+\beta_j(-\arg+\arg_{d(j),\tepsilon_j})-\alpha_jq_{s(j),\tepsilon_j}\big)(\tDelta_j)
\end{align*}
\end{lemma}

\begin{proof}
The proof is identical to the evaluations of the diagonal terms in the proof of Lemma \ref{lemma:log_arg_index_removals}, using 
\begin{gather*}
\tilde{\mathcal{C}}(\log^E,g)=\frac{\hat{g}_0}{2}-\frac{g(0)}{2}, \quad \tilde{\mathcal{C}}(\arg^{E},g)=\frac{1}{2}\thlbrt g(0).
\end{gather*}
\end{proof}

\begin{lemma}\label{lemma:L_cancellation}
Given $\mathcal{L}_1,\mathcal{L}_2,\mathcal{L}_3,\tilde{\mathcal{L}}_1,\tilde{\mathcal{L}}_2,\tilde{\mathcal{L}}_3$ as in the previous lemmas, we have $\mathcal{L}_i=\tilde{\mathcal{L}}_i+O(N^{-\kappa})$ for all $i=1,2,3$.
\end{lemma}
\begin{proof} All the proofs simply follow by change of variables and Taylor expansions. We prove the equality of each term in the summations $\mathcal{L}$'s and $\tilde{\mathcal{L}}$'s separately, up to a negligible error. For notational convenience, we omit the indices related to $j$ and define
\begin{align*} 
f^{E}_{\theta}:=\gamma(\log^{E}-\log_{\epsilon\theta}^{E})+\beta(\arg^{E}-\arg_{\epsilon\theta}^{E})+\alpha q_{\epsilon\theta}^{E},\quad \tilde{f}_{\theta}:=\gamma(\log-\log_{\epsilon\theta})|1-e^{i\cdot}|+\beta(\arg-\arg_{\epsilon\theta})+\alpha q_{\epsilon\theta}
\\
g^{E}_{\theta}:=\gamma(\log_{\Delta\theta}^{E}-\log_{\epsilon\theta}^{E})+\beta(\arg_{\Delta\theta}^{E}-\arg_{\epsilon\theta}^{E})+\alpha q_{\epsilon\theta}^{E},\quad \tilde{g}_{\theta}:=\gamma(\log_{\Delta\theta}-\log_{\epsilon\theta})|1-e^{i\cdot}|+\beta(\arg_{\Delta\theta}-\arg_{\epsilon\theta})+\alpha q_{\epsilon\theta}
\end{align*}
where all the $\arg$ regularizations are either left or right and as usual, $f^{E}=f^{E}_{1}$, $f=f^{0}$ (similarly for $g$).\\

\noindent\textit{Proof of $\mathcal{L}_1=\tilde{\mathcal{L}}_1+O(N^{-\kappa})$:} For single time expressions, $\mathcal{C}(\cdot,\cdot)$ and $\tilde{\mathcal{C}}(\cdot,\cdot)$ coincides with the $H^{1/2}$-inner product on $(-2,2)$ and $\T$ respectively, and the following integral representations are well known (e.g. see (2.23) and (3.4) in \cite{forrester2023review}):
\begin{align*}
\mathcal{C}(f,g)=\frac{1}{4}\sum_{k=0}^\infty k\hat{f}_k\hat{g}_k=\frac{1}{4\pi^2}\int_{-2}^2\int_{-2}^2 \frac{(f(x)-f(y))(g(x)-g(y))}{(x-y)^2} \frac{4-xy}{\sqrt{4-x^2}\sqrt{4-y^2}}\D x\D y 
\end{align*} 
and 
\begin{align*}
\tilde{\mathcal{C}}(f,g)=\sum_{k=-\infty}^\infty |k| \hat{f}_k\hat{g}_{-k}=\frac{1}{16\pi^2}\int_{-\pi}^{\pi}\int_{-\pi}^{\pi} \frac{(f(x)-f(y))(g(x)-g(y))}{\sin^2(\frac{x-y}{2})}\D x\D y.
\end{align*}
Using these, denoting $\theta=2\pi\rhosc (E)$, we obtain
\begin{align*} 
4\pi^2 \mathcal{C}(f^{E},g^{E}) &=\int_{[E-N^{-1+2\kappa},E+N^{-1+2\kappa}]^2}\frac{(f^{E}(x)-f^{E}(y))(g^{E}(x)-g^{E}(y))}{(x-y)^2} \big(1+O((x-y)^2)\big)\D x\D y+O(N^{-\kappa/2})
\\
&=\int_{[-\theta N^{-1+2\kappa},\theta N^{-1+2\kappa}]^2}\frac{(f_{\theta}(x)-f_{\theta}(y))(g_{\theta}(x)-g_{\theta}(y))}{(x-y)^2} \D x\D y+O(N^{-\kappa/2})
\\
&=\frac{1}{4}\int_{[-\theta N^{-1+2\kappa},\theta N^{-1+2\kappa}]^2}\frac{(f_{\theta}(x)-f_{\theta}(y))(g_{\theta}(x)-g_{\theta}(y))}{\sin^2\frac{x-y}{2}} \D x\D y+O(N^{-\kappa/2})
\\
&=\frac{1}{4}\int_{[-\theta N^{-1+2\kappa},\theta N^{-1+2\kappa}]^2}\frac{(\tilde{f}_{\theta}(x)-\tilde{f}_{\theta}(y))(\tilde{g}_{\theta}(x)-\tilde{g}_{\theta}(y))}{\sin^2\frac{x-y}{2}} \D x\D y+O(N^{-\kappa/2})
\\
&=\frac{1}{4}\int_{[-\pi.\pi]^2}\frac{(\tilde{f}_{\theta}(x)-\tilde{f}_{\theta}(y))(\tilde{g}_{\theta}(x)-\tilde{g}_{\theta}(y))}{\sin^2\frac{x-y}{2}} \D x\D y+O(N^{-\kappa})=4\pi^2  \mathcal{C}(\tilde{f}_{\theta},\tilde{g}_{\theta})
\end{align*}
as desired.\\

\noindent\textit{Proof of $\mathcal{L}_2=\tilde{\mathcal{L}}_2+O(N^{-\kappa})$:} This calculation is more straightforward,
\begin{align*} 
-&\sqrt{4-E^2}\cdot \vhlbrt g^{E}(E)=\sqrt{4-E^2}\cdot\dashint\frac{g^E(y)}{E-y}\rhoarcsin (y)\D y=-\sqrt{4-E^2}\cdot\dashint_{-2\epsilon}^{2\epsilon}\frac{g(y)}{y}\rhoarcsin (y+E)\D y
\\
&=-\dashint_{-2\epsilon}^{2\epsilon}\frac{g(y)}{y}\frac{1}{\pi}\D y+O(N^{-\kappa})=-\dashint\frac{\tilde{g}_{\theta}(y)}{y}\frac{1}{\pi}\D y+O(N^{-\kappa})=\dashint\frac{\tilde{g}_\theta(y)}{\tan\frac{0-y}{2}}\frac{1}{2\pi}\D y+O(N^{-\kappa})=\tilde{\uhlbrt}(\tilde{g}_{\theta})(0)+O(N^{-\kappa}).
\end{align*}

\noindent\textit{Proof of $\mathcal{L}_3=\tilde{\mathcal{L}}_3+O(N^{-\kappa})$:} Similar to the previous one, we obtain
\begin{align*} 
-\sqrt{4-E^2}\cdot \vhlbrt f^{E}(E+\Delta)=\sqrt{4-E^2}\cdot\dashint_{-2\epsilon}^{2\epsilon}\frac{f(y)}{\Delta-y}\rhoarcsin (y+E)\D y=\dashint_{-2\epsilon}^{2\epsilon}\frac{f(y)}{\Delta-y}\frac{1}{\pi}\D y+O(N^{-\kappa})
\\
=\int_{-2\theta\epsilon}^{2\theta\epsilon}\frac{\tilde{f}_{\theta}(y)}{\Delta\theta-y}\frac{1}{\pi}\D y+O(N^{-\kappa})=\int_{-2\theta\epsilon}^{2\theta\epsilon}\frac{\tilde{f}_{\theta}(y)}{\tan\frac{\Delta\theta-y}{2}}\frac{1}{2\pi}\D y+O(N^{-\kappa})=\tilde{\uhlbrt}(\tilde{f}_{\theta})(\Delta\theta)+O(N^{-\kappa})
\end{align*}
where, denoting $\rho=\rhoarcsin $, we have used the following non-trivial observation in the second equality
\begin{align*} 
\dashint_{-2\epsilon}^{2\epsilon}\frac{f(y)}{\Delta-y}(\rho(y+E)-\rho(E))\D y=&O\Big(\epsilon\int_{-\Delta}^{\Delta/2}\frac{|f(y)|}{\Delta}\D y+\epsilon\int_{[-2\epsilon,-\Delta]\cup [3\Delta/2,2\epsilon]}\frac{\log N}{|y|}\D y\Big)
\\
&+\int_{-\Delta/2}^{\Delta/2}\frac{f(\Delta-y)(\rho(\Delta-y+E)-\rho(E))-f(\Delta+y)(\rho(\Delta+y+E)-\rho(E))}{y}\D y
\\
=&O(N^{-\kappa})+\int_{-\Delta/2}^{\Delta/2}\big(\Delta\cdot \|f'\|_{L^\infty((\Delta/2,3\Delta/2))}+O(\log N)\big)\D y=O(N^{-\kappa}).
\end{align*}
\end{proof}

\subsection{Gaussian field approximation}\label{app:gaussian_field_approx}

In the proof of Theorem \ref{thm:GMC}, we followed treatment from \cite{berestycki2018random} and every step except the Gaussian field approximation is almost line to line the same. For the Gaussian field approximation, the involvement of time introduces a technical difficulty that we preferred to address in Appendix (cf. proofs of \cite[Lemma 4.1]{webb2015characteristic} and \cite[Proposition 2.9]{berestycki2018random}).

\begin{lemma}\label{lemma:gaussian_field_approx}
Define $\mu_{N,\gamma}^{(M)}$ and $\mu_{\gamma}^{(M)}$ as in the proof of Theorem \ref{thm:GMC}. Then, for every $\gamma\geq 0$, compactly supported continuous function $\psi:\R\times(-2,2)\to\R$ and $M\in\N$, $\mu_{N,\gamma}^{(M)}(\psi)$ converges in distribution to $\mu_{\gamma}^{(M)}(\psi)$ as $N$ goes to infinity.
\end{lemma}

\begin{proof}
Fix $\psi$ and $M\in\N$ and say $\psi$ is supported in $(-T,T)\times(-2,2)$. Define processes $(Y_N)_t$ and $Z_t$ by
\begin{align*} 
(Y_N)_t&:=\int_{(-2,2)}\psi(t,x)\exp\Big(\gamma\sum_{n=1}^{M}\frac{\tilde{T}_n(x)}{n}S_N(-2\breve{T}_n)(H_t)-\frac{\gamma^2}{2}\sum_{n=1}^{M}\frac{\tilde{T}_n(x)^2}{n})\Big)\D x
\\
Z_t&:=\int_{(-2,2)}\psi(t,x)\exp\Big(\gamma\sum_{n=1}^{M}\frac{\tilde{T}_n(x)}{n} (A_n)_{nt}-\frac{\gamma^2}{2}\sum_{n=1}^{M}\frac{\tilde{T}_n(x)^2}{n}\Big)\D x
\end{align*}
Then $\mu_{\gamma}^{(M)}(\psi)=\int_{(-T,T)}Z_t\D t$ and by Theorem \ref{thm:FH}, 
\begin{align*} 
\mu_{N,\gamma}^{(M)}(\psi)&=\int_{\R\times(-2,2)}\psi(t,x) \frac{\exp\Big(\sum_{n=1}^{M}\Tr\Big(\frac{-2\gamma\cdot\tilde{T}_n(x)}{n}\cdot \breve{T}_n\Big)(H_t)\Big)}{\E{\exp\Big(\sum_{n=1}^{M}\Tr\Big(\frac{-2\gamma\cdot\tilde{T}_n(x)}{n}\cdot \breve{T}_n\Big)(H_t)\Big)}}\D t\D x=(1+O(N^{-\delta}))\int_{(-T,T)}(Y_N)_t\D t
\end{align*}
for some $\delta>0$.

The convergence of $\int_{(-T,T)}(Y_N)_t\D t$ to $\int_{(-T,T)}Z_t\D t$ can be concluded assuming the following three claims: For every $\varepsilon>0$, there exists a $K\in\N$ and $N_0\in\N$ such that letting $t_0=-T,t_1=-T+\frac{2T}{K},\dots,t_K=T$ we have
\begin{itemize}
\item[(a)] $\sum_{i=0}^{K-1}(Y_N)_{t_i}(t_{i+1}-t_i)\xrightarrow[N\to\infty]{d}\sum_{i=0}^{K-1}Z_{t_i}(t_{i+1}-t_i)$,
\item[(b)] $\Prob\Big(\Big|\int_{-T}^{T}Z_t\D t-\sum_{i=0}^{K-1}Z_{t_i}(t_{i+1}-t_i)\Big|<\varepsilon\Big)>1-\varepsilon$ and 
\item[(c)] $\Prob\Big(\Big|\int_{-T}^{T}(Y_N)_t\D t-\sum_{i=0}^{K-1}(Y_N)_{t_i}(t_{i+1}-t_i)\Big|<\varepsilon\Big)>1-\varepsilon$ for all $N\geq N_0$.
\end{itemize}
Because on the event when the last two terms hold, for any $a\in\R$ and $\varepsilon>0$, there exists an $N_{a,\varepsilon}$ such that for all $N\geq N_{a,\varepsilon}$,
\begin{align*} 
\Prob(\int_{-T}^{T}(Y_N)_t\D t\leq a)\leq \varepsilon+\Prob(\sum_{i=0}^{K-1}(Y_N)_{t_i}(t_{i+1}-t_i)\leq a+\varepsilon)\leq 2\varepsilon+\Prob(\sum_{i=0}^{K-1}Z_{t_i}(t_{i+1}-t_i)\leq a+\varepsilon)
\\
\leq 3\varepsilon+\Prob(\int_{-T}^{T}Z_t\D t\leq a+2\varepsilon)
\end{align*}
Hence, $\limsup_{N}\Prob(\int_{-T}^{T}(Y_N)_t\D t\leq a)\leq 3\varepsilon+\Prob(\int_{-T}^{T}Z_t\D t\leq a+2\varepsilon)\xrightarrow[\varepsilon\to0]{}\Prob(\int_{-T}^{T}Z_t\D t\leq a)$ and similarly $\liminf$ is lower bounded by the same quantity. Now, we can move on to  proving the three claims. \\

\noindent\textit{Proof of \textnormal{(a)}:} For any $K$, we already know that 
\begin{align} \label{eqn:conv_in_fdd}
\big(S_N(-2\breve{T}_n)(H_{t_i})\big)_{\substack{n=1,\dots M\\i=0,\dots,K-1}}\xrightarrow[N\to\infty]{d} \big((A_n)_{nt_i}\big)_{\substack{ n=1,\dots M\\i=0,\dots,K-1}}
\end{align}
which can be seen by the explicit covariance kernel given in \cite{unterberger2018global} as we have discussed in Section \ref{subsec:multitime_loop_asymp_exploration}, or directly by Theorem \ref{thm:FH} via the convergence of Laplace transform. Define the function $F:\R^{M}\to\R$ by 
\begin{align*} 
F(x_1,\dots,x_M):=\int_{(-2,2)}\psi(t,x)\exp\Big(\gamma\sum_{n=1}^{M}\frac{\tilde{T}_n(x)}{n}x_n-\frac{\gamma^2}{2}\sum_{n=1}^{M}\frac{\tilde{T}_n(x)^2}{n}\Big)\D x
\end{align*}
which is continuous by the dominated convergence theorem. Define $G:\R^{MK}\to\R$ by
\begin{align*} 
G((x_{i,n})_{\substack{n=1,\dots,M \\ i=0,\dots,K-1}}):=\sum_{i=0}^{K-1} F(x_{i,1},x_{i,2},\dots,x_{i,M})(t_{i+1}-t_i)
\end{align*}
which is continuous as well. By \eqref{eqn:conv_in_fdd},
\begin{align*} 
G\big(\big(S_N(-2\breve{T}_n)(H_{t_i})\big)_{\substack{n=1,\dots M\\i=0,\dots,K-1}}\big)\xrightarrow[N\to\infty]{d}G\big(\big((A_n)_{nt_i}\big)_{\substack{ n=1,\dots M\\i=0,\dots,K-1}}\big)
\end{align*}
which is what we needed.\\

\noindent\textit{Proof of \textnormal{(b)}:} It is a simple corollary of convergence of Riemann integral. Define $S_{K}:=\{|\sum_{i=0}^{K-1}Z_{t_i}(t_{i+1}-t_i)-\int_{-T}^{T}Z_t\D t|<\varepsilon\}$. Then $\Prob(\liminf_{K} S_K)=1$, so there exists a $K_0$ such that for all $K\geq K_0$ we have $\Prob( S_K)>1-\varepsilon$.\\

\noindent\textit{Proof of \textnormal{(c)}:} Fix a large $K$ depending on $\varepsilon$ which is to be determined later. Due to the convergence of $\gamma\sum_{n=1}^{M}\frac{\tilde{T}_n(x)}{n}S_N(-2\breve{T}_n(H_{t_i}))$ to a Gaussian with mean zero and variance of order one for all $i=0,\dots,K-1$, by the Gaussian tail bound there exists an $N_1$ such that for all $N\geq N_1$:
\begin{align*} 
\Prob\Big(\underbrace{\bigcap_{i=0}^{K-1}\bigcap_{n=1}^{M}\Big\{S_N(-2\breve{T}_n(H_{t_i}))\leq (\log K)^{1/2}\log\log K\Big\}}_{=:S_1}\Big)\geq 1-\frac{1}{K^{10}}
\end{align*}
where $K$ is taken to be sufficiently large depending only on $\gamma$ and $M$. Also, assume that there exists an $N_2$ such that for all $N\geq N_2$:
\begin{align}\label{eqn:tr_change_poly}
\Prob\Big(\underbrace{\bigcap_{i=0}^{K-1}\bigcap_{n=1}^{M}\Big\{\sup_{t\in[t_i,t_{i+1}]} \Big|S_N(-2\breve{T}_n(H_{t}))-S_N(-2\breve{T}_n(H_{t_i}))\Big|<\frac{1}{K^{1/4}}\Big\}}_{=:S_2}\Big)\geq 1-\frac{1}{K^{10}}.
\end{align}
On the set $S_1\cap S_2$ we get
\begin{align*} 
\Big|\int_{-T}^{T}(Y_N&)_t\D t-\sum_{i=0}^{K-1}(Y_N)_{t_i}(t_{i+1}-t_i)\Big|
\\
\leq\sum_{i=0}^{K-1} \Big|&\int_{t_i}^{t_{i+1}} \int_{(-2,2)} \big(\psi(t_i,x)+O(\frac{1}{K})\big)\exp\Big(\gamma\sum_{n=1}^{M}\frac{\tilde{T}_n(x)}{n}S_N(-2\breve{T}_n)(H_{t_i})-\frac{\gamma^2}{2}\sum_{n=1}^{M}\frac{\tilde{T}_n(x)^2}{n})\Big)(1+O(\frac{1}{K^{1/3}}))\D x \D t
\\
&-\int_{t_i}^{t_{i+1}} \int_{(-2,2)} \psi(t_i,x)\exp\Big(\gamma\sum_{n=1}^{M}\frac{\tilde{T}_n(x)}{n}S_N(-2\breve{T}_n)(H_{t_i})-\frac{\gamma^2}{2}\sum_{n=1}^{M}\frac{\tilde{T}_n(x)^2}{n})\Big)\D x \D t\Big|
\\
&\hspace{-1cm}=\sum_{i=0}^{K-1}\int_{t_i}^{t_{i+1}} \int_{(-2,2)} O(\frac{1}{K^{1/4}})e^{(\log K)^{2/3}}\D x\D t<\frac{1}{K^{1/5}}
\end{align*}
when $K$ is sufficiently large. Thus, taking $K$ sufficiently large concludes the proof. 

Finally, we prove equation \eqref{eqn:tr_change_poly}. By the union bound it suffices to prove 
\begin{align*} 
\Prob\Big(\sup_{t\in[0,\frac{2T}{K}]} \Big|S_N(-2\breve{T}_n(H_{t}))-S_N(-2\breve{T}_n(H_{0}))\Big|<\frac{1}{K^{1/4}}\Big)\geq 1-\frac{1}{K^{20}}
\end{align*}
for any $n=1,\dots,M$. We can assume that the eigenvalues are in the rigidity set $\tilde{\mathcal{G}}$ for $t\in[-T,T]$ by Proposition \ref{prop:rig_until_logN2}. When the eigenvalues are on the rigidity set, $\breve{T}_n$ will be evaluated only as the Chebyshev polynomial of order $n$. So, we just need to control the change of $\Tr(H_t^n)$ on a time interval of order $\frac{1}{K}$. In fact, we will prove that typical scale for that change is $O(K^{-1/2})$ using Dyson Brownian motion and an induction-like argument over $n$.

We know that $\Tr(H_0^n)$ converges to a Gaussian distribution with $\asymp 1$ variance and mean $\int_{(-2,2)}x^n\rhosc (x)\D x=C_n$ where $C_n:={n\choose n/2}\frac{1}{\frac{n}{2}+1}$ if $n$ is even and $0$ if $n$ is odd. That means, for sufficiently large $N$ we have $\Prob(|\Tr H_0^n-NC_n|>\log K)<\frac{1}{K^{\log\log K}}$ for every $n=1,\dots,M$. Thus, we can assume $|\Tr H_0^n-NC_n|<\log K$ for every $n=1,\dots,M$ from now on. Adding up Dyson Brownian motion \eqref{eqn:DB_motion} for all the eigenvalues we obtain:
\begin{align*} 
\D(\Tr H_t)=\sum_{i=1}^{N}\frac{1}{\sqrt{N}}\D (B_i)_t-\frac{\Tr H_t}{2}\D t.
\end{align*}
Applying BDG to the first term and Gr\"onwall's Lemma for the second term we obtain
\begin{align*} 
\Prob(\sup_{t\in[0,\frac{2T}{K}]} |\Tr H_t-\Tr H_0|<\frac{1}{K^{1/3}})\geq 1-e^{-K^{1/4}}.
\end{align*}
Assume that, for some $m>1$, $\sup_{t\in[0,\frac{2T}{K}]} |\Tr H_t^j-\Tr H_0^j|<\frac{1}{K^{1/3}}$ for every $j=1,\dots,m-1$ with such overwhelming probability and prove that it holds for $m^{th}$ power of trace as well. Note that by our assumptions we have that $|\Tr H_t^j-NC_j|<\log K$ as well, on the set with overwhelming probability for $j=1,\dots,m-1$. Writing the SDE for $\Tr H_t^m$ using Dyson Brownian motion (we drop the time indices in $B_i$'s and $\lambda_i$'s for convenience) we obtain
\begin{align*} 
\D(\Tr H_t^m)=\sum_{i=1}^{N}m\lambda_i^{m-1}\D \lambda_i=\sum_{i=1}^{N}\frac{m\lambda_i^{m-1}}{\sqrt{N}}\D B_i+m\Big(\frac{1}{N}\sum_{i<j}\frac{\lambda_i^{m-1}-\lambda_j^{m-1}}{\lambda_i-\lambda_j}-\frac{1}{2}\sum_{i=1}^{N}\lambda_i^{m}\Big)\D t.
\end{align*}
For the martingale term we still have the same type of bound by BDG. Using the induction hypothesis, the term in the parentheses can be written as
\begin{align*} 
\frac{1}{N}\sum_{i<j}&\frac{\lambda_i^{m-1}-\lambda_j^{m-1}}{\lambda_i-\lambda_j}-\frac{1}{2}\sum_{i=1}^{N}\lambda_i^{m}=\frac{N-1}{N}\sum_{i=1}^{N}\lambda_i^{m-2}+\frac{1}{2N}\sum_{k=1}^{m-3}\Big(\sum_{i\neq j}\lambda_i^{k}\lambda_j^{m-2-k}\Big)-\frac{1}{2}\sum_{i=1}^{N}\lambda_i^{m}
\\
&=O((\log K)^2)+\frac{N-1}{N}(NC_{m-2})+\frac{1}{2N}\sum_{k=1}^{m-3}\big(NC_{k}NC_{m-2-k}-NC_{m-2}\big)-\frac{1}{2}(\Tr H_t^m-\Tr H_0^m)-\frac{1}{2}NC_{m}
\\
&=O((\log K)^2)+\frac{N}{2}\sum_{k=0}^{m-2}C_{k}C_{m-2-k}-\frac{N}{2}C_m-\frac{1}{2}(\Tr H_t^m-\Tr H_0^m)=O((\log K)^2)-\frac{1}{2}(\Tr H_t^m-\Tr H_0^m)
\end{align*}
where we have used the recurrence relations of Catalan numbers. Thus, Gr\"onwall's inequality gives the required upper bound for $\sup_{t\in[0,\frac{2T}{K}]}|\Tr H_t^m-\Tr H_0^m|$.
\end{proof}

{\small
\bibliographystyle{alphaabbr}
\bibliography{references}
}

\end{document}